\documentclass[11pt,reqno]{amsart} 
\pdfoutput=1
\usepackage{CLY_Macros}

\begin{document}
\title{Asymptotic Stability of Moving Kinks with an Internal Mode}
\author[G. Chen]{Gong Chen}
\author[Y. Li]{Yongming Li}
\author[Z. Yang]{Zhaojie Yang}

\email{gc@math.gatech.edu}
\email{yli3349@gatech.edu}
\email{zyang768@gatech.edu}
\address{School of Mathematics, Georgia Institute of Technology, Atlanta, GA 30332, USA}

\thanks{G.C. was partially supported by NSF Grants DMS-2350301 and CAREER-DMS-2540992, by the Simons Foundation MP-TSM-00002258, and by the AMS Stefan Bergman Fellowship.}
\begin{abstract}

We prove asymptotic stability of moving kinks for one-dimensional
scalar field equations with one internal mode and no threshold
resonance, under a transformed repulsivity condition  introduced in {Kowalczyk, Martel, Mu\~noz, and Van Den Bosch~}~\cite{kowalczyk2021sufficient} and a Fermi golden
rule at the second harmonic. Near every fixed subluminal velocity,
arbitrary sufficiently small perturbations in $H^1\times L^2$ converge
locally in the energy space to a kink with a limiting velocity.
The principal
difficulty is the coupling between  the damping of the internal mode and modulations: the internal oscillation enters the modulation equations
quadratically, whereas the damping mechanism provides only quartic
time integrability. Inspired by {Martel}~\cite{martel2024asymptotic},  we derive a velocity normal form from conserved
momentum and incorporate the corrected velocity into the virial
functional, cancelling the quadratic modulation terms that obstruct
the local energy estimates. Matrix Darboux transformations, virial
identities at two spatial scales, and a localized Fermi golden rule
then yield simultaneous spacetime control of the internal mode and
radiation, together with convergence of the velocity. 

We apply our criterion to the double sine-Gordon equation and suitable perturbations of the $\phi^4$ model. Through a perturbative analysis of threshold resonances, we determine the full-space discrete spectrum and establish the absence of threshold resonances. We further verify the repulsivity of the twice-transformed potential and the Fermi golden rule.

\end{abstract}

\maketitle 

\vspace{-2em}
\tableofcontents

\section{Introduction}
\subsection{Scalar field equation}

We consider scalar field theories on the line governed by the Lagrangian action functionals 
\begin{equation} \label{eq:intro_lagrangian}
    \iint_{\bbR^{1+1}} \biggl( \frac12 (\pt \phi)^2 - \frac12 (\px \phi)^2 - W(\phi) \biggr) \, \ud x \, \ud t,
\end{equation}
where $\phi \colon \bbR \times \bbR \rightarrow \bbR$ is the scalar field and $W \in \calC^\infty\big(\bbR;[0,\infty)\big)$ 
is a double-well self-interaction potential. We assume that $\zeta_-<\zeta_+$ are consecutive nondegenerate minima satisfying
\begin{equation}\label{eq:vacua}
 \begin{gathered}
 W(\zeta_\pm)=W'(\zeta_\pm)=0,\qquad
 W''(\zeta_\pm)=\omega_\pm^2>0,\\
 W(s)>0\qquad\text{for }\zeta_-<s<\zeta_+.
 \end{gathered}
\end{equation}
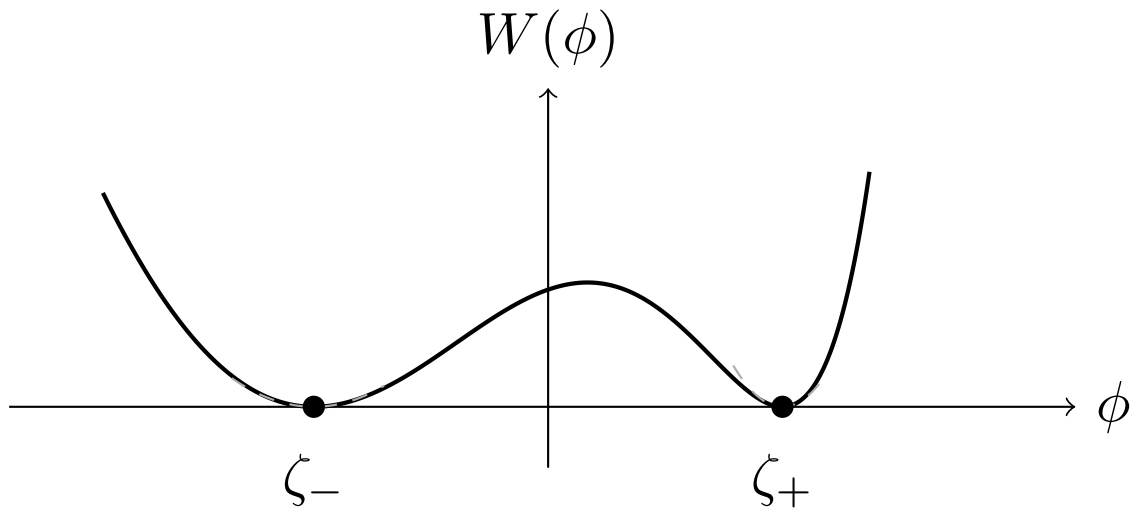
\begin{figure}[t]
    \begin{center}
\begin{tikzpicture}[scale=1.55]
  \draw[->, thin] (-2.30,0) -- (2.25,0) node[right] {$\phi$};
  \draw[->, thin] (0,-0.26) -- (0,1.36) node[above] {$W(\phi)$};
  \draw[thick, smooth, samples=250, domain=-1.9:1.375]
      plot (\x,{0.5*(1-\x*\x)^2*exp(0.69314718056*\x)});
  \draw[thin, dashed, gray!65, samples=45, domain=-1.35:-0.70]
      plot (\x,{(\x+1)^2});
  \draw[thin, dashed, gray!65, samples=45, domain=0.79:1.19]
      plot (\x,{4*(\x-1)^2});
  \fill (-1,0) circle (1.35pt);
  \fill (1,0) circle (1.35pt);
  \node[below=5pt] at (-1,0) {$\zeta_-$};
  \node[below=5pt] at (1,0) {$\zeta_+$};
\end{tikzpicture}
\end{center}
\caption{Double-well potential}
\end{figure}
  
The Euler-Lagrange equations associated with \eqref{eq:intro_lagrangian} are the wave-type equations
\begin{equation}\label{eqn:scalar_field_equation}
\pt^2\phi -\px^2 \phi + W'(\phi) = 0.
\end{equation}

Writing $\bmphi=(\phi,\phi_t)$, the conserved energy and momentum are
\begin{equation}\label{eq:conservation}
 \mathcal E[\bmphi]=\frac12\int_\bbR
 \bigl(\phi_t^2+\phi_x^2+2W(\phi)\bigr)\,\ud x,
 \qquad
 \mathcal P[\bmphi]=\int_\bbR\phi_t\phi_x\, \ud x .
\end{equation}
Setting $\bmH=(H,0)$, we define the affine energy space in the kink sector by
\begin{equation}\label{eq:intro-energy-space}
 \begin{aligned}
 \bmE_{\bmH}
 :=\bmH+\bigl(H^1(\bbR)\times L^2(\bbR)\bigr)
 =\bigl\{(\phi_1,\phi_2):\phi_1-H\in H^1(\bbR),\quad \phi_2\in L^2(\bbR)\bigr\}.
 \end{aligned}
\end{equation}
Up to translation, the increasing kink solution to \eqref{eqn:scalar_field_equation},  $H$,  is determined by
\[
 H''=W'(H),\qquad H(\pm\infty)=\zeta_\pm,
 \qquad H'=\sqrt{2W(H)}>0.
\]
In particular, $H-\zeta_\pm$ and the derivatives of $H$ decay exponentially at the corresponding infinities.

The invariance of \eqref{eqn:scalar_field_equation} under Lorentz transformations and under spatial translations leads to the following two-parameter family of moving kink solutions: for $c\in(-1,1)$ and $y_0 \in \bbR$, set
\[
 \gamma=\gamma(c)=(1-c^2)^{-1/2},\qquad
 H_c(x)=H(\gamma x),\qquad
 \bmH_c=(H_c,-cH_c').
\]
Then $\bmH_c(x-ct-y_0)$ is an exact moving kink.

In this paper, we are interested in  the long-time behavior of solutions to \eqref{eqn:scalar_field_equation} close to a moving kink.   
The variational characterization of kinks gives orbital stability: a small perturbation remains close to the family of translated and boosted kinks. This does not determine whether the perturbation leaves the vicinity of the kink, whether a localized internal oscillation survives, or whether the velocity approaches a limit.

We address the relaxation of the internal mode and the convergence of the velocity together.  The present work concerns the case of one internal mode and no threshold resonance, without any symmetry condition on the perturbation. At the linear level, the internal mode is a localized periodic oscillation and cannot decay. At the nonlinear level, its second harmonic can radiate into the continuous spectrum. For a moving kink, however, this transfer of energy is coupled to translation and velocity modulation. The internal mode oscillation contributes to the acceleration at quadratic order, whereas the radiation damping mechanism naturally controls its fourth power in time. 

We combine the two-scale virial estimates of Kowalczyk, Martel, Mu\~noz, and Van Den Bosch~\cite{kowalczyk2021sufficient} and the localized internal mode analysis of Kowalczyk and Martel~\cite{kowalczyk2022kink}.

The additional step is a velocity normal form derived from conserved momentum and its incorporation into the coefficients of the moving kink virial functional, which is inspired by \cite{martel2024asymptotic}. See \cite{cuccagna2025cubic,rialland2025internal} for similar normal form computations in the context of NLS. In this way, the nonintegrable terms are removed before the local energy and internal mode estimates are closed. In particular, we obtain the asymptotic stability of a moving kink measured in the localized energy norm under small energy perturbations without symmetry assumptions.

\subsection{Main result}
Before stating our main results, we introduce the necessary notation and assumptions.

 We use the equivalent energy norm
\begin{equation}\label{eq:norm}
\|\bmv\|_c^2=\gamma^{-1}\|v_1'\|_{L^2_x(\bbR)}^2+\gamma\|v_1\|_{L^2_x(\bbR)}^2
  +\gamma\|v_2+cv_1'\|_{L^2_x(\bbR)}^2.
\end{equation}
The notation $\|\bmv\|_{c,R}$ means that the three $L^2$-norms are restricted to $|x|<R$. Equivalence with the usual $H^1 \times L^2$ energy norm is uniform when $c$ ranges over a compact subset of $(-1,1)$.

The spectral assumptions are most simply expressed in the rest frame. Let $L_0 = -\px^2 + W''(H)$ be the linearized operator around the static kink. We make the following spectral assumptions 
\begin{equation}\label{eq:spectrum}
\left\{\begin{aligned}
    &\sigma(L_0) = \sigma_{\mathrm{p}}(L_0) \cup \sigma_{\mathrm{ac}}(L_0), \qquad \sigma_{\mathrm{p}}(L_0) = \{0,\lambda^2\}, \qquad \sigma_{\mathrm{ac}}(L_0) = [\omega^2,\infty),\\
&\text{where} \quad \omega=\mathrm{min} \{{\omega_-,\omega_+\}\quad\text{and}\quad\frac{\omega}{2}<\lambda<\omega.}
\end{aligned}\right.
\end{equation}
Both eigenvalues $0$ and $\lambda^2$ are simple. The kernel is generated by $H'$, and we choose a real normalized eigenfunction $Y$ such that $L_0Y=\lambda^2Y$. 

To state a sufficient condition for local decay, we remove the two eigenfunctions by successive Darboux transformations:
\begin{equation}\label{eq:darboux}
 \begin{gathered}
 U_0=\partial_x-\frac{H''}{H'},\qquad
 L_0=U_0^*U_0,\qquad L_1=U_0U_0^*,\\
 Y_1=U_0Y>0,\qquad U_1=\partial_x-\frac{Y_1'}{Y_1},\\
 L_1=\lambda^2+U_1^*U_1,\qquad
 L_2=\lambda^2+U_1U_1^*=-\partial_x^2+P_2.
 \end{gathered}
\end{equation}

After fixing the center of $H(x)$, we can choose $Y(x)$ such that $Y(0)=0$ and take $Y_1$ to be a positive ground state of $L_1$. These factorizations give
\[
 U_1U_0L_0=L_2U_1U_0.
\]
We assume that $P_2$ is nonconstant and repulsive about some $x_0\in\bbR$:
\begin{equation}\label{eq:repulsivity}
 (x-x_0)P_2'(x)\leq0\quad\text{for all }x\in\bbR,
 \qquad P_2'\not\equiv0.
\end{equation}
This condition excludes a bounded threshold state for $L_2$, and hence a threshold resonance for $L_0$.

The remaining assumption concerns nonlinear coupling to the continuum. We require a real bounded solution $g$ of $L_0g=4\lambda^2g$ such that
\begin{equation}\label{eq:fgr}
 \Gamma=-\frac12\int_\bbR W'''(H)Y^2g\,\ud x \ne0.
\end{equation}
Changing the sign of $g$, we take $\Gamma>0$, which is a convenient form of the Fermi golden rule. 
The strict inequality $2\lambda>\omega$ places the second harmonic in the continuous spectrum. 

\begin{theorem}\label{main_thm1} Assume \eqref{eq:vacua}, \eqref{eq:spectrum}, \eqref{eq:repulsivity}, and \eqref{eq:fgr}. For any $c_0 \in (-1,1)$, there exists $\delta_0>0$ such that for any $\bmphi_0 \in \bmE_{\bmH}$ with 
\begin{equation}
    \| \bmphi_0- \bmH_{c_0} \|_{c_0} \leq \delta_0,
\end{equation}
there exist $c_-,c_+ \in (-1,1)$ and a $\calC^1$ function $y: \bbR \rightarrow \bbR$ with $\lim_{t \rightarrow \pm \infty} \dot{y}(t) = c_\pm$ such that the global solution $\bmphi$ of \eqref{eqn:scalar_field_equation} with $\bmphi(0) = \bmphi_0$ satisfies for any $R>0$
\begin{equation}
    \lim_{t \rightarrow \pm \infty} \| \bmphi(t,\cdot + y(t) ) - \bmH_{c_\pm} \|_{c_\pm,R} = 0.
\end{equation}
\end{theorem}

\begin{remark}
The proof of Theorem~\ref{main_thm1} builds on the general virial framework of \cite{kowalczyk2021sufficient}, and relies on establishing an integrated spacetime estimate for both the internal mode and the radiation. For suitable modulation parameters $c(t)$ and $y(t)$, let $a(t)$ denote the internal mode amplitude and let $\bmnu(t)$ denote the radiation component. The integrated-in-time local energy decay estimate \eqref{eqn:final_estimate},
\begin{equation*}
    \int_0^\infty \left( \big|a(t)\big|^4 +\left \| \sech\left(\tfrac{\kappa \gamma(t) }{40}\cdot\right)\bmnu\big(t,\cdot+y(t)\big) \right\|_{c(t)}^2 \right)\,\ud t \leq C \delta^2 < \infty,
\end{equation*}
combined with the differential estimates in the proof of Proposition~\ref{prop:main}, yields decay for the internal mode and local decay for the radiation.  Here $\kappa := \sqrt{\omega^2-\lambda^2}> 0$ denotes the size of the spectral gap. The velocity normal form derived in Lemma~\ref{lemma:normal_form_c} and the modulation equations imply the limits $\dot{y}(t)\rightarrow c_+$ and $c(t)\rightarrow c_+$. Reversing time implies the corresponding conclusions as $t \rightarrow -\infty$.
\end{remark}

\begin{remark}We illustrate the applications of our theorem with the double sine-Gordon equation and a suitable perturbation of the $\phi^4$ model; see Section~\ref{subsec:models} for details.
\end{remark}

\subsection{Literature}
The study of the asymptotic stability of solitons is a rich and vast subject that we
do not attempt to comprehensively review here. Instead, we focus on the significant advancements
made in recent years regarding the asymptotic stability of kinks in scalar field theories on the line
and of solitary waves in nonlinear Schr\"odinger models in one space dimension.   For more complete references and historical introductions, we refer to the surveys \cite{germain2024review,kowalczyk2017asymptotic,martel2024survey}.

The orbital stability of kinks in relativistic scalar field theories on the line is classical by the variational theory of Henry, Perez and Wreszinski~\cite{henry1982stability}.  The variational approach does not distinguish radiation from a positive eigenmode of the linearized operator. The latter generates a time-periodic solution of the linear equation. Its relaxation, if it occurs, is necessarily nonlinear: higher harmonics of the internal oscillation may enter the continuum and transfer energy to spatially propagating waves. This picture goes back to the analysis of wobbling kinks~\cite{segur1983wobbling} and the general theory of nonlinear resonances and radiation damping~\cite{sigal1993instability,soffer1999resonances}. In particular, a nonzero Fermi golden rule measures this coupling between the internal oscillation and the continuous spectrum. 

The asymptotic stability problem establishes a stronger property: small perturbations, in fact, disperse over time. In order to establish the asymptotic stability, one has to understand the coupling between the internal mode and radiation and derive the decay mechanism.  
The earliest kink asymptotic stability results relevant here were obtained by Kopylova and Komech~\cite{kopylova2011kink,kopylova2011moving}. They treated, respectively, an internal mode in a parity-restricted setting and nonsymmetric perturbations of moving kinks when the internal spectrum is absent. Their weighted dispersive analysis assumes, among other properties, sufficiently flat nonlinearities near the limiting vacua. The latter condition removes the persistent quadratic far-field terms that are present for the usual polynomial double-well potentials.  But neither covers a moving kink with an internal mode under arbitrary small, unweighted energy perturbations.

For the $\phi^4$ kink, the seminal work by Kowalczyk, Martel and Mu\~noz~\cite{kowalczyk2017kink} proved asymptotic stability for odd perturbations in the energy space.  The translation mode of the kink is even, as is its bounded threshold resonance, whereas its internal mode is odd. Oddness thus removes the two former obstructions while leaving precisely the localized oscillation that must be damped. The argument combines a virial estimate for radiation with a nonlinear functional detecting transfer at the second harmonic. Snelson~\cite{Snelson_2018_VariableSpeedKink} extended
this odd-perturbation theory to stationary $\phi^4$ kinks
with a small spatial variation of the propagation speed.
For related results on the absence of small odd breathers,
see~\cite{kowalczyk2017breathers}.

A different obstruction, namely the motion of the kink itself, is treated in the work of Kowalczyk, Martel,
Mu\~noz and Van Den Bosch~\cite{kowalczyk2021sufficient}. They establish a sufficient condition for
asymptotic stability of moving kinks under
general energy space perturbations, without any symmetry restriction. Their proof uses a Darboux transformation
to remove the translation mode and reveal a repulsive
potential. The resulting coercivity is combined with a
virial estimate in the original variables, and
weighted estimates relate the two systems to obtain
local energy decay. Their criterion excludes internal
modes, whose interaction with velocity modulation is
a central difficulty in the present work.

Kowalczyk and Martel~\cite{kowalczyk2022kink} developed the complementary theory for a kink with one internal mode under odd perturbations. A second Darboux transformation removes this mode, while a localized Fermi golden rule controls its interaction with the continuum. The criterion is robust under suitable perturbations of the potential and includes the $\phi^4$ equation in the odd class. Cuccagna and Maeda~\cite{cuccagna2023kink} subsequently allowed configurations with several internal modes, still under odd perturbations; they explicitly identify the extension to moving kinks as a separate problem. For small nonlinear Klein--Gordon solutions with a potential, Cuccagna, Maeda and Scrobogna~\cite{cuccagna2023stabilization} combine the transformed virial framework with refined profiles to analyze the discrete dynamics.

For the related problem of unstable solitons for the one-dimensional nonlinear Klein--Gordon equation 
\begin{equation}
    \pt^2\phi - \px^2\phi + \phi = \phi^p,
\end{equation}
Kowalczyk, Martel and
Mu\~noz~\cite{kowalczyk2022jems} proved conditional asymptotic stability in the energy space for even solutions on a Lipschitz center hypersurface, in the regime $p>3$ where no internal modes are present. The virial and Fermi golden rule treatment of an internal mode for the quadratic Klein--Gordon soliton under even symmetry is another related example~\cite{li2023quadratic}. Cuccagna, Maeda,
Murgante and Scrobogna~\cite{Cuccagna_2024_CenterHypersurface}
extended this result to $5/3<p\leq 2$.
For $p=3$, the presence of the even threshold resonance poses significant challenges; 
Palacios and Pusateri~\cite{palacios2026local} obtain time-weighted local energy control over exponentially long times for even perturbations on its center-stable manifold. They also considered an application to odd perturbations of the sine--Gordon kink.

Thus, the existing energy-space results provide two different ingredients: The freedom of translation and velocity in~\cite{kowalczyk2021sufficient}, and radiative damping of bound states in~\cite{kowalczyk2022kink,cuccagna2023kink}. The present setting requires both ingredients without a parity decomposition. 
The closest analogy for the treatment of the internal mode is Martel's result~\cite{martel2024asymptotic} for the focusing--focusing cubic--quintic Schr\"odinger equation
\begin{equation}\label{eq:cq}
 i\psi_t+\psi_{xx}+|\psi|^2\psi+|\psi|^4\psi=0.
\end{equation}
It concerns small standing solitary waves under \emph{even} energy perturbations. A discrete eigenvalue bifurcates from the resonance at the edge of the essential spectrum of the cubic problem. The proof combines successive transformations, a localized Fermi golden rule and virial estimates at two scales; frequency and phase modulation and a normal-form correction enter the analysis. The corresponding focusing--defocusing equation, for which there is no internal mode in the relevant regime, was treated earlier in~\cite{martel2023asymptotic}. The distinction between these two spectral situations is explained further in Martel's survey~\cite{martel2024survey}. Rialland~\cite{rialland2024absence,rialland2025internal} obtained corresponding results for broader near-cubic nonlinearities. The creation or disappearance of a discrete state at a threshold is related to the weak-coupling and edge-bifurcation mechanisms in~\cite{simon1976bound,coles2016edge}. Other near-power NLS results of Cuccagna and Maeda~\cite{cuccagna2025cubic,cuccagna2025quadratic} emphasize that removing parity restrictions and treating an internal mode are distinct spectral tasks; in particular, their near-cubic theorem is restricted to even solutions.

Martel's work \cite{martel2024asymptotic}  provides a principal inspiration for our treatment of the internal mode. Its adaptation to moving kinks, however, is not straightforward. Without symmetry restrictions, the translation--velocity directions must be treated together with the internal mode in a Lorentz-modulated system. Moreover, the quadratic coefficient $W'''(H)$ need not decay at spatial infinity, and the moving-wave virial does not directly control the full localized energy. These features require combining the internal mode analysis with the matrix Darboux transformations and transfer estimates developed from the moving kink framework of \cite{kowalczyk2021sufficient}. A further difficulty is the coupling between internal mode damping and velocity modulation. The key additional step is to incorporate a momentum-based velocity correction inspired by \cite{martel2024asymptotic} into the virial functional itself, allowing the modulation and local decay estimates to be closed together.

A complementary approach to kink stability exploits weighted dispersive estimates, normal forms and distorted Fourier analysis. For quadratic {Klein--Gordon} equations with a potential, see Shatah's normal-form method~\cite{shatah1985normal}, the work of Germain and Pusateri~\cite{germain2022quadratic}, and the asymptotic analyses of Lindblad, L\"uhrmann and Soffer~\cite{lindblad2021asymptotics} and Lindblad, L\"uhrmann, Schlag and Soffer~\cite{lindblad2023modified}. Delort and Masmoudi~\cite{delort2022kink} obtain long-time dispersive information near a kink under their additional localization and regularity assumptions. A broader account of approaches in one dimension is given in~\cite{germain2024review}. Such methods resolve detailed radiation asymptotics, including threshold resonant corrections that cannot generally be read from a local-energy inequality. They also indicate why low-dimensional quadratic terms and spectral thresholds must be handled carefully. 
For equations with an internal mode in other settings, L\'eger and
Pusateri~\cite{leger2021internal} analyze radiation damping in the
three-dimensional quadratic Klein--Gordon setting. Chen, Diebou,
Kairzhan, and Pusateri~\cite{chen2026global} treat a one-dimensional
model with localized quadratic and nonlocalized quartic terms,
obtaining internal mode asymptotics and radiation scattering for
a class of weighted Sobolev data. Factorization and radiation damping also underlie the asymptotic
stability result of L\"uhrmann et al.~\cite{luhrmann2026asymptoticstabilitydegreeonevortexstability}
for the degree-one abelian Yang--Mills--Higgs vortex at self-dual
coupling under small equivariant perturbations in weighted
Sobolev spaces. Their companion papers
establish the spectral and Fermi golden rule
conditions~\cite{luhrmann2026asymptoticstabilitydegreeonevortex} and the dispersive
and local energy estimates~\cite{luhrmann2026asymptotic} used
in the nonlinear analysis. 

The sine--Gordon equation presents the sharpest contrast. The linearized operator around the kink has no internal eigenmode, but it has a threshold resonance; the completely integrable evolution also supports breathers and kink--breather configurations, see \cite{chen2020long,koch2026asymptotic}. Unrestricted local asymptotic stability in the energy space is therefore false at this endpoint. Alejo, Mu\~noz and Palacios~\cite{alejo2023sinegordon} identify sets of energy data on which asymptotic stability does hold through B\"acklund transformations. L\"uhrmann and Schlag~\cite{luhrmann2023sinegordon} prove asymptotic stability under odd perturbations. Chen and L\"uhrmann~\cite{chen2024sinegordon} establish asymptotic stability without symmetry under small perturbations in weighted Sobolev spaces, using a distorted-Fourier and space--time resonance approach; the threshold resonance of the moving linearization and quadratic null structures play essential roles. Their treatment of the moving translation and boost parameters is closely related in spirit to our problem, but the spectral obstacle and the initial-data topology differ: we assume the absence of a threshold resonance, retain one internal mode, and work in the unweighted energy space. These assumptions are justified for the perturbations of the sine-Gordon equation and $\phi^4$ model.  Our asymptotic stability results in particular \emph{rule out} the existence of wobbling kinks in our model. 
Finally, we mention that nonpersistence of breather families under perturbation~\cite{denzler1993nonpersistence} and the small-energy breather analysis of Kowalczyk and Martel~\cite{kowalczyk2026breathers} also caution against extrapolating estimates uniformly to the integrable limit.

\subsection{Perturbations of sine--Gordon and \texorpdfstring{$\phi^4$}{phi4}}\label{subsec:models}

To apply Theorem~\ref{main_thm1} to scalar field models, we construct perturbations of the sine--Gordon equation and the $\phi^4$ model such that all the hypotheses of our main theorem hold. 

The sine--Gordon and $\phi^4$ equations provide two different starting points for the one-internal mode problem. The sine--Gordon kink has a threshold resonance, but no internal mode; a perturbation can therefore
create the oscillation whose damping must be proved. The $\phi^4$ kink already has an internal mode, together with an additional threshold resonance. In this case the relevant perturbation preserves the internal
mode while removing the threshold obstruction.

\medskip\noindent
\emph{The double sine--Gordon model.}
The double sine--Gordon family itself has long been studied in connection with kink interactions~\cite{campbell1986kink}. 
Following~\cite{campbell1986kink},  we take
\begin{equation}\label{eq:dsg}
 W_{\mathrm{dsG}}(\phi)
 =\frac{4}{1+|4\eta|}
 \left[1+\cos\frac{\phi}{2}+\eta(1-\cos\phi)\right],
 \qquad |\eta|\ll1.
\end{equation}
Here an additive constant has been chosen so that the consecutive vacua $\zeta_\pm=\pm2\pi$ have zero potential energy. This does not change the equation or its linearization. The vacuum masses satisfy
\[
 W_{\mathrm{dsG}}''(\zeta_\pm)
 =\frac{1+4\eta}{1+|4\eta|};
\]
in particular, they are equal to one in the positive-parameter regime. At $\eta=0$, the linearized operator around the static kink, $H(x) = 4\arctan(\sinh(x))$, coincides with that of the sine--Gordon model, and is given by
\begin{equation}\label{eq:sglinear}
 L_0=-\partial_x^2+1-2\sech^2(x).
\end{equation}
It has an even translational mode and an odd threshold resonance $R(x)=\tanh(x)$ with $L_0R = R$. For sufficiently small $\eta>0$, the threshold resonance 
bifurcates into a simple internal eigenvalue, with
\begin{equation}\label{eq:dsgbifurcation}
\lambda^2 = 1 - \frac{256}{9}\eta^2 + O(\eta^3).
\end{equation}
Also see \cite[Figure 16]{campbell1986kink} for the numerical bifurcation graph of the eigenvalue.

The opposite sign, $\eta < 0$, removes the resonance without creating an internal mode and 
the moving kink is asymptotically stable by \cite[Theorem~10]{kowalczyk2021sufficient}. 

Our interest is in the positive branch. Although the internal mode lies below the continuous spectrum, its second harmonic lies inside it. The spectral bifurcation alone does not determine the resulting energy transfer: the limiting sine--Gordon coupling vanishes. The model
calculation must therefore follow the bifurcating mode and the continuum state together, with their normalizations fixed, to identify the first
nonzero term in the Fermi golden rule. In 
Proposition~\ref{prop:dsg-main}, it is shown that
\begin{align}\label{eq:intro-dsg-fgr}
 \Gamma_{\mathrm{dsG}}(\eta)
   =\Gamma_0\eta+O(\eta^2),\qquad \Gamma_0
   =8\pi\sech\!\left(\frac{\sqrt3\pi}{2}\right) \approx3.294635.
\end{align} 
For the sine--Gordon model, this coefficient vanishes at $\eta=0$. This is a key null condition for the quadratic term that allows normal form transformations in the study of the asymptotic stability of the sine-Gordon kink; see \cite{luhrmann2023sinegordon,chen2024sinegordon}. For the positive $\eta$ branch, after the translation and internal modes are removed, repulsivity is checked for the second Darboux transformed potential.

Then our results can be applied to the positive branch. Together with the asymptotic stability results for the negative branch from \cite{kowalczyk2021sufficient}, it follows that, at least along this family parametrized by small $\eta$, the existence of wobbling kinks is tied to the threshold resonance at the bottom of the continuous spectrum of the linearized operator around the sine-Gordon kink.

\medskip\noindent
\emph{A perturbation of the $\phi^4$ model.}
The second family is
\begin{equation}\label{eq:sextic}
 W_{\mathrm{perturbed}\,\phi^4}(\phi)
 =\frac{(1-\phi^2)^2}{8(1-\eta)}
        \bigl(1-\eta\chi(\phi)\phi^2\bigr),
 \qquad 0<\eta<\eta_0<\frac18.
\end{equation}
Here $\chi\in \calC^\infty(\bbR;[0,1])$ is even, equals $1$ on
$[-3/2,3/2]$, and vanishes outside $(-2,2)$. The cutoff here ensures that the double-well potential is nonnegative. In this model, the vacua remain at $\pm1$, with vacuum mass one. At $\eta=0$, the linearized operator around the static kink, $H(x) = \tanh(x/2)$, coincides with that of the $\phi^4$ model, and is given by 
\begin{equation}\label{eq:phi4}
  L_0=-\partial_x^2+1-\frac32\sech^2(x/2).
\end{equation}
Besides the translational mode, $L_0$ has an odd internal mode 
\begin{equation}
    Y(x) = \tfrac{\sqrt{3}}{2}\tanh(x/2)\sech(x/2), \quad (L_0 - \tfrac34)Y = 0,
\end{equation}
and an even threshold resonance
\begin{equation}\label{eq:phi4resonance}
 R(x)=1-\tfrac32\sech^2(x/2),\qquad (L_0-1)R=0.
\end{equation}
Thus odd perturbations retain the shape mode but exclude both translation and the threshold resonance, as in~\cite{kowalczyk2017kink,kowalczyk2022kink}.
For the stability problem for the $\phi^4$ kink without symmetry, the even translational mode and threshold resonance must also be treated.

In the convention~\eqref{eq:sextic}, the one-internal mode branch is
$\eta>0$. The internal eigenvalue persists and satisfies
\begin{equation}\label{eq:phi4eigenvalues}
 \lambda^2=\frac34+\frac35\eta+O(\eta^2),
\end{equation}
whereas the threshold resonance is removed without producing another
internal eigenvalue. The sign is important: small negative $\eta$ creates an additional even internal mode, so this is not the one-internal-mode regime.
For positive $\eta$, the expression~\eqref{eq:sextic} is used on a
neighborhood of $[-1,1]$ and extended smoothly and nonnegatively outside
that neighborhood. This extension leaves the kink, the linearized
operator and the coupling coefficients along the kink unchanged.  By a perturbative computation of the Fermi golden rule, we show in Proposition~\ref{prop:phi4-main} that
\begin{align}
\Gamma_{\mathrm{perturbed}\,\phi^4}(\eta)
   =\Gamma_{\phi^4}+O(\eta),\qquad
 \Gamma_{\phi^4}
   =\frac{27\sqrt2\pi}{8\sinh(\sqrt2\pi)}
     \approx0.352784.\label{eq:intro-phi4-fgr}
\end{align}
For the negative branch, the threshold resonance \eqref{eq:phi4resonance} bifurcates into an additional internal mode. This will require a generalization of the Fermi golden rule argumnent.

\subsection{Notation}
In this subsection, we specify our conventions used in this paper.

For complex-valued vector functions $\bmf = (f_1,f_2)$ and $\bmg=(g_1,g_2)$,  we use the complex Hermitian product 
\begin{equation}
    \langle \bmf,\bmg\rangle_{\bbC} := \int f_1 \barg_1 + \int f_2 \barg_2,
\end{equation}
and the real inner product 
\begin{equation}
    \langle \bmf,\bmg\rangle := \Re\langle \bmf,\bmg\rangle_{\bbC} = \Re \int f_1 \barg_1 + \Re \int f_2 \barg_2.
\end{equation}
The symplectic form for complex-valued vector functions $\bmf,\bmg$ is given by 
\begin{equation}
    \bmomega(\bmf,\bmg) := \int \bmf \cdot \bmJ \bmg = \int f_1 g_2 - f_2 g_1,
\end{equation}
where
\begin{equation}
	\bmJ := \begin{pmatrix}
	    0 & 1 \\ -1 & 0 
	\end{pmatrix}.
\end{equation}
%
For operators $A$ and $B$, their commutator is given by 
\begin{equation}
    [A,B] := AB - BA.
\end{equation}
Throughout this paper, we will interchangeably use  $\|\cdot\| \equiv  \|\cdot\|_{L_x^2(\bbR)}$.
%
All implicit constants $C$ may depend on the double-well potential $W(\phi)$, the spectral parameters $\omega^2$, $\lambda^2$, or the initial velocity $c_0$. As usual, we write $X \lesssim Y$ if $X \leq CY$ for non-negative numbers $X,Y$, and we use the notation $X \ll Y$ to indicate that the implicit constant should be regarded as small. We also use the Japanese bracket notation,
\begin{equation}
    \jx := (1+x^2)^{\frac12}.
\end{equation}

\subsection{Organization of the paper}
The proof of Theorem~\ref{main_thm1} occupies Sections 2-9 and rests on spacetime virial estimates for the radiation and the internal mode components. Section~2 reviews the spectral theory of the linearized operator, modulation theory, and orbital stability around moving kinks. In Section~3, we pass the spectral decomposition of the perturbation to coordinates adapted to the reference kink $H(x)$ (which were introduced in \cite[\S4.1]{kowalczyk2021sufficient}), and derive the corresponding evolution equations for the radiation and the internal mode components. We also briefly describe the virial method and introduce weights and multipliers used for the virial functionals. In Section~4, we construct a velocity normal form using the conserved momentum and establish a local energy estimate for the radiation at large spatial scales. For this estimate, we incorporate the corrected velocity into the coefficients of a modified virial functional to cancel the worst terms. Next, Section~5 uses the Fermi golden rule to obtain quartic time integrability of the internal mode amplitude. Section~6 constructs the matrix Darboux transformations and derives the equations for the transformed and regularized radiation variables. Section~7 establishes the weighted coercivity estimates needed to transfer control back to the original radiation variables. Then in Section~8, we exploit the repulsivity of the twice-transformed potential in a virial argument at a smaller spatial scale. Lastly, in Section~9, we combine all of these virial estimates to complete the proof of Theorem~\ref{main_thm1}, including convergence of the velocity. 

Sections 10 and 11 verify the hypotheses of the main theorem for the double sine-Gordon model and a perturbation of the $\phi^4$ model, respectively. By perturbation theory, we explicitly determine the spectrum of the linearized operators for these models and verify the repulsivity hypotheses and the Fermi golden rule constants.

\subsection{Acknowledgements}
We would like to thank Claudio Mu\~noz for enlightening discussions at the early stage of this work, and we thank Jonas L\"uhrmann for encouraging comments.

While completing this work, we learned that Michał Kowalczyk and Yvan Martel were independently working on a similar problem, see \cite{kowalczyk-martel26}.  Their work covers moving kinks with multiple internal modes. We thank Yvan Martel for sharing their preprint and communication.

\textbf{AI statement:} ChatGPT 6.0 Astra and ChatGPT 6.1 Sol were used to aid our computations and to improve writing. The authors have reviewed all mathematical content and take full responsibility for the manuscript.
\section{Modulation and orbital stability around the kink}
In this section, we develop the modulation theory for perturbations of moving kinks. We begin by reviewing the decay and spectral properties of the kink and its internal mode, then discuss the Lorentz-invariant norms (introduced in \cite[\S 2.6]{kowalczyk2021sufficient}) and the orbital stability of moving kinks. Finally, we spectrally decompose the perturbation into radiation plus internal mode components and derive their evolution equations.
\subsection{Regularity of the kink and the internal mode eigenfunction}
First, we recall the following lemma from \cite{kowalczyk2021sufficient} concerning the existence and the regularity of the kink.
\begin{lemma}[Lemma~1.1 of \cite{kowalczyk2021sufficient}]
Under assumption \eqref{eq:vacua}, there exists a solution $H(x)$ of class $\calC^\infty$ of the equation 
\begin{equation}
    \left\{ \begin{aligned}
        H''(x) &= W'(H(x)), \quad x \in \bbR,\\
        \lim_{x \rightarrow \pm \infty} H(x) &= \zeta_\pm,
    \end{aligned}\right.
\end{equation}
unique up to translation, which satisfies $H'>0$ on $\bbR$. Moreover, there exists $C>0$ such that for any $1 \leq k \leq 10$,
\begin{equation}\label{eqn:decay_of_H}
    |H(x)-\zeta_\pm| \leq C e^{\mp \omega x}, \quad |H^{(k)}(x)| \leq C e^{-\omega |x|} .
\end{equation}
\end{lemma}
Next, we prove exponential decay rates for the internal mode eigenfunction $Y$ satisfying $L_0 Y = \lambda^2 Y$. We introduce  the positive number 
\begin{equation}
    \kappa := \sqrt{\omega^2 - \lambda^2},
\end{equation}
which measures the size of the spectral gap. The following lemma is a standard Agmon pointwise exponential decay estimate. See, for example, \cite{Agmon1982Lectures}.
\begin{lemma}
For every integer $k \geq 0$, there exists a $C(k,\kappa)>0$ such that\footnote{The number $\tfrac{11}{12}$ is not special and it can be replaced by any number less than 1.}
\begin{equation}\label{eqn:decay_of_Y}
    |Y^{(k)}(x)| \leq C e^{-\frac{11}{12}\kappa|x|}.
\end{equation}
\end{lemma}
\begin{proof}
We recall that $Y'' = (W''(H)-\lambda^2)Y$. Since $W''(H)$ is bounded and $Y \in L^2$, we have $Y \in H^2$. Moreover, since $\lim_{x \rightarrow \pm \infty}W''(H(x)) =\omega_\pm^2 \geq \omega^2 > \lambda^2$, there exists a large $R>0$ such that 
\begin{equation*}
    W''(H(x)) - \lambda^2 - \left(\tfrac{11\kappa}{12}\right)^2 \geq 0, \quad \pm x \geq R.
\end{equation*}
Hence, for $x \in [R,\infty)$,
\begin{equation*}
    \big(-\px^2 + W''(H) - \lambda^2\big) e^{-\frac{11}{12}\kappa(x-R)} \geq 0.
\end{equation*}
We now argue by a maximum principle. Let $C_0 \geq |Y(R)|$ be a large constant and let $F_\pm(x):= \pm Y(x) - C_0 e^{-\frac{11}{12}\kappa(x-R)}$. Using the eigenvalue equation for $Y$ and the preceding inequality, it holds on  $[R,\infty)$ that 
\begin{equation*}
    \big(-\px^2 + W''(H)-\lambda^2\big)F_\pm \leq 0.
\end{equation*}
By construction, we have $F_\pm(R)\leq 0$ and $F_\pm(x)\rightarrow 0$ as $x \rightarrow \infty$. Now suppose that $F_\pm$ is positive at some point on $[R,\infty)$. Then, it must have a local positive maximum at some $x_* >R$, which means that $F_\pm''(x_*)\leq 0$ and $F_\pm(x_*)>0$. Consequently,
\begin{equation*}
    \big(-\px^2 + W''(H) - \lambda^2\big)F_\pm(x_*) \geq \left(\tfrac{11\kappa}{12}\right)^2 F_\pm(x_*) >0.
\end{equation*}
However, this contradicts the previous inequality. Thus, we have $F_\pm \leq 0$ on $[R,\infty)$, which proves that $|Y(x)| \leq C_0 e^{-\frac{11}{12}\kappa (x-R)}$ for all $x \geq R$. After repeating the same argument for $x \leq -R$, we may find a larger constant $C_1>0$ such that $|Y(x)|\leq C_1 e^{-\frac{11}{12}\kappa |x|}$ for all $x \in \bbR$. Using the eigenvalue equation then gives $|Y''(x)|+|Y'(x)| \lesssim e^{-\frac{11}{12}\kappa |x|}.$
For any $k \geq 3$, we have 
\begin{align*}
    Y^{(k)} = \big(W''(H)-\lambda^2 \big)Y^{(k-2)} + \sum_{j=1}^{k-2} \binom{k-2}{j} \px^j \big (W''(H)\big) Y^{(k-2-j)}.
\end{align*}
From \eqref{eqn:decay_of_H} and the fact that $H,W\in \calC^\infty$, we conclude the assertion \eqref{eqn:decay_of_Y} by induction.
\end{proof}

\subsection{Spectrum of the linearized matrix operator around the boosted kink}
In this subsection, we review the spectral theory of the linearized operator around the boosted kink, following \cite{kowalczyk2021sufficient,chen2024sinegordon}. Throughout this subsection, we fix the parameters $c \in (-1,1)$, $\gamma := (1-c^2)^{-\frac12} \in [1,\infty)$, and $y \in \bbR$. We adopt the vectorized notation from \cite{kowalczyk2021sufficient}.

For the Lorentz boosted kink, we write 
\begin{equation}
\begin{aligned}
    H_{c}(x) &:= H(\gamma x) ,\quad \bmH_c(x) := \begin{pmatrix}
        H_c(x) \\ -c H_c'(x)
    \end{pmatrix},\\
    H_{c,y}(x) &:= H_c(x-y), \quad \bmH_{c,y}(x) := \bmH_{c}(x-y),
\end{aligned}    
\end{equation}
with $H_{c,y}'(x)= \gamma H'\big(\gamma(x-y)\big)$. Using the same subscript convention, we also write 
\begin{equation}
    Y_{c,y}(x) := Y\big(\gamma(x-y)\big).
\end{equation}
Recall that the linearized matrix operator around the boosted kink is given by
\begin{equation}
\bmJ\bmL_{c,y} = \begin{pmatrix}
c \px & 1\\
-L_{c,y} & c \px 
\end{pmatrix},
\end{equation}
where\footnote{In this notation, we have $L_{0,0} \equiv L_0 = L := -\px^2 + W''(H)$ which is the linearized operator around the static kink.}
\begin{equation}
    \bmL_{c,y} := \begin{pmatrix}
L_{c,y} & - c \px \\ c \px & 1        
    \end{pmatrix}, \qquad  L_{c,y} := -\px^2 + W''(H_{c,y}).
\end{equation}
In the following proposition, we describe the spectrum of the linearized operator.

\begin{proposition}\label{prop:spectrum_for_JL}
Define the following (generalized) eigenfunctions
\begin{align}
\bmT_{c,y} &:= -\py \bmH_{c,y} = \begin{pmatrix}
    H_{c,y}'(x)\\
    -c H_{c,y}''(x)
\end{pmatrix}, \label{eqn:definition_eigfn_1}\\
\bmD_{c,y} &:= \partial_c \bmH_{c,y} = \begin{pmatrix}
c \gamma^2(x-y) H_{c,y}'(x)\\
-\gamma^2 H_{c,y}'(x) - c^2 \gamma^2 (x-y) H_{c,y}''(x)
\end{pmatrix}, \label{eqn:definition_eigfn_2}\\
\bmY_{c,y}^+ &:= e^{-ic\gamma \lambda(x-y)}\begin{pmatrix} Y_{c,y}(x)\\ i\lambda \gamma Y_{c,y}(x)  - {c} Y_{c,y}'(x)\end{pmatrix}, \label{eqn:definition_eigfn_3}\\
\bmY_{c,y}^- &:= \overline{\bmY_{c,y}^+}= e^{ic\gamma \lambda(x-y)}\begin{pmatrix}  Y_{c,y}(x)\\ -i\lambda \gamma Y_{c,y}(x)  - {c} Y_{c,y}'(x)\end{pmatrix}.\label{eqn:definition_eigfn_4}
\end{align}
The following holds under the spectral assumption \eqref{eq:spectrum} and the repulsivity assumption \eqref{eq:repulsivity}.
\begin{enumerate}
\item The essential spectrum of $\bmJ \bmL_{c,y}$ is $(-i\infty, -i \omega\gamma^{-1}] \cup [i \omega\gamma^{-1},i\infty)$ without embedded eigenvalues.
\item The discrete spectrum of $\bmJ \bmL_{c,y}$ is $\{0, -i \lambda \gamma^{-1},i \lambda \gamma^{-1} \}$, and the edges of the essential spectrum $\pm i \omega\gamma^{-1}$ are not threshold resonances. The geometric multiplicity of the zero eigenvalue is one and its algebraic multiplicity is two. Moreover, we have
    \begin{align}
        \bmJ \bmL_{c,y} \bmT_{c,y} &=\bm0, \qquad \bmJ \bmL_{c,y} \bmD_{c,y} = - \bmT_{c,y},\label{eqn:eig_eqn-zero}\\
        \bmJ \bmL_{c,y} \bmY_{c,y}^\pm &= \pm i \lambda \gamma^{-1}\bmY_{c,y}^\pm.\label{eqn:eig_eqn-Y_cy}
    \end{align}
\end{enumerate}
\end{proposition}
\begin{proof}
In what follows, we show that the spectrum of $\bmJ\bmL_{c,y}$ can be computed directly from the spectrum of the scalar operator $L_0 = -\px^2 + W''(H)$ by rescaling, and we refer the reader to \cite[Lemmas~4.3-4.5]{chen2024sinegordon} for the assertion about the essential spectrum using Weyl's criterion as well as the absence of embedded eigenvalues. 

Consider the eigenvalue equation $    \bmJ \bmL_{c,y} \bmf = \mu \bm f$ for some $\mu \in \bbC$. The components of the vector $\bmf=(f_1,f_2)$ satisfy $c f_1' + f_2 = \mu f_1$ and $-L_{c,y}f_1 + c f_2' = \mu f_2$. Substituting $f_2 = \mu f_1 - c f_1'$ into the second equation, we obtain  $L_{c,y} f_1+c^2f''_1 - 2 c \mu f_1' + \mu^2 f_1 = 0$. Since $L_{c,y}$ is a rescaling of $L_0$ and $L_0 \geq 0$ is a non-negative operator from \eqref{eq:spectrum}, we must have $\mu \in i \bbR$. We now make a change of variables by setting $g(\gamma(x-y)) = e^{c \gamma^2 \mu(x-y)}f_1(x)$ which yields 
\begin{equation*}
    L_0 g = - \gamma^2 \mu^2 g.
\end{equation*}

For the discrete spectrum of $\bmJ\bmL_{c,y}$, it follows that $-\gamma^2\mu^2 \in \{0,\lambda^2\}$ which is equivalent to $\mu \in \{0,-i\lambda \gamma^{-1},i\lambda\gamma^{-1}\}$. In particular, the change of variables allows us to determine the expressions for the eigenfunctions  \eqref{eqn:definition_eigfn_3} and \eqref{eqn:definition_eigfn_4} at the eigenvalues  $\pm i \lambda \gamma^{-1}$. At $\mu=0$, we differentiate $\bmH_{c,y}$ with respect to $y$ and $c$ to obtain \eqref{eqn:definition_eigfn_1} and \eqref{eqn:definition_eigfn_2} respectively, and compute \eqref{eqn:eig_eqn-zero} directly. Moreover, we note that $\langle \bmD_{c,y},\bmJ \bmT_{c,y}\rangle_{\bbC} = \gamma^3 \| H'\|^2 \neq 0$. If the equation $\bmJ \bmL_{c,y} \bmh = a\bmD_{c,y} + b \bmT_{c,y}$ holds for some  $\bmh \in L^2 \times L^2$ and $a,b \in \bbC$, then testing it with $\bmJ \bmT_{c,y}$ gives 
\begin{equation*}
    0= \langle \bmL_{c,y} \bmh, \bmT_{c,y}\rangle_{\bbC} = a \gamma^3 \| H' \|^2,
\end{equation*}
so $a \equiv 0$. It follows from $\bmJ \bmL_{c,y} \bmD_{c,y} = -\bmT_{c,y}$ that $\ker((\bmJ \bmL_{c,y})^2) = \mathrm{span}\{\bmT_{c,y},\bmD_{c,y}\}$ and $\ker((\bmJ \bmL_{c,y})^3)=\ker((\bmJ \bmL_{c,y})^2)$.
Hence, zero has geometric multiplicity one and algebraic multiplicity two. 

Finally, we verify that there are no threshold resonances at the edges of the essential spectrum. From the change of variables as before, it suffices to show that the equation $L_0 g = \omega^2 g$ has no bounded solution. Suppose for contradiction that there is such a function $g$. Then, using Darboux transform $U_1U_0 L_0 = L_2 U_1 U_0$, we have $(L_2 - \omega^2)h=0$ where $h := U_1U_0 g$. Under the repulsivity assumption \eqref{eq:repulsivity}, we have $P_2 \geq \omega^2$ and $P_2-\omega^2 \not\equiv 0$. Now multiply $(L_2-\omega^2 h)=0$ by $\chi_R^2h$ where $\chi_R(x) = \chi(x/R)$, and $\chi(\cdot)$ is a smooth cut-off function. Integration by parts and H\"older's inequality yield
\begin{align*}
    \int \big|(\chi_R h)'\big|^2 + \int (P_2 - \omega^2)\big|\chi_R h\big|^2 = \int\big|\chi_R'\big|^2 |h|^2 \leq C R^{-1} \| h\|_{L^\infty}^2.
\end{align*}
Letting $R\rightarrow \infty$ shows that $h$ is identically constant and using $P_2 - \omega^2 \not \equiv 0$ then implies $h = U_1 U_0 g \equiv 0$. Hence, we can write $g = a Y + bH'$ for some $a,b \in \bbC$. The equation $L_0 g = \omega^2 g$ then implies that $0=(\lambda^2 - \omega^2)Y-b \omega^2H'$. Since $H'$ and $Y$ are linearly independent, we must have $a\equiv b \equiv  0$, so $g\equiv 0$. Thus, we conclude the proof of the proposition.
\end{proof}
From here onward, we normalize the eigenfunction $Y$ so that $\|Y\|_{L^2}=1$. We record the following symplectic pairings
\begin{align}
\bmomega(Y_{c,y}^+,Y_{c,y}^-) &= - 2i \lambda \| Y \|_{L^2}^2 = -2i \lambda  \label{eqn:identity_Y_JY_+},\\
\bmomega(Y_{c,y}^-,Y_{c,y}^+) &= 2i \lambda \| Y \|_{L^2}^2 = 2i \lambda \label{eqn:identity_Y_JY_-}, \\
\bmomega(\bmD_{c,y},\bmT_{c,y})&= -\bmomega(\bmT_{c,y},\bmD_{c,y})= \gamma^3 \| H'\|_{L^2}^2.
\end{align}
In view of Proposition~\ref{prop:spectrum_for_JL}, we define the following symplectic orthogonal projections 
\begin{equation}\label{eqn:definition_spectral_projectors}
    \begin{aligned}
P_c^0 \bmu &:= -\frac{\bmomega(\bmu,\bmD_{c,y}) }{\gamma^3 \| H'\|_{L^2}^2} \bmT_{c,y} + \frac{\bmomega(\bmu,\bmT_{c,y}) }{\gamma^3 \| H'\|_{L^2}^2} \bmD_{c,y},\\
P_c^\mathrm{i} \bmu &:=  -\frac{\bmomega(\bmu,\bmY_{c,y}^-)}{2i\lambda }\bmY_{c,y}^+ + \frac{\bmomega(\bmu,\bmY_{c,y}^+)}{2i\lambda }\bmY_{c,y}^-, \\
P_c^\mathrm{e} \bmu &:= \bmu - P_c^0\bmu - P_c^\mathrm{i}\bmu.
\end{aligned}
\end{equation}
These projections will be used for the spectral decomposition of the perturbation after we recall the modulation theory in Section~\ref{Sec:modulation}.

\subsection{Lorentz invariant norm} In this subsection, we introduce the Lorentz invariant norms used in \cite{kowalczyk2021sufficient} to measure the size of a perturbation relative to a boosted kink with velocity $c \in (-1,1)$. 
For $\bmv = (v_1,v_2) \in H^1(\bbR) \times L^2(\bbR)$ and $c \in (-1,1)$, we set
\begin{equation}\label{eqn:lorentz_norm_c}
    \| \bmv \|_{c}^2 = \gamma^{-1} \| \px v_1 \|_{L^2}^2 + \gamma \| v_1 \|_{L^2}^2 + \gamma \| v_2 + c \px v_1 \|_{L^2}^2.
\end{equation}
The connection between the Lorentz invariant norm \eqref{eqn:lorentz_norm_c} and the quadratic form induced by the flat linearized operator can be seen by the following observation
\begin{equation*}
\| \bmv \|_{c}^2=    \gamma \left \langle \begin{pmatrix}
        -\px^2 + 1 & -c \px \\ c \px & 1
    \end{pmatrix}\bmv,\bmv \right \rangle .
\end{equation*}
Moreover, for every fixed $c \in (-1,1)$, the Lorentz invariant norm is equivalent to the energy norm:
\begin{equation}
\frac{1}{\sqrt{2\gamma}}\| \bmv \|_c \leq     \| \bmv \|_{H^1 \times L^2} \leq \sqrt{2\gamma} \| \bmv \|_{c}.
\end{equation}
Consequently, we have $\|v_1\|_{L^\infty} \leq \| \bmv \|_c$. For any $R>0$, we define the local version of the Lorentz invariant norm
\begin{equation}\label{eqn:lorentz_norm_c_R}
    \| \bmv \|_{c,R}^2 = \gamma^{-1} \| \px v_1 \|_{L^2(|x|<R)}^2 + \gamma \| v_1 \|_{L^2(|x|<R)}^2 + \gamma \| v_2 + c \px v_1 \|_{L^2(|x|<R)}^2.
\end{equation}
Finally, we recall the following stability lemma from \cite[Lemma~2.7]{kowalczyk2021sufficient} which will be used in the proof of the main theorem.
\begin{lemma}\label{lemma:closeness_stability}
    There exist constants $C,\delta_0>0$ such that the following holds. Let $c_0 \in (-1,1)$ and $y_0 \in \bbR$. Set $\gamma_0 = (1-c_0^2)^{-\frac12}$. For any $c \in (-1,1)$ and $y \in \bbR$ satisfying
    \begin{equation}\label{eqn:smallness_c_c_0}
        \gamma_0^2 |c-c_0| + \gamma_0 |y-y_0| =: \delta \leq \delta_0,
    \end{equation}
    it holds
    \begin{enumerate}
        \item Close kinks comparison estimate:
        \begin{equation}\label{eqn:close_norm_kink_estimate}
            \big\| \bmH_{c_0,y_0} - \bmH_{c,y} \big\|_{c_0} \leq C \delta.
        \end{equation}
        \item Close norms equivalence: 
        \begin{equation}\label{eqn:close_lorentz_norm_equivalence}
C^{-1}\| \bmv \|_{c_0} \leq \| \bmv \|_{c} \leq C \| \bmv \|_{c_0}, \quad \forall \bmv \in H^1 \times L^2.
        \end{equation}
    \end{enumerate}
\end{lemma}
\begin{proof}
See Lemma~2.7 of \cite{kowalczyk2021sufficient} for details of the first estimate \eqref{eqn:close_norm_kink_estimate}. To prove \eqref{eqn:close_lorentz_norm_equivalence}, we first note that
\begin{align*}
    \| \bmv \|_c^2 - \| \bmv \|_{c_0}^2 = (\gamma \gamma_0^{-1} -1) \| \bmv \|_{c_0}^2 + 2 \gamma(c-c_0)\langle v_1',v_2\rangle.
\end{align*}
Under the smallness condition \eqref{eqn:smallness_c_c_0}, we have 
\begin{equation*}
    \frac12 \gamma_0 < \gamma \leq 2 \gamma_0, \quad |\gamma \gamma_0^{-1} -1 | \leq C \delta, \quad |\langle v_1',v_2\rangle| \leq \|\bmv\|_{H^1 \times L^2}^2 \leq 2 \gamma_0 \| \bmv\|_{c_0}^2,
\end{equation*}
which implies that $\big|\| \bmv \|_c^2 - \| \bmv \|_{c_0}^2\big|\leq C \delta \|\bmv \|_{c_0}^2$. Hence, for $\delta_0$ small enough, we obtain \eqref{eqn:close_lorentz_norm_equivalence}.
\end{proof}
\subsection{Orbital stability, modulation, and spectral decomposition}\label{Sec:modulation}
In this subsection, we recall orbital stability results and modulation theory for the moving kink, first established in \cite{henry1982stability}. We quote \cite[Theorem 1]{kowalczyk2021sufficient} below for the orbital stability of kinks.

\begin{theorem}[Theorem~1 of \cite{kowalczyk2021sufficient}]\label{theorem:orbital_stability}
There exist $\delta_*>0$ and $C_*>0$ such that for any $c_0 \in (-1,1)$ and for any $\bmphi_0 \in \bmE_{\bmH}$ with 
\begin{equation}
\delta := \left \| \bmphi_0 - \bmH_{c_0} \right\|_{c_0} \leq \delta_*,
\end{equation}
there exists a unique global solution $\bmphi \in \calC(\bbR,\bmE_{\bmH})$ of \eqref{eqn:scalar_field_equation} with $\bmphi(0) = \bmphi_0$. Moreover, for all $t \in \bbR$
\begin{equation}
    \inf_{y \in \bbR} \left \| \bmphi(t,\cdot + y) - \bmH_{c_0} \right\|_{c_0} \leq C_* \delta. 
\end{equation}
\end{theorem}
The orbital stability of moving kinks follows from standard local well-posedness theory, conservation laws, and modulation around the moving kink. We refer to \S3.2 of \cite{kowalczyk2021sufficient} for more details on the proof of Theorem~\ref{theorem:orbital_stability}. In the following lemma, we recall the modulation around the moving kink.
\begin{lemma}\label{lemma:modulation_lemma}
There exist $C,\delta_1>0$ such that for any $c_0 \in (-1,1)$, $y_0 \in \bbR$, $\delta \in (0,\delta_1)$, if $\bmphi = (\phi_1,\phi_2)$ is a solution of \eqref{eqn:scalar_field_equation} on $[T_1,T_2]$ satisfying
\begin{equation}
    \sup_{t \in [T_1,T_2]}  \inf_{y_0\in \bbR} \| \bmphi(t) - \bmH_{c_0,y_0}\|_{c_0}  \leq \delta, 
\end{equation}
then there exist unique functions $c:[T_1,T_2]\rightarrow (-1,1)$ and $y:[T_1,T_2] \rightarrow \bbR$ of class $\calC^1$, such that decomposing 
\begin{equation}\label{eqn:scalar_field_decomposition}
    \bmphi(t,x') = \bmH_{c(t),y(t)}(x') + \bmu(t,x'),
\end{equation}
the following properties hold, for all $t \in [T_1,T_2]$,
\begin{enumerate}
    \item Smallness:
    \begin{equation}\label{eqn:mod_stability}
        \gamma_0^2 |c(t) - c_0| + \| \bmu(t)\|_{c_0} \leq C \delta , 
    \end{equation}
\item Symplectic orthogonality
\begin{equation}\label{eqn:orthogonality_bmu}
\bmomega\big(\bmu(t),\bmT_{c(t),y(t)}\big) = \bmomega\big(\bmu(t),\bmD_{c(t),y(t)}\big) =  0,
\end{equation}
\item Perturbation equation for $\bmu = (u_1,u_2)$:
\begin{equation}
    \pt \bmu = (\bmJ \bmL_{c,y} - I c\partial_{x'}) \bmu + (\dot{y}-c)\bmT_{c,y} - \dot{c} \bmD_{c,y} +\bmR,
\end{equation}
where 
\begin{equation}
    \bmR = - \begin{pmatrix}
        0 \\ W'(H_{c,y}+u_1) - W'(H_{c,y}) - W''(H_{c,y})u_1
    \end{pmatrix}
\end{equation}
satisfies $|\bmR| \leq C u_1^2$.
\item Modulation equations for $(c,y)$:
\begin{equation}
\begin{bmatrix}
    \gamma^3 \big\|H' \big\|_{L^2}^2 + \bmomega(\bmu,\bmA_{c,y}) & - \bmomega( \bmu, \bmB_{c,y}) \\ \bmomega(\bmu,\bmT_{c,y}') & \gamma^3 \big\|H' \big\|_{L^2}^2 - \bmomega( \bmu, \bmA_{c,y})  
\end{bmatrix}   \begin{bmatrix}
        \dot{y}-c\\ \dot{c} 
    \end{bmatrix} = \begin{bmatrix}
        \bmomega( \bmR, \bmD_{c,y}) \\ \bmomega(\bmR, \bmT_{c,y})
    \end{bmatrix},
\end{equation}
where $\bmA_{c,y} = \pc \bmT_{c,y}$ and $\bmB_{c,y} = \pc \bmD_{c,y}$.
\item Comparison estimate
\begin{equation}\label{eqn:modulation_comparison}
    |c(t) - c_0| \leq \frac{1}{4}\gamma_0^{-2}, \quad \frac{1}{2}\gamma_0 \leq \gamma(t) \leq 2 \gamma_0.
\end{equation}
\item Estimate for modulation parameters
\begin{equation}\label{eqn:mod_estimates}
    \gamma_0 |(\dot{y}-c)(t)| + \gamma_0^2 |\dot{c}(t)| \leq C \int u_1(x)^2 e^{-\frac12 \omega \gamma_0|x-y|} \ud x.
\end{equation}
\end{enumerate}
\end{lemma}
\begin{proof}
    See the proofs of Lemma 2.11 of \cite{kowalczyk2021sufficient}, and Proposition~8.1 of \cite{chen2024sinegordon}.
\end{proof}

Combining Theorem~\ref{theorem:orbital_stability} with Lemma~\ref{lemma:modulation_lemma} yields a global modulation decomposition of the scalar field into a moving kink $\bmH_{c(t),y(t)}$ and the perturbation $\bmu(t)$. More precisely, for all times $t \in \bbR$, we have 
\begin{equation*}
    \bmphi(t) = \bmH_{c(t),y(t)} + \bmu(t),
\end{equation*}
with $\bmu(t)$ satisfying \eqref{eqn:orthogonality_bmu} and with continuously differentiable velocity $c:\bbR \rightarrow (-1,1)$ and center parameter $y:\bbR \rightarrow \bbR$. We further spectrally decompose the perturbation $\bmu(t)$ relative to the time-dependent linearized operator $\bmJ \bmL_{c(t),y(t)}$. Notice that \eqref{eqn:orthogonality_bmu} already implies that $P_{c(t)}^{\mathrm{0}}\bmu(t) = 0$. We define the variables 
\begin{equation}
\begin{aligned}
    P_{c(t)}^{\mathrm{e}}\bmu(t) &=: \bmnu(t),\\
    P_{c(t)}^{\mathrm{i}}\bmu(t) &=: a_+(t) \bmY_{c(t),y(t)}^+ + a_-(t)\bmY_{c(t),y(t)}^-,
\end{aligned}
\end{equation}
where 
\begin{equation}
    a_\pm(t) := \mp \frac{\bmomega\left(\bmu(t),\bmY_{c(t),y(t)}^\mp\right)}{2i \lambda}.
\end{equation}
Since $\bmu(t)$ is real-valued, we have $a_+(t) = \overline{a_-(t)}$. Therefore, it suffices to use the notation $a(t) := a_+(t)$ to keep track of the internal mode amplitude throughout the paper. Hence, we have 
\begin{equation}\label{eqn:spectral_decomposition_u}
    \bmu(t) = \bmnu(t) + 2\Re\left( a(t) \bmY_{c(t),y(t)}^+\right).
\end{equation} Using the perturbation equation for $\bmu$, we derive the following system of equations for the variables 
\begin{equation}\label{eqn:system_moving_frame}
\left\{\begin{aligned}
    \pt{\bmnu}(t) &= (\bmJ \bmL_{c,y} - I c \partial_{x})\bmnu + \dot{P_c^{\mathrm{e}}}(\bmu) + P_c^\mathrm{e}\big((\dot{y}-c)\bmT_{c,y} - \dot{c} \bmD_{c,y} +\bmR \big)    \\
    \dot{a}(t) &= i \lambda \gamma^{-1} a - \frac{1}{2 i \lambda } \Big( \bmomega\big(\bmR,\bmY_{c,y}^-\big) + \bmomega\big(\bmu, \dot{\bmY}_{c,y}^-\big)\Big)     
\end{aligned}
\right.,
\end{equation}
where we use the following notation 
\begin{align}
    \dot{\bmY}_{c,y}^- &:= \dot{c} \pc\bmY_{c,y}^- + (\dot{y}-c) \py \bmY_{c,y}^-,\\
    \dot{P_c^{\mathrm{e}}} &:= \dot{c} \pc P_c^{\mathrm{e}} + (\dot{y}-c) \py P_c^{\mathrm{e}}.
\end{align}
From \eqref{eqn:mod_estimates}, \eqref{eqn:spectral_decomposition_u}, and Lemma~\ref{lemma:closeness_stability}, we obtain the following modulation estimates 
\begin{equation}
    \gamma_0^2 |\dot{y}-c| + \gamma_0^3|\dot{c}| \lesssim \gamma_0 \int |\nu_1(t,x)|^2 e^{-\frac12 \omega \gamma_0 |x-y|}\,\ud x  + |a|^2.
\end{equation}
\section{Setting up the analysis in the fixed frame}
In this section, we express the perturbation $\bmu(t)$ in new coordinates adapted to the reference kink $\bmH=(H,0)$, which were introduced in \cite{kowalczyk2021sufficient}. This change of variables transfers the spectral decomposition \eqref{eqn:spectral_decomposition_u} to the fixed reference frame and allows us to analyze the perturbations using the natural $H^1 \times L^2$ energy norm. We then recall general virial identities and virial notations which will be used in the subsequent sections. 

\subsection{Coordinates adapted to the reference kink}
For every modulation parameters $(c(t),y(t))$, we define the following transformation 
\begin{equation}\label{eq:Ttransform}
\calT_{c,y} \bm f(t,x') := \begin{pmatrix}
    f_1 \big(t,\gamma(x'-y)\big) \\ f_2\big(t,\gamma(x'-y)\big) - c \gamma \px f_1 \big(t,\gamma(x'-y)\big)
\end{pmatrix},
\end{equation}
and its inverse by
\begin{equation}
    \calT_{c,y}^{-1}\bmg(t,x) := \begin{pmatrix}
        g_1\big(t,\frac{x}{\gamma}+y\big) \\ g_2\big(t,\frac{x}{\gamma}+y\big) + c \partial_{x'} g_1\big(t,\frac{x}{\gamma}+y\big)
    \end{pmatrix},
\end{equation}
where the spatial variables $x$ and $x'$ are related by 
\begin{equation}
    x= \gamma(t)(x'-y(t)) \Longleftrightarrow x' = \frac{x}{\gamma(t)}+y(t).
\end{equation}
Observe then that $\bmH(x) = \calT_{c,y}^{-1}\bmH_{c,y}(x)$. We also pullback the radiation variable $\bmnu$ by introducing the new radiation variable $\bmz = (z_1,z_2)$ by 
\begin{equation}
    \bmz(t,x) := \calT_{c,y}^{-1}\bmnu(t,x) \Longleftrightarrow 
    \begin{cases}
        \nu_1(t,x') = z_1\big(t,\gamma(x'-y)\big),\\
        \nu_2(t,x') = z_2\big(t,\gamma(x'-y)\big) - c \gamma \px z_1 \big(t,\gamma(x'-y)\big).
    \end{cases}
\end{equation}
Equivalently, we have\footnote{A key feature of introducing the map $\calT_{c,y}$ is that the directional derivative $(\px z_1 + c \gamma z_2)$ will replace the usual derivative $\px z_1$ in the general virial computations. See Section~\ref{sec:general_virial_computations}. 
}
\begin{equation}
    \begin{cases}
        (\px \nu_1 + c\nu_2)\big(t, \frac{x}{\gamma(t)} + y(t)\big) = \frac{1}{\gamma(t)} (\px z_1 + c \gamma z_2)(t,x),\\
        (\nu_2+c \px \nu_1) \big(t, \frac{x}{\gamma(t)} + y(t)\big) = z_2(t,x).
    \end{cases}
\end{equation}
For all times $t \in \bbR$, we have the following norm equality between the two perturbation variables
\begin{equation}\label{eqn:norm_comparisons_nu_z}
\begin{split}
    \| \bmnu(t) \|_{c(t)}^2 &= \| z_1(t) \|_{H^1}^2 + \| z_2(t) \|_{L^2}^2,\\
    \| \bmnu(t,\cdot+y(t)) \|_{c(t),\gamma(t)^{-1}R}^2 &= \| z_1(t) \|_{H^1(|x|<R)}^2 + \| z_2(t) \|_{L^2(|x|<R)}^2, \quad R>0.
\end{split}
\end{equation}

For each $(c,y)$, the map $\calT_{c,y}$ is a bounded linear symplectomorphism on $H^1 \times L^2$, equipped with the pulled-back symplectic form
\begin{equation}\label{eqn:definition_pulled-back_symplectic_form}
\bmomega_\ast(\bmf,\bmg) := (\calT_{c,y})^\ast \bmomega(\bmf,\bmg) =  \bmomega(\calT_{c,y}\bmf,\calT_{c,y}\bmg) =\int  \bmf \cdot \tilde{\bmJ}_c\bmg,
\end{equation}
where
\begin{equation}
    \tilde{\bmJ_c} := \calT_{c,y}^* \bmJ \calT_{c,y} = \frac{1}{\gamma} \begin{pmatrix}
        -2 c \gamma \px & 1 \\ -1 & 0
    \end{pmatrix}.
\end{equation}
The linearized matrix operator acting on $\bmz=(z_1,z_2)$ is then denoted by 
\begin{equation}
    \calH_c = \calT_{c,y}^{-1} \bmJ \bmL_{c,y} \calT_{c,y} = \begin{pmatrix}
        0 & 1 \\ -L_0 & 2 c \gamma \px
    \end{pmatrix}.
\end{equation}
Note that the operator $\calH_c$ explicitly depends only on $H(x)$ and $c$. The spectral theory for $\calH_c$ can be transferred directly from Proposition~\ref{prop:spectrum_for_JL} via the map $\calT_{c,y}$. Let us define the following (generalized) eigenfunctions for $\calH_c$ by
\begin{align}
\Phi_c^{(0)}(x) &:= \begin{pmatrix}
\gamma H' \\ 0 
\end{pmatrix} \equiv \calT_{c,y}^{-1} \bmT_{c,y} \label{eqn:eig_eqn_Phi_0},\\
\Phi_c^{(1)}(x) &:= \begin{pmatrix}
-c \gamma^2 x H' \\ \gamma H'
\end{pmatrix} \equiv -\calT_{c,y}^{-1} \bmD_{c,y} \label{eqn:eig_eqn_Phi_1},\\
\Phi_c^\pm(x) &:= e^{\mp ic \lambda x}\begin{pmatrix}
Y \\ \pm i \lambda \gamma^{-1} Y
\end{pmatrix} \equiv \calT_{c,y}^{-1} \bmY_{c,y}^\pm \label{eqn:eig_eqn_Phi_c_pm}.
\end{align}
It is straightforward to check from Proposition~\ref{prop:spectrum_for_JL} that 
\begin{equation}\label{eqn:eig_eqn-Phi_c}
\calH_c \Phi_c^{(0)} = \bm{0},\quad \calH_c \Phi_c^{(1)} = \Phi_c^{(0)}, \quad \calH_c \Phi_c^\pm = \pm i \lambda \gamma^{-1} \Phi_c^\pm.
\end{equation}
Next, we define the symplectic orthogonal projections relative to $\calH_c$ by the following pulled-back relations
\begin{equation*}
\tPcO:= \calT_{c,y}^{-1} P_c^{\mathrm{0}} \calT_{c,y}, \quad \tPci:= \calT_{c,y}^{-1} P_c^{\mathrm{i}} \calT_{c,y}, \quad \tPce:= \calT_{c,y}^{-1} P_c^{\mathrm{e}} \calT_{c,y}.
\end{equation*}
From the spectral decomposition \eqref{eqn:spectral_decomposition_u}, we define the internal mode component in the transformed coordinates by
\begin{equation}
    \bmm(t,x)= \begin{pmatrix}
        m_1(t,x)\\m_2(t,x)
    \end{pmatrix} := \tPci \calT_{c,y}^{-1}\bmu = a(t) \Phi_c^{+}(x)+\bara(t) \Phi_c^{-}(x)\label{eqn:definition_m}.
\end{equation}
In the transformed coordinates, the pulled-back scalar field decomposition \eqref{eqn:scalar_field_decomposition} now reads
\begin{equation}\label{eqn:bmpsi_definition}
\bmpsi(t,x)  := \calT_{c(t),y(t)}^{-1}\bmphi(t,x) = \bmH(x) + \bmz(t,x) + \bmm(t,x).
\end{equation}

\subsection{Evolution equations for the radiation and internal mode components}
As in \cite{kowalczyk2021sufficient}, we introduce the modulation operators
\begin{equation}
    \bmOmega(\bmf) := \begin{pmatrix}
        \Omega_1(\bmf)\\
        \Omega_2(\bmf)
    \end{pmatrix}, \qquad \left\{ \begin{aligned}
        &\Omega_1(\bmf) := (\dot{y}-c) \gamma \px f_1 - \dot{c} c \gamma^2 \Lambda f_1 ,\\ 
        &\Omega_2(\bmf) := (\dot{y}-c)\gamma \px f_2 - \dot{c} c \gamma^2 \Lambda f_2 + \dot{c} \gamma \px f_1,
    \end{aligned}\right.
\end{equation}
where $\Lambda$ is the dilation operator
\begin{equation}
    \Lambda := x \px.
\end{equation}
In particular, we note that 
\begin{equation}\label{eqn:bmOmega(bmH)_identity}
    \bmOmega(\bmH) = (\dot{y}-c) \Phi_c^{(0)}+ \dot{c} \Phi_c^{(1)}.
\end{equation}
To derive the evolution equations for $\bmz = (z_1,z_2)$, we first express the scalar field equations in the coordinates adapted to the moving kink in the following lemma. 
\begin{lemma}
The scalar field $\bmphi = (\phi_1,\phi_2)$ solves the system
\begin{equation}
    \begin{cases}
        \pt \phi_1 = \phi_2, \\
        \pt \phi_2 = \partial_{x'}^2 \phi_1 - W'(\phi_1),
    \end{cases}
\end{equation}
if and only if the transformed field $\bmpsi(t,x) := \calT_{c(t),y(t)}^{-1}\bmphi(t,x)=\big(\psi_1(t,x),\psi_2(t,x)\big)$ solves
\begin{equation}\label{eqn:equation_for_psi12}
    \begin{cases}
        \pt \psi_1 = \psi_2 +  \Omega_1(\bmpsi),\\
        \pt \psi_2 = \px^2 \psi_1 - W'(\psi_1) + 2 c \gamma \px \psi_2 + \Omega_2(\bmpsi).
    \end{cases}
\end{equation}
\end{lemma}
\begin{proof}
The equivalence is a direct computation using the product and chain rules.
\end{proof}
Following \cite{martel2024asymptotic}, we introduce the notation
\begin{equation}
    \bmv^\perp := (I - \tPci)\bmv = \bmv - \frac{\langle \bmv,\Psi_c^+\rangle_{\bbC}}{2i\lambda  }\Phi_c^+ + \frac{\langle \bmv,\Psi_c^-\rangle_{\bbC}}{2i\lambda  }\Phi_c^-,
\end{equation}
where we have expressed the symplectic pairings in terms of the complex Hermitian product with the adjoint eigenfunctions, which are defined by 
\begin{equation}
\Psi_c^\pm(x)  :=  - \tilde{\bmJ}_c \Phi_c^\pm =  e^{\mp i c \lambda x}\begin{pmatrix}
2 c  Y'(x) \mp i \lambda (1+c^2)Y(x) \\ \gamma^{-1}Y(x)
\end{pmatrix}\label{eqn:eig_eqn-Psi_c}.
\end{equation}
We can now state the evolution equations for the radiation and internal mode components.
\begin{proposition}\label{prop:equations_for_bmz}
We have 
\begin{align}
\pt\bmz &= \calH_c \bmz +   \bmOmega(\bmH)+ \bmOmega(\bmz)^\perp + \bmOmega(\bmm)^\perp +\calN^\perp- \calR^\perp - \calS\label{eqn:equation_for_z_original},\\
\dot{a} &=  i\lambda \gamma^{-1}a + \frac{1}{2 i \lambda}\Big( \langle  \bmz, \dot{c} \partial_c \Psi_{c}^+\rangle_{\bbC}  + \langle \calN + \bmOmega(\bmpsi)  - \calR, \Psi_c^+\rangle_{\bbC} \Big)\label{eqn:equation_a_+},
\end{align}
where
\begin{align}
\calN &= (0,N_2), \quad N_2 := -[W'(H+z_1+m_1) - W'(H) - W''(H)(z_1+m_1)],\\
\calR&= \dot{c}\left( a \partial_c  \Phi_c^+ + \bara \partial_c  \Phi_c^-\right) \label{eqn:definition_calR},\\
\calS &= \dot{c} \left(\frac{\langle  \bmz,  \partial_c \Psi_{c}^+\rangle_{\bbC}}{2i\lambda}\Phi_c^+ - \frac{\langle  \bmz,  \partial_c \Psi_{c}^-\rangle_{\bbC}}{2 i \lambda}\Phi_c^-\right) \label{eqn:definition_calS}.
\end{align}

\end{proposition}
\begin{proof}
From \eqref{eqn:eig_eqn-Phi_c} and \eqref{eqn:definition_m}, we have
\begin{align*}
    \pt \bmm = \dot{a}\Phi_c^+ + \dot{\bara}\Phi_c^- + a \dot{c} \pc\Phi_c^+ + \bara \dot{c} \pc \Phi_c^- = \dot{a}\Phi_c^+ + \dot{\bara}\Phi_c^- + \calR,
\end{align*}
and
\begin{align*}
    \calH_c \bmm = i \gamma^{-1} \lambda a \Phi_c^+ -i \gamma^{-1} \lambda \bara \Phi_c^-.
\end{align*}
Substituting the transformed field $\bmpsi = \bmH + \bmz + \bmm$ into \eqref{eqn:equation_for_psi12}, and using $H''=W'(H)$, we obtain 
\begin{equation}\label{eqn:proof_eqn_bmz}
\pt\bmz + \dot{a}\Phi_c^+ + \dot{\bara}\Phi_c^- + \calR = \calH_c \bmz+ \bmOmega(\bmH) + \bmOmega(\bmz) + \bmOmega(\bmm) + i  \lambda \gamma^{-1}  a \Phi_c^+ - i  \lambda \gamma^{-1} \bara \Phi_c^- + \calN.
\end{equation}
The symplectic orthogonality conditions imply that $\bmomega_\ast(\bmz,\Phi_c^\pm) = -\langle \bmz,\Psi_c^\mp\rangle_\bbC = 0$, and so time differentiating $\langle \bmz,\Psi_c^\pm\rangle_\bbC = 0$ gives 
\begin{align*}
    \langle \pt \bmz ,\Psi_c^\pm\rangle_{\bbC} = - \langle \bmz,\dot{c} \pc \Psi_c^\pm\rangle_{\bbC}.
\end{align*}
Moreover, since $\calH_c^*\Psi_c^+ = -i\lambda \gamma^{-1}\Psi_c^+$, we have $\langle \calH_c\bmz,\Psi_c^+\rangle_{\bbC}=0$. Testing \eqref{eqn:proof_eqn_bmz} against $\Psi_c^+$ and using \eqref{eqn:identity_Y_JY_-}, \eqref{eqn:identity_Y_JY_+}, and \eqref{eqn:definition_pulled-back_symplectic_form}, we obtain the equation \eqref{eqn:equation_a_+} for $a(t)$. Next, using the equation for $a$ \eqref{eqn:equation_a_+} and its conjugate, we find the following projection identity
\begin{equation}
    \tPci\big(\bmOmega(\bmz)+\bmOmega(\bmm) + \calN - \calR \big) + \calS = \big(\dot{a}-i \lambda \gamma^{-1} a\big)\Phi_c^+ + \big(\dot{\bara}+i  \lambda \gamma^{-1} \bara \big)\Phi_c^-.
\end{equation}
From \eqref{eqn:bmOmega(bmH)_identity} and the fact that $\tPci \Phi_c^{(0)}=\tPci \Phi_c^{(1)}= 0$, we have 
\begin{equation}
    \bmOmega(\bmH)^\perp = \bmOmega(\bmH).
\end{equation}
Thus, we obtain the evolution equation \eqref{eqn:equation_for_z_original} from \eqref{eqn:proof_eqn_bmz} and the preceding identities.
\end{proof}

\subsection{General virial computations}\label{sec:general_virial_computations}
In this section, we briefly describe the virial identity for the moving kink. Following the notation introduced by Martel \cite{martel2024asymptotic},  we define the following  operator
\begin{equation}
    \Theta[\varphi] = 2 \varphi \px + \varphi',
\end{equation}
which satisfies  the integration by parts property
\begin{equation}
    \int (\Theta[\varphi]f)g = - \int f (\Theta[\varphi]g), \qquad \int (\Theta[\varphi]f)f = 0.
\end{equation}
If $\varphi$ is a bounded cutoff approximation to the multiplication operator $Xf = xf$, then $\Theta[\varphi]$ is a localized version of the dilation operator $\Theta[x]= X\px + \px X\equiv 2 x \px + 1$. It is illustrative to first consider the following static Klein--Gordon system
\begin{equation*}
\left\{    \begin{aligned}
        \dot{g}_1 &= g_2,\\
        \dot{g}_2 &= \px^2 g_1 - P g_1,
    \end{aligned}\right.
\end{equation*}
and the virial functional 
\begin{align*}
    \calA = \langle \Theta[\varphi]g_1,g_2\rangle,
\end{align*}
where $\varphi:\bbR \rightarrow \bbR$ is a bounded increasing function. Setting $\zeta = \sqrt{\varphi'}$ and $h_1 = \zeta g_1$, we obtain from integration by parts that
\begin{equation}
    -\dot{\calA} = 2 \int (\px h_1)^2 + \int \ln(\zeta)'' h_1^2 - \int h_1^2 P' \frac{\varphi}{\varphi'}.
\end{equation}
Thus, if $\calA$ is well-defined and bounded, and if $x P' \leq 0$, then $\int (\px h_1)^2$ will be time-integrable. 

In the moving kink setting, the transport term  $2c \gamma \px z_2$ in  \eqref{eqn:equation_for_z_original} requires counter-terms for the standard virial $\langle\Theta[\varphi]z_1,z_2\rangle$. The following general virial identity, established in \cite[Lemma~4.2]{kowalczyk2021sufficient}, accounts for the transport term through the added virial functionals $\calA_2$ and $\calA_3$ as a modification. 
\begin{lemma}[Lemma~4.2 of \cite{kowalczyk2021sufficient}]\label{lemma:general_virial}

Consider the following system for $\bma = (a_1,a_2)$,
\begin{equation}
\left\{\begin{aligned}
\dot{a}_1 &= a_2 + \Omega_1(\bma) + \alpha_1 \\
\dot{a}_2 &= \px^2 a_1 - f(x,a_1) + 2 c \gamma \px a_2 + \Omega_2(\bma) + \alpha_2
\end{aligned}\right.,
\end{equation}
where $f(x,a_1)$ satisfies $f(x,0) = 0$, and $F(x,a_1) = \int_0^{a_1} f(x,b)\ud b. $ For $\varphi \in \calC^3$, let 
\begin{equation}\label{eqn:general_virial_functional}
    \calA = \calA_1 + c \gamma (\calA_2 + \calA_3),
\end{equation}
where
\begin{align}
\calA_1 = \calA_1[\varphi](a_1,a_2) &:= \langle \Theta[\varphi]a_1,a_2\rangle = \int (2 \varphi \px a_1 + \varphi' a_1) a_2,\\
\calA_2 = \calA_2[\varphi](a_1,a_2) &:= -2 \int (\px a_1)^2 \varphi, \\
\calA_3 = \calA_3[\varphi](a_1,a_2) &:= \int [a_2^2 + (\px a_1)^2 + 2F(a_1)]\varphi .
\end{align}
Then\footnote{The notation $\px F(x,a)$ means differentiation with respect to the first argument only throughout.} 
\begin{enumerate}
\item It holds
\begin{equation}
\begin{aligned}
\dot{\mathcal{A}}&= -2 \int \left( \partial_x a_1 + c\gamma a_2 \right)^2 \varphi' 
- 2 \int \left( \partial_x a_1 + c\gamma a_2 \right) a_1 \varphi'' 
- \frac{1}{2} \int a_1^2 \varphi''' \\
&\quad + \int \Big( 2F(a_1) - a_1 f(x, a_1) \Big) \varphi' 
+ 2 \int \varphi (\partial_x F)(a_1) \\
&\quad + (\dot{y} - c)\gamma \mathcal{I} + \dot{c}\gamma \mathcal{J} + \calM+ \dot{c}\gamma^3 (\mathcal{A}_2 + \mathcal{A}_3),
\end{aligned}
\end{equation}
where
\begin{equation*}
\begin{aligned}
\mathcal{I} &= \int \left[ c\gamma (\partial_x a_1)^2 - c\gamma a_2^2 - 2 (\partial_x a_1) a_2 \right] \varphi' 
- \int a_1 a_2 \varphi'' - 2c\gamma \int \left[ F(a_1) \varphi' + (\partial_x F)(a_1) \varphi \right],
\end{aligned}
\end{equation*}
\begin{equation*}
\begin{aligned}
\mathcal{J} &= 2 \gamma^2 \int (\partial_x a_1)^2 \varphi 
+ c\gamma \left[ -c\gamma (\partial_x a_1)^2 + 2(\partial_x a_1) a_2 + c\gamma a_2^2 \right] (x\varphi)' \\
&\quad - \frac{1}{2} \int a_1^2 \varphi'' 
+ c\gamma \int a_1 a_2 (x\varphi')' 
+ 2c^2\gamma^2 \int \left[ (x\varphi)' F(a_1) + x\varphi (\partial_x F)(a_1) \right],
\end{aligned}
\end{equation*}
and
\begin{equation*}
\begin{aligned}
\calM &= \int (2\varphi \partial_x \alpha_1 + \varphi' \alpha_1) a_2 
+ \int (2\varphi \partial_x a_1 + \varphi' a_1) \alpha_2 \\
&\quad + 2c\gamma \int \left[ -\varphi (\partial_x \alpha_1)(\partial_x a_1) + \varphi \alpha_1 f(a_1) \right] 
+ 2c\gamma \int \varphi \alpha_2 a_2.
\end{aligned}
\end{equation*}

\item Assuming that $\varphi' > 0$ and setting
\begin{equation}
\zeta = \sqrt{\varphi'}, \quad b_1 = a_1 \zeta, \quad b_2 = a_2 \zeta,
\end{equation}
it holds
\begin{equation}\label{eqn:general_virial}
\begin{aligned}
\dot{\mathcal{A}} &= -2 \int \left( \partial_x b_1 + c\gamma b_2 \right)^2 
- \int b_1^2 \left( \frac{\zeta''}{\zeta} - \frac{(\zeta')^2}{\zeta^2} \right) \\
&\quad + \int \Big( 2F(a_1) - a_1 f(x, a_1) \Big) \varphi' 
+ 2 \int \varphi (\partial_x F)(a_1) \\
&\quad + (\dot{y} - c)\gamma \mathcal{I} 
+ \dot{c}\gamma \mathcal{J} 
+ \mathcal{M} 
+ \dot{c}\gamma^3 (\mathcal{A}_2 + \mathcal{A}_3).
\end{aligned}
\end{equation}
\end{enumerate}
\end{lemma}
\begin{remark}\label{remark:virial_modification}
In Section~\ref{sec:modified_virial}, we further modify \eqref{eqn:general_virial_functional} by introducing suitable time-dependent coefficients for each $\calA_j$. See Lemma~\ref{lemma:coefficients_dj}. This modification cancels the term $2 \gamma^2 \int (\partial_x a_1)^2 \varphi $ in $\calJ$, whose weight lacks the required spatial localization. Moreover, the modification is required to improve the decay of the remainder term $\dot{c}\gamma^3(\calA_2 + \calA_3)$ because $|\dot{c}|$ decays only like $|a|^2$, but a remainder of order $|a|^4$ is required to close the virial arguments. Our construction of the modified virial functional in the proof of Lemma~\ref{lemma:HE_key_estimate1} is inspired by \cite[Lemma~18]{martel2024asymptotic}.    
\end{remark}

\subsection{Notation for virial arguments}
In this subsection we introduce the various cut-off functions and weight functions used in standard virial computations. We first recall the Darboux-transformed potential $P_2$ defined in \eqref{eq:darboux}. By the repulsivity assumption \eqref{eq:repulsivity}, since $P_2' \not\equiv 0$, there exists $C_1>0$ and $x_1,x_2 \in \bbR$, with $x_1 < x_2$ such that 
\[|P_2'(x)| \geq \frac{1}{C_1} \quad \forall x \in [x_1,x_2].\]
Let $\barx = \frac{1}{2}(x_1 + x_2)$ and let $\eta:\bbR \rightarrow [-1,1]$ be a smooth, non-decreasing function satisfying 
\begin{equation}
    \eta(x) = \begin{cases}
        -1 &\quad x < x_1\\
        1 &\quad x > x_2
    \end{cases}, \quad \text{and} \quad \eta(\barx) = 0.
\end{equation}
For any $K\geq 1$, we define 
\begin{equation}
\zeta_K(x) = \exp\left(-\frac{\kappa}{K}(x-\barx)\eta(x)\right),
\end{equation}
where we recall that $\kappa := \sqrt{\omega^2 - \lambda^2}>0$. Note that $\zeta_K(x) = \exp(-\frac{\kappa}{K}|x-\barx|)$ for any $x \in \bbR \setminus [x_1,x_2]$.
For any $A\geq10$ and for any $B\geq 10$, we define 
\begin{equation}
    \varphi_A(x) = \int_0^x \zeta_A^2(y) \ud y,
\end{equation}
and
\begin{equation}
    \varphi_B(x) = \int_{x_0}^x \zeta_B^2(y) \ud y,
\end{equation}
where $x_0$ is introduced in \eqref{eq:repulsivity}. Since $|P_2'|>0$ on $[x_1,x_2]$, we have $x_0 \notin [x_1,x_2]$. Thus, there exists $C_2>0$ such that for all $B\geq 1$ and for all $x \in [x_1,x_2]$, we have 
\begin{align*}
    |\varphi_B(x)| \geq \left| \int_{x_0}^x \zeta_{B=1}^2(s) \,\ud s\right| \geq \frac{1}{C_2}.
\end{align*}
In particular, we have for any $B\geq 1$ that
\begin{equation}
    \varphi_B P_2' \leq -\frac{1}{C_1C_2}\mathbbm{1}_{[x_1,x_2]}, \quad \text{on $\bbR$}.
\end{equation}

Next, we consider a smooth, even bump function $\chi:\bbR \rightarrow [0,1]$ satisfying
\begin{equation}
    \chi(x) = \begin{cases}
        1, \quad x \in [-1,1]\\
        0,\quad x \in (-\infty,-2]\cup[2,\infty)    \end{cases}, \quad \chi'(x)\leq 0,\quad x \in [0,\infty),
\end{equation}
and set $\chi_K(x) = \chi(x/K)$, for $K\geq 1$. We define the following compactly supported multiplier
\begin{equation}
    \psi_{A,B}(x) = \chi_A^2(x) \varphi_B(x),
\end{equation}
and use the following fixed weight function
\begin{equation}\label{eqn:definition_rho}
\rho(x) = \sech\left(\frac{\kappa }{40}x\right).
\end{equation}
In the next lemma, we collect from \cite[\S4.4]{kowalczyk2021sufficient} the basic identities and estimates related to the above weight functions.

\begin{lemma} 
\begin{enumerate}
    \item Algebraic identities for $\zeta_K$:
    \begin{align*}
\frac{\zeta_K'}{\zeta_K} &= -\frac{\kappa}{K} \left( \eta(x) + (x - \bar{x}) \eta'(x) \right), \\
\frac{\zeta_K''}{\zeta_K} &= \left( \frac{\zeta_K'}{\zeta_K} \right)^2 - \frac{\kappa}{K} \left( 2\eta'(x) + (x - \bar{x}) \eta''(x) \right),\\
\ln(\zeta_K)''&=\frac{\zeta_K''}{\zeta_K} - \left( \frac{\zeta_K'}{\zeta_K} \right)^2 = -\frac{\kappa}{K} \left( 2\eta'(x) + (x - \bar{x}) \eta''(x) \right).
\end{align*}
\item 
There exists $C\geq 1$ such that for any $K\geq 1$, on $\bbR$
\begin{equation}\label{eqn:bound_on_zeta_K}
\begin{aligned}
\frac{1}{C} \operatorname{sech}\left(\frac{\kappa}{K} x\right) &\leq \zeta_K(x) \le C \operatorname{sech}\left(\frac{\kappa}{K} x\right), \\
\left|\zeta_K'(x)\right| &\leq C K^{-1} \operatorname{sech}\left(\frac{\kappa}{K} x\right), \\
\left|\zeta_K''(x)\right| &\leq C K^{-2} \operatorname{sech}\left(\frac{\kappa}{K} x\right) + C K^{-1} \mathbbm{1}_{[x_1, x_2]}(x), \\
\left|x \zeta_K(x)\right| &\leq C K, \quad \left|x \zeta_K'(x)\right| \le C,\\
|\ln(\zeta_K)''|=\left |\frac{\zeta_K''}{\zeta_K} - \frac{(\zeta_K')^2}{\zeta_K^2} \right | &\leq \frac{C}{K} \mathbbm{1}_{[x_1,x_2]}(x) \leq \frac{C}{K}\rho(x)^8.
\end{aligned}
\end{equation}

\item For any $A\geq 10$, \begin{equation}\label{eqn:bound_on_phi_A}
\begin{cases}
0 < \varphi_A'(x) \le 1, \quad |\varphi_A(x)| \le |x|, \quad |\varphi_A(x)| \le CA, \\
|\varphi_A(x) H^{(k)}(x)| + |\varphi_A(x) \Lambda^k H(x)| \leq C \rho^{40}, \quad \text{for } 1 \leq k \leq 7,\\
|\varphi_A(x) Y^{(k)}(x)| + |\varphi_A(x) \Lambda^k Y(x)| \leq C \rho^{30}, \quad \text{for } 0 \leq k \leq 7.
\end{cases}
\end{equation}

\item For any $B \geq 10$ \begin{equation}\label{eqn:bound_on_phi_B}
\begin{cases}
|\varphi_B(x)| \leq CB, \quad |\psi_{A,B}(x)| \leq CB \mathbbm{1}_{|x| < 2A}, \\
|\psi_{A,B}'(x)| \leq C \frac{B}{A} \mathbbm{1}_{A < |x| < 2A}(x) + C \mathrm{sech}^2\left( \frac{\kappa}{B} x \right) \mathbbm{1}_{|x| < 2A}(x).
\end{cases} 
\end{equation}

\item For $K \geq 10$,
\begin{equation}\label{eqn:comparison_rho_zeta_K}
    \rho^4 \leq C \zeta_K.
\end{equation}
For $K\geq 40$,
\begin{equation}\label{eqn:comparison_rho_zeta_K_40}
    \rho \leq C \zeta_K.
\end{equation}

\item For $\theta > \frac12$, $A\geq 10$,
\begin{equation}\label{eqn:bound_varphi'}
\varphi_A' =    \zeta_A^2 \leq \zeta_{\theta A}, \quad |x \varphi_A'| = |x|\zeta_A^2 \leq C_\theta  A \zeta_{\theta A},
\end{equation}
\begin{equation}\label{eqn:bound_varphi''}
    |\varphi_A''| = \left| \frac{2\kappa}{A} [\eta(x) + (x-\barx)\eta'(x)]\right|\zeta_A^2 \leq  A^{-1} \zeta_{\theta A}, \qquad |x\varphi_A''| \leq C_\theta  \zeta_{\theta A}.
\end{equation}

\end{enumerate}
\end{lemma}
\begin{proof}
Items (1)--(4) can be deduced directly from the definitions introduced. Moreover, see Lemma~1 of \cite{kowalczyk2022jems} for final estimate in \eqref{eqn:bound_on_zeta_K}, $|\ln(\zeta_K)''|\leq \frac{C}{K} \mathbbm{1}_{[x_1,x_2]}(x) \leq \frac{C}{K}\rho(x)^8$. The bounds \eqref{eqn:comparison_rho_zeta_K} and \eqref{eqn:comparison_rho_zeta_K_40} can be deduced from the first estimate in \eqref{eqn:bound_on_zeta_K} and \eqref{eqn:definition_rho}. For the final item, we note 
\begin{equation*}
    \frac{\zeta_A^2}{\zeta_{\theta A}} = \exp\left(-\frac{\kappa}{A}\Big( 2 - \frac{1}{\theta}\Big) (x-\barx) \eta(x)\right)\leq 1,
\end{equation*}
whenever $\theta > \frac12$. This proves $\varphi_A' \leq \zeta_{\theta A}$. Moreover, for any $x \in \bbR\setminus[x_1,x_2]$, it holds that $(x-\barx) \eta(x) = |x-\barx|$ which implies that $|x| \lesssim 1+(x-\barx)\eta(x)$ for all $x \in \bbR$. Hence,
\begin{equation*}
    |x| \frac{\zeta_A^2}{\zeta_{\theta A}} \lesssim \big(1 + (x-\barx) \eta(x)\big) \exp\left(-\frac{\kappa}{A}\Big( 2 - \frac{1}{\theta}\Big) (x-\barx) \eta(x)\right) \leq 1 + \frac{A}{\kappa(2- \frac{1}{\theta})}\leq C_\theta A.
\end{equation*}
This proves \eqref{eqn:bound_varphi'}. The final estimate \eqref{eqn:bound_varphi''} then follows from \eqref{eqn:bound_varphi'}.
\end{proof}

In the following sections, we use the multiplier $\varphi_A$ for the radiation component $\bmz$, and $\psi_{A,B}$ for the regularized transformed radiation component $\bmf$ at the scales $A$ and $B$ respectively. Moreover, the parameter $\delta>0$ is used to control the size of the initial perturbation, and the parameter $0<\eps\ll1$ is used in the regularization of the Darboux transformed variable $\bmg$, defined in \eqref{eqn:def_g}.

In the proof of Theorem~\ref{main_thm1}, we first fix $B$ sufficiently large, then choose $A \gg B$ and set $\eps= A^{-\frac14}$, taking $A$ large enough to satisfy all subsequent technical smallness requirements on $\eps$. Then, we take $\delta>0$ small enough depending on $A$. In the next sections, we will assume the following scaling regime
\begin{equation}\label{main_scaling}
1 \ll B  \ll A , \quad 0 < \delta \leq \frac{1}{A^4}, \quad \eps = {A^{-\frac14}}.
\end{equation}
\section{Estimate of the radiation component at large scale}
In this section, we prove the first key proposition which captures the integrated-in-time local energy norm of the radiation component.
\begin{proposition}\label{prop:high_energy_estimate}
Assuming \eqref{main_scaling}, the following bound holds for all $10 \leq K \leq \frac{A}{2}$,
\begin{equation}\label{eqn:estimate_HE}
\int_0^T \left( \| \zeta_K \px z_1 \|^2 + \| \zeta_K z_1 \|^2 + \| \zeta_Kz_2 \|^2\right)  \ud t \leq CK^2A \delta^2 + CK^2\int_0^T \left( |a|^4 + \| \rho^4 z_1 \|^2 \right) \ud t.
\end{equation}

\end{proposition}

\subsection{Preparations}
We first regroup the terms in the equations for $\bmz$ from Proposition~\ref{prop:equations_for_bmz} by isolating $W'(H+z_1)-W'(H)$ from the nonlinear term $\calN^\perp$. We write  
\begin{equation}\label{eqn:z_1_virialform}
\left\{
\begin{aligned}
\dot{z}_1 &= z_2 + \Omega_1(\bmz)  + \beta_1 \\
\dot{z}_2 &= \px^2 z_1 - f(x,z_1) 
+ 2 c \gamma \px z_2 + \Omega_2(\bmz) + \beta_2 + N_{\bmm}
\end{aligned}\right.,
\end{equation}
where we set
\begin{align}
&f(x,z_1) := W'(H+z_1)-W'(H), \qquad F(x,z_1) := W(H+z_1) - W(H) - W'(H)z_1\label{eqn:definition_fxz_1},\\
&N_{\bmm} :=  -[W'(H+z_1+m_1) - W'(H+z_1) - W''(H)m_1],\label{eqn:N_bmm}\\
&\bmbeta \equiv (\beta_1,\beta_2) := \bmOmega(\bmH)+ \big(\bmOmega(\bmz)^\perp - \bmOmega(\bmz)\big) + \bmOmega(\bmm)^\perp - \tilP_c^i(\calN) - \calR^\perp -\calS.\label{eqn:definition_beta}
\end{align}

From \eqref{eqn:mod_stability}, \eqref{eqn:spectral_decomposition_u}, and $\|\bmnu\|_c^2 = \|\bmz \|_{H^1 \times L^2}^2$, we record the following smallness estimate for all $t \geq 0$
\begin{equation}\label{eqn:stability_smallness}
    \gamma_0^2 |c(t)- c_0| + \| \bmz(t) \|_{H^1 \times L^2} + |a(t)| \leq C \delta.
\end{equation}
From \eqref{eqn:mod_estimates}, we have the following modulation estimates for all $t \geq 0$
\begin{equation}\label{eqn:estimate_modulation_crude}
    \gamma_0^2 | \dot{y}(t) - c(t)| + \gamma_0^3 |\dot{c}(t)| \leq C \big( \| \rho^4 z_1 \|^2 + |a(t)|^2\big).
\end{equation}
The following lemma records the decay rates for the eigenfunctions of the linearized operator $\calH_c$. \begin{lemma}\label{lemma:preparatory_bounds_on_eigfn}
For all integers $k \geq 0$, we have
\begin{equation}
\big|\px^k \Phi_c^{(0)}\big| + \big|\px^k \Phi_c^{(1)}\big|\lesssim_k \rho^{12}, \label{eqn:estimate_on_Phi_c_0}
\end{equation}
\begin{equation}\label{eqn:estimate_on_Phi_c_pm}
\big|\px^k \Phi_c^\pm\big| + \big|\px^k \Psi_c^\pm\big| + \big|\px^k \pc \Phi_c^\pm\big| + \big|\px^k \pc \Psi_c^\pm\big| \lesssim_k \rho^{12}.
\end{equation}
\end{lemma}
\begin{proof}
These bounds are a direct consequence of the estimates \eqref{eqn:decay_of_H} and \eqref{eqn:decay_of_Y}.
\end{proof}

We begin by estimating the localized error terms in \eqref{eqn:z_1_virialform}.
\begin{lemma}\label{lemma:bound_on_beta_terms} For $0 \leq k \leq 2$, we have 
\begin{equation}\label{eqn:bound_on_beta}
\big\| \rho^{-6} \px^k \beta_1 \big\| + \big\|  \rho^{-6} \px^k \beta_2 \big\| \leq C_k\big( |a|^2 + \big\| \rho^4 z_1 \big\|^2\big). 
\end{equation}
Moreover, it holds that 
\begin{equation}\label{eqn:bound_on_Nmm}
    \big\| \rho^{-6} N_{\bmm} \big\| \leq C \big( |a|^2 + \big\| \rho^4 z_1 \big\|^2\big).
\end{equation}
\end{lemma}
\begin{proof}
In what follows, we show that every term in \eqref{eqn:definition_beta} satisfies the same bound as  \eqref{eqn:bound_on_beta}. First, from \eqref{eqn:bmOmega(bmH)_identity}, the bounds \eqref{eqn:estimate_modulation_crude} and \eqref{eqn:estimate_on_Phi_c_0} imply the estimate \eqref{eqn:bound_on_beta} for $\bmOmega(\bmH)$. Next, we have  
\begin{align*}
(\Omega(\bmz)^\perp-\Omega(\bmz)) = -\tilP_{c}^i \Omega(\bmz) = -\frac{\langle\bmOmega(\bmz),\Psi_c^+\rangle_{\bbC}}{2 i \lambda }\Phi_c^+ + \frac{\langle \bmOmega(\bmz),\Psi_c^-\rangle_{\bbC}}{2 i \lambda}\Phi_c^-. 
\end{align*}
Hence, it suffices to bound $|\langle \bmOmega(\bmz),\Psi_c^\pm\rangle_{\bbC}|$. Integration by parts yields the following expression
\begin{align*}
\langle \bmOmega(\bmz),\Psi_c^\pm\rangle_{\bbC} &= (\dot{y}-c)\gamma \int (z_1 \px \Psi_{c,1}^\mp + z_2 \px \Psi_{c,2}^\mp) \\
&\quad- \dot{c} c \gamma^2 \int [z_1(\Psi_{c,1}^\mp + x \px \Psi_{c,1}^\mp) +z_2 (\Psi_{c,2}^\mp + x \px \Psi_{c,2}^\mp) ]+ \dot{c} \gamma \int z_1 \px \Psi_{c,2}^\mp.
\end{align*}
By the Cauchy-Schwarz inequality, \eqref{eqn:stability_smallness}, and \eqref{eqn:estimate_modulation_crude}, we have 
\begin{align*}
| \langle \bmOmega(\bmz),\Psi_c^\pm\rangle_{\bbC} | \leq C \delta(|a|^2 + \|\rho^4 z_1 \|^2),
\end{align*}
whence \eqref{eqn:bound_on_beta} is satisfied by $\bmOmega(\bmz)^\perp - \bmOmega(\bmz)$  using the bound \eqref{eqn:estimate_on_Phi_c_pm}. Similarly, by Taylor expansion, the nonlinearity satisfies $|N_2|= |\calN| \leq C(|z_1|^2 + |m_1|^2)$. Since 
\[\langle \calN,\Psi_c^\pm\rangle_{\bbC}= \int N_2 \Psi_{c,2}^\mp,\]
we deduce that the localized term $\tilP_c^\mathrm{i}(\calN)$ satisfies \eqref{eqn:bound_on_beta} using  \eqref{eqn:estimate_on_Phi_c_pm}. Finally, note that $\calR$ and $\calS$ are of the form $\dot{c}$ multiplied by a localized function. Hence, $\calR^\perp$ and $\calS^\perp$ can be straightforwardly estimated using \eqref{eqn:estimate_modulation_crude} and \eqref{eqn:estimate_on_Phi_c_pm} to prove that they satisfy \eqref{eqn:bound_on_beta}. 
Finally, by Taylor expansion, we have 
\begin{align*}
-N_{\bmm} = m_1^2 \int_0^1 (1-s)W'''(H+z_1+sm_1) \,\ud s + z_1 m_1 \int_0^1 W'''(H+sz_1) \,\ud s.
\end{align*}
Using the smoothness of $W$ and the smallness estimate \eqref{eqn:stability_smallness}, we have $|N_{\bmm}| \leq C |m_1|(|m_1| + |z_1|)$. Since $m_1 = 2\Re(a \Phi_{c,1}^+)$, we have $|m_1| \lesssim |a| \rho^{12}$ by \eqref{eqn:estimate_on_Phi_c_pm}, which implies \eqref{eqn:bound_on_Nmm}. Thus we conclude the proof.
\end{proof}

\subsection{Transfer estimates}
For $A\geq 10$, we define the localization of $\bmz$ at scale $A$ by the variables
\begin{equation}
w_1 = \zeta_A z_1, \qquad w_2 =\zeta_A z_2.
\end{equation}

In the next several lemmas, we record transfer estimates for the radiation component at the intermediate scale $K\leq A$ which are then passed to the directional derivative $\px w_1 + c \gamma w_2$.
\begin{lemma} Assuming \eqref{main_scaling}, it holds for any $10 \leq K \leq \frac{2A}{3}$ that\footnote{The proof can also be adapted to any $10\leq K \leq \theta A$ where $\theta \in (0,1)$, with the constant $C$ depending on $\theta$.}
\begin{equation}\label{eqn:transfer_estimate_sigmaA_z1}
    \frac{1}{K^2} \int_0^T \|\zeta_K  z_1\|^2 \, \ud t \leq C\delta^2 + C\int_0^T \big(\| \px w_1 + c \gamma w_2 \|^2 + |a|^4 + \|\rho^4 z_1\|^2 \big)\, \ud t.
\end{equation}
\end{lemma}
\begin{proof}
For $10 \leq K \leq 2A/3$, we first observe the elementary inequality 
\begin{equation}\label{eqn:elementary_inequality_K_A}
    \frac{K}{2} \leq \frac{AK}{2(A-K)} \leq \frac{3K}{2}.
\end{equation}
In particular, for this range of $K$, we have
\begin{align*}
    \frac{\zeta_K(x)^2}{\zeta_A(x)^2} = \exp\big[-2\kappa (\tfrac{1}{K} - \tfrac{1}{A}\big) (x-\barx)\eta(x)] \leq C \sech\Big( \tfrac{2\kappa (A-K)}{AK}x\Big),
\end{align*}
whence
\begin{align*}
\int_0^T \int \zeta_K^2 z_1^2  \leq \int_0^T \int w_1^2 \sech\Big( \tfrac{2\kappa (A-K)}{AK}x\Big).
\end{align*}
Now choose $\delta$ sufficiently small so that $|\pt (c\gamma)|= \gamma^3 |\dot{c}| \lesssim \delta^2$ and $\sup_{t}|c\gamma| \cdot \|  w_1 \|^2 \lesssim \delta^2$. We are in a position to apply Lemma~\ref{KMMvdB_lem4.6}  (with $J\leftarrow (x_1,x_2)$, $\omega \leftarrow \kappa$, $K \leftarrow \frac{AK}{2(A-K)}$, $w \leftarrow w_1$ and $\sigma \leftarrow c\gamma$), which yields
\begin{align*}
\int_0^T \int w_1^2 \sech\Big( \tfrac{2\kappa (A-K)}{AK}x\Big)\leq CK^2 \int_0^T \int (\px w_1 + c \gamma \dot{w}_1)^2 + CK \int_0^T \int_{x_1}^{x_2} w_1^2 + C K\delta^2,
\end{align*}
using also the upper bound from \eqref{eqn:elementary_inequality_K_A}. Note that 
\begin{equation*}
    \px w_1 + c \gamma \dot{w}_1 = (\px w_1 + c \gamma w_2) + c \gamma \zeta_A \underbrace{\big( (\dot{y}-c)\gamma \px z_1 - \dot{c}c \gamma^2 x \px z_1\big)}_{=\Omega_1(\bmz)}+ c \gamma \zeta_A \beta_1.
\end{equation*}
Since $|x \zeta_A| \leq CA$ (see \eqref{eqn:bound_on_zeta_K}), we have
\begin{equation*}
    \|\zeta_A \Omega_1(\bmz)\|^2 \leq C \delta^2(|\dot{y}-c| + A |\dot{c}|)^2 \lesssim |a|^4 + \|\rho^4 z_1\|^2.
\end{equation*}
Moreover, using  $\mathbbm{1}_{[x_1,x_2]} \leq C \rho^8$, it holds that
\begin{equation*}
    \int_0^T \int_{x_1}^{x_2} w_1^2 \leq C \int_0^T \int \rho^8 z_1^2.
\end{equation*}
Furthermore, Lemma~\ref{lemma:bound_on_beta_terms} provides that
\begin{equation*}
    \| \zeta_A \beta_1 \| \leq C\big(|a|^2 + \|\rho^4 z_1\|^2\big).
\end{equation*}
Thus, we conclude \eqref{eqn:transfer_estimate_sigmaA_z1} from the preceding estimates.
\end{proof}
The next lemma transfers estimates for the directional derivatives at scale $K$ to scale $A$.
\begin{lemma}
Assuming \eqref{main_scaling}, it holds for any $10 \leq K \leq A$ that
\begin{equation}\label{eqn:HE_transfer_zK1}
\int_0^T    \| \zeta_K (\px z_1 + c \gamma z_2) \|^2 \leq C \int_0^T \| \px w_1 + c \gamma w_2 \|^2 + \frac{C}{A^2}\int_0^T \| \zeta_K z_1\|^2,
\end{equation}
and hence for any $10 \leq K \leq \frac{2A}{3}$ that
\begin{equation}\label{eqn:HE_transfer_zK2}
\int_0^T \| \zeta_K (\px z_1 + c \gamma z_2)\|^2 + \frac{1}{K^2}\int_0^T \| \zeta_K z_1 \|^2 \leq C\delta^2 +  C\int_0^T \big(\| \px w_1 + c \gamma w_2 \|^2 + |a|^4 + \| \rho^4 z_1 \|^2 \big).
\end{equation}
\end{lemma}
\begin{proof}
For $K \leq A$, we have $0< \frac{\zeta_K}{\zeta_A} \leq 1$, and $|\zeta_A'/\zeta_A| \leq C A^{-1}$. Since
\begin{equation*}
    \zeta_K(\px z_1 + c \gamma z_2) = \frac{\zeta_K}{\zeta_A}(\px w_1 + c \gamma w_2) - \frac{\zeta_A'}{\zeta_A}\zeta_K z_1,
\end{equation*}
we obtain $\| \zeta_K (\px z_1 + c \gamma z_2) \|^2 \leq C \| \px w_1 + c \gamma w_2 \|^2 + CA^{-2} \| \zeta_K z_1\|^2$ and we have \eqref{eqn:HE_transfer_zK1}  after integrating over $[0,T]$. Hence, \eqref{eqn:HE_transfer_zK2} follows from \eqref{eqn:transfer_estimate_sigmaA_z1}, \eqref{eqn:HE_transfer_zK1}, and $A^{-2} \leq K^{-2}$.

\end{proof}

The lemma below controls $\zeta_K z_2$ in terms of the directional derivatives. 
\begin{lemma}
Assuming \eqref{main_scaling}, it holds for any $K \geq 10$ that
\begin{equation}\label{eqn:HE_transfer_z2}
    \int_0^T \|\zeta_K z_2 \|^2 \lesssim  \delta^2 +  \int_0^T |a|^4 +  \int_0^T \|\zeta_K (\px z_1+ c \gamma z_2)\|^2 + \int_0^T \|\zeta_K z_1\|^2.
\end{equation}
\end{lemma}
\begin{proof}The proof closely follows Lemma~4.8 of \cite{kowalczyk2021sufficient}. We define the virial functional 
\begin{equation}
    \calH_{\bmz} = \int z_1 z_2 \zeta_K^2.
\end{equation}
We have, using \eqref{eqn:z_1_virialform}, integration by parts, and then rearranging the terms, that
\begin{align*}
\dot{\calH}_{\bmz} &= \int \dot{z}_1 z_2 \zeta_K^2 + \int z_1 \dot{z}_2 \zeta_K^2 \\
&= (1+c^2 \gamma^2) \int z_2^2 \zeta_K^2 - \int (\px z_1 + c \gamma z_2)^2 \zeta_K^2 - \int z_1 f(x,z_1) \zeta_K^2 + \frac{1}{2} \int z_1^2 (\zeta_K^2)''\\
&\quad - 2 c \gamma \int z_1 z_2 (\zeta_K^2)' + \int z_2(\Omega_1(\bmz) + \beta_1) \zeta_K^2 + \int z_1(\Omega_2(\bmz) + \beta_2 + N_{\bmm}) \zeta_K^2  .
\end{align*}

From the identity $1+c^2\gamma^2 = \gamma^2 \geq 1$, we have a lower bound
\begin{equation*}
    (1+c^2 \gamma^2) \int z_2^2 \zeta_K^2 \geq \| \zeta_K z_2 \|^2 .
\end{equation*} 
By Taylor expansion, we have $|z_1 f(x,z_1)| \leq Cz_1^2$ which implies that 
\begin{equation*}
    \int z_1 f(x,z_1) \zeta_K^2  \leq C\| \zeta_K z_1 \|^2.
\end{equation*}
Moreover, using \eqref{eqn:bound_on_zeta_K}, we obtain that
\begin{equation*}
    \left|\frac{1}{2} \int z_1^2 (\zeta_K^2)''\right| \leq  C \| \zeta_K z_1 \|^2,
\end{equation*}
and
\begin{equation*}
    \left| 2 c \gamma \int z_1 z_2 (\zeta_K^2)'\right| \leq \frac{1}{2^4} \| \zeta_K z_2 \|^2 + C\| \zeta_K z_1 \|^2.
\end{equation*}
For the terms involving the modulation operator $\bmOmega(\bmz)$, we integrate by parts which yields
\begin{equation*}
\begin{split}
&\int z_2 \Omega_1(\bmz) \zeta_K^2 + \int z_1 \Omega_2(\bmz) \zeta_K^2  = -(\dot{y}-c)\gamma \int z_1 z_2 (\zeta_K^2)' + \dot{c}c \gamma^2 \int z_1 z_2 (\zeta_K^2 + x (\zeta_K^2)') - \frac{\dot{c}{\gamma}}{2}\int z_1^2 (\zeta_K^2)'.
\end{split}
\end{equation*}
Using $\zeta_K^2 + |(\zeta_K^2)'| + |x (\zeta_K^2)'|\leq C$ from  \eqref{eqn:bound_on_zeta_K}, \eqref{eqn:comparison_rho_zeta_K}, \eqref{eqn:stability_smallness}, and \eqref{eqn:estimate_modulation_crude}, we obtain
\begin{align*}
&\left| \int z_2 \Omega_1(\bmz) \zeta_K^2 + \int z_1 \Omega_2(\bmz) \zeta_K^2 \right| \leq C (|\dot{y}-c| + |\dot{c}|)(\delta \|\zeta_K z_1\| +  \| \zeta_K z_1 \|^2) \leq C \delta \big(|a|^4 + \| \zeta_K z_1 \|^2\big).
\end{align*}
By Lemma~\ref{lemma:bound_on_beta_terms},\eqref{eqn:comparison_rho_zeta_K}, and Young's inequality, 
\begin{equation*}
\left| \int (z_2 \beta_1 + z_1 \beta_2 + z_1N_{\bmm}) \zeta_K^2\right| \leq \frac{1}{2^4}\|\zeta_K z_2 \|^2 +  C (|a|^4  + \|\zeta_K z_1 \|^2).
\end{equation*}
Thus, the preceding bounds yield the estimate 
\begin{equation*}
    \dot{\calH}_{\bmz} \geq \frac{1}{2} \|\zeta_K z_2 \|^2 - \| \zeta_K(\px z_1 + c \gamma z_2)\|^2 - C(\| \zeta_K z_1\|^2 + |a|^4).
\end{equation*}
Integrating in time on $[0,T]$ and using $|\calH_{\bmz}| \leq \|z_1\| \|z_2 \| \leq C \delta^2$, we obtain
\begin{equation}
    \frac{1}{K^2} \int_0^T \| \zeta_K z_2 \|^2  \lesssim \delta^2 + \frac{1}{K^2} \int_0^T ( \| \zeta_K (\px z_1 + c \gamma z_2)\|^2 + \| \zeta_K z_1\|^2 + |a|^4),
\end{equation}
which finishes the proof.
\end{proof}
The previous lemmas combine to yield the following transfer estimate. 
\begin{corollary}\label{corollary:HE_main_transfer}
Assuming \eqref{main_scaling}, it holds for any $10 \leq K \leq \frac{2A}{3}$ that
\begin{equation}\label{eqn:HE_main_transfer}
\int_0^T (\|\zeta_K \px z_1 \|^2 + \|\zeta_K  z_1 \|^2 +  \|\zeta_K  z_2 \|^2 )\lesssim K^2\Big(\delta^2  + \int_0^T \|\px w_1 + c \gamma w_2 \|^2 + \int_0^T (|a|^4 + \|\rho^4 z_1 \|^2) \Big).
\end{equation}
\end{corollary}
\begin{proof}
We use $\|\zeta_K \px z_1 \|^2 \lesssim\|\zeta_K(\px z_1 + c \gamma z_2)\|^2 + \|\zeta_K z_2 \|^2$ and combine the estimates \eqref{eqn:transfer_estimate_sigmaA_z1}, \eqref{eqn:HE_transfer_zK1}, \eqref{eqn:HE_transfer_zK2}, and \eqref{eqn:HE_transfer_z2} to conclude \eqref{eqn:HE_main_transfer}.
\end{proof}

\subsection{Normal forms for the velocity parameter}
In view of Corollary~\ref{corollary:HE_main_transfer}, the next step is to estimate $\| \px w_1 + c \gamma w_2\|_{L_t^2L_x^2}$. However, as we explained in Remark~\ref{remark:virial_modification}, {the general virial functional \eqref{eqn:general_virial}} will require modifications because $|\dot{c}|$ is only bounded by $|a|^2$ and $\|\rho^4 z_1\|^2$; see  \eqref{eqn:estimate_modulation_crude}. In the key lemma below, we construct a normal form correction for the velocity parameter.
\begin{lemma}\label{lemma:normal_form_c} There exists a $\calC^1$ functional $    \Pi[c,\bmz,\bmm]$ satisfying $c + \Pi[c,\bmz,\bmm] \in (-1,1)$,
\begin{equation}\label{eqn:bound_on_Pi}
    |\Pi[c,\bmz,\bmm]| \lesssim |a|^2 + \|\rho^4 z_1\|^2 + \| \rho^4 z_2 \|^2,
\end{equation}
and
\begin{equation}
    \bigg|\dot{c} + \frac{\ud}{\ud t}{\Pi[c,\bmz,\bmm]} \bigg| \lesssim  |a|^4 + \| \rho^4 z_1\|^2 + \| \rho^4 z_2 \|^2.
\end{equation}
\end{lemma}
\begin{proof}
We first recall the  notation  $\bmphi = \calT_{c,y}[\bmpsi]=\calT_{c,y}[\bmH + \bmz + \bmm]$. After expanding the conserved momentum $\calP[\bmphi]$ and using the modulational orthogonality \eqref{eqn:orthogonality_bmu}, we obtain the identity
\begin{align*}
&\calP[\bmphi] = \int \phi_2 \px \phi_1 = \int \psi_2\psi_1' - c \gamma \int (\psi_1')^2 =- c \gamma \| H' \|^2 + \bmP[c,\bmz] + \bmQ[c,\bmz,\bmm],
\end{align*}
where
\begin{align*}
    \bmP[c,\bmz] &:=\int (z_2 - c \gamma z_1') z_1',\\
    \bmQ[c,\bmz,\bmm] &:= \int[m_2 m_1' - c \gamma (m_1')^2 + z_2 m_1' + z_1(2c \gamma m_1'' - m_2')].
\end{align*}
Since $\calP[\bmphi]$ is a conserved quantity, we have, after differentiating in time, 
\begin{equation}
    \frac{\ud}{\ud t} \left( c \gamma - \frac{\bmQ}{\|H'\|^2}\right) = \frac{1}{\|H'\|^2} \frac{\ud}{\ud t} \bmP[c,\bmz].
\end{equation}
After lengthy computation using \eqref{eqn:z_1_virialform}, we obtain that 
\begin{align*}
\frac{\ud}{\ud t} \bmP[c,\bmz] &= \int \dot{z}_2 z_1' + \int z_2 \px \dot{z}_1 - \frac{\ud}{\ud t}{(c\gamma)} \int (z_1')^2 - 2 c \gamma \int z_1' \px \dot{z}_1\\
&= \int(\partial_x^2z_1)(\partial_xz_1)
-\int f(x,z_1)\partial_xz_1 +2c\gamma\int(\partial_xz_2)(\partial_xz_1)
+\int z_2\partial_xz_2
-2c\gamma\int(\partial_xz_1)(\partial_xz_2)\\
&\quad+(\dot y-c)\gamma
\int \Big[(\partial_xz_2)(\partial_xz_1)+ z_2\partial_x^2z_1-2c\gamma (\partial_xz_1)(\partial_x^2z_1)\Big]\\
&\quad-\dot c\,c\gamma^2
\int \Big[  x(\partial_xz_2)(\partial_xz_1) + z_2\partial_xz_1+ xz_2\partial_x^2z_1\Big]\\
&\quad+\dot c\,\gamma
\int(\partial_xz_1)^2
-\dot c\,\gamma^3
\int(\partial_xz_1)^2+2\dot c\,c^2\gamma^3
\int \Big[(\partial_xz_1)^2+ x(\partial_xz_1)(\partial_x^2z_1)\Big]\\
&\quad+\int(\beta_2 + N_{\bmm})\,\partial_xz_1
+\int z_2\partial_x\beta_1
-2c\gamma\int(\partial_xz_1)(\partial_x\beta_1)\\
&= -\int f(x,z_1)\px z_1 + \int (\beta_2 + N_{\bmm}) \px z_1 + \int (z_2 - 2c \gamma \px z_1) \px \beta_1 .
\end{align*}
Note that all terms in the intermediate equality that contain $\dot{y}-c$ and $\dot{c}$, as well as the anti-symmetric terms, either cancel or evaluate to zero after integrating by parts, which can be justified by smooth approximation. 

We first estimate $\frac{\ud}{\ud t} \bmP[c,\bmz]$. Using integration by parts, Taylor expansion, and \eqref{eqn:decay_of_H}, we have
\begin{equation*}
    \left| \int f(x,z_1) \px z_1 \right| = \left| \int H'[W'(H+z_1)-W'(H)-W''(H)z_1]\right| \lesssim \| \rho^4 z_1\|^2,
\end{equation*}
and by \eqref{eqn:bound_on_Nmm}
\begin{align*}
\left| \int N_{\bmm} \px z_1 \right| = \left| \int \int_0^{z_1} \px [W'(H+s+m_1)-W'(H+s)-W''(H)m_1]\,\ud s \right| \lesssim |a|^4 + \| \rho^4 z_1 \|^2.
\end{align*}
Integration by parts, \eqref{eqn:modulation_comparison},  \eqref{eqn:bound_on_beta}, and Young's inequality then imply that 
\begin{align*}
\left| \int \beta_2 \px z_1 \right| + \left|\int (z_2 - 2c \gamma \px z_1) \px \beta_1  \right| &\lesssim \int (|\px \beta_1| + |\px \beta_2|) \big(|z_1| + |z_2|\big )+  \int |c \gamma z_1||\px^2 \beta_1|\\
&\lesssim |a|^4 + \| \rho^4 z_1\|^2 +  \|\rho^4 z_2 \|^2.
\end{align*}
Hence, we obtain the desired quartic bound
\begin{equation}\label{eqn:bound_on_P}
    \left| \frac{\ud}{\ud t} \bmP[c,\bmz] \right| \lesssim |a|^4 + \| \rho^4 z_1\|^2 +  \|\rho^4 z_2 \|^2.
\end{equation}
Now, let $h(x) := \frac{x}{\jx} = \frac{x}{\sqrt{1+x^2}}$ and note that $h'(x) = \frac{1}{(1+x^2)^{3/2}}$. Both $h(x)$ and $h'(x)$ are uniformly bounded by 1 over $\bbR$. 
We define the functional
\begin{equation}\label{eqn:definition_Pi}
    \Pi[c,\bmz,\bmm] := h\left( c \gamma - \frac{\bmQ[c,\bmz,\bmm]}{\|H'\|^2}\right) - h(c\gamma).
\end{equation}
Since $c = h(c\gamma)$, we have $c + \Pi[c,\bmz,\bmm] = h\left( c \gamma - \frac{\bmQ[c,\bmz,\bmm]}{\|H'\|^2}\right) \in (-1,1)$. Moreover, recall that $\bmm = 2\Re(a(t)\Phi_c^+)$. Using the uniform boundedness of $h'$, \eqref{eqn:estimate_on_Phi_c_pm}, and \eqref{eqn:bound_on_P}, we obtain from the mean value theorem that   
\begin{equation*}
    |\Pi[c,\bmz,\bmm]| \leq \max(|h'|) \cdot \frac{|\bmQ[c,\bmz,\bmm]|}{\|H'\|^2} \lesssim |a|^2 + \|\rho^4 z_1 \|^2 +  \|\rho^4 z_2 \|^2,
\end{equation*}
and 
\begin{equation*}
    \left|\frac{\ud}{\ud t}\left(c+ \Pi[c,\bmz,\bmm]\right)\right| \leq h'\left( c \gamma - \frac{\bmQ[c,\bmz,\bmm]}{\|H'\|^2}\right) \frac{1}{\|H'\|^2} \left |\frac{\ud}{\ud t} \bmP[c,\bmz] \right| \lesssim |a|^4 + \|\rho^4 z_1 \|^2 +  \|\rho^4 z_2 \|^2.
\end{equation*}
\end{proof}
In the following lemma, we introduce scalar coefficients which will be used to modify the main virial identity.
\begin{lemma}\label{lemma:coefficients_dj}
Let $d_j=d_j(c,\Pi[c,\bmz,\bmm])$, $j\in\{1,2,3\}$ be scalars given by 
\begin{equation}\label{eqn:d_j_coefficients}
d_1 := \frac{1-c(c+\Pi)}{1-(c+\Pi)^2},\qquad d_2 := c\gamma \frac{1-c(c+\Pi)}{1-(c+\Pi)^2},\qquad d_3 := \frac{c+\Pi}{\gamma[1-(c+\Pi)^2]},
\end{equation}
where $\Pi = \Pi[c,\bmz,\bmm]$ is given by \eqref{eqn:definition_Pi}. Then, we have the following identities 
\begin{align}
d_2 &= c\gamma d_1,\qquad    d_3 - d_2 = \frac{\gamma \Pi}{1-(c+\Pi)^2},\label{eqn: d_2-d_3=2cgd_1}\\
\dot{d}_1 + \dot{c}\gamma d_3 &= \frac{\dot{c}+\dot{\Pi}}{1-(c+\Pi)^2}[-c + 2(c + \Pi)d_1]\label{eqn:cancellation_dj_1},\\
\dot{d}_2 - \dot{c}\gamma d_1 - \dot{c} c \gamma^2 d_2 + \dot{c} c \gamma^2 d_3 &= \frac{\dot{c}+\dot{\Pi}}{1-(c+\Pi)^2}[-c^2 \gamma + 2(c + \Pi) d_2]  \label{eqn:cancellation_dj_2},\\
\dot{d}_3 + \dot{c} c \gamma^2 d_3 &= \frac{\dot{c}+\dot{\Pi}}{1-(c+\Pi)^2}[\gamma^{-1} + 2(c+\Pi)d_3] \label{eqn:cancellation_dj_3}.
\end{align}
Moreover, for sufficiently small $\delta>0$, the following comparison estimates hold
\begin{equation}\label{eqn:comparison_d_1}
    \frac{1}{2} \leq d_1 \leq \frac{3}{2},
\end{equation}
\begin{equation}\label{eqn:comparison_d_23}
|d_2 - c_0 \gamma_0| \leq C\delta, \qquad  |d_3 - c_0 \gamma_0| \leq C\delta .
\end{equation}
\end{lemma}
\begin{proof}
The first identity \eqref{eqn: d_2-d_3=2cgd_1} immediately follows from \eqref{eqn:d_j_coefficients} and $1+c^2\gamma^2 = \gamma^2$. By the product rule and $\pt{\gamma^{-1}} = -\dot{c}c\gamma$, we compute the following derivatives
\begin{align*}
    \dot{d}_1 &= - \frac{\dot{c}(c+\Pi)}{1-(c+\Pi)^2} +  \frac{\dot{c} + \dot{\Pi}}{1-(c+\Pi)^2}[-c+2(c+\Pi)d_1],\\
    \dot{d}_3 &= \frac{\gamma^{-1}(\dot{c} + \dot{\Pi})}{1-(c+\Pi)^2} - \frac{\dot{c} c \gamma(c+\Pi)}{1-(c+\Pi)^2} + \frac{\dot{c} + \dot{\Pi}}{1-(c+\Pi)^2}\cdot 2(c+\Pi)d_3.
\end{align*}
Since $ \gamma d_3 = (c+\Pi)/[1-(c+\Pi)^2]$, the identities \eqref{eqn:cancellation_dj_1} and \eqref{eqn:cancellation_dj_3} follow from the above expressions. Moreover, since $\pt(c\gamma) = \dot{c} \gamma^3$, $\gamma^3 = \gamma + c^2 \gamma^3$, and $d_2 = c \gamma d_1$, we have
\begin{align*}
\dot{d}_2 = \dot{c} \gamma^3 d_1 + c \gamma \dot{d}_1 = \dot{c} \gamma d_1 + \dot{c} c^2 \gamma^3 d_1 + c \gamma \dot{d}_1 = \dot{c} \gamma d_1 + \dot{c} c \gamma^2 d_2 + c \gamma \dot{d}_1.
\end{align*}
Using \eqref{eqn:cancellation_dj_1}, we obtain 
\begin{align*}
    \dot{d}_2 - \dot{c}\gamma d_1 - \dot{c} c \gamma^2 d_2 + \dot{c} c \gamma^2 d_3 = c \gamma (\dot{d}_1 + \dot{c} \gamma d_3)=\frac{\dot{c}+\dot{\Pi}}{1-(c+\Pi)^2}[-c^2 \gamma + 2(c + \Pi) d_2] ,
\end{align*}
which proves \eqref{eqn:cancellation_dj_2}. Next, we show the comparison bounds. By \eqref{eqn:stability_smallness}, \eqref{eqn:estimate_modulation_crude}, and Lemma~\ref{lemma:normal_form_c}, we have $|c-c_0| \lesssim \delta$ and $|\Pi|\lesssim \delta^2$. Hence, for sufficiently small $\delta>0$, we have 
\begin{align*}
    1-(c+\Pi)^2 \geq \frac12 (1-c_0^2) = \frac{1}{2\gamma_0^2}, \qquad |\Pi| \leq \frac{1}{4\gamma_0^2}.
\end{align*}
From these bounds, we infer from \eqref{eqn:d_j_coefficients} that 
\begin{equation*}
|d_1 - 1| \leq \frac{|c+\Pi||\Pi|}{|1-(c+\Pi)^2|}\leq  2 \gamma_0^2 |\Pi|\leq \frac12,
\end{equation*}
which proves \eqref{eqn:comparison_d_1}. Finally, since $\pc(c \gamma) = \gamma^3$, the mean value theorem, together with  \eqref{eqn:modulation_comparison}, gives $|c\gamma - c_0 \gamma_0| \leq C\gamma_0^3|c-c_0| \leq C\gamma_0^3 \delta$. Since $d_2 = c \gamma d_1$, we have $|d_2 - c_0 \gamma_0| \leq |c\gamma| |d_1 - 1| + |c \gamma - c_0 \gamma_0| \lesssim \delta$ which ,together with \eqref{eqn: d_2-d_3=2cgd_1}, implies \eqref{eqn:comparison_d_23}.
\end{proof}
\begin{remark}
The coefficients $d_j$ introduced in the previous lemma will incorporate the corrected velocity $c + \Pi$ into the virial argument in Lemma~\ref{lemma:HE_key_estimate1}. The cancellation identities \eqref{eqn:cancellation_dj_1}--\eqref{eqn:cancellation_dj_3} ensure that the terms arising from differentiation of these coefficients combine into expressions proportional to $\dot{c}+\dot{\Pi}$, which yields an improved control from $|a|^2$ to $|a|^4$. Moreover, the relation $d_2 = c \gamma d_1$ preserves the structure of the term $(\px w_1 + c \gamma w_2)^2$, and that $d_1 \geq \frac12$ ensures coercivity.
\end{remark}

\subsection{Modified virial estimate}\label{sec:modified_virial}

In the next lemma, we modify the virial functional from Lemma~\ref{lemma:general_virial} by incorporating the coefficients $d_j$ and estimate the directional derivative of the radiation component at scale $A$.
\begin{lemma}\label{lemma:HE_key_estimate1}
Assuming \eqref{main_scaling}, it holds for $A\geq 20$ that 
\begin{equation}\label{eqn:HE_key_estimate1}
 \int_0^T \| \px w_1 + c \gamma w_2\|^2 \ud t  \lesssim A \delta^2 + \int_0^T \Big( |a|^4 + \| \rho^4 z_1 \|^2 \Big) \ud t + \int_0^T \int |z_1|^3 \zeta_A^2.
\end{equation}
 
\end{lemma} 
\begin{proof}
We introduce the modified virial functional\footnote{Observe that when $\bmm\equiv 0$, we have $\Pi = 0$, and the coefficients reduce to $d_1 = 1$, $d_2 = d_3 = c \gamma$ and the virial functional $\calZ$ reduces to $\calZ_1 + c \gamma (\calZ_2 + \calZ_3)$ which is the virial functional used in Lemma~4.4 of \cite{kowalczyk2021sufficient}.}
\begin{equation}
    \calZ = d_1 \calZ_1 + d_2 \calZ_2 + d_3 \calZ_3,
\end{equation}
where $d_j = d_j(c,\Pi)$ are defined in \eqref{eqn:d_j_coefficients} and set $\calZ_j := \calA_j[\varphi_A](\bmz)$, for $j \in \{1,2,3\}$, in the notation of Lemma~\ref{lemma:general_virial}. From the system \eqref{eqn:z_1_virialform} and Lemma~\ref{lemma:general_virial}, we obtain
\begin{align*}
\dot{\calZ} &=d_1 \bigg \{ -2 \int \varphi_A' (z_1')^2 + \int \frac12 \varphi_A''' z_1^2 + \int \varphi_A'(2F-z_1f) + 2 \int \varphi_A (\px F) - 4 c \gamma \int z_1'' z_2 \varphi_A\\
&\quad \qquad  -6 c \gamma \int \varphi_A' z_1' z_2 - 2 c \gamma \int \varphi_A'' z_1 z_2 - (\dot{y}-c)\gamma \bigg( 2 \int \varphi_A' z_1' z_2 + \int \varphi_A'' z_1 z_2 \bigg)\\
&\quad \qquad + \dot{c} \bigg[-\gamma \calZ_2 + 2 c \gamma^2 \int x \varphi_A' z_1' z_2 + c \gamma^2 \int \varphi_A' z_1 z_2 + c \gamma^2 \int x \varphi_A'' z_1 z_2 - \frac{1}{2}\gamma \int \varphi_A'' z_1^2 \bigg]\\
&\quad \qquad - \int (2 \varphi_A z_2' + \varphi_A' z_2)\beta_1 + \int(2\varphi_A z_1' + \varphi_A' z_1)(\beta_2 + N_{\bmm}) \bigg\}\\
&+ d_2 \bigg\{4 \int \varphi_A z_1'' z_2 + 4 \int \varphi_A' z_1' z_2 + 2(\dot{y}-c) \gamma \int \varphi_A' (z_1')^2 - \dot{c} c \gamma^2 \bigg[\calZ_2+ 2\int x \varphi_A' (z_1')^2\bigg] -4 \int \varphi_A z_1' \beta_1'\bigg\}\\
&+d_3 \bigg\{ -2 \int \varphi_A' z_1' z_2 - 2 c \gamma \int \varphi_A' z_2^2 - (\dot{y}-c)\gamma \int \varphi_A'\big( (z_1')^2 + z_2^2 + 2F\big) -2 (\dot{y}-c)\gamma \int \varphi_A (\px F)\\
&\quad \qquad +\dot{c} \bigg[2\gamma \int \varphi_A z_1' z_2 + c \gamma^2(\calZ_2 + \calZ_3) + c \gamma^2 \int x\varphi_A' \big( (z_1')^2 + z_2^2 + 2F\big) +2c \gamma^2 \int x \varphi_A (\px F) \bigg]\\
&\quad \qquad + 2 \int \varphi_A \big( z_1' \beta_1' + z_2 (\beta_2 + N_{\bmm}) + f \beta_1 \big)\bigg\}+\dot{d}_1 \calZ_1 + \dot{d}_2 \calZ_2 + \dot{d}_3 \calZ_3 .
\end{align*}
First, integration by parts yields
\begin{align*}
    2 \dot{c} \gamma d_3 \int \varphi_A z_1' z_2 = \dot{c} \gamma d_3 \calZ_1 - \dot{c} \gamma d_3 \int \varphi_A' z_1 z_2.
\end{align*}
At this point, we group the $\dot{c}$ terms multiplied by one of the $\calZ_j$ terms and we define the remainder term 
\begin{align*}
    \calR_{\calZ} := (\dot{d}_1 + \dot{c} \gamma d_3) \calZ_1 + (\dot{d}_2 - \dot{c} \gamma d_1 - \dot{c} c \gamma^2 d_2 + \dot{c} c \gamma^2 d_3) \calZ_2 + (\dot{d}_3 + \dot{c} c \gamma^2 d_3)\calZ_3.
\end{align*}
Notice that the coefficients of $\calZ_j$ in the above expression match the key cancellation identities given in \eqref{eqn:cancellation_dj_1}--\eqref{eqn:cancellation_dj_3}. Hence, we obtain a favorable expression for the remainder term, 
\begin{equation*}
    \calR_{\calZ} \equiv \frac{\dot{c}+\dot{\Pi}}{1-(c+\Pi)^2}\big[-c\mathcal Z_1-c^2\gamma\mathcal Z_2
+\gamma^{-1}\mathcal Z_3+2(c+\Pi)\mathcal Z\big].
\end{equation*}
We now regroup the other terms in $\dot{\calZ}$  and display the directional derivative $(\px w_1 + c \gamma w_2)^2$ more explicitly. We compute that 
\begin{align*}
\dot{\calZ} 
&=-2d_1 \int (\px w_1 + c \gamma w_2)^2 - d_1 \int w_1^2 (\ln\zeta_A)'' + 2(d_2 - d_3) \int [(\px w_1 + c \gamma w_2)w_2 - w_1w_2 \ln(\zeta_A)']\\
&\quad+ d_1 \int [2F(x,z_1) - z_1 f(z_1)]\zeta_A^2 + 2d_1 \int (\px F)(x,z_1) \varphi_A  - (\dot{y}-c)\gamma \calI + \dot{c} \calJ + \calM_1 + \calM_2 + \calR_{\calZ},
\end{align*}
where 
\begin{align*}
\calI &:= 
\int\varphi_A'
\Big[
2d_1(\partial_xz_1)z_2
+(d_3-2d_2)(\partial_xz_1)^2
+d_3\bigl(z_2^2+2F(x,z_1)\bigr)
\Big]\\
&\qquad +d_1\int\varphi_A''z_1z_2
+2d_3\int\varphi_A(\partial_xF)(x,z_1),
 \\
\calJ &:=
c\gamma^2\int x\varphi_A'
\Big[
2d_1(\partial_xz_1)z_2
+(d_3-2d_2)(\partial_xz_1)^2
+d_3\bigl(z_2^2+2F(x,z_1)\bigr)
\Big]\\
&\qquad+ (d_2 -d_3)\gamma\int\varphi_A'z_1z_2
+c\gamma^2d_1\int x\varphi_A''z_1z_2\\
&\qquad -\frac{\gamma d_1}{2}\int\varphi_A''z_1^2
+2c\gamma^2d_3\int x\varphi_A(\partial_xF)(x,z_1),
\\
\calM_1&:=
d_1\int(2\varphi_A\partial_x\beta_1+\varphi_A'\beta_1)z_2+2(d_3-2d_2)\int\varphi_A(\partial_xz_1)(\partial_x\beta_1)
+2d_3\int\varphi_A f(x,z_1)\beta_1,\\
\calM_2&:=
d_1\int(2\varphi_A\partial_xz_1+\varphi_A'z_1)(\beta_2+ N_{\bmm})
+2d_3\int\varphi_Az_2(\beta_2+ N_{\bmm}).
\end{align*} 

First, by $d_1 \geq \frac12$ from \eqref{eqn:comparison_d_1}, we have 
\begin{align*}
    -2d_1\int (\px w_1 + c \gamma w_2)^2\leq -\int (\px w_1 + c \gamma w_2)^2.
\end{align*}
By \eqref{eqn:comparison_rho_zeta_K} and \eqref{eqn:bound_on_zeta_K},
\begin{align*}
    \int |w_1^2 \ln(\zeta_A)'' | = \int z_1^2 | \zeta_A'' \zeta_A - (\zeta_A')^2| \leq C A^{-1} \int z_1^2 \rho^8.
\end{align*}
Next, for a small $\theta \in (0,1)$ to be chosen later, we use \eqref{eqn: d_2-d_3=2cgd_1} and  $\|w_2\| \leq C\delta$ to obtain the bound 
\begin{align*}
&2|d_2-d_3|    \int |(\px w_1 + c \gamma w_2) w_2|\lesssim \delta |\Pi| \| \px w_1 + c \gamma w_2 \| \leq \theta \| \px w_1 + c \gamma w_2\|^2 + C_\theta \delta^2 |\Pi|^2.
\end{align*}
Since $|\ln(\zeta_A)'| \leq CA^{-1}$ and $\zeta_A^2 = \zeta_{A/2}\leq \zeta_{2A/3}$,
\begin{align*}
    \left|2 (d_2-d_3) \int (\ln \zeta_A)'w_1 w_2 \right| \lesssim  \delta A^{-1} |\Pi|  \| \zeta_{2A/3} z_1\|  \leq \theta A^{-2} \| \zeta_{2A/3} z_1 \|^2 + C_\theta \delta^2|\Pi|^2.
\end{align*}
By \eqref{eqn:bound_on_Pi}, we have $\delta^2|\Pi|^2 \lesssim |a|^4 + \|\rho^4 z_1\|^2 + \delta^4 \|\rho^4 z_2\|^2$. These preceding bounds imply that
\begin{align*}
    &\left| 2(d_2 - d_3) \int [(\px w_1 + c \gamma w_2)w_2 - w_1w_2 \ln(\zeta_A)']\right| \\
&\leq \theta \| \px w_1 + c \gamma w_2\|^2 +\theta A^{-2} \| \zeta_{2A/3} z_1 \|^2+C_\theta(|a|^4 + \|\rho^4 z_1\|^2) + C_\theta \delta^4 \|\rho^4 z_2\|^2.
\end{align*}

For the nonlinear terms, we use Taylor expansion, \eqref{eqn:comparison_rho_zeta_K}, and the bound $\| z_1 \|_\infty \leq C\delta$ to obtain
\begin{align*}
   &d_1\int |[2F(z_1)- z_1 f(z_1)]| \zeta_A^2 \\
   &= d_1\int |2 W(H+z_1) - (W'(H+z_1) +W'(H))z_1 - 2W(H)| \zeta_A^2 \leq C \int |z_1|^3 \zeta_A^2, 
\end{align*}
and
\begin{align*}
2d_1    \int |(\px F)(z_1) \varphi_A| = 2d_1\int | [W'(H+z_1) - W'(H) - W''(H)z_1]H' \varphi_A| \leq C \| \rho^4 z_1 \|^2.
\end{align*}

Next, we estimate the quadratic terms $\calI$ and $\calJ$. By \eqref{eqn:bound_varphi'}--\eqref{eqn:bound_varphi''}, we recall the following comparison bounds
\begin{equation}\label{eqn:bound_varphi_A'}
\varphi_A' \leq \zeta_{2A/3}, \quad |x \varphi_A'| \leq CA \zeta_{2A/3}, \quad |\varphi_A''| \leq CA^{-1} \zeta_{2A/3}, \quad |x\varphi_A''| \leq C \zeta_{2A/3}.
\end{equation}

Moreover, using \eqref{eqn:bound_varphi_A'}, \eqref{eqn:decay_of_H} and Taylor expansion, we observe that the $F(x,z_1)$-terms in $\calI$ and $\calJ$ satisfy 
\begin{align*}
\int    |\varphi_A' F(x,z_1)| + |\varphi_A (\px F)(x,z_1)| +|x\varphi_A' F(x,z_1) | + |x \varphi_A (\px F)(x,z_1)|&\leq C A \delta \| \zeta_{2A/3} z_1\|.
\end{align*}
Thus, using the above bound, \eqref{eqn:bound_varphi_A'}, the comparison estimates \eqref{eqn:comparison_d_1}-\eqref{eqn:comparison_d_23}, the Cauchy-Schwarz inequality, and then the smallness estimate \eqref{eqn:stability_smallness}, we derive that
\begin{align*}
|\calI| &\lesssim \delta \big( \| \zeta_{2A/3} z_1 \| + \| \zeta_{2A/3} \px z_1 \| +\| \zeta_{2A/3} z_2 \|\big) ,\\
|\calJ| &\lesssim A\delta \big( \| \zeta_{2A/3} z_1 \| + \| \zeta_{2A/3} \px z_1 \| +\| \zeta_{2A/3} z_2 \|\big) .
\end{align*}
By \eqref{eqn:estimate_modulation_crude}, Young's inequality, and $\delta \leq A^{-4}$, we obtain
\begin{align*}
    |(\dot{y}-c)\gamma \calI| + |\dot{c} \calJ| &\leq C(|a|^2+\|\rho^4 z_1\|^2) A\delta \big( \| \zeta_{2A/3} z_1 \| + \| \zeta_{2A/3} \px z_1 \| +\| \zeta_{2A/3} z_2 \|\big)\\
    &\leq \frac{\theta}{A^2} \big( \| \zeta_{2A/3} \px z_1 \|^2 + \| \zeta_{2A/3}z_1\|^2 + \| \zeta_{2A/3} z_2 \|^2\big) + C_\theta (|a|^4 + \| \rho^4 z_1 \|^2).
\end{align*}

For the localized error terms $\calM_1,\calM_2$, we use Lemma~\ref{lemma:bound_on_beta_terms}, $|\varphi_A|\leq |x|$, and Young's inequality to obtain 
\begin{align*}
    &|\calM_1| + |\calM_2| \\
    &\leq C \big( \| \rho^4 z_1\| + \|\rho^4 \px z_1 \| + \| \rho^4 z_2\|\big)\big(\|\jx  \rho^{-4} \beta_1\|+\|\jx  \rho^{-4} \px\beta_1\|+\|\jx  \rho^{-4} \beta_2\|+\|\jx  \rho^{-4} N_{\bmm}\|\big)\\
    &\leq \theta\big(  \|\rho^4 \px z_1 \|^2 + \| \rho^4 z_2\|^2\big) + C_\theta(|a|^4 + \|\rho^4 z_1\|^2).
\end{align*}

Finally, we note that 
\begin{equation}
    |\calZ_1| +|\calZ_2| +|\calZ_3| +|\calZ| \lesssim A \delta^2.
\end{equation}
Hence, Lemma~\ref{lemma:normal_form_c} gives 
\begin{align*}
    |\calR_{\calZ}| \lesssim A \delta^2|\dot{c}+\dot{\Pi}| \lesssim A \delta^2 (|a|^4 + \|\rho^4 z_1 \|^2 + \| \rho^4 z_2\|^2).
\end{align*}

Combining the preceding estimates, we conclude that 
\begin{equation}
\begin{split}
    \dot{\calZ} &\leq - \frac{3}{4} \| \px w_1 + c \gamma w_2\|^2  + \frac{C\theta}{A^2}\Big(
\|\zeta_{2A/3}\partial_xz_1\|^2
+\|\zeta_{2A/3} z_1\|^2
+\|\zeta_{2A/3} z_2\|^2
\Big)\\
&\qquad +C\theta\Big(\|\rho^4 \partial_xz_1\|^2+\|\rho^4 z_1\|^2+\|\rho^4 z_2\|^2
\Big)+C_\theta\bigl(|a|^4+\|\rho^4z_1\|^2\bigr)
+C\int|z_1|^3\zeta_A^2\\
&\qquad +\bigl(C_\theta\delta^4+CA\delta^2\bigr)\|\rho^4z_2\|^2.
\end{split}
\end{equation}
We integrate over $[0,T]$, use the bound $|\calZ| \leq C A \delta^2$, and then apply \eqref{eqn:HE_main_transfer} twice (with $K=2A/3$ and $K=10$ together with $\rho^4 \leq C\zeta_{K=10}$) to obtain the integrated-in-time estimate 
\begin{equation*}
\begin{split}
&\frac34\int_0^T\|\partial_xw_1+c\gamma w_2\|^2\,dt\\
&\leq CA\delta^2+C_\theta\int_0^T\bigl(|a|^4+\|\rho^4z_1\|^2\bigr)\,\ud t+C\int_0^T\int|z_1|^3\zeta_A^2\\
&\qquad +\frac{C\theta}{A^2} \int_0^T \big( \| \zeta_{2A/3} \px z_1 \|^2 + \| \zeta_{2A/3}z_1\|^2 + \| \zeta_{2A/3} z_2 \|^2\big) + C\theta \int_0^T \Big(\|\rho^4 z_1\|^2+\|\rho^4 \px z_1\|^2+\|\rho^4 z_2\|^2\Big)\\
&\leq CA\delta^2+C_\theta\int_0^T\bigl(|a|^4+\|\rho^4z_1\|^2\bigr)\,\ud t+C\int_0^T\int|z_1|^3\zeta_A^2\\
&\qquad +\bigl(C\theta+C_\theta\delta^4+CA\delta^2\bigr)
\left[\delta^2+\int_0^T\Big(
\|\partial_xw_1+c\gamma w_2\|^2
+|a|^4+\|\rho^4z_1\|^2
\Big)\,\ud t
\right].
\end{split}
\end{equation*}
Now choose $\theta$ small enough so that $C\theta \leq \frac{1}{8}$ and then choose $\delta$ sufficiently small so that $(C\theta+C_\theta\delta^4+CA\delta^2\bigr) \leq \frac12$, which then allows us to absorb the term $\frac12 \| \px w_1 + c \gamma w_2\|_{L_t^2L_x^2}^2$ in parentheses to the left-hand side and conclude \eqref{eqn:HE_key_estimate1}.
\end{proof}

Next, we estimate the cubic terms using a technical lemma from \cite[Lemma~4.5]{kowalczyk2021sufficient}. 
\begin{lemma} Assuming \eqref{main_scaling}, it holds for all $A\geq 20$ that 
\begin{equation}\label{eqn:HE_estimate_cubic}
    \int_0^T \int |z_1|^3 \zeta_A^2 \lesssim A \delta^3 + A^2 \delta \int_0^T \big(\| \px w_1 + c \gamma w_2 \|^2 + |a|^4 + \| \rho^4 z_1 \|^2 \big) \ud t.
\end{equation}
\end{lemma}
\begin{proof}
For sufficiently large $A$, we have $\sup_t |\pt (c\gamma)| \leq C\delta^2 \leq \frac{\kappa}{8A}$.
Since $\zeta_A(x) \simeq \sech(\kappa x /A)$ and $w_1 = \zeta_A z_1$, we may apply Lemma~\ref{lemma45} (with $w \leftarrow w_1$, $\sigma \leftarrow c \gamma$, $K \leftarrow  \frac{\omega}{\kappa}A$) to obtain 
\begin{align*}
&\int_0^T \int |z_1|^3 \zeta_A^2 \leq C \int_0^T \int |w_1|^3 \cosh\big(\frac{\kappa x}{A}\big)\\
&\leq CA \|z_1\|_\infty \bigg[A \int_0^T \int (\px w_1 + c \gamma \dot{w}_1)^2 + \sup_{[0,T]}\|w_1\|^2\bigg]\\
&\leq CA \|z_1\|_\infty \bigg[A \int_0^T \int (\px w_1 + c \gamma w_2)^2 + A\int_0^T \int \zeta_A^2 (\Omega_1(\bmz)+\beta_1)^2 + \sup_{[0,T]}\|w_1\|^2\bigg].
\end{align*}
Using the preparatory estimates from Lemma~\ref{lemma:bound_on_beta_terms}, the estimates$\|w_1\|_{L^\infty([0,T]\times\bbR)}  \leq C\delta$, and $|x\zeta_A|\leq CA$, we then conclude \eqref{eqn:HE_estimate_cubic}.
\end{proof}

We now complete the proof of Proposition~\ref{prop:high_energy_estimate}. Combining \eqref{eqn:HE_key_estimate1} and \eqref{eqn:HE_estimate_cubic}, we obtain the estimate
\begin{equation}\label{eqn:HE_key_estimate2}
 \int_0^T \| \px w_1 + c \gamma w_2\|^2 \ud t  \lesssim A \delta^2 + \int_0^T \Big( |a|^4 + \| \rho^4 z_1 \|^2 \Big) \ud t.
\end{equation}
We then apply \eqref{eqn:HE_main_transfer} with $K\in[10,\frac{A}{2}]$ and use $\rho^4 \leq C\zeta_K$ to conclude \eqref{eqn:estimate_HE} which finishes the proof.

\section{Estimate of the internal mode component} 
In this section, we use the Fermi golden rule to estimate a quartic-in-time bound for the internal mode amplitude. The goal for this section is to prove the following proposition.
\begin{proposition}\label{prop:internal_mode_estimate} Assuming \eqref{main_scaling} and $A^{-\frac12} \ll \Gamma/2$, we have 
    \begin{equation}
        \int_0^T |a|^4 \ud t \lesssim   A \delta^2 + \frac{1}{\sqrt{A}} \int_0^T \| \rho^4 z_1 \|^2 \ud t.
    \end{equation}
\end{proposition}

First, we begin by constructing the adjoint generalized eigenfunction at the resonant double frequency.
\begin{lemma}\label{lemma:generalized_eigenfunction}
    Let $g$ be a bounded solution to 
    \begin{equation*}
        L_0 g = 4 \lambda^2 g.
    \end{equation*}
    Then
    \begin{equation}
        \bmG_c(x) := e^{-2ic \lambda x}\begin{pmatrix}
            2c \gamma (\px g)(x) - 2i \lambda \gamma(1+c^2) g(x) \\ 
            g(x)
        \end{pmatrix}
    \end{equation}
    solves
    \begin{equation}\label{eqn:eigenvalue_Gc}
        \calH_c^* \bmG_c = -2i \lambda \gamma^{-1} \bmG_c.
    \end{equation}
\end{lemma} 
\begin{proof}
It suffices to verify the generalized eigenvalue equation \eqref{eqn:eigenvalue_Gc} directly. The first and second components of the vector $\calH_c^* \bmG_c$ simplify as
\begin{equation*}
\begin{split}
    -L_0 (e^{-2ic \lambda x} g(x))  &= e^{-2ic\lambda x}(-4 c^2\lambda^2 g - 4 i c \lambda g' -L_0 g )=-e^{-2ic\lambda x}(4 c^2\lambda^2 g + 4 i c \lambda g' +4 \lambda^2 g ) \\
    &\quad = -2i \lambda \gamma^{-1} e^{{-2ic\lambda x}}(2 c \gamma g'(x) - 2i \lambda \gamma(1+c^2)g(x)) \equiv  -2i \lambda \gamma^{-1} (\bmG_c(x))_1,
\end{split}
\end{equation*}
and
\begin{equation*}
    \begin{split}
        e^{-2ic\lambda x}\big(2c \gamma g' - 2 i\lambda \gamma (1+c^2)g\big) - 2c \gamma \px(e^{-2ic\lambda x}g) = - 2i\lambda \gamma (1-c^2)g \equiv  -2i \lambda \gamma^{-1} (\bmG_c(x))_2,
    \end{split}
\end{equation*}
where we used $\gamma(1-c^2)= \gamma^{-1}$ in the last line.
\end{proof}
We recall the Taylor expansion of the nonlinear term 
\begin{equation*}
    \calN = \begin{pmatrix}
        0 \\ -\frac{1}{2}W'''(H)(z_1 + m_1)^2 + O((z_1+m_1)^3)
    \end{pmatrix}.
\end{equation*}
Since $a := a_+ = \bar{a}_-$ and $m_1(t,x) = \big(a(t) e^{-ic\lambda x} + \bara(t)e^{ic\lambda x}\big)\cdot Y(x)$, the contribution of the pure quadratic forcing term due to the internal mode in the nonlinearity is given by 
\begin{equation*}
\begin{pmatrix}
        0 \\ -\frac{1}{2}W'''(H)m_1^2
    \end{pmatrix} = \begin{pmatrix}
        0 \\ - \frac{1}{2}W'''(H)Y^2(a^2 e^{-2ic\lambda x} + 2 |a|^2 + \bara^2 e^{2ic \lambda x})
    \end{pmatrix} =: a^2 \bmF_{++} + 2 |a|^2 \bmF_{+-} + \bara^2 \bmF_{--}.
\end{equation*}
Thus, we have the following relation to the Fermi golden rule assumption \eqref{eq:fgr} through the identity
\begin{equation}\label{eqn:definition_FGR_Gamma}
    \langle \bmF_{++},\bmG_{c}\rangle = \int -\frac{1}{2}W'''(H) e^{-2ic\lambda x}Y(x)^2\cdot \overline{e^{-2ic\lambda x}g(x)} \ud x = \int - \frac{1}{2}W'''(H) Y^2 g = \Gamma.
\end{equation}
Inspired by \cite{kowalczyk2022kink} and \cite{martel2024asymptotic}, we introduce a virial functional to estimate the internal mode in the following lemma. 
\begin{lemma}
Define the functional
\begin{equation}\label{eqn:definition_calJ_virial}
    \calJ := \langle \bmz,a^2\bmG_c \chi_A \rangle - \left \langle \frac{i\gamma}{\lambda}\left( |a|^2 \bmF_{+-} + \frac{1}{4} \bara^2 \bmF_{--}\right) , a^2 \bmG_c \chi_A  \right \rangle.
\end{equation}
We have
\begin{equation}\label{eqn:definition_dot_calJ}
\dot{\calJ} = \Gamma |a|^4 + \sum_{k=1}^6 \rmJ_k,
\end{equation}
where
\begin{align*}
&\rmJ_1 :=\left \langle \tPci\big[ \calN + \bmOmega(\bmz) + \bmOmega(\bmm) - \calR\big] + \calS - \bmOmega(\bmH) - a^2 \bmF_{++} ,a^2 \bmG_c(1-\chi_A) \right\rangle ,\\
&\rmJ_2 := \left\langle \calN - a^2 \bmF_{++} -2|a|^2 \bmF_{+-}-\bara^2 \bmF_{--},a^2 \bmG_{c}\chi_A \right\rangle,\\
&\rmJ_3 := \left\langle \bmOmega(\bmz) + \bmOmega(\bmm) ,a^2 \bmG_{c}\chi_A \right\rangle,\\
&\rmJ_4 := -\left\langle 2\chi_A'(\px z_1 + c \gamma z_2) + \chi_A'' z_1 ,a^2  e^{-2ic\lambda x}g(x)\right\rangle,\\
&\rmJ_5 := - \langle \calR,a^2 \bmG_c \chi_A\rangle +  \dot{c} \left\langle \bmz, a^2 (\pc \bmG_{c}) \chi_A\right\rangle + \frac{\dot{c}}{\lambda}\left\langle |a|^2 \pc(\gamma \bmF_{+-}) + \frac14 \bara^2 \pc(\gamma \bmF_{--}),ia^2 \bmG_c \chi_A \right\rangle\\
&\qquad \qquad + \frac{\dot{c}\gamma}{\lambda} \left\langle |a|^2 \bmF_{+-} + \frac14 \bara^2 \bmF_{--},i a^2 (\pc \bmG_{c})\chi_A\right\rangle,\\
&\rmJ_6 := 2 \langle \bmz, a(\dot{a}-i\lambda \gamma^{-1} a) \bmG_{c} \chi_A \rangle + \frac{\gamma}{\lambda} \left\langle (\pt |a|^2) \bmF_{+-} + \frac12 \bara (\dot{\bara} + i \lambda \gamma^{-1} \bara) \bmF_{--},ia^2 \bmG_c \chi_A \right \rangle \\
&\qquad \qquad + \frac{2\gamma}{\lambda} \left\langle |a|^2 \bmF_{+-} + \frac14 \bara^2 \bmF_{--}, ia (\dot{a} -i\lambda \gamma^{-1} a) \bmG_c \chi_A \right \rangle.
\end{align*}

\end{lemma}
\begin{proof}
By the product rule, we find that 
\begin{align*}
\frac{\ud}{\ud t}\langle \bmz,a^2\bmG_c \chi_A \rangle = \langle \pt \bmz, a^2 \bmG_{c} \chi_A\rangle + 2\langle \bmz, a \dot{a}\bmG_{c}\chi_A \rangle + \dot{c}\langle \bmz, a^2 (\pc \bmG_c)\chi_A\rangle.
\end{align*}
Recall from Proposition~\ref{prop:equations_for_bmz} that $\pt{\bmz} = \calH_c \bmz + \bmOmega(\bmH) + \big(\bmOmega(\bmz + \bmm) + \calN - \calR \big)^\perp -\calS$, where we used the notation $\bmh^\perp = \bmh - \tPci\bmh$. Using \eqref{eqn:eigenvalue_Gc}, the pairing of the linear term gives
\begin{align*}
    \langle \calH_c \bmz,a^2 \bmG_c \chi_A\rangle = -2 \langle \bmz,i\lambda \gamma^{-1} a^2 \bmG_c \chi_A\rangle + \langle \bmz,a^2 [\calH_c^*,\chi_A]\bmG_c\rangle.
\end{align*}
By direct computation, the commutator term is $\rmJ_4 = \langle \bmz,a^2 [\calH_c^*,\chi_A]\bmG_c\rangle$. Moreover, Lemma~\ref{lemma:generalized_eigenfunction} implies the orthogonality $\langle \Phi_c^\pm,\bmG_c\rangle =\langle \Phi_c^{(0)},\bmG_c\rangle =\langle \Phi_c^{(1)},\bmG_c\rangle =0$, so we have 
\begin{align*}
    \left \langle \tPci\big[ \calN + \bmOmega(\bmz) + \bmOmega(\bmm) - \calR\big] + \calS - \bmOmega(\bmH),a^2\bmG_c\right\rangle =0,
\end{align*}
which implies
\begin{align*}
&-\left \langle \tPci\big[ \calN + \bmOmega(\bmz) + \bmOmega(\bmm) - \calR\big] + \calS - \bmOmega(\bmH),a^2\bmG_c\chi_A\right\rangle \\
&\qquad = \left \langle \tPci\big[ \calN + \bmOmega(\bmz) + \bmOmega(\bmm) - \calR\big] + \calS - \bmOmega(\bmH),a^2\bmG_c(1-\chi_A)\right\rangle.
\end{align*}
From \eqref{eqn:definition_FGR_Gamma}, 
\begin{align*}
    \langle a^2 \bmF_{++},a^2 \bmG_c \chi_A\rangle = \Gamma |a|^4 - \langle a^2 \bmF_{++},a^2 \bmG_c (1-\chi_A)\rangle.
\end{align*}
Thus, after decomposing the nonlinear term $\calN$ as before, we obtain that 
\begin{align*}
\frac{\ud}{\ud t}\langle \bmz,a^2\bmG_c \chi_A \rangle &= \Gamma |a|^4 +  \sum_{k=1}^4 \rmJ_k + \langle 2|a|^2\bmF_{+-} + \bara^2 \bmF_{--},a^2 \bmG_c \chi_A\rangle - \langle \calR,a^2\bmG_c \chi_A \rangle\\
&\quad+ \dot{c} \langle \bmz,a^2 \pc(\bmG_c)\chi_A\rangle + 2 \langle \bmz, a(\dot{a}-i\lambda \gamma^{-1}a)\bmG_c \chi_A\rangle.
\end{align*}
Next, we differentiate the second term in \eqref{eqn:definition_calJ_virial} with respect to time, which produces a key cancellation against the term $ \langle 2|a|^2\bmF_{+-} + \bara^2 \bmF_{--},a^2 \bmG_c \chi_A\rangle $. By the product rule, 
\begin{align*}
&\frac{\ud}{\ud t}    \left \langle \frac{i\gamma}{\lambda}\left( |a|^2 \bmF_{+-} + \frac{1}{4} \bara^2 \bmF_{--}\right) , a^2 \bmG_c \chi_A  \right \rangle\\
&= \left \langle \frac{i\gamma}{\lambda}\left( \pt(|a|^2) \bmF_{+-} + \frac{1}{2} \bara\, \dot{\bara} \bmF_{--}\right) , a^2 \bmG_c \chi_A  \right \rangle + \dot{c}\left \langle \frac{i}{\lambda}\left( |a|^2 \pc(\gamma\bmF_{+-}) + \frac{1}{4} \bara^2 (\pc \gamma\bmF_{--})\right) , a^2 \bmG_c \chi_A  \right \rangle\\
&\quad + 2 \left \langle \frac{i\gamma}{\lambda}\left( |a|^2 \bmF_{+-} + \frac{1}{4} \bara^2 \bmF_{--}\right) , a \, \dot{a} \bmG_c \chi_A  \right \rangle + \dot{c}  \left \langle \frac{i\gamma}{\lambda}\left( |a|^2 \bmF_{+-} + \frac{1}{4} \bara^2 \bmF_{--}\right) , a^2 (\pc\bmG_c) \chi_A  \right \rangle.
\end{align*}
Using the simple identity $\dot{a} = i\lambda \gamma^{-1} a + (\dot{a}-i\lambda \gamma^{-1}a)$, we have 
\begin{align*}
\left \langle \frac{i\gamma}{2\lambda} \bara\, \dot{\bara} \bmF_{--}, a^2 \bmG_c \chi_A  \right \rangle = \left \langle \frac{1}{2} \bara^2\bmF_{--}, a^2 \bmG_c \chi_A  \right \rangle + \left \langle \frac{i\gamma}{2\lambda} \bara\, (\dot{\bara}+i\lambda \gamma^{-1}\bara) \bmF_{--}, a^2 \bmG_c \chi_A  \right \rangle ,
\end{align*}
and
\begin{align*}
&2 \left \langle \frac{i\gamma}{\lambda}\left( |a|^2 \bmF_{+-} + \frac{1}{4} \bara^2 \bmF_{--}\right) , a \, \dot{a} \bmG_c \chi_A  \right \rangle\\
&\quad =\langle 2|a|^2 \bmF_{+-} + \frac12 \bara^2 \bmF_{--},a^2 \bmG_c \chi_A\rangle  + 2 \left \langle \frac{i\gamma}{\lambda}\left( |a|^2 \bmF_{+-} + \frac{1}{4} \bara^2 \bmF_{--}\right) , a(\dot{a}-i\lambda \gamma^{-1}a)\bmG_c \chi_A  \right \rangle.
\end{align*}
After grouping $\dot{c}$--terms into $\rmJ_5$ and $\dot{a}$--terms into $\rmJ_6$, we then conclude the identity \eqref{eqn:definition_dot_calJ} which finishes the lemma.

\end{proof}

Since $\chi_A$ is compactly supported on $|x| \leq 2A$, we have the following bound
\begin{equation}\label{eqn:bound_chi_A}
    |\chi_A| + A |\chi_A'| + A^2 | \chi_A''| \leq C_\kappa \zeta_{A/2}^2.
\end{equation}
Indeed, by definition of $\chi_A$, we have $|\chi'(x/A)| = A |\chi_A'(x)|$ and $|\chi''(x/A)| = A^2 |\chi_A''(x)|$ which implies the bound
\begin{align*}
|\chi_A| + A |\chi_A'| + A^2 | \chi_A''| \leq \big(1 + \| \chi'\|_{\infty} + \| \chi''\|_{\infty}\big) \mathbbm{1}_{|x|\leq 2A}.    
\end{align*}
On $|x|\leq 2A$, we have $\zeta_{A/2}(x)^2 \simeq \exp( - \frac{10\kappa}{A}|x|) \geq \exp(-20\kappa)$ which implies \eqref{eqn:bound_chi_A}. The following auxiliary lemma will be needed to estimate the error terms. 
\begin{lemma}\label{lemma:auxiliary_lemma_bmGc} There exists a $C \equiv C(c_0)>0$ such that 
\begin{align}
|\bmG_c| + |\px \bmG_{c}| &\leq C \label{eqn:bound_1_bmG_c},\\
| \pc \bmG_c| &\leq C\jx \label{eqn:bound_2_bmG_c},\\
\big\| \px(\bmG_c \chi_A)\big\| &\leq C \sqrt{A} \label{eqn:bound_3_bmG_c},\\
\big\| \px(x \bmG_c \chi_A) \big\| &\leq C A^{\frac32} \label{eqn:bound_4_bmG_c}.
\end{align}
\end{lemma}
\begin{proof}
Since $W''(H)$ is bounded, the equation for $g$ gives $g''=(W''(H)-4\lambda^2)g$ which shows that $\|g''\|_{L^\infty} \leq C \|g\|_{L^\infty} < \infty$. Moreover, Taylor expansion yields
\begin{align*}
    g'(x) = g(x+1) - g(x) - \int_0^1 (1-s) g''(x+s) \,\ud s,
\end{align*}
which implies that $\|g'\|_{L^\infty} \leq C(\|g\|_{L^\infty} + \|g''\|_{L^\infty}) < \infty$. Consequently, we have $g \in W^{2,\infty}(\bbR)$ and $\bmG_c \in W^{1,\infty}(\bbR)$. Moreover, we find that 
\begin{align*}
    \px \bmG_c = -2ic \lambda \bmG_c + e^{-2ic\lambda x}\begin{pmatrix}
        2c \gamma g'' - 2i\lambda \gamma (1+c^2)g'\\
        g'
    \end{pmatrix}.
\end{align*}
Hence, the bounds on $g$ together with \eqref{eqn:modulation_comparison} imply \eqref{eqn:bound_1_bmG_c}. We also compute that 
\begin{align*}
    \pc \bmG_c = -2 i \lambda x \bmG_c + e^{-2ic\lambda x}\begin{pmatrix}
        2 \gamma^3 g' - 2i \lambda c \gamma^3(3-c^2)g\\ 0
    \end{pmatrix},
\end{align*}
which implies \eqref{eqn:bound_2_bmG_c}. The remaining bounds \eqref{eqn:bound_3_bmG_c} and \eqref{eqn:bound_4_bmG_c} now follow from the elementary bounds
\begin{align*}
    \| \chi_A \| \leq C \sqrt{A}, \quad \| x \chi_A \| \leq C A^{\frac32}, \quad \| \chi_A'\| \leq CA^{-\frac12}, \quad \| x \chi_A'\| \leq C \sqrt{A}.
\end{align*}

\end{proof}

In the following lemma, we estimate the error terms.
\begin{lemma}\label{lemma:error_bounds_rmJ}
Assuming \eqref{main_scaling}, it holds for any $A\geq 20$ that
\begin{align*}
&|\rmJ_1| \leq C e^{-\kappa A/5} \big(|a|^4 + \| \rho^4 z_1 \|^2\big),\\
&|\rmJ_2| \leq C \delta \big(|a|^4 + \| \rho^4 z_1 \|^2\big) + C \delta^2 \| \zeta_{A/2} z_1\|^2 ,\\
&|\rmJ_3| \leq C A^{\frac32}\delta \big(|a|^4 + \| \rho^4 z_1 \|^2\big),\\
&|\rmJ_4| \leq C A^{-\frac12} \big(|a|^4 + \| \zeta_{A/2}(\px z_1 + c \gamma z_2) \|^2 + A^{-2} \| \zeta_{A/2} z_1\|^2\big),\\
&|\rmJ_5| \leq CA^{\frac32}\delta \big(|a|^4 + \| \rho^4 z_1 \|^2\big),\\
&|\rmJ_6| \leq C \sqrt{A} \delta \big(|a|^4 + \| \rho^4 z_1 \|^2 +\| \zeta_{A/2} z_1 \|^2 + \| \zeta_{A/2} z_2 \|^2\big).
\end{align*}
\end{lemma}
\begin{proof}
\underline{Estimate of $\rmJ_1$:}
First, we note that 
\begin{align*}
    \int \rho^8 |\calN | \leq C (|a|^2 + \|\rho^4 z_1 \|^2).
\end{align*}
Moreover, replicating the proof of Lemma~\ref{lemma:bound_on_beta_terms}, we obtain 
\begin{align*}
    \left| \tPci [\calN + \bmOmega(\bmz) + \bmOmega(\bmm) - \calR] + \calS - \bmOmega(\bmH) - a^2 \bmF_{++}\right| \leq C \rho^{8} (|a|^2 + \| \rho^4 z_1 \|^2),
\end{align*}
where we used the spatial localization of each term on the left-hand side. Since $1-\chi_A$ is supported on $|x|\geq A$, we have 
\begin{align*}
    |J_1| \leq C|a|^2(|a|^2 + \| \rho^4 z_1 \|^2) \int_{|x|\geq A} \rho^{8} \leq C e^{-\kappa A/5} (|a|^4 + \|\rho^4 z_1 \|^2).
\end{align*}

\underline{Estimate of $\rmJ_2$:} For $\rmJ_2$, by Taylor expansion and \eqref{eqn:stability_smallness}, we have 
\begin{align*}
    \left| \calN - a^2 \bmF_{++} - 2 |a|^2 \bmF_{+-} - \bara^2 \bmF_{--} \right| \leq C(|z_1|^2 + |z_1m_1| + |z_1 + m_1|^3) \leq C(|z_1|^2 + \rho^8|a| |z_1| + |a|^3 \rho^{24}).
\end{align*}
Hence, by \eqref{eqn:bound_chi_A}, we have 
\begin{align*}
    |J_2| \leq C |a|^2 \int \chi_A |z_1|^2 +  C|a|^3 \int \rho^8 |z_1| + C|a|^5 \leq C \delta^2 \| \zeta_{A/2} z_1\|^2 + C \delta (|a|^4 + \| \rho^4 z_1\|^2).    
\end{align*}

\underline{Estimate of $\rmJ_3$:} Integration by parts yields
\begin{align*}
\langle \bmOmega(\bmz),a^2 \bmG_c \chi_A\rangle = -(\dot{y}-c)\gamma \langle \bmz,a^2 \px(\bmG_c \chi_A)\rangle + \dot{c} c \gamma^2 \langle \bmz,a^2\px(x \bmG_c \chi_A)\rangle - \dot{c} \gamma \langle z_1,a^2\px((\bmG_c)_2 \chi_A)\rangle.
\end{align*}
Hence, Lemma~\ref{lemma:auxiliary_lemma_bmGc} implies that 
\begin{align*}
| \langle \bmOmega(\bmz),a^2 \bmG_c \chi_A\rangle | \leq C A^{\frac32} \delta |a|^2(|\dot{y}-c| + |\dot{c}|) \leq C A^{\frac32} \delta (|a|^4+\|\rho^4 z_1\|^2).
\end{align*}
Since $\bmm$ is spatially localized, we obtain the bound
\begin{align*}
| \langle \bmOmega(\bmm),a^2 \bmG_c \chi_A\rangle |  \leq C \delta (|a|^4+\|\rho^4 z_1\|^2).
\end{align*}
Therefore, we conclude the desired estimate for $\rmJ_3$.

\underline{Estimate of $\rmJ_4$:} Using \eqref{eqn:bound_chi_A}, the boundedness of $g$, and the compact support of the cutoffs, we obtain 
\begin{align*}
|J_4| &\leq C A^{-\frac12}|a|^2 \|\zeta_{A/2}(\px z_1 + c \gamma z_2)\| + C A^{-\frac32} |a|^2 \| \zeta_{A/2} z_1\|\\
&\leq CA^{-\frac12}(|a|^4 + \| \zeta_{A/2} (\px z_1 + c \gamma z_2)\|^2 + A^{-2} \| \zeta_{A/2} z_1 \|^2).
\end{align*}
Since $A^{-1} \leq A^{-\frac12}$, this proves the desired estimate for $\rmJ_4$.

\underline{Estimate of $\rmJ_5$:} For $\rmJ_5$, the definition of $\calR$ from \eqref{eqn:definition_calR}, the exponential decay of $\pc \Phi_c^\pm$ from Lemma~\ref{lemma:preparatory_bounds_on_eigfn}, and the bound \eqref{eqn:estimate_modulation_crude} imply that
\begin{align*}
    |\langle \calR, a^2 \bmG_c \chi_A \rangle| \leq C |\dot{c}| |a|^3 \leq C \delta (\|\rho^4 z_1\|^2 + |a|^4).
\end{align*}
Moreover, by the Cauchy-Schwarz inequality, \eqref{eqn:stability_smallness}, and \eqref{eqn:estimate_modulation_crude}, we have
\begin{align*}
|\dot{c} \langle \bmz, a^2 (\pc \bmG_c) \chi_A \rangle| \leq C |\dot{c}|\, \|\bmz\|_{H^1 \times L^2} |a|^2 \| \pc \bmG_c \chi_A \|\leq C A^{\frac32} \delta (|a|^4 + \|\rho^4 z_1\|^2),
\end{align*}
where we used $|\pc \bmG_c| \leq \jx$. Furthermore, since both $\bmF_{+-}$ and $\bmF_{--}$, together with their derivatives, decay exponentially, the last two terms in $\rmJ_5$ can be bounded by 
\[C |\dot{c}| |a|^4 \leq C \delta^2 |a|^4.\]
Combining these bounds imply the desired estimate for $\rmJ_5$.

\underline{Estimate of $\rmJ_6$:} For $\rmJ_6$, we recall from Proposition~\ref{prop:equations_for_bmz} that
\begin{align*}
    \dot{a} - i \lambda \gamma^{-1} a = \frac{1}{2i\lambda } \left[\int \big(\calN + \bmOmega(\bmH +\bmz + \bmm) - \calR) \cdot \Psi_c^{-}  + \dot{c} \int \bmz \cdot \pc \Psi_c^-\right ].
\end{align*}
After integrating by parts on the $\bmOmega(\bmz)$ term, we then use \eqref{eqn:stability_smallness} and \eqref{eqn:estimate_modulation_crude} to obtain the bound 
\begin{equation}\label{eqn:estimate_on_dot_a}
|\dot{a} - i \lambda \gamma^{-1} a| \leq C \int |\calN \cdot \Psi_c^{-}|+ C(|\dot{c}| + |\dot{y}-c|)(\|z_1\| + \|z_2 \| + |a|)\leq C(|a|^2 + \|\rho^4 z_1\|^2).
\end{equation}
The preceding bound, \eqref{eqn:bound_chi_A}, and Lemma~\ref{lemma:auxiliary_lemma_bmGc} hence imply that 
\begin{align*}
\left|2 \langle \bmz, a(\dot{a}-i\lambda \gamma^{-1} a) \bmG_{c} \chi_A \rangle \right|  &\leq C \sqrt{A} |a|(|a|^2 + \|\rho^4 z_1 \|^2)(\| \zeta_{A/2} z_1\| + \| \zeta_{A/2} z_2 \|)\\
&\leq C \sqrt{A} \delta (|a|^4 + \| \rho^4 z_1 \|^2 +\| \zeta_{A/2} z_1 \|^2 + \| \zeta_{A/2} z_2 \|^2).
\end{align*}
Similarly, we have 
\begin{align*}
    |\pt|a|^2| + |\bara(\dot{\bara} + i \lambda \gamma^{-1} \bara)| \leq C |a|(|a|^2 + \|\rho^4 z_1\|^2).
\end{align*}
Therefore, the remaining two terms in $\rmJ_6$ can be bounded by 
\begin{align*}
    C|a|^3(|a|^2 + \|\rho^4 z_1\|^2) \leq C \delta (|a|^4 + \|\rho^4 z_1\|^2),
\end{align*}
using the exponential decay of the forcing terms $\bmF_{+-}$ and $\bmF_{--}$. This completes the estimate for $\rmJ_6$ and the proof of the lemma.
\end{proof}

We are now in the position to finish the proof of Proposition~\ref{prop:internal_mode_estimate}.
\begin{proof}[Proof of Proposition~\ref{prop:internal_mode_estimate}] We first apply Proposition~\ref{prop:high_energy_estimate} (with $K \leftarrow \frac{A}{2}$) to obtain 
\begin{align*}
    \int_0^T (\| \zeta_{A/2} z_1\|^2 + \| \zeta_{A/2} z_2 \|^2)  \leq C A^2 \bigg(A \delta^2 + \int_0^T |a|^4 + \int_0^T \|\rho^4 z_1\|^2\bigg),
\end{align*}
and use the transfer estimates \eqref{eqn:HE_transfer_zK2} and \eqref{eqn:HE_key_estimate2} (both with $K \leftarrow \frac{A}{2}$) to find that 
\begin{align*}
    \int_0^T \Big(\|\zeta_{A/2}(\px z_1 + c \gamma z_2)\|^2 + A^{-2} \|\zeta_{A/2} z_1 \|^2 \Big) \leq C A\delta ^2+ C \int_0^T |a|^4 + C \int_0^T \| \rho^4 z_1 \|^2.
\end{align*}
The error bounds from Lemma~\ref{lemma:error_bounds_rmJ} and the smallness condition $\delta \leq A^{-4}$ now imply that 
\begin{align*}
    \sum_{k=1}^6 \int_0^T  |\rmJ_k| &\leq  C\big(e^{-\frac{\kappa A}{5}} + A^2 \delta^2 + A^{\frac32} \delta + A^{-\frac12} + A^{\frac52}\delta\big) \bigg(A \delta^2 + \int_0^T |a|^4 + \int_0^T \|\rho^4 z_1\|^2\bigg)\\
    &\leq C A \delta^2 + C A^{-\frac12} \bigg( \int_0^T |a|^4 + \int_0^T \| \rho^4 z_1 \|^2\bigg).
\end{align*}
Moreover, the boundedness of $\bmG_c$, the spatial localization of $\bmF_{+-}$ and $\bmF_{--}$, and \eqref{eqn:estimate_modulation_crude} imply that $|\calJ| \leq C A \delta^2$. Therefore, integrating \eqref{eqn:definition_dot_calJ} from $[0,T]$ yields
\begin{align*}
    \Gamma \int_0^T |a|^4 \leq C A \delta^2  + C A^{-\frac12} \bigg( \int_0^T |a|^4 + \int_0^T \| \rho^4 z_1 \|^2\bigg).
\end{align*}
Taking $A$ sufficiently large so that $CA^{-\frac12} \leq \frac12\Gamma$ and choosing $\delta$ sufficiently small as before concludes Proposition~\ref{prop:internal_mode_estimate}.
\end{proof}

\section{Darboux factorization and transformed equations}
In this section we construct matrix Darboux transformations for the linearized operator $\calH_c$ and derive the equations satisfied by the transformed radiation component.  We also compute the commutator terms arising in the derivation of the transformed equations, and recall the standard regularization $X_\eps = (1-\eps\px^2)^{-1}$ used in virial estimates. These preliminaries are needed to establish the coercivity and virial estimates in the next two sections.
\subsection{Matrix Darboux factorization}
In this subsection, we recall the following notation from our introduction for the scalar Darboux transforms applied to the static linearized operator $L_0 = -\px^2 + W''(H)$. We use the notation $Y_0 \equiv H'$ and 
\begin{equation}\label{eqn:def_U0}
    U_0 := Y_0 \cdot \px \cdot Y_0^{-1} = \px - V_0,\qquad V_0 := \frac{Y_0'}{Y_0}.
\end{equation}
Then, we let \begin{equation}
    Y_1 = U_0 Y,
\end{equation}
and
\begin{equation}\label{eqn:def_U1}
    U_1 := Y_1 \cdot \px \cdot Y_1^{-1} = \px - V_1, \qquad V_1 := \frac{Y_1'}{Y_1}.
\end{equation}
We define the potentials $P_j$ through the following scalar Schr\"odinger operators
\begin{align}
    &L_0 := -\px^2 + W''(H) =: -\px^2 + P_0,\\
    &L_1 := -\px^2 + 2 \left(\frac{Y_0'}{Y_0}\right)^2 - \frac{Y_0''}{Y_0} =: - \px^2 + P_1 ,\\
    & L_2 := -\px^2 + \lambda^2 +2\left(\frac{Y_1'}{Y_1}\right)^2 - \frac{Y_1''}{Y_1}=: -\px^2 + P_2.\label{eqn:def_P2}
\end{align}
Moreover, we note that 
\begin{align}
    &\lim_{x \rightarrow \pm \infty}P_0 = \lim_{x \rightarrow \pm \infty}P_1 = \lim_{x \rightarrow \pm \infty}P_2 = \omega^2_\pm,\\
    & L_1 Y_1 = \lambda^2 Y_1, \qquad L_2 U_1 Y_1 = 0.
\end{align}
In the next lemma, we recall some standard identities for the above objects.
\begin{lemma}\label{lemma:scalar_darboux_factorization}
The following identities hold.
\begin{enumerate}
    \item Conjugation identities
\begin{equation}
    L_0 = U_0^* U_0,\quad L_1 = U_0 U_0^* = U_1^* U_1 + \lambda^2, \quad L_2 = U_1 U_1^* + \lambda^2,
\end{equation}
    \begin{equation}\label{eqn:scalar_factorization}
        U_0 L_0 = L_1 U_0,\qquad         U_1 L_1 = L_2 U_1, \qquad U_1 U_0 L_0 = L_2 U_1 U_0 .
    \end{equation}
    \item Commutator identities
    \begin{equation}\label{eqn:commutator_U0_px}
       [U_0,\px] =  U_0 \px - \px U_0 = \frac{1}{2}(L_0 - L_1) = \frac{1}{2}(P_0 - P_1) = V_0',
    \end{equation}
\begin{equation}\label{eqn:commutator_U1_px}
    [U_1,\px] = U_1 \px - \px U_1 =  \frac{1}{2}(L_1 - L_2) = \frac{1}{2}(P_1 - P_2) = V_1'.
\end{equation}
\end{enumerate}
\end{lemma}
We need a matrix version of the Darboux factorization to intertwine the matrix operator $\calH_{c}$. Inspired by the construction in \cite[\S4.4]{kowalczyk2021sufficient}, we define the matrix operators 
\begin{equation}\label{eqn:calD_calH_def}
    \calD_{j,c} := \begin{pmatrix}
        U_j & c \gamma \\ - c\gamma L_j & U_j + 2 c^2 \gamma^2\px
    \end{pmatrix}, \qquad \calH_{j,c} = \begin{pmatrix}
        0 & 1 \\ - L_j & 2 c \gamma \px
    \end{pmatrix}.
\end{equation}
In this notation, we note that $\calH_c \equiv \calH_{0,c}$. The following lemma establishes the intertwining property.
\begin{lemma}
We have 
\begin{equation}\label{eqn:intertwining_property}
    \calD_{1,c} \calD_{0,c} \calH_{0,c} = \calH_{2,c} \calD_{1,c} \calD_{0,c}.
\end{equation}
\end{lemma}
\begin{proof}
It suffices to prove the intertwining identity $\calD_{j,c} \calH_{j,c} = \calH_{j+1,c} \calD_{j,c}$ for $j=0,1$. We have
\begin{equation*}
    \calD_{j,c} \calH_{j,c} = \begin{pmatrix}
        - c \gamma L_j & U_j + 2 c^2 \gamma^2 \px \\
        -U_j L_j - 2 c^2 \gamma^2 \px L_j & - c \gamma L_j + 2 c \gamma U_j \px + 4 c^3 \gamma^3 \px^2
    \end{pmatrix},
\end{equation*}
and on the other hand
\begin{equation*}
    \calH_{j+1,c} \calD_{j,c} = \begin{pmatrix}
        - c \gamma L_j & U_j + 2 c^2 \gamma^2 \px \\ -L_{j+1}U_j - 2 c^2 \gamma^2 \px L_j & - c \gamma L_{j+1} + 2 c \gamma \px U_j + 4 c^3 \gamma^3 \px^2
    \end{pmatrix}.
\end{equation*}
The first row of $\calD_{j,c} \calH_{j,c}$ and $\calH_{j+1,c} \calD_{j,c}$ obviously agree, while the bottom row is also equal using $U_jL_j=L_{j+1}U_j$ and $2[U_j,\px] = L_j - L_{j+1}$ from Lemma~\ref{lemma:scalar_darboux_factorization}.
\end{proof}

\subsection{Transformed variables}

We introduce the following short-hand notation for the iterated Darboux factorization
\begin{equation}\label{eqn:definition_sfABCD}
    \begin{pmatrix}
        \sfA & \sfB \\ \sfC & \sfD
    \end{pmatrix} := \calD_{1,c} \calD_{0,c},
\end{equation}
where
\begin{equation*}
\begin{split}
    &\sfA := U_1 U_0 - c^2 \gamma^2 L_0, \\
    &\sfB := c\gamma (U_1 + U_0) + 2c^3 \gamma^3 \px,\\
    &\sfC := -c \gamma L_1 U_0 - c \gamma U_1 L_0 - 2c^3 \gamma^3 \px L_0,\\
    &\sfD :=-c^2 \gamma^2 L_1 +U_1 U_0+ 2c^2\gamma^2(U_1\px+\px U_0) + 4 c^4 \gamma^4 \px^2. 
\end{split}
\end{equation*}
Using the notation from \eqref{eqn:def_U0}--\eqref{eqn:def_P2}, we rewrite the above expressions as the following differential operators
\begin{equation}\label{eqn:rewritten_sfABCD}
\left\{\begin{aligned}
\sfA &= \gamma^2 \px^2 -(V_0 + V_1) \px + V_0 V_1 - V_0' - c^2 \gamma^2 P_0, \\
\sfB &= 2c \gamma^3 \px - c \gamma (V_0 + V_1), \\
\sfC &= 2 c \gamma^3 \px^3 - c \gamma (V_0 + V_1) \px^2 - 2 c \gamma^3 P_0 \px - c \gamma[2 \gamma^2 P_0' - (V_0 + V_1) P_0] ,\\
\sfD &=  \gamma^4(1+3c^2) \px^2 - \gamma^2(1+c^2)(V_0+V_1) \px + V_0(V_0+V_1) - \gamma^2 P_0.
\end{aligned}
\right.    
\end{equation}
In particular, all coefficients appearing in \eqref{eqn:rewritten_sfABCD} are smooth and bounded. 

For our moving kink setting, we need to commute the Darboux transformation $\calD_{1,c}\calD_{0,c}$ with time differentiation (which produces $\dot{c}\pc$ terms) and with the modulation operator $\bmOmega$ (which is a linear combination of $\px$ and $\Lambda=x\px$) in order to analyse the transformed equations. We record these commutator identities in the following lemma. 
\begin{lemma}\label{lemma:commutator_ABCD}
The following identities hold.
\begin{enumerate}
    \item Commutation with $\pc$:
\begin{align*}
\pc \sfA &= 2c\gamma^4 \px^2 - 2 c \gamma^4 P_0,\\
\pc \sfB &= 2\gamma^5(1+2c^2)\px-\gamma^3(V_0+V_1),\\
\pc \sfC &= 2 \gamma^5(1+2c^2)\px^3 - \gamma^3(V_0+V_1) \px^2 - 2 \gamma^5 (1+2c^2) P_0 \px\\
&\quad - 2 \gamma^5(1+2c^2)P_0' + \gamma^3 (V_0+V_1) P_0,\\
\pc \sfD &= 2c\gamma^6(5+3c^2)\px^2-4c\gamma^4(V_0+V_1)\px-2c\gamma^4P_0.
\end{align*}
\item Commutation with translation: 
\begin{align*}
[\sfA,\px]=&(V_0'+V_1')\px-[V_0(V_0+V_1)]'+\gamma^2P_0',\\
[\sfB,\px]=&c\gamma(V_0'+V_1'),\\
[\sfC,\px]=&c\gamma(V_0'+V_1')\px^2 + 2 c \gamma^3 P_0'\px + c \gamma ( 2 \gamma^2 P_0'' - (V_0'+V_1')P_0 - (V_0+V_1)P_0'),\\
[\sfD,\px]=&\gamma^2(1+c^2)(V_0'+V_1')\px-[V_0(V_0+V_1)]'+\gamma^2P_0'.
\end{align*}
\item Commutation with dilation:
\begin{align*}
[\sfA,\Lambda]=&2\gamma^2\px^2 +[x(V_0'+V_1')-(V_0+V_1)]\px-x[V_0(V_0+V_1)]'+x\gamma^2P_0',\\
[\sfB,\Lambda]=&2c\gamma^3\px+c\gamma x(V_0'+V_1'),\\
[\sfC,\Lambda]=&6 c \gamma^3 \px^3 + c \gamma[x(V_0'+V_1')-2(V_0+V_1)]\px^2 + 2 c \gamma^3(xP_0' - P_0)\px + 2c \gamma^3 x P_0'' \\
&\quad - c \gamma x(V_0' + V_1')P_0 - c \gamma x(V_0 + V_1)P_0',\\
[\sfD,\Lambda]=&2\gamma^4(1+3c^2)\px^2+\gamma^2(1+c^2)[x(V_0'+V_1')-(V_0+V_1)]\px -x[V_0(V_0+V_1)]'+x\gamma^2P_0'.
\end{align*}
\end{enumerate}
\end{lemma}
\begin{proof}
These identities follow directly from \eqref{eqn:rewritten_sfABCD}.
\end{proof}

We now introduce the transformed variables
\begin{equation}\label{eqn:def_g}
    \bmg = (g_1,g_2) := \calD_{1,c} \calD_{0,c} \bmz \in H^{-1} \times H^{-2}.
\end{equation}
Recall from Proposition~\ref{prop:equations_for_bmz} that
\begin{equation*}
    \pt\bmz = \calH_{0,c} \bmz +   \bmOmega(\bmH)+ \bmOmega(\bmz)^\perp + \bmOmega(\bmm)^\perp +\calN^\perp- \calR^\perp - \calS.
\end{equation*}
Since the Darboux factorization annihilates the discrete spectrum of $\calH_c$, we obtain from Lemma~\ref{lemma:annihilation_kernel} 
\begin{equation*}
    \calD_{1,c} \calD_{0,c} \tilde{P}_c^\mathrm{i}=    \calD_{1,c} \calD_{0,c} \bmOmega(\bmH)= \calD_{1,c} \calD_{0,c} \calS = \calD_{1,c} \calD_{0,c} \bmm = \bm0.
\end{equation*}
Moreover, differentiating $\calD_{1,c} \calD_{0,c} \bmm = \bm0$ in time yields the identity
\begin{equation*}
    \dot{c} \pc(\calD_{1,c}\calD_{0,c})\bmm = -\calD_{1,c} \calD_{0,c} \calR.
\end{equation*}
With the above identities and the intertwining property \eqref{eqn:intertwining_property}, we formally arrive at the evolution equation for $\bmg$ by the product rule,
\begin{equation*}
    \pt \bmg = \calH_{2,c} \bmg + \bmOmega(\bmg)+ \dot{c} \pc(\calD_{1,c}\calD_{0,c})(\bmz + \bmm) + [\calD_{1,c}\calD_{0,c},\bmOmega](\bmz +\bmm) + \begin{pmatrix}
        \sfB N_2 \\ \sfD N_2
    \end{pmatrix}.
\end{equation*}

For the virial computations, we express the transformed evolution equations for $\bmg =(g_1,g_2)$ componentwise as the following system
\begin{equation}\label{eqn:transformed_equations_g}
\left\{
\begin{aligned}
&\dot{g}_1 = g_2 + \Omega_1(\bmg)+ \eta_1 + \sfB N_2\\
&\dot{g}_2 = -L_2 g_1 + 2 c \gamma \px g_2 + \Omega_2(\bmg) + \eta_2 +  \sfD N_2
\end{aligned}\right.,
\end{equation}
where $\eta_1$ and $\eta_2$ capture the commutator terms given respectively by 
\begin{equation}\label{eqn:definition_eta1}
\begin{split}
\eta_1 :=& \dot{c}\left\{\pc \sfA - c \gamma^2 [\sfA,\Lambda] + \gamma \sfB \px  \right\}(z_1+m_1)+ \dot{c}\left\{\pc \sfB - c \gamma^2 [\sfB,\Lambda]\right\}(z_2+m_2)\\
&+ (\dot{y}-c)\gamma\big([\sfA,\px](z_1+m_1) + [\sfB,\px](z_2+m_2) \big) ,
\end{split}
\end{equation}
and
\begin{equation}\label{eqn:definition_eta2}
\begin{split}
\eta_2 :=& \dot{c} \left\{\pc \sfC - c \gamma^2[\sfC,\Lambda] + \gamma (\sfD \px - \px \sfA) \right\} (z_1+m_1)\\
&+ \dot{c} \left\{\pc \sfD - c \gamma^2[\sfD,\Lambda] - \gamma \px \sfB\right\}(z_2+m_2)\\
&+(\dot{y}-c)\gamma \big([\sfC,\px](z_1+m_1) + [\sfD,\px](z_2+m_2) \big) .
\end{split}
\end{equation}
Notice that the system for transformed variables \eqref{eqn:transformed_equations_g} now features the repulsive potential $P_2$. To perform virial estimates on $\bmg$, we first need to regularize the variables to the class $H^1 \times L^2$. 
\subsection{Regularized transformed variables}

For fixed $0 < \eps \ll 1$, we set
\begin{equation}
    X_\varepsilon= (1-\varepsilon \px^2)^{-1}, \quad \bmX_{\eps} = \begin{pmatrix}
        X_\eps & 0 \\ 0 & X_\eps
    \end{pmatrix},
\end{equation}
and define the regularized transformed variables 
\begin{equation}
    \bmf = (f_1,f_2) := \bmX_\eps \bmg = (X_\eps g_1, X_\eps g_2) \in H^1 \times L^2.
\end{equation}
Although $X_\eps$ commutes with $\pt$ and $\px$, it does not commute with the dilation operator $\Lambda = x\px$ and with multiplication by the potential $P_2$ or $x$. Using the well-known commutator identity $[A^{-1}, B] = -A^{-1} [A, B] A^{-1}$, we compute that 
\begin{equation}
[X_\eps,x] = 2 \eps X_\eps^2 \px, \qquad [X_\eps,\Lambda] = 2\eps X_\eps^2 \px^2. 
\end{equation}
Thus, we obtain from \eqref{eqn:transformed_equations_g} the equations for the regularized transformed variables 
\begin{equation}\label{eqn:equations_for_f12}
\left\{\begin{aligned}
&\dot{f}_1 = f_2 + \Omega_1(\bmf)+ \xi_1 + X_\eps\sfB N_2\\
&\dot{f}_2 = -L_2 f_1 + 2 c \gamma \px f_2 + \Omega_2(\bmf) + \xi_2 + \sfT  + X_\eps\sfD N_2 
\end{aligned}\right.,
\end{equation}
with the error terms
\begin{equation}\label{eqn:definition_xi_T}
\begin{aligned}
&\xi_1 := X_\eps \eta_1 - 2 \eps \dot{c} c \gamma^2 X_\eps \px^2 f_1,\\
&\xi_2 := X_\eps \eta_2 - 2 \eps \dot{c} c \gamma^2 X_\eps \px^2 f_2, \\
&\sfT:= - \eps X_\eps( P_2' \px f_1 + \px (P_2' f_1)).
\end{aligned}
\end{equation}
In the next sections, we perform virial estimates on the regularized transformed variables $(f_1,f_2)$.
\section{Coercivity of the transformed components}
In this section we establish recovery estimates needed to transfer control of the transformed variables $\bmf=(f_1,f_2)$ back to the original radiation component $\bmz=(z_1,z_2)$. It is crucial here to use the orthogonality relations satisfied by $\bmz$, and ensure that the regularities of $\bmz$ and $\bmf$ are compatible. The goal is to prove the following proposition.
\begin{proposition}\label{prop:coercivity_estimate} 
Assume \eqref{main_scaling}. We have 
\begin{equation}
    \int_0^T \big\| \rho^4 z_1 \big\|^2 \ud t  \lesssim \int_0^T \int \big(f_1^2 + (\px f_1)^2 + f_2^2\big)\rho^2.
\end{equation}
\end{proposition}
We first show that the matrix Darboux transformation annihilates the discrete spectrum of $\calH_c$. 
\begin{lemma}\label{lemma:annihilation_kernel}
It holds that 
\begin{equation}\label{eqn:kernel_darboux_transform}
\ker(\calD_{1,c} \calD_{0,c} ) = \mathrm{span}\left\{\Phi_{c}^{(0)}, \Phi_{c}^{(1)},\Phi_{c}^+,\Phi_{c}^- \right\}.
\end{equation}
\end{lemma}
\begin{proof}
Since $\calD_{0,c}$ and $\calD_{1,c}$ are second-order differential operators, the kernel of the composition $\calD_{1,c}\calD_{0,c}$ has dimension at most four. Hence, it suffices to prove the reverse inclusion of \eqref{eqn:kernel_darboux_transform}. From $U_0 H' = L_0 H' = 0$, $U_0(xH')=H'$ and $L_0(xH')=-2H''$, we compute that
\begin{equation*}
\calD_{0,c} \begin{pmatrix}
    H' \\  0 
\end{pmatrix} = \begin{pmatrix}
    0 \\ 0 
\end{pmatrix},\qquad    \calD_{0,c} \begin{pmatrix}
    -c\gamma xH' \\  H' 
\end{pmatrix} = \begin{pmatrix}
    -c \gamma U_0(xH') + c \gamma H'\\
    c^2 \gamma^2 L_0 (xH') + U_0H' + 2 c^2 \gamma^2 \px H'
\end{pmatrix} = \begin{pmatrix}
    0 \\ 0 
\end{pmatrix}.
\end{equation*}
Using $U_1U_0 Y = 0$ from Lemma~\ref{lemma:scalar_darboux_factorization}, we patiently compute that
\begin{equation*}
    \calD_{1,c} \calD_{0,c} \Phi_{c}^\pm = \calD_{1,c}\left( e^{\mp i c \lambda x} \begin{pmatrix}
        U_0 Y \\ \pm i \lambda \gamma^{-1} U_0 Y
    \end{pmatrix} \right) = e^{\mp i c \lambda x} \begin{pmatrix}
        U_1 U_0 Y \\ \pm i \lambda \gamma^{-1} U_1 U_0 Y
    \end{pmatrix}= \begin{pmatrix}
        0 \\ 0
    \end{pmatrix}.
\end{equation*}
\end{proof}
Next, we need to find right-inverse formulas for the matrix operator $\calD_{j,c}$ to obtain coercivity. These recovery formulas are analogous to the scalar versions used in \cite[\S5]{li2023quadratic}. See also \cite{luhrmann2023sinegordon}, \cite[Lemma~6]{kowalczyk2022jems}, \cite[Lemma~2]{kowalczyk2022kink}. For $j\in\{0,1\}$ we define the operators
\begin{equation}\label{eqn:recovery_def}
    (\calR_{j,c}\bmg)(x) :=  Y_j(x) \int_0^x Y_j(y)^{-1} \calE_{j,c}(x-y) \begin{pmatrix}
        g_1(y) \\ \gamma^{-2} g_2(y) + c \gamma^{-1} (U_j^*g_1)(y)
    \end{pmatrix}\,\ud y,
\end{equation}
where
\begin{equation}\label{eqn:calE0}
    \calE_{0,c}(x) := \exp\left[x \begin{pmatrix}
        0 & -c \gamma \\ 0 & 0 
    \end{pmatrix}\right]=\begin{pmatrix}
        1 & - c \gamma x \\ 0 & 1
    \end{pmatrix},
\end{equation}
and
\begin{equation}\label{eqn:calE1}
    \calE_{1,c}(x) :=\exp\left[x\begin{pmatrix}
        0 & - c \gamma \\ \lambda^2 c \gamma^{-1} & 0 
    \end{pmatrix}\right] =\begin{pmatrix}
        \cos(c \lambda x) & - \frac{\gamma}{\lambda} \sin(c \lambda x) \\ \frac{\lambda}{\gamma}\sin(c \lambda x) & \cos(c \lambda x)
    \end{pmatrix}.
\end{equation}
The next lemma shows that $\calR_{j,c}$ are right-inverses of $\calD_{j,c}$, and its proof explains how \eqref{eqn:recovery_def} is derived.
\begin{lemma}
For $j\in\{0,1\}$, we have 
\begin{equation}\label{eqn:recovery_identity}
\calD_{j,c}\calR_{j,c} = I,
\end{equation}
and
\begin{equation}\label{eqn:recovery_darboux_identity}
\calR_{j,c}\calD_{j,c}\bmh(x) = \bmh(x) - Y_j^{-1}(0)Y_j(x) \calE_{j,c}(x) \bmh(0).
\end{equation}
\end{lemma}
\begin{proof}
We solve the equation $\calD_{j,c}\bmf = \bmg$ with $\bmf(0) =0$ using the variation-of-constants formula. First, by Lemma~\ref{lemma:scalar_darboux_factorization}, we have the identity
\begin{equation*}
    \begin{pmatrix}
        1 & 0 \\ c \gamma^{-1} U_j^* & \gamma^{-2}
    \end{pmatrix}\calD_{j,c} = \begin{pmatrix}
        U_j & c \gamma \\ - c \gamma^{-1}(L_j - U_j^* U_j) & U_j
    \end{pmatrix} = \begin{pmatrix}
        U_j & c \gamma \\ - c \gamma^{-1}\lambda_j^2 & U_j
    \end{pmatrix} ,
\end{equation*}
where we used $L_j - U_j^* U_j = \lambda_j^2$ with
\[\lambda_j^2 := \begin{cases}
    0,& j=0 \\ \lambda^2,&j=1
\end{cases}.\]
Applying this identity to the equation $\calD_{j,c}\bmf = \bmg$ and using  $U_j f_k = Y_j \px (Y_j^{-1}f_k)$, we obtain
\begin{equation*}
\begin{split}
    &\left\{\begin{aligned} U_j f_1 + c \gamma f_2 &= g_1 \\ U_j f_2 - c \gamma^{-1} \lambda_j^2f_1 &= \gamma^{-2} g_2 + c \gamma^{-1} U_j^* g_1 \end{aligned}\right. \\
    &\Longleftrightarrow\left\{\begin{aligned}
    \px (Y_j^{-1}f_1) - Y_j^{-1}(-c\gamma f_2) &= Y_j^{-1}(g_1) \\
    \px (Y_j^{-1}f_2) - Y_j^{-1}(c\gamma^{-1}\lambda_j^2 f_1) &=  Y_j^{-1}(\gamma^{-2}g_2 + c \gamma^{-1} U_j^* g_1)
\end{aligned}\right. .
\end{split}
\end{equation*}

We obtain, for both cases $j\in \{0,1\}$, by using integrating factors \eqref{eqn:calE0}--\eqref{eqn:calE1}, that
\begin{equation*}
    \px \Big( \calE_{j,c}(-x) Y_j^{-1}(x) \bmf(x)\Big) = \calE_{j,c}(-x) Y_j^{-1}(x) \begin{pmatrix}
        g_1(x) \\ \gamma^{-2} g_2(x) + c \gamma^{-1}(U_j^* g_1)(x)
    \end{pmatrix}.
\end{equation*}
Integrating from $0$ to $x$, using $\bmf(0) = 0$, gives the definition \eqref{eqn:recovery_def} and also shows that $\bmf = \calR_{j,c}\bmg$, which implies \eqref{eqn:recovery_identity}. Next we show \eqref{eqn:recovery_darboux_identity}. Since $\calD_{j,c}\bmh = U_j \bmh - \calE_{j,c}'(0) \bmh$ and $U_j \bmh = Y_j \px (Y_j^{-1}\bmh)$, we have 
\[\py \Big(Y_j^{-1}(y)\calE_{j,c}(-y)\bmh(y)\Big) = Y_j^{-1}(y)\calE_{j,c}(-y)[U_j \bmh - \calE_{j,c}'(0)\bmh]=Y_j^{-1}(y) \calE_{j,c}(-y)\calD_{j,c}\bmh.\]
Thus,
\[\calR_{j,c} \calD_{j,c}\bmh (x) = Y_j(x) \calE_{j,c}(x) \int_0^x \py \big( Y_{j}^{-1}(y)\calE_{j,c}(-y) \bmh(y)\big)\,\ud y = \bmh(x) - Y_j^{-1}(0)Y_j(x) \calE_{j,c}(x) \bmh(0) .\]
\end{proof}
Next, we determine an identity for $\calR_{0,c}\calR_{1,c}\calD_{1,c}\calD_{0,c}$ for sufficiently smooth functions $\bmh$.
\begin{lemma}\label{lemma:recover_formula}
We have 
\begin{equation}
\calR_{0,c}\calR_{1,c}\calD_{1,c}\calD_{0,c}\bmh(x) = \bmh(x) - \frakc_0[\bmh] \Phi_c^{(0)}(x) - \frakc_1[\bmh] \Phi_c^{(1)}(x) - \frakc_+[\bmh] \Phi_c^{+}(x) - \frakc_-[\bmh] \Phi_c^{-}(x),
\end{equation}
with the coefficients
\begin{equation}
\frakc_0[\bmh] := \frac{h_1(0)}{\gamma H'(0)},\qquad \frakc_1[\bmh] := \frac{h_2(0)}{\gamma H'(0)},\qquad \frakc_\pm[\bmh] := \frac{(\calD_{0,c}\bmh)_1(0) \mp i \gamma \lambda^{-1} (\calD_{0,c} \bmh )_2 (0)}{2Y'(0)},
\end{equation}
where $Y(x)$ is chosen so that $Y(0) = 0$ and $Y_1(0) = Y'(0) > 0$.
\end{lemma}
\begin{proof}
By iterating \eqref{eqn:recovery_darboux_identity}, we have 
\begin{align*}
\big(\calR_{0,c}\calR_{1,c}\calD_{1,c}\calD_{0,c}\bmh\big)(x) = \bmh(x) - \frac{H'(x)}{H'(0)} \calE_{0,c}(x) \bmh(0) - \bigg(\calR_{0,c} \bigg[ \frac{Y_1}{Y_1(0)}\calE_{1,c} (\calD_{0,c}\bmh) (0)\bigg]\bigg)(x).
\end{align*}
Recalling the definitions \eqref{eqn:eig_eqn_Phi_0}, \eqref{eqn:eig_eqn_Phi_1}, and \eqref{eqn:calE0} simplifies the second term in the above display:
\begin{align*}
\frac{H'(x)}{H'(0)} \calE_{0,c}(x) \bmh(0) = \frac{H'(x)}{H'(0)} \begin{pmatrix}
    h_1(0) - c \gamma x h_2(0) \\ h_2(0)
\end{pmatrix} = \frakc_0[\bmh]\Phi_c^{(0)}(x) + \frakc_1[\bmh]\Phi_{c}^{(1)}(x).
\end{align*}
Next, using Euler's identity $e^{\pm ic\lambda x} = \cos(c \lambda x) \pm i \sin (c \lambda x)$, definition  \eqref{eqn:calE1}, and the following identity
\begin{align*}
    (\calD_{0,c} \Phi_{c}^\pm)(x) = e^{\mp i  c \lambda x} \begin{pmatrix}
        Y_1(x) \\ \pm i \lambda \gamma^{-1} Y_1(x)
    \end{pmatrix},
\end{align*}
we derive that 
\begin{align*}
\frac{Y_1(x)}{Y_1(0)}\calE_{1,c}(x) (\calD_{0,c}\bmh) (0) &=     \frac{Y_1(x)}{Y_1(0)} \begin{pmatrix}
    \cos(c\lambda x) (\calD_{0,c}\bmh)_1(0) - \gamma \lambda^{-1} \sin(c \lambda x) (\calD_{0,c}\bmh)_2(0)\\
    \gamma \lambda^{-1} \sin(c \lambda x)(\calD_{0,c}\bmh)_1(0) + \cos(c \lambda x) (\calD_{0,c}\bmh)_2(0)
\end{pmatrix}\\
&= \frakc_+[\bmh](\calD_{0,c}\Phi_c^+)(x) + \frakc_-[\bmh](\calD_{0,c}\Phi_c^-)(x).
\end{align*}
Since $Y(0)=0$, we have $\Phi_{c}^\pm(0) =\bm0$. Applying \eqref{eqn:recovery_darboux_identity} hence implies that 
\begin{align*}
    \bigg(\calR_{0,c} \bigg[ \frac{Y_1}{Y_1(0)}\calE_{1,c} (\calD_{0,c}\bmh) (0)\bigg]\bigg)(x) = \frakc_+[\bmh]\Phi_c^+(x) + \frakc_-[\bmh] \Phi_{c}^-(x),
\end{align*}
which completes the proof. 

\end{proof}
Consequently, the previous lemmas imply the operator equality
\begin{equation}\label{eqn:iterated_darboux_identity}
    \calD_{1,c}\calD_{0,c}\calR_{0,c}\calR_{1,c}=I, \qquad \tilde{P}_c^{\mathrm{e}} \calR_{0,c}\calR_{1,c}\calD_{1,c}\calD_{0,c}= \tilde{P}_c^{\mathrm{e}}.
\end{equation}
Note that \eqref{eqn:iterated_darboux_identity} holds for all sufficiently smooth functions, and hence for all $L^2$ functions by routine density arguments. In the following lemma, we determine the integral kernel for $\tilde{P}_c^{\mathrm{e}} \calR_{0,c}\calR_{1,c}$.

\begin{lemma}\label{lemma:kernel_of_PceR0cR1c}
Let $\bmh \in H^1 \times L^2$. We have 
\begin{equation}\label{eqn:kernel_of_PceR0cR1c}
    \tilde{P}_c^{\mathrm{e}} \calR_{0,c}\calR_{1,c} \bmh(x) = \tPce \int_{\bbR} \calK(x,y) \bmh(y)\,\ud y,
\end{equation}
where
\begin{equation}\label{eqn:kernel_formula_of_PceR0cR1c}
    \calK(x,y) := \begin{pmatrix}
        0 & 0 \\ - c \gamma^{-1} & 0
    \end{pmatrix}\calK_{1,c}(x,y) + \calK_{0,c}(x,y) \begin{pmatrix}
        0 & 0 \\ - c \gamma^{-1} & 0
    \end{pmatrix} + \int_\bbR \calK_{0,c}(x,s)\calK_{1,c}(s,y)\,\ud s,
\end{equation}
and
\begin{equation}
    \calK_{j,c}(x,y) := \frac{Y_j(x)}{Y_j(y)}\calE_{j,c}(x-y)\begin{pmatrix}
        1+c^2 & 0 \\ -2c\gamma^{-1} V_j(y) & \gamma^{-2}
    \end{pmatrix}  \big[ \mathbbm{1}_{[0,\infty)}(x)\mathbbm{1}_{[0,x]}(y) - \mathbbm{1}_{(-\infty,0)}(x)\mathbbm{1}_{[x,0]}(y)\big],
\end{equation}
for $j \in\{0,1\}$.
\end{lemma}
\begin{proof}
Integrating by parts in \eqref{eqn:recovery_def} using $U_j^*f = -f' - V_j f$ gives
\begin{align*}
    \calR_{j,c}\bmh(x) = \begin{pmatrix}
        0 & 0 \\ - c \gamma^{-1} & 0
    \end{pmatrix} \bmh(x) + \int_\bbR \calK_{j,c}(x,y) \bmh(y) \,\ud y + \frac{c Y_j(x)}{\gamma Y_j(0)}\calE_{j,c}(x) \begin{pmatrix}
        0\\ h_1(0)
    \end{pmatrix},
\end{align*}
for $j \in \{0,1\}$. In particular, we have $(\calR_{1,c}\bmh)_1(0)=0$.
Using $\begin{pmatrix}
        0 & 0 \\ - c \gamma^{-1} & 0
    \end{pmatrix}^2 =0$ and Fubini's theorem, we obtain
\begin{align*}
\calR_{0,c}\calR_{1,c}\bmh(x)= \int_\bbR \calK(x,y) \bmh(y) \,\ud y + c \gamma^{-1}  \calR_{0,c} \left[\frac{Y_1(\cdot)}{Y_1(0)}\calE_{1,c}(\cdot)\begin{pmatrix}
    0 \\ h_1(0)
\end{pmatrix} \right](x).
\end{align*}
The last term in the above identity belongs to $\ker(\calD_{1,c}\calD_{0,c})$ using $\calD_{0,c}\calR_{0,c} = I$ and $\calD_{1,c}(Y_1\calE_{1,c})=0$. Hence, we obtain \eqref{eqn:kernel_of_PceR0cR1c} after applying the projection $\tPce$ to the above identity.
\end{proof}
In the next lemma, we establish coercivity for $\bmz$ under the orthogonality conditions.
\begin{lemma} For any $0 < \eps \leq 1$ and for all $\bmz \in H^1 \times L^2$ satisfying $\bmz = \tilP_c^{\mathrm{e}}\bmz$, it holds that
\begin{equation}\label{eqn:estimate_z_f}
    \| \rho z_1 \|  \lesssim \| \rho f_1 \| + \eps \| \rho \px f_1 \| + \| \rho f_2 \|,
\end{equation}
where $\bmf = \bmX_\eps \calD_{1,c}\calD_{0,c} \bmz$. 
\end{lemma}
\begin{proof} We follow the proof of Lemma~5.3 of \cite{li2023quadratic}. 
Let $\calK(x,y)$ be the integral kernel appearing in Lemma~\ref{lemma:kernel_of_PceR0cR1c}. Integrating twice by parts using the integral kernel expressions yields 
\begin{equation}\label{eqn:proof_kernel_identity1}
\begin{split}
    &\tPce \calR_{0,c}\calR_{1,c} \px^2 \bmf \\
    &= \tPce\left[\begin{pmatrix}
        0 \\ -2c \gamma^{-1} \px f_1
    \end{pmatrix}  + \begin{pmatrix}
        1 + 3c^2 & 0 \\ -c \gamma^{-1}(3+c^2)(V_0+V_1) & \gamma^{-2}
    \end{pmatrix}\bmf + \int_\bbR \py^2 \calK(\cdot,y)\bmf(y)\,\ud y\right].
\end{split}
\end{equation}
From the boundedness of $V_j$, from the linear growth of $\calE_j$, and from the decay of $Y_j$ and its derivatives, we find that
\begin{equation*}
    \sum_{j=0}^2 | \py^j \calK(x,y) | \lesssim \jap{x-y}^2 e^{-\frac{\kappa}{2}|x-y|}.
\end{equation*}
Since $\rho(x) \rho^{-1}(y) \leq \exp(\frac{\kappa}{40}|x-y|)$, Schur's test implies that 
\begin{equation}\label{eqn:proof_kernel_identity2}
    \left \| \rho \int_\bbR \py^j \calK(\cdot,y)\bmf(y) \,\ud y \right\|_{L^2 \times L^2} \lesssim \| \rho f_1 \| + \| \rho f_2 \|, \qquad 0 \leq j \leq 2.
\end{equation}
Moreover, since the discrete eigenfunctions of $\calH_c$ are exponentially localized, it follows from Cauchy-Schwarz inequality that 
\begin{equation*}
    \| \rho \tPce \bmh \|_{L^2 \times L^2} \lesssim\| \rho \bmh \|_{L^2 \times L^2}, \quad \forall \bmh \in L^2 \times L^2.
\end{equation*}
Hence, \eqref{eqn:proof_kernel_identity1} and \eqref{eqn:proof_kernel_identity2} imply that
\begin{equation*}
    \| \rho \tPce \calR_{0,c} \calR_{1,c} \bm f \|_{L^2 \times L^2} \lesssim \| \rho f_1 \| + \| \rho f_2 \|,
\end{equation*}
and
\begin{equation*}
    \| \rho \tPce \calR_{0,c} \calR_{1,c} \px^2 \bm f \|_{L^2 \times L^2} \lesssim \| \rho f_1 \| + \| \rho \px f_1 \| + \| \rho f_2 \|.
\end{equation*}
Since $\bmz = \tilP_c^{\mathrm{e}}\bmz$, the recovery formula \eqref{eqn:iterated_darboux_identity} yields
\begin{equation*}
    \bmz = \tilP_c^{\mathrm{e}}\calR_{0,c} \calR_{1,c} \bmf - \eps \tilP_c^{\mathrm{e}}\calR_{0,c} \calR_{1,c} \px^2 \bmf.
\end{equation*}
Thus, the estimate \eqref{eqn:estimate_z_f} now follows from the preceding bounds.
\end{proof}
Integrating \eqref{eqn:estimate_z_f} in time and using $\rho^4 \leq \rho$ concludes the proof of Proposition \ref{prop:coercivity_estimate}.

\section{Estimate on the transformed components}
In this section we prove Proposition~\ref{prop:low_energy_estimate}, which controls the integrated-in-time local energy norms of the regularized transformed variables. We first collect weighted regularization and transfer estimates and establish bounds for the directional derivative components. We then exploit the crucial repulsivity of the potential $P_2$ in a modified virial argument at the smaller scale $B$, which finishes the last key proposition.
\begin{proposition}\label{prop:low_energy_estimate}
Assuming \eqref{main_scaling}, it holds that 
\begin{equation}
\int_0^T \int \big(f_1^2 + (\px f_1)^2 + f_2^2 \big)\rho^2  \lesssim   A \delta^2 +  \int_0^T |a|^4  + \frac{1}{\sqrt{A}} \int_0^T \big\| \rho^4 z_1 \big\|^2 .
    \end{equation}
\end{proposition}

\subsection{Technical estimates}
We begin with the following technical estimates for the regularization operator. 
\begin{lemma}\label{lemma:commutator_Xeps_lemma}
The following estimates hold.
\begin{enumerate}
\item For any $\varepsilon \in (0, 1)$ and $g \in L^2$,
\begin{equation}\label{eqn:bound_X_eps_basic}
\begin{split}
&\|X_\varepsilon g\|_{L^2} \le \|g\|_{L^2}, \quad \|\partial_x X_\varepsilon g\|_{L^2} \le \varepsilon^{-\frac{1}{2}} \|g\|_{L^2},\\
& \|\partial_x^2 X_\varepsilon g\|_{L^2} \le \varepsilon^{-1} \|g\|_{L^2}, \quad \|\partial_x X_\varepsilon^{\frac12} g\|_{L^2} \le \varepsilon^{-\frac{1}{2}} \|g\|_{L^2}.
\end{split}
\end{equation}
\item There exist $\varepsilon_1 > 0$ and $C > 0$ such that for any $\varepsilon \in (0, \varepsilon_1)$, $K \ge 1$ and $g \in L^2$,
\begin{equation}\label{eqn:bound_X_eps_commute_sech}
\left\| \mathrm{sech}\left(\frac{\kappa}{K} x\right) X_\varepsilon g \right\|_{L^2} \le C \left\| X_\varepsilon \left[ \mathrm{sech}\left(\frac{\kappa}{K} x\right) g \right] \right\|_{L^2},
\end{equation}
and
\begin{equation}\label{eqn:bound_X_eps_conjugate_cosh}
\left\| \cosh\left(\frac{\kappa}{K} x\right) X_\varepsilon \left[ \mathrm{sech}\left(\frac{\kappa}{K} x\right) g \right] \right\|_{L^2} \le C \|X_\varepsilon g\|_{L^2}.
\end{equation}
\item  For any $\varepsilon \in (0, 1)$, $P \in W^{1,\infty}(\bbR)$, and $g \in L^2$,
\begin{equation}\label{eqn:bound_X_eps_multiplication_basic}
\begin{split}
& \| [X_\eps,P]g \|_{L^2} \lesssim \sqrt{\eps} \| P'\|_{L^\infty} \| g \|_{L^2}.
\end{split}
\end{equation}
\end{enumerate}
\end{lemma}
\begin{proof}
See Lemma~4.7 from \cite{kowalczyk2021sufficient} for (1) and (2). For (3), the commutator identity gives
\begin{align*}
    [X_\eps,P] = -X_\eps[1-\eps \px^2,P]X_\eps = \eps \px X_\eps P' X_\eps + \eps X_\eps P' \px X_\eps.
\end{align*}
Thus \eqref{eqn:bound_X_eps_multiplication_basic} now follows from \eqref{eqn:bound_X_eps_basic}.
\end{proof}
In particular, from \eqref{eqn:bound_X_eps_commute_sech}, \eqref{eqn:bound_X_eps_conjugate_cosh} with $K=40$ and \eqref{eqn:comparison_rho_zeta_K_40}, we obtain that 
\begin{align}
    \| \rho X_\eps g \| \leq C\| X[\rho g] \| \leq C \| \rho g \|,\\
    \| \rho^{-1} X_\eps [ \rho g] \| \leq C \| X_\eps g \| \leq C \| g \|.
\end{align}
In the next lemma, we record weighted transfer estimates relating $\bmf$ to $\bmz$ which will be used in the virial estimates for the transformed components.
\begin{lemma}\label{eqn:transfer_fg_z}
There exists $C>0$ such that for any $\eps \in (0,\eps_1)$ and for any $K\geq 10$, the following estimates hold
\begin{align}
 \sqrt{\eps} \| f_1 \|_{L^2} + \eps \| \px f_1 \|_{L^2} + \eps \| f_2 \|_{L^2} &\leq C \| \bmz \|_{H^1 \times L^2}\label{eqn:transfer_f_z},\\
 \sqrt{\eps} \| \rho f_1 \|_{L^2} + \eps \|\rho \px f_1 \|_{L^2} + \eps \|\rho f_2 \|_{L^2} &\leq C \big(\| \rho z_1\|_{L^2} + \| \rho \px z_1 \|_{L^2} + \| \rho z_2 \|_{L^2} \big)\label{eqn:transfer_rho_K_f_z},\\
 \sqrt{\eps} \| \zeta_K f_1 \|_{L^2} + \eps \| \zeta_K \px f_1 \|_{L^2} + \eps \| \zeta_K f_2 \|_{L^2}  &\leq C\big(\| \zeta_K z_1\|_{L^2} + \| \zeta_K \px z_1 \|_{L^2} + \| \zeta_K z_2 \|_{L^2} \big).\label{eqn:transfer_zeta_K_f_z}
\end{align}
\end{lemma}
\begin{proof}
Let $\frakm \in \{1,\rho,\zeta_K\}$. Since $\bmf = \bmX_{\eps} \calD_{1,c}\calD_{0,c}\bmz$, we obtain from Lemma~\ref{lemma:commutator_Xeps_lemma} and \eqref{eqn:rewritten_sfABCD} that 
\begin{align*}
\| \frakm f_1 \| &\lesssim \eps^{-\frac12} \big( \| \frakm z_1 \| + \| \frakm \px z_1 \| + \| \frakm  z_2 \|\big),\\
\| \frakm \px f_1 \| + \| \frakm f_2 \| &\lesssim \eps^{-1} \big( \| \frakm z_1 \| + \| \frakm \px z_1 \| + \| \frakm  z_2 \|\big),
\end{align*}
which proves the estimates \eqref{eqn:transfer_f_z}--\eqref{eqn:transfer_zeta_K_f_z} after multiplying them by $\eps^{\frac12}$ and $\eps$.
\end{proof}

\subsection{Deriving the equations for the directional derivatives and their bounds}
For nonzero velocity $c$, the transport operator $2c\gamma \px$ requires additional virial arguments for both $\bmz$ and $\bmf$ to account for the directional derivative components. Following \cite{kowalczyk2021sufficient}, we derive the evolution systems associated with the directional derivative $k_1 := \px z_1 + c \gamma z_2$ and its regularized version. We then combine weighted regularization estimates with the radiation bounds to obtain spacetime control over the directional derivative and its regularized components. The estimates are then used in the next subsection to control cutoff errors for $\bmf$.

We define the variables $(k_1,k_2) \in L^2 \times H^{-1}$ formally by
\begin{equation}\label{eqn:definition_of_k12}
    \left\{\begin{aligned}
            &k_1 := \px z_1 + c \gamma z_2,\\
    &k_2 := \px z_2 + c \gamma[ \px^2 z_1 - f(x,z_1)+ 2 c \gamma \px z_2].
    \end{aligned}
    \right.
\end{equation}
From \eqref{eqn:z_1_virialform} and $\pt(c \gamma) = \dot{c}\gamma^3$, we formally compute that
\begin{equation*}
\begin{split}
\dot{k}_1 
&= \px(z_2 + \Omega_1(\bmz) + \beta_1) + \dot{c}\gamma^3 z_2 + c \gamma [\px^2 z_1 - f(x,z_1) + 2 c \gamma \px z_2 + \Omega_2(\bmz) + \beta_2 + N_{\bmm}].
\end{split}
\end{equation*}
Using $[\px,\Lambda] = \px$, we note that
\begin{equation*}
\begin{split}
\px \Omega_1(\bmz) + c \gamma \Omega_2(\bmz) &= (\dot{y}-c)\gamma( \px^2 z_1+ c \gamma \px z_2)- \dot{c} c\gamma^2 (\px z_1 + \Lambda \px z_1 + c \gamma \Lambda z_2) + \dot{c} c \gamma^2 \px z_1 \\
&= (\dot{y}-c) \gamma \px k_1  - \dot{c} c \gamma^2 \Lambda k_1  \equiv \Omega_1(\bmk).
\end{split}
\end{equation*}
Moreover, since $\pt(c \gamma) = \dot{c}\gamma^3$ and $\pt(1+2c^2\gamma^2) = 4 \dot{c} c \gamma^4$, we obtain that 
\begin{align*}
    \dot{k}_2 = c \gamma [\px^2 \dot{z}_1 - W''(H+z_1) \dot{z}_1] + (1+2c^2 \gamma^2) \px \dot{z}_2 + \dot{c} \gamma^3[\px^2 z_1 - f(x,z_1) + 4 c \gamma \px z_2].
\end{align*}
Let us keep the final term as an error term. After substituting \eqref{eqn:z_1_virialform}, we find that 
\begin{align*}
    (1+2c^2 \gamma^2)\px[\px^2 z_1 - f(x,z_1) + 2c \gamma \px z_2] + c \gamma [\px^2 z_2 - W''(H+z_1)z_2]\\
    = \px^2 k_1 - W''(H)k_1 + 2 c \gamma \px k_2 - [W''(H+z_1) - W''(H)](k_1 + H'),
\end{align*}
and 
\begin{align*}
    (1+2c^2\gamma^2) \px \Omega_2(\bmz) + c \gamma[\px^2 \Omega_1(\bmz) - W''(H+z_1)\Omega_1(\bmz)] + 2 \dot{c} c \gamma^4 \px z_2 \\
    = \Omega_2(\bmk) + \big[c\gamma^2(\dot y-c) - \dot{ c}c^2\gamma^3x\big]\bigl[W''(H+z_1)-W''(H)\bigr]H',
\end{align*}
where we used $[\px^2,\Lambda] = 2 \px^2$ and $\gamma^2 = 1 + c^2 \gamma^2$. Therefore, we obtain the evolution system for $(k_1,k_2)$:
\begin{equation}
\left\{
\begin{aligned}
&\dot{k}_1 = k_2 + \Omega_1(\bmk)+ \mu_1 \\
&\dot{k}_2 = \px^2 k_1 - W''(H) k_1 + 2 c \gamma \px k_2 + \Omega_2(\bmk) + \mu_2
\end{aligned}\right.,
\end{equation}
where
\begin{equation}\label{eqn:definition_mu1}
\begin{aligned}
\mu_1 := \dot{c} \gamma^3 z_2 + \px \beta_1 + c \gamma (\beta_2+ N_{\bmm}),
\end{aligned}
\end{equation}
and
\begin{equation}\label{eqn:definition_mu2}
\begin{aligned}
\mu_2 :=&\dot{c}\gamma^3
\bigl[\partial_x^2z_1- f(x,z_1)+2c\gamma\partial_xz_2\bigr]-\bigl[W''(H+z_1)-W''(H)\bigr]k_1\\
&-\bigl[1-c\gamma^2(\dot y-c)+\dot c\,c^2\gamma^3x\bigr]
  \bigl[W''(H+z_1)-W''(H)\bigr]H'\\
&+c\gamma\bigl[\partial_x^2\beta_1
               -W''(H+z_1)\beta_1\bigr]
 +(1+2c^2\gamma^2)\partial_x(\beta_2+ N_{\bmm}).
\end{aligned}
\end{equation}
Next, we define the regularized directional derivative variables $(j_1,j_2) \in H^{2} \times H^{1}$ by
\begin{align}
j_1 := X_\eps k_1,\qquad j_2 := X_\eps k_2.
\end{align}
From the equations for $(k_1,k_2)$ we obtain
\begin{equation}
\left\{ \begin{aligned}
\dot{j}_1 &= j_2 + \Omega_1(\bmj) + \nu_1 \\
\dot{j}_2 &= \px^2 j_1 - W''(H)j_1 + 2 c \gamma \px j_2 + \Omega_2(\bmj) + \nu_2 
\end{aligned}\right.,
\end{equation}
where
\begin{equation}
\begin{aligned}
&\nu_1 := X_\eps \mu_1 - 2 \eps \dot{c} c \gamma^2 X_\eps \px^2 j_1, \\
&\nu_2 := X_\eps \mu_2 - 2 \eps \dot{c} c \gamma^2 X_\eps \px^2 j_2 - \eps X_\eps[P_0' \px j_1 + \px(P_0' j_1)].
\end{aligned}
\end{equation} 

We now estimate $\bmj=(j_1,j_2)$. We start with a transfer bound for $j_1$ in terms of $k_1$.
\begin{lemma} There exists $\eps_2>0$ such that for any $\eps \in (0,\eps_2)$ and for any $K\geq 10$, 
\begin{equation}\label{eqn:bound_on_j1_by_k1}
    \int [j_1^2 + \eps(\px j_1)^2 + \eps^2 (\px^2 j_1)^2]\zeta_K^2 \lesssim \int k_1^2 \zeta_K^2. 
\end{equation}
\end{lemma}
\begin{proof}
Since $k_1 = j_1 - \eps \px^2 j_1$, integration by parts gives 
\begin{equation*}
    \int \zeta_K^2 k_1^2 = \int [j_1^2 + 2 \eps (\px j_1)^2 + \eps^2 (\px^2 j_1)^2]\zeta_K^2 - \eps \int j_1^2 (\zeta_K^2)''.
\end{equation*}
From \eqref{eqn:bound_on_zeta_K}, using $|(\zeta_K^2)''| \leq C K^{-1} \zeta_K^2$, we then obtain the bound
\begin{align*}
      \Big(1-  \frac{C\eps}{K} \Big)\int j_1^2 \zeta_K^2 + \int [2 \eps (\px j_1)^2 + \eps^2 (\px^2 j_1)^2]\zeta_K^2 \leq \int  k_1^2\zeta_K^2.
\end{align*}
Choosing $\eps$ small enough completes the proof. 
\end{proof}
Next, we prove a transfer estimate for $j_2$.
\begin{lemma}
Assume \eqref{main_scaling} and $\eps \leq \min(\eps_1,\eps_2)$. For any $K\geq 10$, it holds that 
\begin{equation}\label{eqn:bound_on_j2_by_k1_a_z1}
\int_0^T \int \big[ \eps j_2^2 + \eps^2(\px j_2)^2] \zeta_K^2  \leq C \delta^2 + C\int_0^T \| \zeta_K k_1 \|^2 + C \int_0^T \big( |a|^4 + \| \rho^4 z_1 \|^2\big).
\end{equation}    
\end{lemma}

\begin{proof}
The proof closely follows Lemma~4.9 of \cite{kowalczyk2021sufficient}. We define the virial functional 
\begin{equation}
\calK_{k_1,j_2} = \int k_1 j_2 \zeta_K^2.
\end{equation}
Note that $\calK_{k_1,j_2}$ is well-defined and satisfies $|\calK_{k_1,j_2}| \leq C \eps^{-\frac12} \delta^2$. Using the equations for $\bmk$ and $\bmj$ and integrating by parts, we compute that 
\begin{align*}
\eps\dot{\calK}_{k_1,j_2}=&\eps\int\bigl[j_2^2+\varepsilon(\partial_xj_2)^2\bigr]\zeta_K^2-\eps\int\bigl[(\partial_xj_1)^2
+\varepsilon(\partial_x^2j_1)^2+W''(H)k_1j_1\bigr]\zeta_K^2\\
&-\frac{\eps^2}{2}\int j_2^2(\zeta_K^2)''
+\frac{\eps}{2}\int j_1^2(\zeta_K^2)''+2\eps c\gamma\int k_1\partial_xj_2\,\zeta_K^2\\
&-\eps^2\int k_1X_\varepsilon
\bigl[(W''(H))'\partial_xj_1+\partial_x((W''(H))'j_1)\bigr]\zeta_K^2\\
&+\eps\int\bigl[\Omega_1(\bm k)j_2
+k_1\Omega_2(\bm j)\bigr]\zeta_K^2 - 2\eps^2 \dot{c}c\gamma^2
\int k_1X_\varepsilon\partial_x^2j_2\,\zeta_K^2\\
&+\eps\int\bigl[\mu_1j_2+k_1X_\varepsilon\mu_2\bigr]\zeta_K^2.
\end{align*}

Applying $|(\zeta_K^2)''| \leq C K^{-1} \zeta_K^2$ and Cauchy-Schwarz inequality, and then \eqref{eqn:bound_on_j1_by_k1}, yields
\begin{align*}
&\int\big|\eps(\partial_xj_1)^2
+\eps^2(\partial_x^2j_1)^2\big|\zeta_K^2 + \int |\eps W''(H)k_1j_1|\zeta_K^2 +\frac{\eps^2}{2}\int |j_2^2(\zeta_K^2)''|
+\frac{\eps}{2}\int |j_1^2(\zeta_K^2)''|\\
&\leq  \eps \| \zeta_K \px j_1 \|^2 + \eps^2 \|\zeta_K \px^2 j_1 \|^2+ C\eps \| \zeta_K k_1\|^2 + C\eps \| \zeta_K j_1 \|^2+ C \eps^2 K^{-1} \| \zeta_K j_2\|^2 + C \eps K^{-1} \| \zeta_K j_1\|^2 \\
&\leq C \eps^2 K^{-1} \| \zeta_K j_2\|^2 + C \| \zeta_K k_1 \|^2.
\end{align*}
Moreover, by Young's inequality, 
\begin{align*}
\left| 2 \eps c \gamma \int k_1 (\px j_2) \zeta_K^2 \right| \leq \frac{\eps^2}{8} \| \zeta_K \px j_2\|^2 + C \| \zeta_K k_1 \|^2.
\end{align*}
Next, since $[X_\eps,P_0]k_1 = \eps X_\eps(P_0' \px j_1 + \px(P_0' j_1))$ where $P_0 = W''(H)$,  the bound \eqref{eqn:bound_X_eps_multiplication_basic} gives 
\begin{align*}
&\eps^2\left|\int k_1X_\varepsilon
\bigl[(W''(H))'\partial_xj_1+\partial_x((W''(H))'j_1)\bigr]\zeta_K^2\right|=\eps \left| \int k_1 [X_\eps,P_0]k_1 \zeta_K^2\right| \leq C \eps^{\frac32} \| \zeta_K k_1\|^2.
\end{align*}
For the modulation operator terms, integration by parts gives
\begin{align*}
\int\bigl[\Omega_1(\bm k)j_2
+k_1\Omega_2(\bm j)\bigr]\zeta_K^2 = -(\dot{y}-c)\gamma \int k_1 j_2 (\zeta_K^2)' + \dot{c} c \gamma^2 \int k_1 j_2 [\zeta_K^2 + x(\zeta_K^2)'] + \dot{c} \gamma \int k_1 \px j_1 \zeta_K^2.
\end{align*}
Then, we use $|(\zeta_K^2)'| \leq C \zeta_K^2$ and $|x(\zeta_K^2)'| \leq C \zeta_K$ from \eqref{eqn:bound_on_zeta_K} with $\|k_1 \| \lesssim \delta$ to bound 
\begin{align*}
&\eps \left| \int\bigl[\Omega_1(\bm k)j_2
+k_1\Omega_2(\bm j)\bigr]\zeta_K^2 \right| \leq C \eps (|\dot{y}-c|+|\dot{c}|)\| k_1 \| \| \zeta_K j_2 \|\\
&\leq \frac{\eps}{16}\|\zeta_K j_2\|^2 + C \eps \delta^2(|\dot{y}-c|^2 + |\dot{c}|^2)
\leq \frac{\eps}{16}\| \zeta_K j_2 \|^2 + C|a|^4 + C \|\rho^4 z_1 \|^2.
\end{align*}
Next, using $\eps X_\eps \px^2 = X_\eps -1$ and \eqref{eqn:bound_X_eps_basic}, we have
\begin{align*}
\left|2\eps^2 \dot{c}c\gamma^2
\int k_1X_\varepsilon\partial_x^2j_2\,\zeta_K^2\right| \leq C \eps |\dot{c}c \gamma^2|\left( \left| \int k_1 X_\eps j_2 \zeta_K^2 \right| +  \left| \int k_1 j_2 \zeta_K^2 \right| \right) \leq \frac{\eps}{16}\| \zeta_K j_2\|^2 + C \| \zeta_K k_1 \|^2.
\end{align*}
Finally, using \eqref{eqn:stability_smallness}, \eqref{eqn:estimate_modulation_crude}, \eqref{eqn:bound_on_beta}, \eqref{eqn:bound_on_Nmm}, and \eqref{eqn:bound_X_eps_basic}, we obtain the following bounds for the source terms \eqref{eqn:definition_mu1}-\eqref{eqn:definition_mu2}
\begin{align*}
\|\zeta_K\mu_1\|&\lesssim\delta|\dot c|+|a|^2+\|\rho^4z_1\| \lesssim |a|^2 + \| \rho^4 z_1 \|,\\
\varepsilon\|\zeta_KX_\varepsilon\mu_2\|&\lesssim\delta|\dot c|+\varepsilon\delta\|\zeta_K k_1\|+|a|^2+\|\rho^4z_1\|  \lesssim \|\zeta_K k_1\|+|a|^2+\|\rho^4z_1\|,
\end{align*}
whence
\begin{align*}
\left|\eps\int\bigl[\mu_1j_2+k_1X_\varepsilon\mu_2\bigr]\zeta_K^2 \right| \leq \frac{\eps}{16} \| \zeta_K j_2 \|^2 + C \| \zeta_K k_1 \|^2 + C |a|^4 + C \|\rho^4 z_1 \|^2.
\end{align*}
Combining all these estimates, we have proved that 
\begin{align*}
    \eps\dot{\calK}_{k_1,j_2} \geq \frac{1}{2}\int \big( \eps j_2^2 + \eps^2 (\px j_2)^2 \big) \zeta_K^2 - C \|\zeta_K k_1 \|^2 - C |a|^4 - C \|\rho^4 z_1 \|^2.
\end{align*}
Integrating in time on $[0,T]$ and using the bound $\eps|\calK_{k_1,j_2}|\leq C \delta^2$ proves the lemma. 
\end{proof}
Combining the preceding two lemmas with the directional derivative estimate yields the following corollary.
\begin{corollary}\label{corollary:directional_derivative} Assume that \eqref{main_scaling} holds, $\eps \leq \min(\eps_1,\eps_2)$ and $A$ is sufficiently large. For any $10 \leq K \leq \frac{A}{2}$, it holds that 
    \begin{equation}\label{eqn:estimate_directional_j_z}
        \int_0^T \int \big( j_1^2 + \eps (\px j_1)^2 + \eps^2 (\px^2 j_1)^2 + \eps j_2^2 + \eps^2 (\px j_2)^2 \big) \zeta_K^2\lesssim A\delta^2 + \int_0^T |a|^4 + \int_0^T \| \rho^4 z_1 \|^2.
    \end{equation}
\end{corollary}
\begin{proof}
Recall that $k_1 = \px z_1 + c \gamma z_2$. From \eqref{eqn:HE_transfer_zK2} and \eqref{eqn:HE_key_estimate2},
\begin{align*}
    \int_0^T \| \zeta_K k_1 \|^2 \lesssim  \delta^2 + \int_0^T \Big( \| \px w_1 + c \gamma w_2 \|^2 + |a|^4 + \| \rho^4 z_1 \|^2\Big)\lesssim A \delta^2 + \int_0^T |a|^4 + \int_0^T \| \rho^4 z_1 \|^2.
\end{align*}
After integrating \eqref{eqn:bound_on_j1_by_k1} in time and adding \eqref{eqn:bound_on_j2_by_k1_a_z1}, we obtain \eqref{eqn:estimate_directional_j_z}.

\end{proof}

\subsection{Bounds for the transformed problem}

We define the localization of the regularized transformed variables $\bmf$ at scale $B$:
\begin{equation}
    h_1 := \chi_A \zeta_B f_1 ,\qquad h_2 := \chi_A \zeta_B f_2.
\end{equation}

 Following \cite{kowalczyk2021sufficient}, we define the potential 
\begin{align*}
    Z_B := \frac12 \ln(\zeta_B)'' - \frac{\varphi_B P_2'}{2\zeta_B^2}.
\end{align*}
By repulsivity of $P_2$, one has 
\begin{equation*}
    - \frac{\varphi_B P_2'}{2\zeta_B^2} \geq \frac{1}{2C_0} \mathbbm{1}_{[x_1,x_2]}.
\end{equation*}
Moreover, from \eqref{eqn:bound_on_zeta_K},
\begin{align*}
    |\ln(\zeta_B)''| \leq C B^{-1} \mathbbm{1}_{[x_1,x_2]}.
\end{align*}
Hence, for all $B \geq B_1 := \max(80,2C_0C)$, we have 
\begin{equation}\label{eqn:lower_bound_Z_B}
    Z_B \geq \frac{1}{4C_0} \mathbbm{1}_{[x_1,x_2]}.
\end{equation}
We fix $B = B_1$ and let  $A \geq 2B$ be sufficiently large so that $\eps = A^{-\frac14} \leq \min(\eps_1,\eps_2)$ from here onward.

\begin{lemma} Assuming \eqref{main_scaling}, it holds that 
\begin{align}
 \| (\px f_1 + c \gamma f_2) \chi_A \rho \| &\lesssim \| \px h_1 + c \gamma h_2 \| + B^{-1} \| f_1 \rho \| \label{eqn:bound_directional_f_chiArho},\\
 \int_{|x|< A} (\px f_1 + c \gamma f_2)^2 \rho^2 &\lesssim \int (\px h_1 + c \gamma h_2)^2 + B^{-2}\int_{|x|<A} f_1^2 \rho^2\label{eqn:bound_scaleA_directional_f_chiArho},\\
\int_{|x|<A} f_1^2 \rho^2 & \lesssim \int h_1^2 \rho \label{eqn:bound_scaleA_f1_h1}.
\end{align}
Moreover, on the cutoff region $|x|<2A$, it holds that 
\begin{equation}\label{eqn:bound_scale2A_directional_f_chiArho}
\begin{split}
 \int_{|x|< 2A} (\px f_1 + c \gamma f_2)^2 &\lesssim \int [(\px^2 j_1)^2 +(\px j_1)^2 + (\px j_2)^2 + j_2^2 + k_1^2] \zeta_A^4 \\
 &\qquad + \| \rho^4 z_1 \|^2 + \delta^2 \int[(\px z_1)^2 + z_1^2]\zeta_A^4,
\end{split}
\end{equation}
and
\begin{equation}\label{eqn:bound_scale2A_f1}
    \int_{|x|<2A} f_1^2 \lesssim \eps^{-1} \int [(\px z_1)^2 + z_2^2]\zeta_A^4 + \int z_1^2 \zeta_A^4.
\end{equation}
\end{lemma}

\begin{proof}
First, we note that
\begin{equation}\label{eqn:identity_directional_h_to_f}
    \px h_1 + c \gamma h_2 = \chi_A \zeta_B (\px f_1 + c \gamma f_2) + (\chi_A \zeta_B)' f_1.
\end{equation}
Moreover, the weight estimates in \eqref{eqn:bound_on_zeta_K} give, for $B \geq 80$,
\begin{equation}\label{eqn:bound_on_zeta_B}
    |(\chi_A \zeta_B)'| \leq  C B^{-1} \zeta_B, \qquad \rho\leq C\zeta_B^2 \leq C \zeta_B.
\end{equation}
Hence, by multiplying \eqref{eqn:identity_directional_h_to_f} by the weight $\rho \zeta_B^{-1}$, we obtain the estimate \eqref{eqn:bound_directional_f_chiArho} using the above bounds. Next, on the interval $|x|<A$, we have $\chi_A \equiv 1$, and the identities
\begin{align*}
    \rho(\px f_1 + c \gamma f_2) &= \frac{\rho}{\zeta_B}(\px h_1 + c \gamma h_2) - \frac{\zeta_B'\rho}{\zeta_B} f_1, \qquad |x|<A,\\
    \rho f_1 &= \frac{\rho}{\zeta_B} h_1, \qquad |x|<A.
\end{align*}
Hence, the estimates \eqref{eqn:bound_scaleA_directional_f_chiArho} and \eqref{eqn:bound_scaleA_f1_h1} now follow from the above identities and \eqref{eqn:bound_on_zeta_B}. 

Next, using the definition \eqref{eqn:definition_of_k12}, and the short-hand notation from \eqref{eqn:definition_sfABCD}, we compute that 
\begin{equation}\label{eqn:relation_directional_f_to_j}
\begin{split}
\px f_1 + c \gamma f_2 &=X_\eps[(\sfA- (V_0'+V_1'))k_1 + \sfB k_2 + ([V_0(V_0+V_1)]'-\gamma^2P_0')z_1 + c \gamma (\sfC + \sfB L_0) z_1\\
&\qquad \qquad + c \gamma (\sfD - \sfA - 2 c \gamma \sfB \px) z_2 + c \gamma \sfB \big(W'(H+z_1) - W'(H) - P_0z_1\big)]\\
&=X_\eps \left[ \big(\sfA - (V_0'+V_1')\big)k_1 + \sfB k_2 + \big([V_0(v_0+v_1)]' - \gamma^2 P_0'\big)z_1 \right. \\
&\qquad\qquad +\left. c \gamma \sfB \big(W'(H+z_1) - W'(H) - P_0z_1\big)\right],
\end{split}
\end{equation}
where we used $\sfC +\sfB L_0 = 0$ and $\sfD - \sfA - 2 c \gamma \sfB \px = 0$.

Since $k_1 = j_1 - \eps \px^2 j_1$, from \eqref{eqn:rewritten_sfABCD} we compute that 
\begin{align*}
X_\eps \big(\sfA - (V_0'+V_1')\big)k_1 &= \gamma^2 \px^2 j_1 - (V_0 + V_1) \px j_1 + X_\eps \big( V_0 (V_0 + V_1) - \gamma^2 P_0 - V_0' V_1'\big) k_1 \\
&\qquad - \eps X_\eps \big( V_0' + V_1') \px^2 j_1 - \eps \px X_\eps \big(V_0' + V_1') \px j_1.
\end{align*}
Note that all of the above coefficients are bounded. We obtain from \eqref{eqn:bound_X_eps_basic} the bound 
\begin{align*}
    \left\| \zeta_A^2 X_\eps \big(\sfA - (V_0'+V_1')\big)k_1 \right\| \lesssim \| \zeta_A^2 \px^2 j_1 \| + \| \zeta_A^2 \px j_1 \| + \| \zeta_A^2 k_1 \|.
\end{align*}
Similarly, using $X_\eps \px k_2 = \px j_2$, we have 
\begin{align*}
    \| \zeta_A^2 X_\eps \sfB k_2 \| \lesssim \| \zeta_A^2 \px j_2 \| + \| \zeta_A^2 j_2 \|.
\end{align*}
Next, from Lemma~\ref{lemma:scalar_darboux_factorization}, we have 
\begin{align*}
    \left| [V_0(V_0+V_1)]' - \gamma^2 P_0' \right| \leq C \rho^4,
\end{align*}
whence
\begin{align*}
    \left\| \zeta_A^2 X_\eps \big([V_0(V_0+V_1)]' - \gamma^2 P_0'\big)z_1\right\|\lesssim \| \rho^4 z_1 \|.
\end{align*}
Finally, by Taylor expansion, we have
\begin{align*}
    |W'(H+z_1) - W'(H) - P_0z_1| &\leq C |z_1|^2,\\
    |\px \big( W'(H+z_1) - W'(H) - P_0z_1\big)| &\leq C( |z_1||\px z_1| + |H'||z_1|^2),
\end{align*}
Using the expression for $\sfB$ in \eqref{eqn:rewritten_sfABCD} and the bounds \eqref{eqn:stability_smallness} \eqref{eqn:bound_X_eps_basic}, we obtain 
\begin{align*}
    \Big\| \zeta_A^2 X_ \eps  \sfB \big(W'(H+z_1) - W'(H) - P_0z_1\big) \Big\| \lesssim \delta( \| \zeta_A^2 \px z_1 \| + \| \zeta_A^2 z_1 \|).
\end{align*}
Thus, using the preceding bounds, we obtain that 
\begin{align*}
    &\int_{|x|<2A} (\px f_1 + c \gamma f_2)^2  \lesssim \int \zeta_A^4 (\px f_1 + c \gamma f_2)^2\\
    &\quad \lesssim \int [(\px^2 j_1)^2 +(\px j_1)^2 + (\px j_2)^2 + j_2^2 + k_1^2] \zeta_A^4 + \| \rho^4 z_1 \|^2 + \delta^2 \int[(\px z_1)^2 + z_1^2]\zeta_A^4,
\end{align*}
which is \eqref{eqn:bound_scale2A_directional_f_chiArho}. Lastly, noting that 
\begin{align*}
    f_1 = \gamma^2 \px X_\eps(\px z_1 + 2c \gamma z_2) + X_\eps[(V_0(V_0+V_1)-\gamma^2 P_0)z_1 -(V_0 + V_1)(\px z_1 + c \gamma z_2) ] ,
\end{align*}
we have 
\begin{align*}
    \| \zeta_A^2 f_1 \| \lesssim \eps^{-\frac12} \big( \| \zeta_A^2 \px z_1 \| + \| \zeta_A^2 z_2 \|\big) + \| \zeta_A^2 z_1 \|,
\end{align*}
whence \eqref{eqn:bound_scale2A_f1}.
\end{proof}
We now combine the repulsivity of $P_2$ with the preceding transfer estimates to establish a virial estimate for the localized transformed variables.
\begin{lemma}\label{lemma:LE_main_lemma} Assume \eqref{main_scaling}. For any $\eta \in (0,1)$, it holds that
\begin{equation}\label{eqn:estimate_from_calF}
\begin{split}
\int_0^T  \left[ \| \px h_1 + c \gamma h_2\|^2 + \int_{x_1}^{x_2} h_1^2 \right] &\leq C_{\eta,B}\bigg[A \delta^2 + \int_0^T|a|^4\bigg] + \frac{C_B}{\sqrt{A}} \int_0^T \| \rho^4 z_1 \|^2    \\
&\quad+ C_B(\eta + A^{-\frac18}\big)\int_0^T \int [f_1^2 + (\px f_1)^2 + f_2^2]\rho^2.
\end{split}
\end{equation}
\end{lemma}
\begin{proof}
We again use the modified virial functional
\begin{align*}
    \calF := d_1\calF_1 + d_2 \calF_2 + d_3 \calF_3,
\end{align*}
where $d_j$ are the coefficients from Lemma~\ref{lemma:coefficients_dj} and $\calF_j := \calA_j[\psi_{A,B}](f_1,f_2)$ in the notation of Lemma~\ref{lemma:general_virial}. Note that
\begin{equation}\label{eqn:bound_on_calF}
    |\calF_1| + |\calF_2| + |\calF_3| + |\calF| \leq C \eps^{-2} \delta^2 \leq C A^{\frac12} \delta^2.
\end{equation}

As in \cite[\S4.11.2]{kowalczyk2021sufficient}, we organize the terms of $\dot{\calF}$ as 
\begin{align*}
    \dot{\calF}=& -2d_1 \int \Big((\px h_1 + c \gamma h_2)^2 +Z_B h_1^2\Big) + 2(d_2-d_3)\int \psi_{A,B}' (\px f_1 + c \gamma f_2)f_2\\
    &\qquad +\calE_{\sfT} + \calE_{\sfS} + \calE_{\xi} + d_1 \theta - (\dot{y} - c) \gamma \calI_{\bmf} + \dot{c} \gamma \calJ_{\bmf} + \calR_{\calF},
\end{align*}
where we define the following error terms (arising from \eqref{eqn:definition_xi_T}) by
\begin{align*}
    \calE_{\sfT} = \calE[0,\sfT], \quad \calE_{\sfS} = \calE[X_\eps \sfB N_2, X_\eps \sfD N_2], \quad \calE_\xi = \calE[\xi_1,\xi_2]
\end{align*}
with the common notation
\begin{align*}
\calE[\alpha_1,\alpha_2] &:= d_1 \int (2 \psi_{A,B} \px \alpha_1 + \psi_{A,B}' \alpha_1) f_2 +d_1 \int \big[2\psi_{A,B}(\px f_1 + c \gamma f_2) + \psi_{A,B}' f_1\big]\alpha_2\\
&\quad + 2 c \gamma d_1 \int \psi_{A,B}\big[P_2 f_1 \alpha_1 - (\px f_1) (\px \alpha_1)\big]\\
&\quad + 2(d_3 - d_2)\int \psi_{A,B}\big[(\px f_1)(\px \alpha_1) + P_2 f_1 \alpha_1 + f_2 \alpha_2\big].
\end{align*}
We decompose $\theta$ as $\theta = \theta_1+\theta_2+\theta_3$, where the cut-off terms are given by
\begin{align*}
    \theta_1 &:= -2 \int (\px f_1 + c \gamma f_2)^2 (\chi_A^2)' \varphi_B,\\
    \theta_2 &:= -2 \int (\px f_1 + c \gamma f_2)f_1[(\chi_A^2)'\varphi_B]',\\
    \theta_3 &:= -\int f_1^2\Big[\frac12 (\chi_A^2)''' \varphi_B + (3 \chi_A \chi_A'' + (\chi_A')^2) \zeta_B^2 + (\chi_A^2)' \zeta_B \zeta_B' \Big].
\end{align*}
The modulation error terms are given by 
\begin{align*}
    \calI_{\bmf}  &:= \int \psi_{A,B}' \big[2d_1 (\px f_1) f_2 + (d_3 - 2d_2)(\px f_1)^2 + d_3(f_2^2 + P_2 f_1^2)\big] \\
    &\qquad + d_1 \int \psi_{A,B}'' f_1 f_2 + d_3 \int \psi_{A,B} P_2' f_1^2,\\
    \calJ_{\bmf} &:= c \gamma \int x\psi_{A,B}' \big[ 2d_1 (\px f_1) f_2 + (d_3 - 2d_2) (\px f_1)^2 + d_3(f_2^2 + P_2 f_1^2) \big]\\
    &\qquad + (d_2 - d_3) \int \psi_{A,B}' f_1f_2 + c \gamma d_1 \int x \psi_{A,B}'' f_1 f_2 \\
    &\qquad - \frac{d_1}{2} \int \psi_{A,B}'' f_1^2 + c \gamma d_3 \int x \psi_{A,B} P_2' f_1^2,
\end{align*}
and the remainder term is given by 
\begin{align*}
\calR_{\calF} = \frac{\dot{c}+\dot{\Pi}}{1-(c+\Pi)^2}\big[ -c \calF_1 - c^2 \gamma \calF_2 + \gamma^{-1} \calF_3 + 2(c+\Pi)\calF\big].
\end{align*}

In what follows, we perform an integrated-in-time estimate for the terms in $\dot{\calF}$, except for the first term, which will be retained for a lower bound. Throughout, we will frequently use that $\psi_{A,B}$ and its derivatives are supported on $|x|<2A$, and we will use the following basic bound from \eqref{eqn:bound_on_phi_B}, 
\begin{equation}\label{eqn:bound_basic_psiAB}
    |\psi_{A,B}| + |\psi_{A,B}'| \leq C_B \mathbbm{1}_{|x|<2A} \lesssim \zeta_A^4 = \zeta_{A/2}^2.
\end{equation}

\underline{Estimate of the cross term}: From \eqref{eqn:stability_smallness}, \eqref{eqn:bound_on_Pi}, and \eqref{eqn: d_2-d_3=2cgd_1}, we note that $|d_3 - d_2| \lesssim \delta^2$. Using \eqref{eqn:bound_basic_psiAB} and \eqref{eqn:transfer_zeta_K_f_z},
\begin{align*}
&\left| 2(d_2-d_3) \int \psi_{A,B}' (\px f_1 + c \gamma f_2) f_2 \right| \leq C \delta^2 \int_{|x|<2A}\big((\px f_1)^2 + f_2^2 \big) \\
&\leq C \delta^2 \eps^{-2} \int \big( z_1^2 + (\px z_1)^2 + z_2^2\big) \zeta_A^4=C \delta^2 \eps^{-2} \int \big( z_1^2 + (\px z_1)^2 + z_2^2\big) \zeta_{A/2}^2.
\end{align*}
Integrating in time and using Proposition~\ref{prop:high_energy_estimate} (with $K \leftarrow \frac{A}{2}$) gives 
\begin{align}\label{eqn:auxiliary_estimate_cross}
    \int_0^T \left| 2(d_2-d_3) \int \psi_{A,B}' (\px f_1 + c \gamma f_2) f_2 \right| \leq C_B \delta^2 \eps^{-2} A^2\left[A \delta^2 + \int_0^T \big( |a|^4 + \| \rho^4 z_1 \|^2 \big)\right].
\end{align}

\underline{Estimate of the commutator term}: We recall from \eqref{eqn:definition_xi_T} that $\sfT = -\eps X_\eps \big(P_2' \px f_1 + \px(P_2' f_1) \big)$,  which gives 
\begin{align*}
    \calE_{\sfT} = \calE[0,\sfT] = d_1 \int \big[2 \psi_{A,B}(\px f_1 + c \gamma f_2) + \psi_{A,B}' f_1 \big] \sfT + 2 (d_3 - d_2) \int \psi_{A,B} f_2 \sfT.
\end{align*}
The commutator bounds \eqref{eqn:bound_X_eps_basic}, \eqref{eqn:bound_X_eps_conjugate_cosh}, \eqref{eqn:bound_X_eps_multiplication_basic}, together with the  decay bound $|P_2'| \leq C \rho^2$, yield
\begin{equation}\label{eqn:estimate_for_sfT}
        \| \rho^{-1} \sfT \| \lesssim \eps \| \rho^{-1} P_2' \px f_1 \| + \sqrt{\eps} \| \rho^{-1} P_2' f_1 \| \lesssim \eps \| \rho \px f_1 \| + \sqrt{\eps} \| \rho f_1\|.
\end{equation}
Using the cut-off bounds $|\psi_{A,B}| \leq C B\chi_A^2$, $|\psi_{A,B}'| \leq C \chi_A$ from \eqref{eqn:bound_on_phi_B}, the transfer estimate \eqref{eqn:bound_directional_f_chiArho}, and Cauchy-Schwarz inequality, we obtain 
\begin{align*}
|\calE_{\sfT}| &\leq C_B \big[\| \px h_1 + c \gamma h_2 \| + \| \rho f_1 \| + |d_3 - d_2| \| \rho f_2 \|\big]\cdot  \big[ \eps \| \rho \px f_1 \| + \sqrt{\eps} \| \rho f_1\|\big]\\
&\leq \frac{1}{2^3}\| \px h_1 + c \gamma h_2 \|^2 + C_B \big[ \eps^{\frac12} \| \rho f_1\|^2 + \eps^{\frac32} \| \rho \px f_1\|^2 + \eps^{\frac12} \delta^4 \| \rho f_2\|^2\big].
\end{align*}
Using $\eps = A^{-\frac14}$ and integrating the above estimate yields 
\begin{equation}\label{eqn:auxiliary_estimate_commutator_T}
\begin{split}
\int_0^T |\calE_{\sfT} | \leq \frac14 \int_0^T \| \px h_1 + c \gamma h_2 \|^2 + C_B A^{-\frac18}\bigg[  \int_0^T \big( f_1^2 + (\px f_1)^2 + f_2^2) \rho^2 \bigg].
\end{split}
\end{equation}

\underline{Estimate of the nonlinear terms}: For any fixed $\eta\in (0,1)$, we claim that 
\begin{equation}\label{eqn:auxiliary_estimate_nonlinear_term}
\begin{split}
\int_0^T |\calE_{\sfS}| &\leq \eta \int_0^T \int \big[ f_1^2 + (\px f_1)^2 + f_2^2]\rho^2 + C_B C_\eta\int_0^T |a|^4 \\
&\qquad + C_B A^{-\frac{13}{8}}\bigg[A \delta^2 + \int_0^T \big(|a|^4 + \|\rho^4 z_1 \|^2 \big)\bigg].
\end{split}
\end{equation}
First, as in \cite{kowalczyk2022kink}, we separate the pure internal mode forcing terms from the nonlinearity by writing $N_2 = N_{m_1} + \tilN$, with
\begin{align*}
    N_{m_1} &:= -[W'(H + m_1) - W'(H) - W''(H)m_1],\\
    \tilN &:=  -[W'(H+z_1+m_1)-W'(H+m_1)-W''(H)z_1].
\end{align*}
Therefore, we can split the nonlinear error term as 
\begin{equation*}
    \calE_{\sfS} = \calE[X_\eps \sfB N_{m_1}, X_\eps \sfD N_{m_1}] + \calE[X_\eps \sfB \tilN, X_\eps \sfD \tilN].
\end{equation*}
By Taylor expansion, \eqref{eqn:stability_smallness}, and \eqref{eqn:estimate_on_Phi_c_pm}, we have the pointwise bounds
\begin{align*}
|\tilN| + |\px \tilN| &\lesssim \delta(|z_1| + |\px z_1|),\\
|\px^j N_{m_1}| &\lesssim |a|^2 \rho^{8}, \quad 0 \leq j \leq 5.
\end{align*}

Using the expressions for $\sfB$ and $\sfD$ from \eqref{eqn:rewritten_sfABCD}, we have 
\begin{align*}
    \| \zeta_A^2 X_\eps \sfB \tilN \| + \sqrt{\eps} \| \zeta_A^2 \px X_\eps \sfB \tilN \| + \sqrt{\eps} \| \zeta_A^2 X_\eps \sfD \tilN \| \lesssim \delta \left(\|\zeta_A^2 z_1 \| + \| \zeta_A^2 \px z_1 \|\right).
\end{align*}
By \eqref{eqn:bound_basic_psiAB} and Cauchy-Schwarz inequality, and then by \eqref{eqn:transfer_zeta_K_f_z} (with $K \leftarrow \frac{A}{2}$),
\begin{align*}
    |\calE[X_\eps\sfB\tilN,X_\eps\sfD\tilN]| &\leq C_B \big(\|\zeta_A^2 f_1\| + \|\zeta_A^2 \px f_1\| +\|\zeta_A^2 f_2\|\big) \big(\| \zeta_A^2 X_\eps \sfB \tilN \| + \| \zeta_A^2 \px X_\eps \sfB \tilN \| + \| \zeta_A^2 X_\eps \sfD \tilN \|\big)\\
    &\lesssim \eps^{-\frac32} \delta \left(\|\zeta_A^2 z_1 \|^2 + \| \zeta_A^2 \px z_1 \|^2 + \| \zeta_A^2 z_2\|^2\right).
\end{align*}
Since $\zeta_A^4 = \zeta_{A/2}^2$, $\eps = A^{-\frac14}$, and $\delta\leq A^{-4}$, we may apply Proposition~\ref{prop:high_energy_estimate} (with $K \leftarrow \frac{A}{2}$) to bound
\begin{equation*}
    \int_0^T |\calE[X_\eps\sfB\tilN,X_\eps\sfD\tilN]| \leq C_B \eps^{-\frac{3}{2}} \delta A^2 \bigg[ A \delta^2 + \int_0^T \big( |a|^4 + \| \rho^4 z_1\|^2\big)\bigg] \leq C_B  A^{-\frac{13}{8}}\bigg[ A \delta^2 + \int_0^T \big( |a|^4 + \| \rho^4 z_1\|^2\big)\bigg].
\end{equation*}
For the pure internal mode forcing terms, we have
\begin{align*}
    \| \rho^{-1} X_ \eps \sfB N_{m_1} \| + \| \rho^{-1} \px X_ \eps \sfB N_{m_1} \| + \| \rho^{-1} X_ \eps \sfD N_{m_1} \|  \lesssim |a|^2,
\end{align*}
which implies 
\begin{align*}
\int_0^T |\calE[X_\eps\sfB N_{m_1},X_\eps\sfD N_{m_1}]| &\leq C_B \int_0^T |a|^2 \big( \| \rho f_1 \| + \| \rho \px f_1 \| + \| \rho f_2 \|\big)\\
&\leq \eta \int_0^T \int \big( f_1^2 + (\px f_1)^2 + f_2^2\big) \rho^2 + C_B C_\eta \int_0^T |a|^4.
\end{align*}
This concludes \eqref{eqn:auxiliary_estimate_nonlinear_term}.

\underline{Estimate of the cut-off terms}: We show that
\begin{equation}\label{eqn:estimate_theta_error_terms}
    \int_0^T \big(|\theta_1|+|\theta_2|+|\theta_3|\big) \lesssim A^{-\frac12} \left[A \delta^2 + \int_0^T \big(|a|^4 + \| \rho^4 z_1 \|^2\big) \right].
\end{equation}
Indeed, to estimate $\theta_1$, we use $|(\chi_A^2)'\varphi_B| \leq C BA^{-1} \mathbbm{1}_{|x|<2A}$,  \eqref{eqn:bound_scale2A_directional_f_chiArho}, $\zeta_A^4 = \zeta_{A/2}^2$, together with \eqref{eqn:estimate_directional_j_z} and Proposition~\ref{prop:high_energy_estimate} (with $K \leftarrow \frac{A}{2}$), to obtain
\begin{align*}
    &\int_0^T |\theta_1| \leq CB A^{-1}\int_0^T\int_{|x|<2A} (\px f_1 + c \gamma f_2)^2 \\
    &\leq CBA^{-1} \bigg[\int_0^T \int\big[(\px^2 j_1)^2 + (\px j_1)^2 + (\px j_2)^2 + j_2^2 + k_1^2  \big]\zeta_{A}^4 +  \int_0^T \| \rho^4 z_1\|^2 + \delta^2 \int_0^T \big[z_1^2+ (\px z_1)^2 ]\zeta_A^4\bigg]\\
    &\leq CBA^{-1} \eps^{-2}\bigg[ A \delta^2 + \int_0^T\big( |a|^4 + \| \rho^4 z_1\|^2\big) \bigg].
\end{align*}
This implies the estimate \eqref{eqn:estimate_theta_error_terms} for $\theta_1$ using the scaling $\eps=A^{-\frac14}$ from \eqref{main_scaling}.

Similarly, \eqref{eqn:bound_scale2A_f1} and Proposition~\ref{prop:high_energy_estimate} (with $K \leftarrow \frac{A}{2})$ yield
\begin{equation}\label{eqn:estimate_auxiliary_f1}
    \int_0^T \int_{|x|<2A} f_1^2 \leq C \eps^{-1} A^2\bigg[ A \delta^2 + \int_0^T \big(|a|^4 + \|\rho^4 z_1\|^2\big)\bigg].
\end{equation}
Noting that $\chi_A'$ and its derivatives are supported on $\{A<|x|<2A\}$, we have
\begin{equation*}
    \zeta_B^2(x) \leq C e^{-\frac{2\kappa A}{B}}, \qquad |\zeta_B(x) \zeta_B'(x)| \leq C B^{-1}  e^{-\frac{2\kappa A}{B}} , \qquad A<|x|<2A,
\end{equation*}
whence
\begin{align*}
    \left|\frac12 (\chi_A^2)''' \varphi_B + (3 \chi_A \chi_A'' + (\chi_A')^2) \zeta_B^2 + (\chi_A^2)' \zeta_B \zeta_B'\right| \lesssim  BA^{-3} + (A^{-2}+A^{-1}B^{-1})   e^{-\frac{2\kappa A}{B}}   \lesssim_\kappa B A^{-3}.
\end{align*}
Thus, we have 
\begin{align*}
    \int_0^T |\theta_3| \lesssim BA^{-3}\int_0^T \int_{|x|<2A} f_1^2 \lesssim A^{-\frac34}\bigg[ A \delta^2 + \int_0^T \big(|a|^4 + \|\rho^4 z_1\|^2\big)],
\end{align*}
which is better than \eqref{eqn:estimate_theta_error_terms}. Finally, by the Cauchy-Schwarz inequality and \eqref{eqn:estimate_auxiliary_f1}, we obtain 
\begin{align*}
    \int_0^T |\theta_2| \leq C B A^{-1} \int_0^T \int_{|x|<2A} (\px f_1 + c \gamma f_2)^2 + CB A^{-3} \int_0^T \int_{|x|<2A} f_1^2\\
    \lesssim A^{-\frac12} \left[A \delta^2 + \int_0^T \big(|a|^4 + \| \rho^4 z_1 \|^2\big) \right].
\end{align*}
In conclusion, we have shown \eqref{eqn:estimate_theta_error_terms}.

\underline{Estimate of the modulation terms}: Here we estimate the terms $\calE_\xi$, $\calI_{\bmf}$, and $\calJ_{\bmf}$. First, after recalling the definition of $\xi_1,\xi_2$ from \eqref{eqn:definition_xi_T} and the definition of $\eta_1,\eta_2$ from \eqref{eqn:definition_eta1}-\eqref{eqn:definition_eta2}, we apply the commutator estimates from Lemma~\ref{lemma:commutator_ABCD} and Lemma~\ref{lemma:commutator_Xeps_lemma} with \eqref{eqn:transfer_f_z} to obtain 
\begin{align*}
    &\| \xi_1 \| + \| \px \xi_1 \| + \| \xi_2 \| \\
    &    \lesssim \| X_\eps \eta_1 \| + \| \px X_\eps \eta_1 \| + \|X_\eps \eta_2\| + \eps |\dot{c}| \big( \| X_\eps \px^2 f_1\|+ \| \px X_\eps \px^2 f_1\| + \| X_\eps \px^2 f_2\|\big)\\
    &\lesssim \big(\eps^{-1}|\dot{y} -c| + \eps^{-\frac12}|\dot{c}|\big)\|\bmz\|_{H^1 \times L^2} + \big(|\dot{y} -c| + |\dot{c}|\big)|a|+\eps^{-1}|\dot{c}| \| \bmz \|_{H^1 \times L^2}.
\end{align*}
Hence, using the Cauchy-Schwarz inequality and the stability estimate \eqref{eqn:stability_smallness}, we obtain
\begin{equation}\label{eqn:estimate_on_xi12}
\| \xi_1 \| + \| \px \xi_1 \| + \| \xi_2 \| \lesssim \eps^{-1} \delta \big( |a|^2 + \|\rho^4 z_1 \|^2\big),
\end{equation}
then by \eqref{eqn:bound_basic_psiAB} and \eqref{eqn:transfer_zeta_K_f_z},
\begin{align*}
    |\calE_\xi| &\leq C_B \eps^{-1}\delta\big(|a|^2 + \|\rho^4 z_1\|^2 \big) \big(\|\zeta_{A/2} f_1\| + \| \zeta_{A/2} \px f_1\| + \| \zeta_{A/2} f_2 \|\big)\\
    &\leq C_B \eps^{-2}\delta\big(|a|^2 + \|\rho^4 z_1\|^2 \big) \big(\|\zeta_{A/2} z_1\| + \| \zeta_{A/2} \px z_1\| + \| \zeta_{A/2} z_2 \|\big).
\end{align*}
We integrate the above estimate in time and use Proposition~\ref{prop:high_energy_estimate} (with $K\leftarrow \frac{A}{2}$) to obtain 
\begin{equation}\label{eqn:auxiliary_estimate_cal_xi}
\begin{aligned}
\int_0^T |\calE_\xi| &\leq C_B \eps^{-2} \delta A^2\bigg[ A \delta^2 + \int_0^T \big( |a|^4 + \|\rho^4 z_1 \|^2\big) \bigg] \\
&\leq  C_B A^{-\frac32}\bigg[ A \delta^2 + \int_0^T \big( |a|^4 + \|\rho^4 z_1 \|^2\big) \bigg] .
\end{aligned}
\end{equation}
Next, we consider $\calI_{\bmf}$ and $\calJ_{\bmf}$, and use the fact that the weights are bounded by $CB \mathbbm{1}_{|x|<2A}$ to derive that
\begin{equation*}
    |\calI_{\bmf}| + |\calJ_{\bmf}| \leq C_B \int_{|x|<2A} \big(f_1^2 + (\px f_1)^2 + f_2^2\big)\leq C_B \int \zeta_{A/2}^2\big(f_1^2 + (\px f_1)^2 + f_2^2\big).
\end{equation*}
Thus, by \eqref{eqn:estimate_modulation_crude} and by applying \eqref{eqn:transfer_zeta_K_f_z} and Proposition~\ref{prop:high_energy_estimate} (with $K\leftarrow \frac{A}{2}$), we conclude that 
\begin{equation}\label{eqn:auxiliary_estimate_modulation}
\begin{aligned}
\int_0^T \big(|(\dot{y}-c)\gamma\calI_{\bmf}| + |\dot{c} \gamma \calJ_{\bmf}|\big) &\leq C_B \eps^{-2} \delta^2 A^2 \bigg[A \delta^2 + \int_0^T \big( |a|^4 + \|\rho^4 z_1 \|^2\big) \bigg]\\
&\leq C_B A^{-\frac{11}{2}}\bigg[A \delta^2 + \int_0^T \big( |a|^4 + \|\rho^4 z_1 \|^2\big) \bigg].   
\end{aligned}
\end{equation}

\underline{Estimate of the remainder term and conclusion}: Finally, we bound the remainder term $\calR_{\calF}$. By Lemma~\ref{lemma:normal_form_c} and \eqref{eqn:bound_on_calF}, we have
\begin{align*}
    |\calR_{\calF}| \leq C_B \eps^{-2} \delta^2 \big(|a|^4 + \|\rho^4 z_1 \|^2 + \|\rho^4 z_2\|^2\big).
\end{align*}
Since $\rho^4 \leq C \zeta_{10}$, Proposition~\ref{prop:high_energy_estimate} (with $K \leftarrow 10$) yields 
\begin{align*}
    \int_0^T \big(\|\rho^4 z_1 \|^2 + \|\rho^4 z_2\|^2\big) \leq C\int_0^T \big(\|\zeta_{10} z_1 \|^2 + \|\zeta_{10} z_2\|^2\big) \leq CA \delta^2 + C \int_0^T\big( |a|^4 + \|\rho^4 z_1 \|^2\big).
\end{align*}
Therefore, we have
\begin{equation}\label{eqn:auxiliary_estimate_remainder}
    \int_0^T |\calR_{\calF}| \leq C_B \eps^{-2} \delta^2 \bigg[A \delta^2 + \int_0^T \big( |a|^4 + \|\rho^4 z_1 \|^2\big) \bigg] \leq C_B A^{-\frac{15}{2}} \bigg[A \delta^2 + \int_0^T \big( |a|^4 + \|\rho^4 z_1 \|^2\big) \bigg].
\end{equation}

In conclusion, using the lower bounds \eqref{eqn:lower_bound_Z_B}, $d_1 \geq \frac12$, and \eqref{eqn:bound_on_calF}, we obtain the integrated-in-time estimate \eqref{eqn:estimate_from_calF} from the preceding bounds \eqref{eqn:auxiliary_estimate_cross}, \eqref{eqn:auxiliary_estimate_commutator_T}, \eqref{eqn:auxiliary_estimate_nonlinear_term}, \eqref{eqn:estimate_theta_error_terms}, \eqref{eqn:auxiliary_estimate_cal_xi}, \eqref{eqn:auxiliary_estimate_modulation}, and \eqref{eqn:auxiliary_estimate_remainder}.

\end{proof}

Next, we transfer the estimates on $\px f_1$ and $f_2$.
\begin{lemma}\label{lemma:LE_transfer_f2}
It holds that \begin{equation}\label{eqn:LE_estimate_from_calT}
\begin{split}
\int_0^T  \big(\|\rho f_2\|^2 +\| \rho \px f_1 \|^2\big)\rho^2&\lesssim A \delta^2 + \int_0^T \big( \| \rho f_1 \|^2 + \| \rho (\px f_1 + c \gamma f_2)\|^2 + |a|^4\big)\\
&\qquad + A^{-2} \bigg[ A \delta^2 + \int_0^T \big( |a|^4 + \| \rho^4 z_1 \|^2 \big) \bigg].
\end{split}
\end{equation}
\end{lemma}

\begin{proof} For this transfer estimate, we use the virial functional
\begin{equation}
    \calT_{\bmf} = \int f_1 f_2 \rho^2,
\end{equation}
which is well-defined. By \eqref{eqn:transfer_f_z} and \eqref{eqn:stability_smallness}, we have 
\begin{equation}\label{eqn:estimate_for_calT}
    |\calT_{\bmf}| \leq C \| f_1 \| \| f_2 \| \leq C \eps^{-\frac32} \delta^2 \leq C A \delta^2.
\end{equation}
Using \eqref{eqn:equations_for_f12}, we compute that 
\begin{align*}
\dot{\calT}_{\bmf} =&\, (1+c^2\gamma^2) \int f_2^2 \rho^2 - \int (\px f_1 + c \gamma f_2)^2 \rho^2 -2c \gamma \int f_1 f_2 (\rho^2)' + \frac12 \int f_1^2 (\rho^2)''- \int P_2 f_1^2 \rho^2\\
&+ \int [f_2 \Omega_1(\bmf) + f_1 \Omega_2(\bmf)]\rho^2 + \int [f_2 \xi_1 + f_1 (\xi_2 + \sfT)]\rho^2+ \int [f_2 X_\eps \sfB N_2 + f_1 X_\eps \sfD N_2]\rho^2.
\end{align*}
We observe that 
\begin{align*}
    \left| 2 c \gamma \int f_1 f_2 (\rho^2)' \right| \leq \frac{1}{2^4} \int f_2^2 \rho^2 + C \int f_1^2 \rho^2,
\end{align*}
and 
\begin{align*}
    \left| \int f_1^2 (\rho^2)'' \right| + \left| \int f_1^2 P_2 \rho^2 \right| \leq C \int f_1^2 \rho^2,
\end{align*}
using $(\rho^2)'' \leq C \rho^2$ and $\|P_2\|_{\infty} \leq C$. Next, we integrate by parts and write 
\begin{align*}
\int\bigl[f_2\Omega_1(f)+f_1\Omega_2(f)\bigr]\rho^2=-(\dot y-c)\gamma\int f_1f_2(\rho^2)'+\dot{c}c\gamma^2\int f_1f_2(x\rho^2)'-\frac{\dot{c}\gamma}{2}\int f_1^2(\rho^2)'.
\end{align*}
Therefore,  using $|(\rho^2)'| \lesssim \rho^2$, $|(x\rho^2)'| \lesssim \rho$, we obtain from \eqref{eqn:stability_smallness}, \eqref{eqn:estimate_modulation_crude} that
\begin{align*}
\left| \int\bigl[f_2\Omega_1(f)+f_1\Omega_2(f)\bigr]\rho^2 \right| \lesssim \delta^2 \int (f_1^2 + f_2^2)\rho.
\end{align*}
Since $\rho^{1/2} \lesssim \zeta_{80}$, we obtain from \eqref{eqn:transfer_zeta_K_f_z} and Proposition~\ref{prop:high_energy_estimate} (with $K \leftarrow 80$) that 
\begin{align*}
\int_0^T \left| \int\bigl[f_2\Omega_1(f)+f_1\Omega_2(f)\bigr]\rho^2 \right| &\lesssim \eps^{-2}\delta^2 \int (z_1^2 + (\px z_1)^2 + z_2^2)\zeta_{80}^2\\
&\lesssim A^{-\frac{15}{2}} \bigg[ A \delta^2 + \int_0^T \big( |a|^4 + \|\rho^4 z_1 \|^2\big)\bigg].
\end{align*}
Next, from \eqref{eqn:estimate_on_xi12},
\begin{align*}
    \left|\int [f_2 \xi_1 + f_1 \xi_2]\rho^2 \right| &\leq C \eps^{-1} \delta \big(|a|^2 + \|\rho^4 z_1\|^2 \big) \big(\| \rho f_1\| + \|\rho f_2\|\big)\\
    &\leq \frac{1}{2^4}\| \rho f_2 \|^2 + C \| \rho f_1 \|^2 + C \eps^{-2} \delta^2 \big( |a|^4 + \|\rho^4 z_1\|^2\big),
\end{align*}
and from \eqref{eqn:estimate_for_sfT},
\begin{align*}
    \left|\int  f_1 \sfT\rho^2 \right| &\leq \frac{1}{2^4}\| \rho f_2 \|^2 +C \| \rho (\px f_1 + c \gamma f_2)\|^2 + C \| \rho f_1 \|^2.
\end{align*}
Finally, as in the proof of the previous lemma, we have the following bounds for the nonlinear terms 
\begin{align*}
    \| \rho X_\eps \sfB N_{m_1} \| + \| \rho X_\eps \sfD N_{m_1} \| &\lesssim |a|^2,\\
    \| \rho X_\eps \sfB \tilN  \| + \sqrt{\eps}\| \rho X_\eps \sfD \tilN  \|  &\lesssim \delta(\|\rho z_1\| + \| \rho \px z_1 \|),
\end{align*}
where we split $N_2 = N_{m_1} + \tilN$. Since $\rho\lesssim \zeta_{40}$, we obtain from Proposition~\ref{prop:high_energy_estimate} (with $K \leftarrow 40$) that 
\begin{align*}
\int_0^T \left|\int [f_2 X_\eps \sfB N_2 + f_1 X_\eps \sfD N_2]\rho^2 \right| &\leq \int_0^T\bigg[ \frac{1}{2^4} \int \rho^2 f_2^2 + C \int \rho^2 f_1^2 + C |a|^4 \bigg]\\
&\quad + C A^{-\frac{15}{2}} \bigg[A \delta^2 + \int_0^T \big(|a|^4 + \|\rho^4 z_1\|^2\big)\bigg].
\end{align*}
Therefore, using \eqref{eqn:estimate_for_calT}, we obtain the integrated-in-time estimate \eqref{eqn:LE_estimate_from_calT} from the preceding bounds.
\end{proof}

In view of the two previous lemmas, we are ready to conclude the proof of Proposition~\ref{prop:low_energy_estimate}.
\begin{proof}[Proof of Proposition~\ref{prop:low_energy_estimate}]
Following \cite[\S4.11.11]{kowalczyk2021sufficient}, we decompose the following term
\begin{align*}
    \int_0^T \int \big( (\px f_1 + c \gamma f_2)^2 + f_1^2\big) \rho^2 &\leq C \rho(A) \int_0^T \int \big( (\px f_1 + c \gamma f_2)^2 + f_1^2\big) \rho\\
    &\quad + C \int_0^T \int_{|x|<A} \big( (\px f_1 + c \gamma f_2)^2 + f_1^2 \big) \rho^2.
\end{align*}
Since $\rho^{\frac12} \simeq \zeta_{80}$ and $\eps^{-2} = A^{\frac12}$, Proposition~\ref{prop:high_energy_estimate} (with $K \leftarrow 80$) and \eqref{eqn:transfer_zeta_K_f_z} yield
\begin{align*}
    \rho(A) \int_0^T \int \big( (\px f_1 + c \gamma f_2)^2 + f_1^2\big) \rho &\lesssim \eps^{-2} \rho(A) \bigg[ A\delta^2 + \int_0^T \big( |a|^4 + \| \rho^4 z_1 \|^2\big)\bigg]\\
    &\lesssim A\delta^2 + \frac{1}{\sqrt{A}}\int_0^T \big( |a|^4 + \| \rho^4 z_1 \|^2\big)\bigg].
\end{align*}

For $B\geq 80$, we have $\frac{\rho^2}{\zeta_B^2} \leq C \rho$, so we obtain 
\begin{align*}
 \int_{|x|<A} \big( (\px f_1 + c \gamma f_2)^2 + f_1^2 \big) \rho^2 \leq \int \Big((\px h_1 + c \gamma h_2)^2 + h_1^2 \rho\Big).
\end{align*}
Applying Lemma~\ref{KMMvdB_lem4.6} (with $w \leftarrow h_1$, $\sigma \leftarrow c \gamma$, $J\leftarrow (x_1,x_2)$, $K \leftarrow 40$, $\omega \leftarrow \kappa$) yields 
\begin{align*}
    \int_0^T \int h_1^2 \rho \lesssim \int_0^T \bigg[ \int (\px h_1 + c \gamma h_2)^2 + \int (\dot{h}_1 - h_2)^2 + \int_{x_1}^{x_2} h_1^2\bigg] + \eps^{-1} \delta^2.
\end{align*}
Moreover,
\begin{align*}
    \dot{h}_1 - h_2 = \chi_A \zeta_B \big( \Omega_1(\bmf) + \xi_1 + X_\eps \sfB N_2 \big).
\end{align*}
After integration by parts, it holds that 
\begin{align*}
    \| \chi_A \zeta_B \Omega_1(\bmf)\| \leq C_B \big(|\dot{y}-c| + |\dot{c}|\big) \| \px f_1 \| \lesssim \eps^{-1} \delta \big(|a|^2 +\|\rho^4 z_1\|\big).
\end{align*}
For the nonlinear terms, we have
\begin{align*}
\| \chi_A \zeta_B  X_\eps \sfB N_2\| \leq \| \chi_A \zeta_B  X_\eps \sfB \tilN\|  +\| \chi_A \zeta_B  X_\eps \sfB N_{m_1}\| \lesssim |a|^2 + \delta \big(\|\zeta_B z_1 \| + \| \zeta_B \px z_1 \|\big).
\end{align*}
Therefore, using these two bounds and \eqref{eqn:estimate_on_xi12}, we obtain from Proposition~\ref{prop:high_energy_estimate} (with $K \leftarrow B$ fixed) that 
\begin{align*}
    \int_0^T \int (\dot{h}_1 - h_2)^2 \lesssim \int_0^T |a|^4 + \eps^{-2} \delta^2 \bigg[ A \delta^2 + \int_0^T\big(|a|^4 + \| \rho^4 z_1 \|^2 \big)\bigg].
\end{align*}
Using $\eps = A^{-\frac14}$, $\delta \leq A^{-4}$, the bound
\[\|\rho \px f_1\|^2 \leq 2 \| \rho (\px f_1 + c \gamma f_2) \|^2 + 2 c^2 \gamma^2 \| \rho f_2 \|^2,\]
and Lemma~\ref{lemma:LE_transfer_f2}, we conclude from the preceding estimates that 
\begin{align*}
    \int_0^T \int \big( f_1^2 + (\px f_1)^2 + f_2^2\big) \rho^2 &\leq C \int_0^T \bigg[ \int (\px h_1 + c \gamma h_2)^2 + \int_{x_1}^{x_2} h_1^2\bigg] \\
    &\quad + CA \delta^2 + C\int_0^T |a|^4 + \frac{C}{\sqrt{A}}\int_0^T \| \rho^4 z_1 \|^2.
\end{align*}

We now use Lemma~\ref{lemma:LE_main_lemma} (with $B\geq B_1$ fixed, with $A$ sufficiently large, then with $\eta$ chosen sufficiently small so that $CC_B(\eta+ A^{-\frac18})< \frac12$) to absorb the term $\|\px h_1 + c \gamma h_2\|_{L_t^2L_x^2}^2 + \| h_1 \|_{L_t^2L_x^2(x_1,x_2)}^2$ and conclude that 
\begin{align*}
\int_0^T \int \big( f_1^2 + (\px f_1)^2 + f_2^2\big) \rho^2 \leq CA \delta^2 + C\int_0^T |a|^4 + \frac{C}{\sqrt{A}}\int_0^T \| \rho^4 z_1 \|^2.
\end{align*}
This finishes the proof of Proposition~\ref{prop:low_energy_estimate}.
\end{proof}

\section{Closing the final estimates}
In this section, we combine the virial and coercivity estimates to obtain a global spacetime bound for the internal mode and the radiation component. Finally, we use the velocity normal form and the modulation equations to complete the proof of the main theorem.
\begin{proposition}\label{prop:main}
For $\delta>0$ sufficiently small, it holds that 
    \begin{equation}\label{eqn:final_estimate}
        \int_0^\infty \big(|a(t)|^4 + \| \rho z_1 \|^2 + \| \rho \px z_1 \|^2 + \| \rho z_2 \|^2 \big)\ud t \leq C \delta^2 < \infty.
    \end{equation}
\end{proposition}

\begin{proof}[Proof of Proposition~\ref{prop:main} and Theorem~\ref{main_thm1}]
Let
\begin{equation}
\calL := \| \rho z_1 \|^2 + \| \rho \px z_1 \|^2 + \| \rho z_2 \|^2 + |a|^4.
\end{equation}
We claim that 
\begin{equation}\label{eqn:estimate_dot_calM}
    |\dot{\calL}| \lesssim \calL.
\end{equation}

Differentiating $\calL$ in time using the system \eqref{eqn:z_1_virialform} and integration by parts yields
\begin{align*}
    \dot{\calL} &= 2 \int \rho^2 (z_1 - f(x,z_1)) z_2 - 2 \int (\rho^2)' (\px z_1)z_2 - 2 c \gamma \int (\rho^2)' z_2^2\\
    &\quad + 2 \int \rho^2 \big( (\px z_1) (\px \beta_1) + z_1 \beta_1 + z_2(\beta_2 + N_{\bmm})\big) + 2 \dot{c} \gamma \int \rho^2(\px z_1)z_2\\
    &\quad + \dot{c} c \gamma^2 \int \Big( \big(x (\rho^2)' - \rho^2\big) (\px z_1)^2 + \big( \rho^2 + x (\rho^2)'\big) (z_1^2+z_2^2\Big)\\
    &\quad - (\dot{y}-c)\gamma \int (\rho^2)' \big( z_1^2 + (\px z_1)^2 + z_2^2 \big) + \frac{\ud}{\ud t}\big( |a|^4\big).
\end{align*}
Notice that the weights satisfy
\begin{equation*}
    |x(\rho^2)'| \leq C \rho, \qquad | (\rho^2)'| \leq C \rho^2.
\end{equation*}
From $|f(x,z_1)|\leq C |z_1|$ and Cauchy-Schwarz inequality, we have
\begin{align*}
    \left| 2 \int \rho^2 (z_1 - f(x,z_1)) z_2 - 2 \int (\rho^2)' (\px z_1)z_2 - 2 c \gamma \int (\rho^2)' z_2^2 \right| \leq C \calL.
\end{align*}
For the source terms,  Lemma~\ref{lemma:bound_on_beta_terms} implies that 
\begin{align*}
    \left| 2 \int \rho^2 \big( (\px z_1) (\px \beta_1) + z_1 \beta_1 + z_2(\beta_2 + N_{\bmm})\big) + 2 \dot{c} \gamma \int \rho^2(\px z_1)z_2\right| \lesssim  \big(|a|^2 + \|\rho^4 z_1 \|^2\big)  \calL^{\frac12} \lesssim \calL.
\end{align*}
By \eqref{eqn:stability_smallness}, \eqref{eqn:estimate_modulation_crude}, \eqref{eqn:estimate_on_dot_a}, and Cauchy-Schwarz inequality, we obtain
\begin{align*}
&\left| \dot{c} c \gamma^2 \int \Big( \big(x (\rho^2)' - \rho^2\big) (\px z_1)^2 + \big( \rho^2 + x (\rho^2)'\big) (z_1^2+z_2^2\Big)\right| \\
&\lesssim |\dot{c}| \int \rho \big(z_1^2 + (\px z_1)^2 + z_2\big)^2 \lesssim \big(|a|^2 + \|\rho^4 z_1 \|^2\big) \delta \calL^{\frac12} \lesssim  \calL.
\end{align*}
Lastly, the bound \eqref{eqn:estimate_on_dot_a} implies that 
\begin{align*}
    \left| \frac{\ud}{\ud t}\big( |a|^4\big)\right| \leq 4 |a|^3 |\dot{a} - i \lambda \gamma^{-1} a| \lesssim |a|^3 \big(|a|^2 + \|\rho^4 z_1 \|^2\big) \lesssim \calL.
\end{align*}
Thus, the claim \eqref{eqn:estimate_dot_calM} has been proven.

We now let $A\gg B \gg 1$ be sufficiently large, and fix $\eps = A^{-\frac14}$. Combining Proposition~\ref{prop:coercivity_estimate} and Proposition~\ref{prop:low_energy_estimate} yields
\begin{equation*}
\int_0^T \| \rho^4 z_1 \|^2 \ud t  \lesssim A \delta^2 + \int_0^T |a|^4 \ud t + \frac{1}{\sqrt{A}} 
\int_0^T \| \rho^4 z_1 \|^2 \ud t .
\end{equation*}
This and Proposition~\ref{prop:internal_mode_estimate} then imply that
\begin{equation*}
    \int_0^T \| \rho^4 z_1 \|^2 \ud t  + \int_0^T |a|^4 \ud t \lesssim A \delta^2. 
\end{equation*}
Fixing such $A$ and taking $\delta \leq A^{-4}$ sufficiently small, we obtain from \eqref{eqn:comparison_rho_zeta_K_40} and \eqref{eqn:estimate_HE} with $K=40$ the estimate 
\begin{equation*}
\int_0^T \calL(t) \,\ud t=  \int_0^T\big( |a(t)|^4 + \| \rho z_1 \|^2 + \| \rho \px z_1 \|^2 + \| \rho z_2 \|^2  \big)\ud t \leq C A\delta^2 \leq C \delta^2.
\end{equation*}
Letting $T \rightarrow +\infty$ completes the proof of  Proposition~\ref{prop:main}. From \eqref{eqn:final_estimate}, there exists a sequence $T_n \rightarrow \infty$ such that 
\begin{equation}
    \lim_{n \rightarrow \infty} \calL(T_n) = \lim_{n \rightarrow \infty} \Big( |a(T_n)|^4 + \|\rho \px z_1(T_n)\|_{L^2}^2 + \|\rho z_1(T_n)\|_{L^2}^2 + \|\rho z_2(T_n)\|_{L^2}^2 \Big)= 0.
\end{equation}
Let $t>0$, and let $n\in \bbN$ be large enough so that $T_n > t$. Integrating \eqref{eqn:estimate_dot_calM} on $[t,T_n]$ and then passing to the limit as $n \rightarrow \infty$, we obtain
\begin{equation*}
    \calL(t) \lesssim \int_t^\infty \calL(s) \,\ud s.
\end{equation*}
Thus, we conclude that $\lim_{t \rightarrow \infty} \calL(t) = 0$.

Finally, by Lemma~\ref{lemma:normal_form_c}, we have $\lim_{t\rightarrow \infty}\Pi[c,\bmz,\bmm]=0$ and $|\dot{c} + \dot{\Pi}[c,\bmz,\bmm]| \lesssim \calL$. It follows from \eqref{eqn:final_estimate} that $\lim_{t \rightarrow \infty} |a|=0$ and there exists $c_+ \in (-1,1)$ such that $\lim_{t\rightarrow \infty}c = c_+$. Moreover, by \eqref{eqn:mod_stability} and \eqref{eqn:estimate_modulation_crude} 
\begin{equation*}
\lim_{t \rightarrow \infty} \dot{y} = c_+, \quad     \gamma_0^2 |\dot{y}-c_0| + \gamma_0^2 |c_+ - c_0| \leq C \delta.
\end{equation*}
By \eqref{eqn:norm_comparisons_nu_z}, for every fixed $R>0$,
\begin{align*}
    \|\bmnu(t,\cdot + y(t))\|_{c(t),R}^2  = \| \bmz(t)\|_{H^1\times L^2(|x|< \gamma(t) R)}^2 \lesssim_R \calL(t).
\end{align*}
Hence, the spectral decomposition \eqref{eqn:spectral_decomposition_u}, the comparison estimate $\gamma(t) \leq 2 \gamma_0$, and Lemma~\ref{lemma:closeness_stability} give 
\begin{align*}
    \| \bmphi(t,\cdot+y(t)) - \bmH_{c_+}\|_{c_+,R} \lesssim_R \calL(t)^{\frac12} + |a(t)| + |c(t)-c_+| \longrightarrow 0, \qquad \text{as $t \longrightarrow \infty$}.
\end{align*}
Applying the same argument by reversing time then completes the proof of Theorem~\ref{main_thm1}. 
\end{proof}

\section{Double sine-Gordon model}
To illustrate the application of our main result, we verify the assumptions of Theorem~\ref{main_thm1} for the
positive-parameter branch of \eqref{eq:dsg}. 
We consider the double sine--Gordon equation
\begin{align}
\label{dsg:dsg}
\begin{gathered}
    \phi_{tt}-\phi_{xx}+W_\eta'(\phi)=0,\\
    W_\eta(\phi):=W_{\mathrm{dsG}}(\phi)
    =\frac{4}{1+4\eta}
       \left[1+\cos\left(\frac{\phi}{2}\right)
                    +\eta(1-\cos\phi)\right], \eta>0.
\end{gathered}
\end{align}
 The key point is that under the double sine-Gordon perturbation, the threshold resonance from the linearized operator for the sine-Gordon equation produces an internal mode.
 
The main result of this section is the following:
\begin{proposition}\label{prop:dsg-main}
For the double sine-Gordon model \eqref{dsg:dsg}, the following statements hold for small $\eta>0$. 
\begin{enumerate}
    \item The kink is given by
    \begin{equation}\label{dsg:exact-kink}
        H_\eta(x) =4\arctan\left(\frac{\sinh (x)}
                                {\sqrt{1+4\eta}}\right).
    \end{equation}
    \item The linearized operator about $H_\eta$ is
    \begin{equation}
        L_{0,\eta} = -\px^2 + P_{0,\eta}(x), \qquad P_{0,\eta}(x) =1 -\frac{2}{4\eta+\cosh^2(x)}
        -\frac{32\eta\sinh^2(x)}
                    {\big(4\eta+\cosh^2(x)\big)^2}.
    \end{equation}
    \item The eigenvalue $\lambda_\eta^2$ admits the asymptotic expansion
    \begin{equation}\label{dsg:internal mode-goal}
        \lambda_\eta^2 =1-\frac{256}{9}\eta^2 +O(\eta^3).
    \end{equation}
This eigenvalue bifurcates from the threshold resonance of the linearized operator for the sine-Gordon model ($\eta=0$). Its full spectrum is $\sigma(L_{0,\eta}) = \{0,\lambda_\eta^2\} \cup [1,\infty)$.
    
    \item The Fermi golden rule coefficient satisfies 
    \begin{equation}\label{dsg:FGR-main-result}
        \Gamma_\eta= \Gamma_0\eta
          +O(\eta^2), \qquad \text{where} \qquad \Gamma_0 := 8\pi\sech\left(\frac{\sqrt3\pi}{2}\right).
    \end{equation}
    \item The second Darboux-transformed potential $P_{2,\eta}$, defined in \eqref{dsg:P2-potential}, is not identically zero and satisfies the repulsivity condition $xP_{2,\eta}'(x) \leq0$ for all $x \in \bbR$, with strict inequality for $x \neq 0$.
\end{enumerate}
\end{proposition}
The proof of Proposition~\ref{prop:dsg-main} is divided into the following subsections.

\subsection{Kink and linearized operator}\label{dsg:kink}
From \eqref{dsg:dsg}, the double angle formulas for the sines and cosines yield the equivalent expression
\begin{align}\label{eqn:W_eta_rewritten}
    W_\eta(\phi)
       =\frac{8}{1+4\eta}\cos^2\left(\frac{\phi}{4}\right)
          \left[1+4\eta\sin^2\left(\frac{\phi}{4}\right)\right],\qquad
    W_\eta(\pm2\pi)=0.
\end{align}
Thus the double-well potential is nonnegative and is strictly positive between the adjacent vacuum states $-2\pi$ and $2\pi$.  
The derivatives of \eqref{dsg:dsg} are
\begin{align}
W_\eta'(\phi)
        &=-\frac{2}{1+4\eta}\sin\left(\frac{\phi}{2}\right)
           +\frac{4\eta}{1+4\eta}\sin(\phi),\notag\\
W_\eta''(\phi)
        &=-\frac{1}{1+4\eta}\cos\left(\frac{\phi}{2}\right)
           +\frac{4\eta}{1+4\eta}\cos(\phi),\notag\\
W_\eta^{(3)}(\phi)
        &=\frac{1}{2(1+4\eta)}\sin\left(\frac{\phi}{2}\right)
           -\frac{4\eta}{1+4\eta}\sin(\phi),\label{dsg:U-third}\\
W_\eta^{(4)}(\phi)
        &=\frac{1}{4(1+4\eta)}\cos\left(\frac{\phi}{2}\right)
           -\frac{4\eta}{1+4\eta}\cos(\phi).\notag
\end{align}
In particular,
\begin{align*}
W_\eta''(-2\pi)&=W_\eta''(2\pi)=1.
\end{align*}
Thus $\omega_- = \omega_+ = \omega = 1$, and both the Klein--Gordon mass and the  bottom of the continuous
spectrum of the linearized spatial operator are exactly $1$.

\begin{figure}[t]
 \centering
 \includegraphics[width=\textwidth]{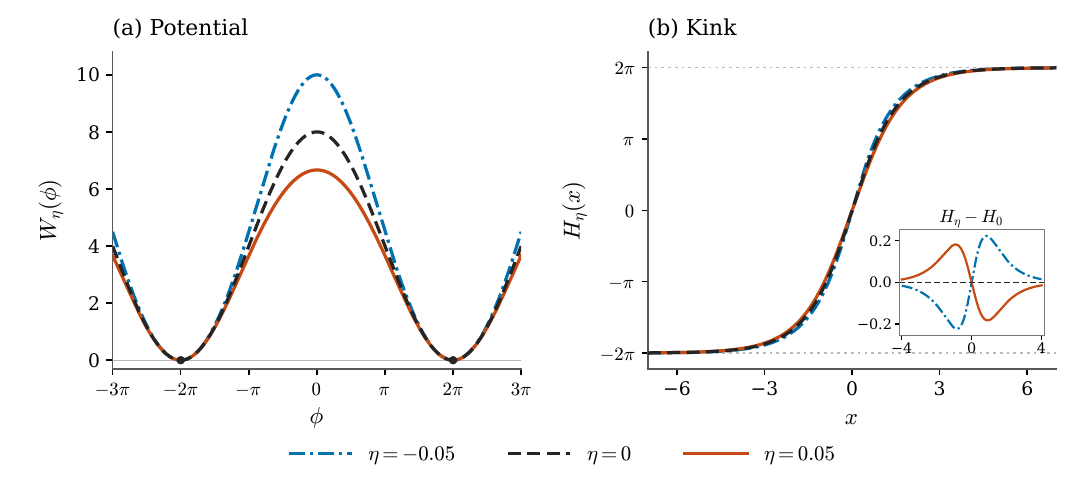}
 \caption{Double sine--Gordon: the potential near the adjacent vacua
 $\pm2\pi$ (left) and the centered kink (right), for
 $\eta=-0.05,0,0.05$. The potential is $4\pi$-periodic, with
 $W_\eta''(\pm2\pi)=1$. The inset shows $H_\eta-H_0$.}
 \label{fig:dsg-potential-kink}
\end{figure}

The kink profile $H_\eta$ is an increasing static solution of
\begin{equation}\label{dsg:kink-ode}
\left\{ \begin{aligned}
-H_\eta''+W_\eta'(H_\eta)&=0\\
   \lim_{x\rightarrow \pm \infty}H_\eta(x)&=\pm 2\pi.
\end{aligned}\right.
\end{equation}
We now solve the above ODE and prove \eqref{dsg:exact-kink}. First, multiplying \eqref{dsg:kink-ode} by $H_\eta'$ and using the limits at spatial infinity gives 
\begin{align*}
\frac12(H_\eta')^2
       &=W_\eta(H_\eta),
\end{align*}
and we obtain from \eqref{eqn:W_eta_rewritten} the identity
\begin{align*}
H_\eta'
       &=\frac{4}{\sqrt{1+4\eta}}
           \cos\left(\frac{H_\eta}{4}\right)
           \sqrt{1+4\eta
                       \sin^2\left(\frac{H_\eta}{4}\right)}.
\end{align*}
With the substitution $v=-\sin(H_\eta/4)$ , we obtain
\begin{align*}
v'&=-\frac{1}{\sqrt{1+4\eta}}
                (1-v^2)\sqrt{1+4\eta v^2},
    \qquad v(0)=0.
\end{align*}
Solving the above ODE using separation of variables yields
\begin{align*}
v(x)&=-\frac{\sinh (x)}
    {\sqrt{4\eta+\cosh^2(x)}}.
\end{align*}
Then the explicit formula \eqref{dsg:exact-kink} follows directly. Figure~\ref{fig:dsg-potential-kink} illustrates the double-well potential and the centered kink.
The sine--Gordon model and a negative parameter value are included
for comparison.

\noindent\textit{The linearized operator.}  Using the kink formula \eqref{dsg:exact-kink}, we define the static linearized operator 
\begin{align*}
L_{0,\eta}&:=-\partial_x^2+1+Q_{\eta}(x),\\
Q_{\eta}(x)&:=W_\eta''(H_\eta(x))-1=-\frac{2
         \big(1+4\eta+(1+16\eta)\sinh^2(x)\big)}
         {(4\eta+\cosh^2(x))^2}.
\end{align*}
Equivalently, we have 
\begin{align*}
Q_{\eta}(x)&=-\frac{2}{4\eta+\cosh^2(x)}-\frac{32\eta\sinh^2(x)}{\big(4\eta+\cosh^2(x)\big)^2}.
\end{align*}
Since $|Q_{\eta}(x)| \leq C e^{-2|x|}$, it follows that $\sigma_{ac}(L_{0,\eta})=[1,\infty).$ Differentiating the kink equation gives the translation mode
\begin{align}
\label{dsg:translation-mode}
L_{0,\eta}Y_{0,\eta}&=0,\qquad Y_{0,\eta}(x):=H_\eta'(x)=\frac{4\sqrt{1+4\eta}\cosh (x)}{4\eta+\cosh^2(x)}.
\end{align}
Since $Y_{0,\eta}>0$, the Sturm--Liouville theorem implies that
zero is a simple eigenvalue and $L_{0,\eta}\geq0$.

In the sine--Gordon limit $\eta=0$, we write
\begin{align*}
H_{\mathrm{SG}}(x)&:=H_\eta(x)|_{\eta=0} =4\arctan(\sinh x),
\end{align*}
and
\begin{align*}
L_{\mathrm{SG}}:=L_{0,\eta}|_{\eta=0}&=-\partial_x^2+1-2\sech^2(x).
\end{align*}
The linearized operator $L_{\mathrm{SG}}$ features the first member of the family of P\"oschl-Teller potentials \cite{PoschlTeller}, and its spectrum is explicit. The function $Y_{\mathrm{SG}}(x):=\tanh (x)$ is an odd threshold resonance,
\begin{align*}
L_{\mathrm{SG}}Y_{\mathrm{SG}}&=Y_{\mathrm{SG}},\qquad
    Y_{\mathrm{SG}}\in L_x^\infty(\mathbb R)\setminus L^2_x(\mathbb R).
\end{align*}
We will prove in the next subsection that for sufficiently small $\eta>0$, the threshold resonance $Y_{\mathrm{SG}}(x)$ bifurcates into
a simple odd internal mode $Y_\eta(x)$, with eigenvalue
$\lambda_\eta^2$,
\begin{align}
\label{dsg:internal mode}
L_{0,\eta}Y_\eta&=\lambda_\eta^2Y_\eta,
    \qquad 0<\lambda_\eta^2<1.
\end{align}

\subsection{Bifurcation of the threshold resonance into an internal mode eigenvalue}\label{dsg:darboux}
Since $H_\eta'(x)>0$, we define its logarithmic derivative
$V_{0,\eta}$ and the Darboux transform $U_{0,\eta}$ by
\begin{align}\label{dsg:Q-explicit}
V_{0,\eta}(x)&:=\frac{H_\eta''(x)}{H_\eta'(x)}=\tanh (x)
      \frac{4\eta-\cosh^2(x)}
      {4\eta+\cosh^2(x)},
    \qquad
    U_{0,\eta}
    :=H_\eta'\,\partial_x\,(H_\eta')^{-1}
    =\partial_x-V_{0,\eta}(x).
\end{align}
Hence,
\begin{align*}
L_{0,\eta}
    &=U_{0,\eta}^*U_{0,\eta}.
\end{align*}
Moreover, the supersymmetric partner operator is given by 
\begin{align*}
L_{1,\eta} &:=U_{0,\eta}U_{0,\eta}^* =-\partial_x^2+1+Q_{1,\eta}(x),
\end{align*}
where 
$$
Q_{1,\eta}(x)
    =-\frac{8\eta}
    {\cosh^2(x) \left(4\eta+\cosh^2(x)\right)}.$$
We also have the intertwining relations
\begin{align}
\label{dsg:first-intertwining}
U_{0,\eta}L_{0,\eta}
      &=L_{1,\eta} U_{0,\eta},
    \qquad
    L_{0,\eta} U_{0,\eta}^*
      =U_{0,\eta}^*L_{1,\eta}.
\end{align}
Since $U_{0,\eta}H_\eta'=0$, the first Darboux transform removes the translation eigenfunction.

We now prove the eigenvalue expansion in \eqref{dsg:internal mode-goal} using the Birman-Schwinger principle. Our proof follows the outline of \cite[Lemma~2]{martel2024asymptotic}. Define 
\begin{align}
\label{dsg:kappa-definition}
\kappa_\eta^2&:=1-\lambda_\eta^2.
\end{align}
By \eqref{dsg:first-intertwining}, every positive discrete eigenvalue of $L_{0,\eta}$ corresponds, with the same multiplicity, to a discrete eigenvalue of $L_{1,\eta}$. More precisely, if
$L_{0,\eta} Y_\eta
=\lambda_\eta^2Y_\eta$, then
$Y_{1,\eta}=U_{0,\eta}Y_\eta$ satisfies
\begin{align}
\label{dsg:Darboux-reduced-eigenvalue}
\left(-\partial_x^2+Q_{1,\eta}\right)Y_{1,\eta}
    &=-\kappa_\eta^2Y_{1,\eta}.
\end{align}
Conversely, for $\kappa_\eta^2<1$, a nonzero solution of
\eqref{dsg:Darboux-reduced-eigenvalue} gives an eigenfunction of
$L_{0,\eta}$ through
\begin{align}
\label{dsg:recover-Y-from-Z}
Y_\eta
    &=\frac{1}{1-\kappa_\eta^2}
      U_{0,\eta}^*Y_{1,\eta}.
\end{align}
To make the Taylor expansion more explicit, we factor out $Q_{1,\eta}(x)=\eta P_\eta(x)$, where
\begin{align}
\label{dsg:P-epsilon-definition}
    P_\eta(x) :=-\frac{8}{\cosh^2(x)\left(4\eta+\cosh^2(x)\right)}.
\end{align}
Taylor expansion of $P_\eta$ about $\eta=0$ gives 
\begin{align}\label{dsg:P_eta}
P_\eta = P_0+\eta P_1+O(\eta^2),\qquad P_0(x):=-8\sech^4(x),\qquad P_1(x):=32\sech^6(x).
\end{align}
The remainder $O(\eta^2)$ satisfies the following weighted estimate
\begin{align*}
\int_{\mathbb R}(1+x^2)
    \left|P_\eta(x)-P_0(x)-\eta P_1(x)\right|\,\ud x
    &\leq C\eta^2.
\end{align*}
Set $S_\eta:=\sqrt{|P_\eta|}$. Since $P_\eta=-S_\eta^2$, the eigenvalue equation \eqref{dsg:Darboux-reduced-eigenvalue} is equivalent to
\begin{align*}
(-\px^2 + \kappa_\eta^2) Y_{1,\eta}
    &=\eta S_\eta^2Y_{1,\eta}.
\end{align*}
For $\kappa>0$, recall that the Green's function is given by 
$
(-\px^2 + \kappa^2)^{-1}(x,y)
    =\frac{e^{-\kappa|x-y|}}{2\kappa}.
$ 
Given $\kappa>0$, setting
$\Phi_\eta:=S_\eta Y_{1,\eta}$ turns the reduced
spectral problem into the Birman--Schwinger equation
\begin{align}
\label{dsg:Birman-Schwinger-equation}
\Phi_\eta
    &=\eta\mathcal K_{\eta,\kappa}\Phi_\eta,
    \qquad
    \mathcal K_{\eta,\kappa}
    :=S_\eta (-\px^2 + \kappa^2)^{-1}S_\eta.
\end{align}
The operator $\mathcal K_{\eta,\kappa}$ is positive and compact on
$L^2_x(\mathbb R)$, and the Birman--Schwinger principle preserves
multiplicity of eigenvalues. To isolate the singular part of the
free resolvent, we expand the Green's function
\begin{align}
\label{dsg:resolvent-decomposition}
\frac{e^{-\kappa|x-y|}}{2\kappa}
    &=\frac{1}{2\kappa}+N_\kappa(x,y),
    \qquad
    N_\kappa(x,y)
    :=\frac{e^{-\kappa|x-y|}-1}{2\kappa}.
\end{align}
Thus the integral kernel of $\mathcal K_{\eta,\kappa}$ is expanded as 
\begin{align}
\label{dsg:BS-decomposition}
\mathcal K_{\eta,\kappa}(x,y) &=\frac{S_\eta(x)S_\eta(y)}{2\kappa} e^{-\kappa|x-y|} =\frac{S_\eta(x)S_\eta(y)}{2\kappa} +M_{\eta,\kappa}(x,y),
\end{align}
where the remainder kernel is
\begin{align*}
M_{\eta,\kappa}(x,y)
    &=S_\eta(x)N_\kappa(x,y)S_\eta(y).
\end{align*}
The elementary bound
\begin{align*}
|N_\kappa(x,y)|&\leq\frac12|x-y|
\end{align*}
and the exponential decay of $S_\eta$ show that
 $M_{\eta,\kappa}$ is Hilbert-Schmidt on $L^2_x(\mathbb R)$ for small
$|\eta|$ and all $\kappa\geq0$.  Furthermore,
\begin{align*}
\|M_{\eta,\kappa}-M_{0,0}\|_{L^2_{x,y}(\mathbb R^2)}&\longrightarrow0
    \quad\text{as }(\eta,\kappa)\longrightarrow(0,0),
\end{align*}
where
\begin{align*}
M_{0,0}(x,y)
    &=-\frac12 S_0(x) |x-y| S_0(y).
\end{align*}
For sufficiently small $\eta\geq0$, the operator
$I-\eta M_{\eta,\kappa}$ is therefore invertible, uniformly
for $\kappa >0$.  Substituting \eqref{dsg:BS-decomposition} into
\eqref{dsg:Birman-Schwinger-equation} gives
\begin{align*}
(I-\eta M_{\eta,\kappa})\Phi_\eta
    &=\frac{\eta}{2\kappa}S_\eta
      \langle \Phi_\eta,S_\eta\rangle.
\end{align*}
Applying $(I-\eta M_{\eta,\kappa})^{-1}$ and then taking the
inner product with $S_\eta$ shows that a nontrivial solution $\Phi_\eta$ exists
if and only if
\begin{align}
\label{dsg:kappa-scalar-equation}
2\kappa
    &=\eta\mathcal A(\eta,\kappa),
    \qquad
    \mathcal A(\eta,\kappa)
    :=\left\langle (I-\eta M_{\eta,\kappa})^{-1}
      S_\eta,S_\eta\right\rangle,
\end{align}
and the resulting eigenspace is one-dimensional.
Define the following functional
\begin{align*}
F(\eta,\kappa)
&:=2\kappa-\eta\mathcal A(\eta,\kappa).
\end{align*}
Then $F(\eta,\kappa)$ is continuously differentiable near
$(0,0)$, with
$F(0,0)=0$ and $\partial_\kappa F(0,0)=2$. Hence the implicit function theorem gives a unique solution for \eqref{dsg:Birman-Schwinger-equation} for  $\kappa=\kappa_\eta$ near $(\eta,\kappa)=(0,0)$. To obtain its
leading-order expansion, first observe from
\eqref{dsg:kappa-scalar-equation} and the boundedness of $\mathcal A$
that $\kappa_\eta=O(\eta)$. We directly compute the trace of $\calA(\eta,\kappa)$ at $(0,0)$,
\begin{align*}
\mathcal A(0,0)&=\|S_0\|_{L^2_x}^2 =8\int_{\bbR}\sech^4(x)\,\ud x =\frac{32}{3}.
\end{align*}  The continuous differentiability
of $\mathcal A$ then gives
\begin{align*}
\mathcal A(\eta,\kappa_\eta)
    &=\mathcal A(0,0)
     +O(\eta+|\kappa_\eta|)
    =\frac{32}{3}+O(\eta).
\end{align*}
Substitution into \eqref{dsg:kappa-scalar-equation} yields
$2\kappa_\eta
=\eta(32/3+O(\eta))$, and therefore
\begin{align*}
\kappa_\eta
    &=\frac{16}{3}\eta+O(\eta^2).
\end{align*}
Hence
\begin{align*}
\lambda_\eta^2
    &=1-\kappa_\eta^2
    =1-\frac{256}{9}\eta^2
      +O(\eta^3),
\end{align*}
which proves expansion \eqref{dsg:internal mode-goal}. The uniqueness of the solution $\Phi_\eta$ of \eqref{dsg:Birman-Schwinger-equation} (or equivalently \eqref{dsg:kappa-scalar-equation}) also follows from the implicit function theorem using the uniform bound for $\|(I-\eta M_{\eta,\kappa})^{-1}\|_{L^2_x(\mathbb R)\rightarrow L^2_x(\mathbb R)}$ and the uniform bound for $\mathcal A(\eta,\kappa)$. Since $Y_{1,\eta}$ is positive and even, by \eqref{dsg:recover-Y-from-Z}, the
corresponding eigenfunction $Y_\eta$ of $L_{0,\eta}$ is simple and odd.

\subsection{Expansion of the internal mode eigenfunction}\label{dsg:eigenfunction}
Before we compute the Fermi golden rule coefficient, we need to  expand the internal eigenfunction $Y_{\eta}$ in terms of $\eta$.
For all sufficiently small $\eta>0$, let $Y_{1,\eta}$ be the positive even eigenfunction from  \eqref{dsg:Darboux-reduced-eigenvalue}, and we consider the following normalization
\begin{align}\label{dsg:internal-normalization}
-\frac{\eta}{2\kappa_\eta}
       \int_{\mathbb R}P_\eta Y_{1,\eta}\,\ud x&=1,
    \qquad
    Y_\eta
       =(\lambda_\eta^2)^{-1}U_{0,\eta}^*Y_{1,\eta}.
\end{align}
We claim that 
\begin{align}
\label{dsg:Y1-expansion-final}
\begin{aligned}
Y_\eta(x)=\tanh (x)+\eta\biggl[\frac{16}{3}\tanh (x)-4\tanh (x)\sech^2(x) -\frac{16}{3}\tanh (x)\log(2\cosh (x))
          \biggr]+R_\eta^Y(x),
\end{aligned}
\end{align}
where the remainder satisfies
\begin{align*}
|R_\eta^Y(x)|&\leq C\eta^2\langle x\rangle^2,
    \qquad
    |(R_\eta^Y)'(x)|\leq C\eta^2\langle x\rangle,
\end{align*}
for some absolute constant $C$.

\begin{proof}[Proof of \eqref{dsg:Y1-expansion-final}]
We rewrite $(-\partial_x^2+\kappa_\eta^2)Y_{1,\eta}
=-\eta P_\eta Y_{1,\eta}$ as a resolvent equation
\begin{align*}
Y_{1,\eta}(x)
       &=-\frac{\eta}{2\kappa_\eta}
          \int_{\mathbb R}e^{-\kappa_\eta|x-y|}
             P_\eta(y)Y_{1,\eta}(y)\,\ud y.
\end{align*}
Since  $-\eta P_\eta Y_{1,\eta}/
(2\kappa_\eta)$ is positive and has unit integral by our normalization \eqref{dsg:internal-normalization}, we obtain the preliminary bounds
\begin{align*}
0&<Y_{1,\eta}\leq1,\qquad
    |Y_{1,\eta}'|\leq\kappa_\eta\leq C\eta.
\end{align*}
Splitting the resolvent kernel as in
\eqref{dsg:resolvent-decomposition} and using the normalization gives
\begin{align*}
Y_{1,\eta}
       &=1-\eta N_{\kappa_\eta}
                         (P_\eta Y_{1,\eta}).
\end{align*}
Since $|N_\kappa(x,y)|\leq |x-y|/2$, we obtain
\begin{align*}
|Y_{1,\eta}(x)-1|
       &\leq\frac{\eta}{2}
            \int_{\mathbb R}|x-y|\,|P_\eta(y)|\,\ud y
       \leq C\eta\langle x\rangle.
\end{align*}
Let
\begin{align*}
\delta_\eta&:=P_\eta Y_{1,\eta}-P_0 =(P_\eta-P_0)Y_{1,\eta}+P_0(Y_{1,\eta}-1).
\end{align*}
From the preceding bounds on $Y_{1,\eta}$, we obtain 
\begin{align*}
|\delta_\eta(x)|
       &\leq C\eta\langle x\rangle e^{-C|x|}.
\end{align*}
We now have the following expansion for $Y_{1,\eta}$, with 
\begin{align*}
Y_{1,\eta}=1-\eta N_{\kappa_\eta}
                           (P_\eta Y_{1,\eta})=1-\eta N_{\kappa_\eta}
                           (P_0+\delta_\eta)=1-\eta A+R_\eta,
\end{align*}
where  $N_0(x,y)=-\frac12|x-y|$, and
\begin{align*}
A(x)&:=(N_0P_0)(x)
       =-\frac12\int_{\mathbb R}|x-y|P_0(y)\,\ud y ,\\
R_\eta(x) &:=-\eta\bigl[(N_{\kappa_\eta}-N_0)P_0\bigr](x) -\eta(N_{\kappa_\eta}\delta_\eta)(x) =Y_{1,\eta}(x)-1+\eta A(x).\end{align*}
Direct computation using \eqref{dsg:P_eta} gives
\begin{align*}
A(x)&=\frac{16}{3}\log(2\cosh x)-\frac43\sech^2(x).
\end{align*}
Next, we estimate $R_\eta$, beginning with $R_\eta(0)$.
Since
$|e^{-u}-1+u|\leq u^2/2$ for $u\geq0$, we have
\begin{align*}
Y_{1,\eta}(0) &=1+\frac{\eta}{2} \int_{\mathbb R}|y|P_\eta(y)Y_{1,\eta}(y)\,\ud y +O(\eta\kappa_\eta) =1-\eta A(0)+O(\eta^2).
\end{align*}
Thus $|R_\eta(0)|\leq C\eta^2$. The remaining estimates follow directly from the differential equation.
Since $A''=-P_0$ and
\begin{align*}
R_\eta''
       &=\kappa_\eta^2Y_{1,\eta}
           +\eta\delta_\eta,
\end{align*}
we obtain $|R_\eta''|\leq C\eta^2$. Moreover, $R_\eta'(0)=0$ by evenness. Integrating this bound from $0$ gives
\begin{align*}
|R_\eta(x)|&\leq C\eta^2\langle x\rangle^2,
    \quad
    |R_\eta'(x)|\leq C\eta^2\langle x\rangle,
    \quad
    |R_\eta''(x)|\leq C\eta^2.
\end{align*}
We now transfer the preceding bounds to $Y_\eta$ from $Y_{1,\eta}$ using Darboux transform. By \eqref{dsg:first-intertwining},
\begin{align}
\label{dsg:Z-eigenvalue-equation}
L_{1,\eta} Y_{1,\eta}
      &=\lambda_\eta^2Y_{1,\eta},
    \qquad
    Y_\eta
      =\frac{1}{\lambda_\eta^2}
       U_{0,\eta}^*Y_{1,\eta}.
\end{align}
Rewriting the logarithmic derivative from \eqref{dsg:Q-explicit} as
\begin{align*}
V_{0,\eta}(x) &=-\tanh (x) +\frac{8\eta\tanh (x)\sech^2(x)}{1+4\eta\sech^2(x)}
\end{align*}
gives the Taylor expansion of $V_{0,\eta}$ at $\eta=0$ with
\begin{align}\label{dsg:V_eta_expansion}
V_{0,\eta}(x)&=-\tanh x+\eta V^{(1)}(x)+R_\eta^V(x),
\end{align}
where
\begin{align*}
V^{(1)}(x)&=8\tanh (x)\sech^2(x),
\end{align*}
and a remainder that satisfies
\begin{align*}
|R_\eta^V(x)|+|\px R_\eta^{V}(x)|&\leq C\eta^2.
\end{align*}

Applying $U_{0,\eta}^*=-\partial_x-V_{0,\eta}$ now gives
\begin{align*}
U_{0,\eta}^*Y_{1,\eta}&=\tanh (x)+\eta(A'-A\tanh (x)-V^{(1)})+E_\eta,
\end{align*}
with
\begin{align*}
E_\eta
       &=-R_\eta'+R_\eta\tanh x
          -\eta V^{(1)}(Y_{1,\eta}-1)
          -R_\eta^V Y_{1,\eta}
\end{align*}
satisfying the error bounds 
\begin{align*}
|E_\eta(x)|&\leq C\eta^2\langle x\rangle^2,
    \qquad |E_\eta'(x)|\leq C\eta^2\langle x\rangle.
\end{align*}
Substituting the preceding expansion of
$U_{0,\eta}^*Y_{1,\eta}$ and using
$(\lambda_\eta^2)^{-1}=1+O(\eta^2)$, we obtain the expansion
\begin{align}\label{dsg:Y_eta_expansion}
Y_\eta(x)&=\tanh (x)+\eta B(x)+R_\eta^Y(x),
\end{align}
with
\begin{align*}
B(x)&=A'(x)-A(x)\tanh (x)-V^{(1)}(x)=\frac{16}{3}\tanh (x)-4\tanh (x)\sech^2(x) -\frac{16}{3}\tanh (x)\log(2\cosh x),\end{align*}
and
\begin{align*}
R_\eta^Y(x)
       &=\frac{E_\eta(x)}{\lambda_\eta^2}
        +\left(\frac{1}{\lambda_\eta^2}-1\right)
                   \bigl[\tanh x+\eta B(x)\bigr].
\end{align*}
Since $(\lambda_\eta^2)^{-1}$ is uniformly bounded and
$(\lambda_\eta^2)^{-1}-1=O(\eta^2)$, the bounds above for
$E_\eta,E_\eta'$ and the explicit formula for $B$ yield
\begin{align*}
|R_\eta^Y(x)|&\leq C\eta^2\langle x\rangle^2,
    \qquad
    |(R_\eta^Y)'(x)|\leq C\eta^2\langle x\rangle,
\end{align*}
which proves \eqref{dsg:Y1-expansion-final}.
\end{proof}

\subsection{Verifying the Fermi golden rule}\label{dsg:FGR}

We retain the normalization of $Y_\eta$ from \eqref{dsg:internal-normalization} and take $\eta>0$ sufficiently small throughout.  At the continuum energy $4\lambda_\eta^2$, the eigenvalue expansion gives
\begin{align*}
4\lambda_\eta^2-1
       &=3-\frac{1024}{9}\eta^2+O(\eta^3),
\end{align*}
and hence
\begin{align*}
\sqrt{4\lambda_\eta^2-1}
       &=\sqrt3+O(\eta^2).
\end{align*}
Let $g_\eta$ be the solution of
\begin{align}
L_{0,\eta} g_\eta&=4\lambda_\eta^2g_\eta,
    \qquad
g_\eta(0)=0,\qquad
    g_\eta'(0)=4\lambda_\eta^2.\label{dsg:continuum-normalization}
\end{align}
We then define the Fermi golden rule constant
\begin{align*}
\Gamma_\eta
       &:=-\frac12\int_{\mathbb R}
          W_\eta^{(3)}(H_\eta(x))
          (Y_\eta(x))^2g_\eta(x)\,\ud x.
\end{align*}
In the remaining part of this subsection, we compute $\Gamma_\eta$.
\subsubsection{Constructing the generalized eigenfunction}\label{dsg:continuum}

Let $h_\eta$ be the solution of
\begin{equation}\label{dsg:FGR-partner-IVP}
\left\{\begin{aligned}
    -h_\eta''+\eta P_\eta h_\eta
        =(4\lambda_\eta^2-1)h_\eta,\\
    h_\eta(0)=1,\qquad h_\eta'(0)=0.
\end{aligned}\right.
\end{equation}
Since $P_\eta$ is smooth and bounded, standard ODE theory guarantees a unique solution $h_\eta$. The solution $h_\eta$ is even, and  satisfies $
g_\eta=U_{0,\eta}^*h_\eta.$  
By taking the adjoint of $U_{0,\eta}$, we obtain 
\begin{align*}
\Gamma_\eta &=-\frac12\int_{\mathbb R} W_\eta^{(3)}(H_\eta)Y_\eta^2 U_{0,\eta}^*h_\eta\,\ud x =-\frac12\int_{\mathbb R} U_{0,\eta}\!\left[ W_\eta^{(3)}(H_\eta)Y_\eta^2\right] h_\eta\,\ud x.
\end{align*}
We now estimate $h_\eta$.  For $x\geq0$, we formulate \eqref{dsg:FGR-partner-IVP} as a Volterra integral equation
\begin{align}
\label{dsg:FGR-Volterra}
\begin{aligned}
    h_\eta(x)=\cos\!\left(\sqrt{4\lambda_\eta^2-1}\,x\right)+\frac{\eta}{\sqrt{4\lambda_\eta^2-1}}
            \int_0^x\sin\!\left(
                  \sqrt{4\lambda_\eta^2-1}\,(x-y)\right)
                  P_\eta(y)h_\eta(y)\,\ud y.
\end{aligned}
\end{align}
The uniform $L^1$ bound for $P_\eta$ and Gronwall's inequality then imply that 
\begin{align*}
|h_\eta(x)|
       &\leq\exp\left(
                    \frac{\eta}{\sqrt{4\lambda_\eta^2-1}}
                    \int_0^x|P_\eta(y)|\,\ud y\right)
       \leq C.
\end{align*}
Differentiating \eqref{dsg:FGR-Volterra} also gives
\begin{align*}
h_\eta'(x)=-\sqrt{4\lambda_\eta^2-1}\,
              \sin\!\left(\sqrt{4\lambda_\eta^2-1}\,x\right)+\eta\int_0^x\cos\!\left(
                 \sqrt{4\lambda_\eta^2-1}\,(x-y)\right)
                        P_\eta(y)h_\eta(y)\,\ud y.
\end{align*}
The preceding identities imply the following estimates
\begin{align}
\label{dsg:FGR-h-bounds}
\begin{aligned}
    \|h_\eta\|_{L^\infty_x(\mathbb R)}+\|h_\eta'\|_{L^\infty_x(\mathbb R)}&\leq C,\\
    \left|h_\eta(x)-\cos\!\left(
         \sqrt{4\lambda_\eta^2-1}\,x\right)\right|
        &\leq C\eta,\\
    \left|h_\eta'(x)
         +\sqrt{4\lambda_\eta^2-1}\,
          \sin\!\left(
          \sqrt{4\lambda_\eta^2-1}\,x\right)\right|
        &\leq C\eta.
\end{aligned}
\end{align}
Let $h_0(x)=\cos(\sqrt3x)$. Using
$|\sqrt{4\lambda_\eta^2-1}-\sqrt3|
\leq C\eta^2$, we obtain the 
weighted estimate
\begin{align}
\label{dsg:FGR-h-comparison}
|h_\eta(x)-h_0(x)|
      +|h_\eta'(x)-h_0'(x)|
       &\leq C\eta\langle x\rangle.
\end{align}

\subsubsection{Expansion of the localized forcing}\label{dsg:forcing}

Define
\begin{align*}
G_\eta
       &:=U_{0,\eta}\!\left[
          W_\eta^{(3)}(H_\eta)Y_\eta^2\right],
\end{align*}
so that
\begin{align}
\label{dsg:FGR-transformed-amplitude}
\Gamma_\eta
       &=-\frac12\int_{\mathbb R}G_\eta(x)h_\eta(x)\,\ud x.
\end{align}
We now expand $G_\eta(x)$. From \eqref{dsg:V_eta_expansion} and \eqref{dsg:Y_eta_expansion}, we have
\begin{align*}
|Y_\eta(x)-\tanh (x)| +|(Y_\eta)'(x)-\sech^2(x)|  &\leq C\eta\langle x\rangle.
\end{align*}
Next, we expand $W_\eta^{(3)}(H_\eta)$. A Taylor expansion of the kink formula \eqref{dsg:exact-kink} gives
\begin{align}
\label{dsg:K-epsilon-expansion}
H_\eta(x)=4\arctan(\sinh x)-8\eta\sech (x)\tanh (x)+R_\eta^H(x),
\end{align}
where the remainder satisfies
\begin{align*}
\big|R_\eta^H(x)\big|+\big|(R_\eta^H)'(x)\big| &\leq C\eta^2e^{-|x|/2}.
\end{align*}

Using \eqref{dsg:U-third} and \eqref{dsg:K-epsilon-expansion} gives
\begin{align*}
W_\eta^{(3)}(H_\eta(x))&=\sech (x)\tanh (x) +\big(\sech (x)\tanh (x)(14-36\sech^2(x))\big)\eta +R_\eta^W,
\end{align*}
where the remainder satisfies
\begin{align*}
|R_\eta^W(x)|+|(R_\eta^W)'(x)|
       &\leq C\eta^2e^{-C|x|},
\end{align*}
for some fixed $c>0$. Now write $F_\eta:=W_\eta^{(3)}(H_\eta(x))(Y_\eta)^2$. Combining the expansion of $W_\eta^{(3)}(H_\eta(x))$ with that of $Y_\eta$, we obtain
\begin{align*}
F_\eta&=F_0+\eta F_1+R_\eta^F,
\end{align*}
with $
F_0(x)=\sech (x)\tanh^3(x)$, 
and
\begin{align*}
F_1(x) =\sech (x)\tanh^3(x)\left(\frac{74}{3}-44\sech^2(x)-\frac{32}{3}\log(2\cosh x)\right),
\end{align*}
and  a remainder term that satisfies $|R_\eta^F(x)|+|(R_\eta^F)'(x)|\leq C\eta^2e^{-C|x|}.$

Finally, from $G_\eta=U_{0,\eta}F_\eta$ 
we obtain the expansion 
\begin{align}
\label{dsg:FGR-G-expansion}
\begin{aligned}
    G_\eta&=G_0+\eta G_1+R_\eta^G,\\
    G_0&=(\partial_x+\tanh x)F_0=3\big(\sech^3(x)-\sech^5(x)\big),\\
    G_1&=(\partial_x+\tanh x)F_1-V^{(1)}F_0\\
       &=-\frac{32}{3}\sech (x)+\frac{526}{3}\sech^3(x)
          -\frac{1130}{3}\sech^5(x)\\
       &\quad+212\sech^7(x)-32\big(\sech^3(x)-\sech^5(x)\big)\log(2\cosh x).
\end{aligned}
\end{align}
The remainder satisfies $
|R_\eta^G(x)| \leq C\eta^2e^{-C|x|}$.

\subsubsection{Computation of the FGR constant}
To evaluate \eqref{dsg:FGR-transformed-amplitude}, we use a localized solution of
an inhomogeneous free equation to account for the
first-order correction to $h_\eta$. Let
\begin{align*}
\chi(x)&:=-\frac14\sech^3(x).
\end{align*}
Then
\begin{align*}
(-\partial_x^2-3)\chi&=3(\sech^3x-\sech^5x)=G_0.
\end{align*}
Using integration by parts and the equation \eqref{dsg:FGR-partner-IVP}, we obtain that 
\begin{align*}
\Gamma_\eta
       &=-\frac12\int_{\mathbb R}(G_\eta-G_0)h_\eta\,\ud x
          -\frac12\int_{\mathbb R}G_0h_\eta\,\ud x\\
&= -\frac12\int_{\mathbb R}(G_\eta-G_0)h_\eta\,\ud x -\frac12 \int_{\mathbb R}(-\chi''-3\chi)h_\eta\,\ud x\\
&= -\frac12\int_{\mathbb R}(G_\eta-G_0)h_\eta\,\ud x -\frac12  \int_{\mathbb R}
            \chi(4\lambda_\eta^2-4-\eta P_\eta)
            h_\eta\,\ud x\\
& =-\frac12\int_{\mathbb R}
          \left[G_\eta-G_0
            +\chi(4\lambda_\eta^2-4-\eta P_\eta)\right]
                  h_\eta\,\ud x.
\end{align*}
From \eqref{dsg:FGR-G-expansion} we obtain an expansion for the terms in the bracketed parenthesis for the above expression,
\begin{align*}
G_\eta-G_0
       +\chi(4\lambda_\eta^2-4-\eta P_\eta)
       &=\eta\calG+ \calR_\eta,
\end{align*}
where
\begin{align}
\label{dsg:FGR-H-definition}
\calG&:=G_1-\chi P_0,
\end{align}
and the remainder is defined exactly by
\begin{align*}
\calR_\eta
       &=R_\eta^G
          +\chi(4\lambda_\eta^2-4)
          -\eta\chi(P_\eta-P_0).
\end{align*}
Recalling that
\begin{align*}
|P_\eta(x)-P_0(x)|
       &\leq C\eta\langle x\rangle e^{-C|x|},
\end{align*}
we obtain
\begin{align*}
\|\calR_\eta\|_{L_x^1} &\leq\|R_\eta^G\|_{L_x^1} +|4\lambda_\eta^2-4|\,\|\chi\|_{L^1} +\eta\int_{\mathbb R} |\chi(x)|\,|P_\eta(x)-P_0(x)|\,\ud x\leq C\eta^2.
\end{align*}
Consequently, inserting the bounds for $h_\eta$ from \eqref{dsg:FGR-h-bounds} and
\eqref{dsg:FGR-h-comparison} now give
\begin{align}
\label{dsg:FGR-first-coefficient}
\Gamma_\eta &=-\frac{\eta}{2}\int_{\bbR}\calG(x)\cos\big(\sqrt3x\big)\,\ud x +O(\eta^2).
\end{align}
After inserting \eqref{dsg:FGR-G-expansion}, $P_0(x) = -8\sech^4 x$, and $\chi(x) = - \frac14 \sech^3 x$ into \eqref{dsg:FGR-H-definition}, we have
\begin{align}
\label{dsg:FGR-H-explicit}
\begin{aligned}
\calG(x)=&-\frac{32}{3}\sech (x)+\frac{526}{3}\sech^3(x)-\frac{1130}{3}\sech^5(x)+210\sech^7(x)\\
          &-32\big(\sech^3(x)-\sech^5(x)\big)\log(2\cosh x).
\end{aligned}
\end{align}
In the proof of Lemma~\ref{lemma:fourier-fgr-values}, we show that 
\begin{align*}
\int_{\mathbb R}\calG(x)\cos\big(\sqrt3x\big)\,\ud x
       &=-16\pi\sech\left(\frac{\sqrt3\pi}{2}\right),
\end{align*}
and hence
\begin{align*}
\Gamma_\eta = 8 \pi \sech\left(\frac{\sqrt3\pi}{2}\right)\eta + O(\eta^2).
\end{align*}
Together with \eqref{dsg:FGR-first-coefficient}, this proves
\eqref{dsg:FGR-main-result}. In particular, after decreasing
$\eta_0$ if necessary,
\begin{align*}
\Gamma_\eta
       &\geq4\pi\sech\left(\frac{\sqrt3\pi}{2}\right)\eta>0,
    \qquad 0<\eta\leq\eta_0.
\end{align*}

\subsection{Repulsivity of the second Darboux transformed potential}
\label{dsg:repulsivity}
\begin{figure}[t]
 \centering
 \includegraphics[width=\textwidth]{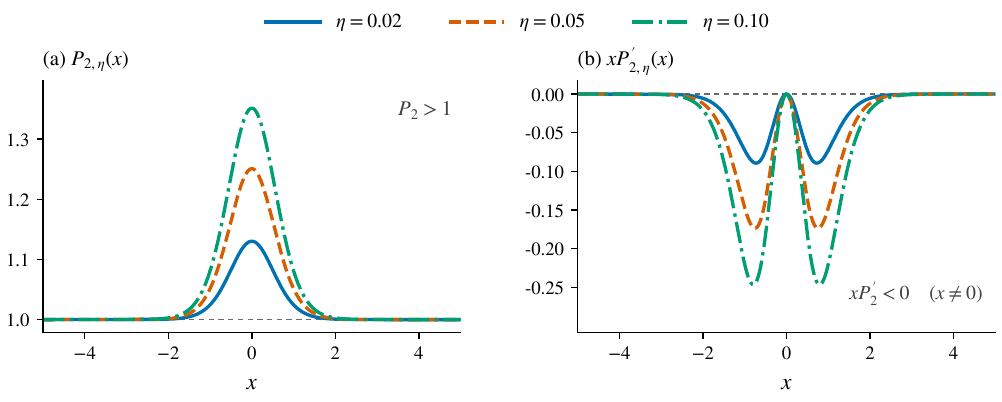}
 \caption{Double sine--Gordon: $P_{2,\eta}$ (left)
 and $xP_{2,\eta}'$ (right), computed from the positive even
 eigenfunction $Y_{1,\eta}$ of $L_{1,\eta}$, for
 $\eta=0.02,0.05,0.10$. The horizontal dashed lines are
 $P_{2,0}\equiv1$ and $xP_{2,0}'=0$. The graphs illustrate $P_{2,\eta}>1$
 and $xP_{2,\eta}'<0$ for $x\ne0$.}
 \label{fig:dsg-repulsive-P2}
\end{figure}

In this final subsection, we prove that the second Darboux transformed potential is repulsive. 
Figure~\ref{fig:dsg-repulsive-P2} illustrates the repulsivity conditions for the Darboux transformed
potential $P_{2,\eta}$. 

First, we define the logarithmic derivative of $Y_{1,\eta}$ and the second Darboux transform $U_{1,\eta}$ by
\begin{align*}
V_{1,\eta}(x):=\frac{Y_{1,\eta}'(x)}{Y_{1,\eta}(x)}, \qquad U_{1,\eta}
      :=Y_{1,\eta}\,\partial_x\,Y_{1,\eta}^{-1}
        =\partial_x-V_{1,\eta}(x).
\end{align*}
From \eqref{dsg:Z-eigenvalue-equation}, we have 
\begin{align}
\label{dsg:R-riccati}
V_{1,\eta}'+V_{1,\eta}^2
    &=\frac{Y_{1,\eta}''}{Y_{1,\eta}}
    =1+Q_{1,\eta}
      -\lambda_\eta^2.
\end{align}
Consequently,
\begin{align*}
L_{1,\eta}-\lambda_\eta^2
    &=U_{1,\eta}^*U_{1,\eta},
\end{align*}
and 
\begin{align*}
L_{2,\eta}
    &:=\lambda_\eta^2
      +U_{1,\eta}U_{1,\eta}^*
    =-\partial_x^2+1+Q_{2,\eta}(x) =: -\px^2 + P_{2,\eta}.
\end{align*}
From \eqref{dsg:R-riccati}, we have 
\begin{align}
\label{dsg:P2-potential}
\begin{aligned}
    Q_{2,\eta}(x)
    &=\lambda_\eta^2-1+2V_{1,\eta}(x)^2
      -\frac{Y_{1,\eta}''(x)}{Y_{1,\eta}(x)}=2\lambda_\eta^2-2
      -Q_{1,\eta}(x)
      +2V_{1,\eta}(x)^2.
\end{aligned}
\end{align}

By symmetry, it suffices to consider $x>0$. From \eqref{dsg:R-riccati}, we obtain 
\begin{equation*}
    V_{1,\eta}'' + 2 V_{1,\eta} V_{1,\eta}' = Q_{1,\eta}'.
\end{equation*}
By $Y_{1,\eta}' = V_{1,\eta} Y_{1,\eta}$, we have 
\begin{equation}
    (Y_{1,\eta}^2 V_{1,\eta}')' = Y_{1,\eta}^2 Q_{1,\eta}'.
\end{equation}
Since $Y_{1,\eta}^2(x)V_{1,\eta}'(x) \rightarrow 0$ as $x \rightarrow \infty$, integrating the above differential identity yields
\begin{equation}\label{dsg:integral_identity}
    V_{1,\eta}'(x) = - \int_x^\infty Q_{1,\eta}'(y) \left( \frac{Y_{1,\eta}(y)}{Y_{1,\eta}(x)}\right)^2 \,\ud y < 0.
\end{equation}
Recall that $Y_{1,\eta}$ is even and positive as it is the ground state of $L_{1,\eta}$. This implies that $V_{1,\eta}(0) = \frac{Y_{1,\eta}'(0)}{Y_{1,\eta}(0)}= 0$. Hence \eqref{dsg:integral_identity} implies that $V_{1,\eta}(x)<0$ and $Y_{1,\eta}(x)$ is decreasing for $x>0$. Moreover, 
\begin{equation*}
    0 < -V_{1,\eta}'(x) < \int_{x}^\infty Q_{1,\eta}'(y) \,\ud y = -Q_{1,\eta}(x).
\end{equation*}
From \eqref{dsg:P-epsilon-definition}, we have $-Q_{1,\eta}(y) \leq \frac{8\eta}{1+4\eta}\sech^2(y)$. Since $V_{1,\eta}(0)=0$, integrating the above inequality from $0$ yields 
\begin{equation}
    0 < -V_{1,\eta}(x) < -\int_0^x Q_{1,\eta}(y)\,\ud y \leq   \frac{8\eta}{1+4\eta}\int_0^x\sech^2(y) \,\ud y =  \frac{8\eta}{1+4\eta} \tanh (x).
\end{equation}
On the other hand, from \eqref{dsg:P-epsilon-definition} and the bound $4 \eta + \cosh^2(x) \geq 1 + 4 \eta$, we obtain
\begin{equation}
    - \frac{Q_{1,\eta}'(x)}{Q_{1,\eta}(x)} = 2 \left( 1 + \frac{\cosh^2 (x)}{4\eta + \cosh^2 (x)}\right) \tanh (x) \geq \frac{4(1+2\eta)}{1+4\eta} \tanh (x).
\end{equation}
Hence, combining \eqref{dsg:P2-potential} with the two preceding inequalities yields
\begin{align*}
    P_{2,\eta}'(x) &= -Q_{1,\eta}'(x) + 4 V_{1,\eta}(x)V_{1,\eta}'(x)\\
    &< -Q_{1,\eta}'(x) + \frac{32\eta}{1+4\eta} \big(-Q_{1,\eta}(x)\big) \tanh(x)\\
    &\leq -\frac{24\eta-4}{1+4\eta}Q_{1,\eta}(x)\tanh(x) \leq 0.
\end{align*}
Thus, by symmetry, we conclude that $x P_{2,\eta}'(x) <0$ for all $x \neq 0$ and for all $0<\eta\leq \frac16$. This finishes the proof of Proposition~\ref{prop:dsg-main}.

\section{Perturbation of the $\phi^4$ model}\label{phi4:section}
In this section, we verify the assumptions of Theorem~\ref{main_thm1} for the
perturbed $\phi^4$ model \eqref{eq:sextic}. Unlike the
sine--Gordon operator, the linearized operator for the  $\phi^4$ model already has
an isolated internal mode. We show that this mode persists, while
the threshold resonance resolves into the continuous spectrum. Consider the perturbed $\phi^4$ equation
\begin{equation}\label{phi4:phi4}
\begin{split}
\phi_{tt}&-\phi_{xx}+W_\eta'(\phi)=0,\\
 W_\eta(\phi)
 &:=\frac{(1-\phi^2)^2}{8(1-\eta)}
        \bigl(1-\eta\chi(\phi)\phi^2\bigr),
 \qquad \eta >0.    
\end{split}
\end{equation}
Here $\chi\in \calC^\infty(\bbR;[0,1])$ is even, equals $1$ on
$[-3/2,3/2]$, and vanishes outside $(-2,2)$. The main result of this section is the following proposition.

\begin{figure}[t]
 \centering
 \includegraphics[width=\textwidth]{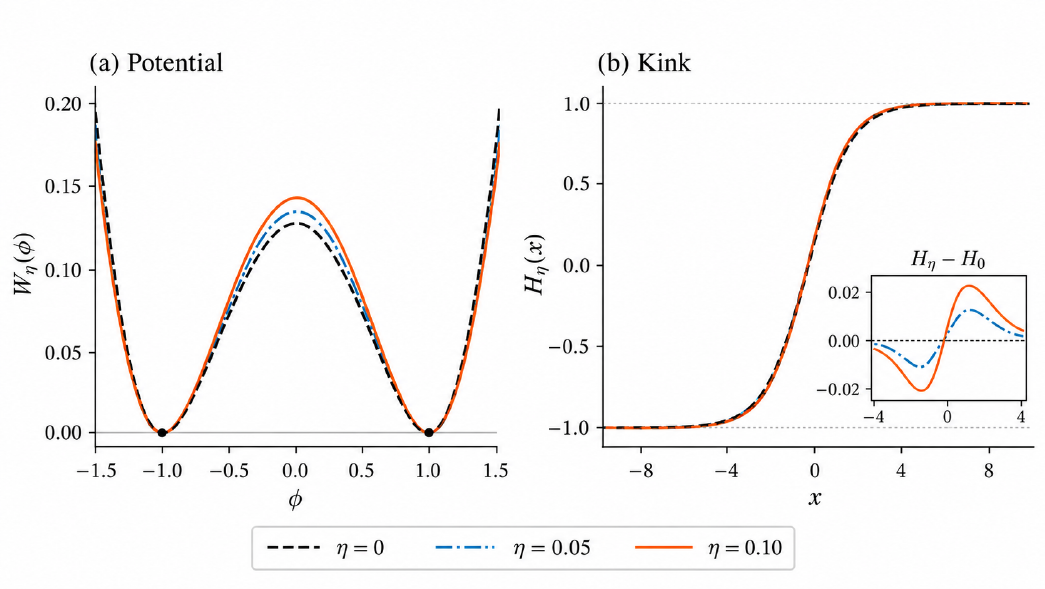}
 \caption{Perturbed $\phi^4$: the double-well potential on $[-3/2,3/2]$ (left)
 and the centered kink (right), for $\eta=0,0.05,0.10$.
 The cutoff is identically one on the displayed potential interval,
 and $W_\eta''(\pm1)=1$. The inset shows $H_\eta-H_0$.} 
 \label{fig:phi4-potential-kink}
\end{figure}

\begin{figure}[t]
 \centering
 \includegraphics[width=\textwidth]{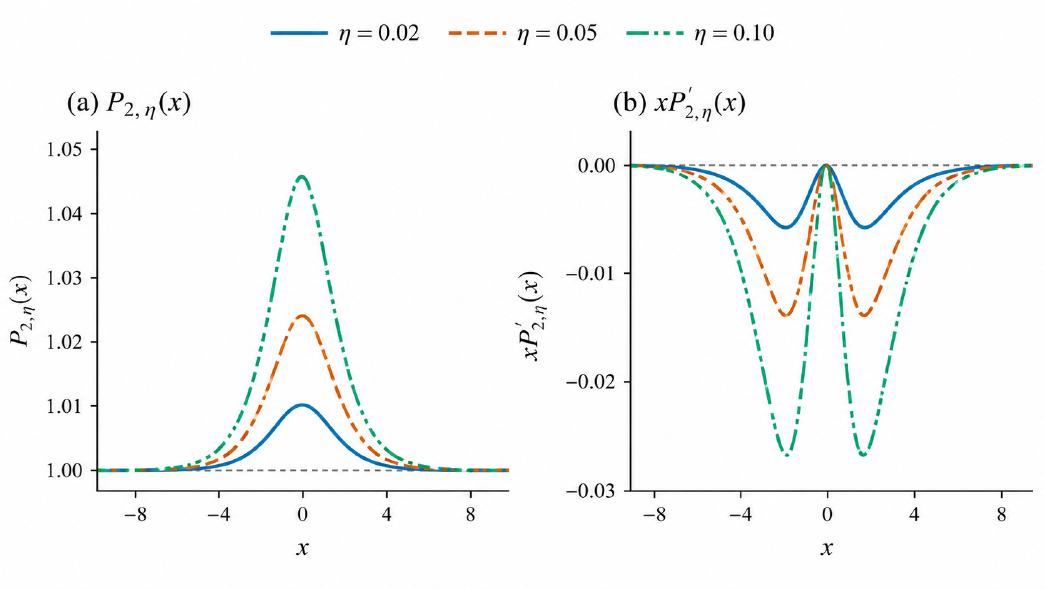}
 \caption{Perturbed $\phi^4$: $P_{2,\eta}$ (left)
 and $xP_{2,\eta}'$ (right), computed from the positive even
 eigenfunction $Y_{1,\eta}$ of $L_{1,\eta}$, for
 $\eta=0.02,0.05,0.10$. The horizontal reference lines are
 $P_{2,0}\equiv1$ and $xP_{2,0}'\equiv0$. The graphs illustrate $P_{2,\eta}>1$
 and $xP_{2,\eta}'<0$ for $x\neq0$.}
 \label{fig:phi4-repulsive-P2}
\end{figure}

\begin{proposition}\label{prop:phi4-main}
For the perturbed $\phi^4$ model above, the following statements hold
for sufficiently small $\eta>0$.
\begin{enumerate}
    \item The centered odd kink is
    \begin{equation}\label{phi4:exact-kink}
        H_\eta(x)=\frac{\sinh(x/2)}{\sqrt{\cosh^2(x/2)-\eta}}.
    \end{equation}
    \item The linearized operator about $H_\eta$ is
    \begin{equation}
    \begin{gathered}
        L_{0,\eta}=-\partial_x^2+P_{0,\eta}(x),\\
        P_{0,\eta}(x)=1+
        \frac{-6(1-3\eta)\sech^2(x/2)-3\eta(3+\eta)\sech^4(x/2)}
             {4(1-\eta\sech^2(x/2))^2}.
    \end{gathered}
    \end{equation}
    \item The operator $L_{0,\eta}$ has a unique simple internal
    eigenvalue $\lambda_\eta^2$, with an odd eigenfunction, and
    \begin{equation}\label{phi4:internal mode-goal}
        \lambda_\eta^2=\frac34+\frac35\eta+O(\eta^2).
    \end{equation}
    Its full spectrum is
    \begin{equation}\label{phi4:complete-spectrum}
        \sigma(L_{0,\eta})=\{0,\lambda_\eta^2\}\cup[1,\infty),
    \end{equation}
    and $\lambda=1$ is not a threshold resonance of $L_{0,\eta}$. 
    \item The second Darboux potential $P_{2,\eta}$ defined in
    \eqref{phi4:second-potential} is even and satisfies
    \[
        xP_{2,\eta}'(x)\leq0\quad(x\in\mathbb R),\qquad
        P_{2,\eta}'\not\equiv0,
    \]
    with strict inequality for $x\ne0$.
    \item With the normalizations \eqref{phi4:internal-normalization}
    and~\eqref{phi4:continuum-normalization},
    the Fermi golden rule coefficient satisfies
    \begin{equation}\label{phi4:FGR-main-result}
        \Gamma_\eta=\Gamma_0+O(\eta),\qquad
        \Gamma_0:=\frac{27\sqrt2\pi}{8\sinh(\sqrt2\pi)}.
    \end{equation}
\end{enumerate}
\end{proposition}
The proof of Proposition~\ref{prop:phi4-main} will be divided into the following subsections, and it will be similar to the proof of Proposition~\ref{prop:dsg-main}. Figure~\ref{fig:phi4-potential-kink} illustrates the double-well potential \eqref{phi4:phi4} and the centered kink in \eqref{phi4:exact-kink}. The $\phi^4$ case ($\eta=0$) is included for comparison. Figure~\ref{fig:phi4-repulsive-P2} illustrates the repulsivity conditions for the transformed potential $P_{2,\eta}$.

\subsection{Kink and linearized operator}\label{phi4:model-section}
Since $1-\eta\chi(\phi)\phi^2\geq1-4\eta>0$, the potential is
smooth, nonnegative, and vanishes precisely at the vacuum states $\phi=\pm1$, and $W_\eta(\phi)$ is strictly positive between these two vacua.
On $|\phi|\leq3/2$, it agrees with the expression 
\begin{align}\label{phi4:pot}
 W_\eta(\phi)=\frac{(1-\phi^2)^2(1-\eta\phi^2)}{8(1-\eta)}.
\end{align}
The cutoff $\chi(\phi)$ in \eqref{phi4:phi4} keeps the double-well potential nonnegative for all $\phi$ without changing any of its derivatives along the kink. Thus, it suffices to use the expression \eqref{phi4:pot} for the rest of the proof.

For $|\phi|\leq3/2$, the derivatives of the potential are
\begin{align*}
 W_\eta'(\phi)
 &=\frac{-(\eta+2)\phi+(2+4\eta)\phi^3-3\eta\phi^5}
              {4(1-\eta)},\\
 W_\eta''(\phi)
 &=\frac{-\eta-2+(6+12\eta)\phi^2-15\eta\phi^4}
              {4(1-\eta)},\\
 W_\eta^{(3)}(\phi)
 &=\frac{(3+6\eta)\phi-15\eta\phi^3}{1-\eta},\\
 W_\eta^{(4)}(\phi)
 &=\frac{3+6\eta-45\eta\phi^2}{1-\eta}.
\end{align*}
In particular, $W_\eta''(-1)=W_\eta''(1)=1.$
Thus $\omega_- = \omega_+ = \omega=1$, as in the double
sine--Gordon case.

The kink profile $H_\eta(x)$ is an increasing static solution of
\begin{equation}\label{phi4:kink-ode}
\left\{ \begin{aligned}
-H_\eta''+W_\eta'(H_\eta)&=0,\\
 \lim_{x\rightarrow \pm\infty}H_\eta(x)&= \pm 1.    
\end{aligned}\right.
\end{equation}
We fix the translation invariance by imposing $H_\eta(0)=0$.
The solution to \eqref{phi4:kink-ode} will be
\begin{align*}
H_\eta(x)&=\frac{\sinh(x/2)}
                 {\sqrt{\cosh^2(x/2)-\eta}}.
\end{align*}
Indeed, from  
$\frac12(H_\eta')^2=W_\eta(H_\eta)$, we have 
\begin{align*}
H_\eta'&=\frac{1}{2\sqrt{1-\eta}}\,
             (1-H_\eta^2)\sqrt{1-\eta H_\eta^2}.
\end{align*}
Substituting $v=H_\eta/\sqrt{1-\eta H_\eta^2}$ gives 
\begin{align*}
v'&=\frac{1-(1-\eta)v^2}{2\sqrt{1-\eta}},\qquad v(0)=0,
\end{align*}
and solving it by separation of variables yields
$$v(x)=(1-\eta)^{-1/2}\tanh(x/2).$$
Substituting back for $H_\eta$ using the hyperbolic identity $1-\tanh^2(x/2)=\sech^2(x/2)$ gives \eqref{phi4:exact-kink}.

Next, writing $\phi(t,x)=H_\eta(x)+u(t,x)$ gives the perturbation equation 
\begin{align*}
 u_{tt}+L_{0,\eta}u
 &=-\bigl[W_\eta'(H_\eta+u)-W_\eta'(H_\eta)-W_\eta''(H_\eta)u\bigr].
\end{align*}
If $\|u\|_{L_x^\infty(\mathbb R)}\leq1/2$, then
$|H_\eta+u|\leq3/2$ and the cutoff $\chi(H_\eta+u)$ is identically one throughout
the Taylor expansion. In this regime the perturbation equation becomes
\begin{align*}
u_{tt}+L_{0,\eta} u
 ={}&-\frac12W_\eta^{(3)}(H_\eta)u^2
     -\frac16W_\eta^{(4)}(H_\eta)u^3+\frac{15\eta H_\eta}{4(1-\eta)}u^4
     +\frac{3\eta}{4(1-\eta)}u^5,
\end{align*}
where the linearized operator is given by 
\begin{equation}\label{phi4:linearized-operator}
\begin{split}
L_{0,\eta}&:=-\partial_x^2+P_{0,\eta}(x),\\
 P_{0,\eta}(x)&:=1+\frac{-6(1-3\eta)\sech^2(x/2)-3\eta(3+\eta)\sech^4(x/2)}
              {4(1-\eta \sech^2(x/2))^2}.    
\end{split}
\end{equation}
In particular, since $P_{0,\eta}-1$ is even and exponentially decaying, we have $\sigma_{\mathrm{ac}}(L_{0,\eta})
 =[1,\infty)$.
 
Differentiating \eqref{phi4:kink-ode} produces the translational mode
\begin{align}\label{phi4:translation-mode}
L_{0,\eta}Y_{0,\eta}&=0,\qquad Y_{0,\eta}:=H_\eta',\qquad
 H_\eta'=\frac{(1-\eta)\sech^2(x/2)}{2(1-\eta \sech^2(x/2))^{3/2}}>0.
\end{align}
Consequently, zero is a simple eigenvalue and $L_{0,\eta}\geq0$.

In the $\phi^4$ case ($\eta=0$), we write 
\begin{align*}
L_{\phi^4}:=L_{0,\eta}|_{\eta=0}&=-\partial_x^2+1-\frac32\sech^2(x/2).
\end{align*}
Note that the potential in the linearized operator $L_{\phi^4}$ is the second member of the P\"oschl-Teller family of potentials \cite{PoschlTeller} and its spectrum is explicit. Besides the translation mode, $L_{\phi^4}$ has the simple odd internal mode
\begin{align}\label{phi4:internal_mode_eigenfunction}
L_{\phi^4}Y_{\phi^4}&=\frac34Y_{\phi^4},
 \qquad Y_{\phi^4}=\frac{\sqrt3}{2}\tanh(x/2) \sech(x/2),
 \qquad \|Y_{\phi^4}\|_{L_x^2}=1,
\end{align}
and an even threshold resonance
\begin{align*}
\begin{gathered}
 L_{\phi^4}R_{\phi^4}=R_{\phi^4},\qquad R_{\phi^4}:=1-\frac32\sech^2(x/2)=\frac{3\tanh^2(x/2)-1}{2}\in L^\infty(\mathbb R)\setminus L^2(\mathbb R).
\end{gathered}
\end{align*}
For sufficiently small positive $\eta$, we show that the odd internal mode persists, while the even threshold resonance bifurcates into the continuous spectrum. We denote the unique internal mode of \eqref{phi4:phi4} by $Y_\eta$ and its
eigenvalue by $\lambda_\eta^2$,
\begin{align}\label{phi4:internal mode}
L_{0,\eta} Y_\eta&=\lambda_\eta^2Y_\eta,\qquad
 0<\lambda_\eta^2<1.
\end{align}
We will prove its existence, uniqueness, and expansion in the next subsection.

\subsection{Expansion of the internal mode eigenvalue and verification of the repulsivity condition}\label{phi4:internal}

Since $H_\eta'(x)>0$, we define the logarithmic derivative of 
\eqref{phi4:translation-mode} by
\begin{align}\label{phi4:log-translation}
V_{0,\eta}:=\frac{H_\eta''}{H_\eta'}
 &=-\tanh(x/2)-\frac{3\eta\tanh(x/2)\sech^2(x/2)}
                         {2(1-\eta\sech^2(x/2))}.
\end{align}
Let $U_{0,\eta}=\partial_x-V_{0,\eta}$ be the first Darboux transform. We have
$L_{0,\eta}=U_{0,\eta}^*U_{0,\eta}$ and
\begin{align*}
L_{1,\eta}:=U_{0,\eta}U_{0,\eta}^*
 &=-\partial_x^2+1+Q_{1,\eta},
\end{align*}
where
\begin{align}\label{phi4:first-potential}
Q_{1,\eta}
 &=-\frac12\sech^2(x/2)+\frac{3\eta}{4}
       \frac{\sech^2(x/2)[2-(1+\eta)\sech^2(x/2)]}
            {(1-\eta\sech^2(x/2))^2}.
\end{align}

We now prove the eigenvalue expansion in \eqref{phi4:internal mode-goal} using the Birman-Schwinger principle, as in Section~\ref{dsg:eigenfunction}. Let 
\begin{align*}
\kappa_\eta^2&:=1-\lambda_\eta^2,\qquad \kappa_\eta>0.
\end{align*}
By the intertwining relation $U_{0,\eta}L_{0,\eta} = L_{1,\eta}U_{0,\eta}$, every positive discrete eigenvalue of $L_{0,\eta}$ corresponds to an eigenvalue of $L_{1,\eta}$.

Suppose for the moment that $Y_\eta$ is a solution found for \eqref{phi4:internal mode}. We set $Y_{1,\eta}=U_{0,\eta}Y_\eta$. By the intertwining relation, the eigenvalue equation \eqref{phi4:internal mode} is equivalent to 
\begin{align}\label{phi4:transformed-eigenvalue}
L_{1,\eta}Y_{1,\eta}&=\lambda_\eta^2Y_{1,\eta}.
\end{align}
Rewriting the above equation using the potential $Q_{1,\eta}$  gives 
\begin{align}\label{phi4:reduced-eigenvalue}
(-\partial_x^2+Q_{1,\eta})Y_{1,\eta}
       &=-\kappa_\eta^2Y_{1,\eta}.
\end{align}
For small $|\eta|\lesssim 1$, we note from \eqref{phi4:first-potential} that $Q_{1,\eta}<0$ everywhere, and so we set $S_\eta=\sqrt{-Q_{1,\eta}}$ and define the following Hilbert-Schmidt operator 
\begin{align*}
\calK_{\eta,\kappa}
 &=S_\eta(-\partial_x^2+\kappa^2)^{-1}S_\eta,\qquad \kappa>0.
\end{align*}
The Birman-Schwinger principle for \eqref{phi4:reduced-eigenvalue} is the equation 
\begin{align}\label{phi4:BS-equivalence}
\calK_{\eta,\kappa_\eta}\Phi_\eta&=\Phi_\eta,\qquad \text{where} \qquad 
\Phi_\eta=S_\eta Y_{1,\eta}.
\end{align}
At $\eta=0$, since $\kappa=\kappa_0=\frac12$, the resolvent $(-\px^2 + \kappa_0^2)^{-1}$ stays bounded on $L^2(\bbR)$. Moreover, recall that $Q_{1,0}(x) =-\frac12\sech^2(x/2)$, and
\begin{align*}
\left(-\partial_x^2+\frac14\right)\sech(x/2)&=\frac12\sech^3(x/2).
\end{align*}
Thus, at $\kappa=\kappa_0=\frac12$, we have 
\begin{align*}
 \calK_{0,\frac12}\Phi_0=\Phi_0, \qquad \text{where} \qquad  \Phi_0(x)&\equiv\frac{\sech^2(x/2)}{\sqrt2}.
\end{align*}
Since $\Phi_0(x)>0$ everywhere and $\calK_{0,\frac12}$ has a strictly positive integral kernel, it follows that $1$ is the simple largest eigenvalue of $\calK_{0,\frac12}$. Let $\mu(\eta,\kappa)$ denote the eigenvalue of the Birman-Schwinger operator $\calK_{\eta,\kappa}$ which is obtained by analytic continuation from $\mu(0,\tfrac12)=1$. 

Indeed, in a neighborhood of $(\eta,\kappa) = (0,\frac12)$, we may differentiate the integral kernel of $\calK_{\eta,\kappa}$ with respect to $\kappa$ to obtain
\begin{equation*}
    \partial_{\kappa} \calK_{\eta,\kappa}(x,y) = - S_\eta(x)  e^{-\kappa|x-y|} \left( \frac{1}{2\kappa^2} + \frac{|x-y|}{2\kappa}\right)S_\eta(y).
\end{equation*}
Using the exponential decay of $S_\eta$, it follows that
\begin{equation*}
    \left| \partial_{\kappa} \calK_{\eta,\kappa}(x,y) \right| \leq C\langle x-y\rangle e^{-C(|x|+|y|)}.
\end{equation*}
The preceding bound shows that $\partial_{\kappa} \calK_{\eta,\kappa}$ is a Hilbert-Schmidt operator, and also justifies differentiating the eigenvalue equation \eqref{phi4:BS-equivalence} near $(\eta,\kappa) = (0,\frac12)$. We compute that 
\begin{align*}
\partial_\kappa\mu(0,\tfrac12) = \frac{\langle \Phi_0,\partial_{\kappa} \calK_{0,\kappa}\big|_{\kappa=\frac12} \Phi_0\rangle}{\| \Phi_0 \|^2}=-\frac{\|\sech(x/2)\|^2}
          {\|\sech^2(x/2)/\sqrt2\|^2}=-3.
\end{align*}
The implicit function theorem therefore gives a unique analytic
branch $\kappa_\eta$ near $\kappa_0 = \frac12$ satisfying $\mu(\eta,\kappa_\eta)=1$, and the corresponding eigenfunction $\Phi_\eta$ is positive and even, with $\kappa_\eta^2 \approx \frac14$ or equivalently, $\lambda_\eta^2 \approx \frac34$. From \eqref{phi4:BS-equivalence} we obtain that $Y_{1,\eta} = (-\partial_{xx} + \kappa_\eta^2)^{-1}(S_\eta\Phi_\eta) = \frac{1}{S_\eta}\Phi_\eta$ and $Y_{1,\eta}>0$ everywhere. Hence, $Y_{1,\eta}$ is the ground state of $L_{1,\eta}$. The intertwining relation $Y_{\eta} = \frac{1}{\lambda_\eta^2} U_{0,\eta}^*Y_{1,\eta}$ then gives
a simple odd internal mode for $L_{0,\eta}$. 

In order to prove the eigenvalue expansion \eqref{phi4:internal mode-goal}, we first expand the linearized operator \eqref{phi4:linearized-operator} as 
\begin{align*}
\begin{gathered}
 L_{0,\eta}=L_{\phi^4}+\eta P+O_{W^{1,\infty}}(\eta^2),
\end{gathered}
\end{align*}
where
\begin{equation*}
    P:=\frac92\sech^2(x/2)-\frac{21}{4}\sech^4(x/2).
\end{equation*}
Differentiating the eigenvalue equation $L_{0,\eta}Y_{\eta}=\lambda_\eta^2 Y_{\eta}$ and using \eqref{phi4:internal_mode_eigenfunction} yields 
\begin{align*}
\left.\frac{\ud}{\ud\eta}\lambda_\eta^2\right|_{\eta=0}=\int_{\mathbb R}P Y_{\phi^4}^2\,\ud x =\int_\bbR \left(\frac{27}{8}\sech^4(x/2) - \frac{117}{16}\sech^6(x/2) + \frac{63}{16}\sech^8(x/2) \right)\,\ud x
 =\frac35.
\end{align*}
Hence, we have proven \eqref{phi4:internal mode-goal}.

Next, we show that there are no other eigenvalues and no threshold resonances for the linearized operator $L_{0,\eta}$, by proving that the second Darboux transformed potential $P_{2,\eta}$ is repulsive. We refer the reader to \cite[Remark~2.1(1)]{kowalczyk2021sufficient} and \cite[Lemma~5.6]{cuccagna2023kink} for why the repulsivity condition on $P_{2,\eta}$ is a sufficient condition for the absence of eigenvalues and threshold resonances for the operator $L_{2,\eta}$. Let us define the logarithmic derivative of $Y_{1,\eta}$ and the second Darboux transform operator by
\begin{align*}
V_{1,\eta}&:=\frac{Y_{1,\eta}'}{Y_{1,\eta}},\qquad
U_{1,\eta}:=\partial_x-V_{1,\eta}.
\end{align*}
We define 
\begin{align*}
L_{2,\eta}&:=\lambda_\eta^2+U_{1,\eta}U_{1,\eta}^*
           =-\partial_x^2+P_{2,\eta},
\end{align*}
with the second Darboux transformed potential given by 
\begin{align}\label{phi4:second-potential}
P_{2,\eta}
 &=1+Q_{1,\eta}-2V_{1,\eta}'
  =2\lambda_\eta^2-1-Q_{1,\eta}+2V_{1,\eta}^2.
\end{align}

Since $V_{1,\eta}^2\rightarrow1-\lambda_\eta^2$ at both spatial infinities, it follows that $P_{2,\eta}\rightarrow1$ as $x \rightarrow \pm \infty$. We claim  for all $x \in \bbR$ and for sufficiently small $\eta>0$ that the following estimate holds
\begin{equation}\label{phi4:main_estimate}
\left| \px^k \left[P_{2,\eta}(x) - 1 - \frac{3\eta}{20}\sech^2(x/2)\big(2 + \sech^2(x/2)\big) \right]\right|\leq C \eta^2 \sech^2(x/2), \quad k \in \{0,1,2\}.
\end{equation}
We first show why the estimate \eqref{phi4:main_estimate} implies repulsivity for $P_{2,\eta}$. Indeed, for small positive $\eta>0$, \eqref{phi4:main_estimate} (with $k=0$) gives 
\begin{equation}
    P_{2,\eta}(x) - 1 \geq \eta\left( \frac{3}{10}-C\eta\right) \sech^2(x/2) > 0.
\end{equation}
Moreover, since the bracketed terms in \eqref{phi4:main_estimate} are even, their first derivative is odd and vanishes at $x=0$. Hence, for $|x|\leq 1$, by Taylor expansion and \eqref{phi4:main_estimate} (with $k=2$), we obtain 
\begin{equation*}
\left| P_{2,\eta}'(x) + \frac{3\eta}{10}\tanh(x/2) \sech^2(x/2)\big(1+\sech^2(x/2)\big)\right| \leq C \eta^2 \int_0^{|x|}1 \leq C \eta^2 |x|.
\end{equation*}
For $|x| \geq 1$, the estimate \eqref{phi4:main_estimate} (with $k=1$) gives 
\begin{equation*}
\left| P_{2,\eta}'(x) + \frac{3\eta}{10}\tanh(x/2) \sech^2(x/2)\big(1+\sech^2(x/2)\big)\right| \leq C \eta^2 \sech^2(x/2) \leq C \eta^2 |\tanh(x/2)|.
\end{equation*}
Hence, for all $x>0$ and for $\eta>0$ sufficiently small, we have 
\begin{equation*}
    P_{2,\eta}'(x) \leq - \eta \left( \frac{3}{10}-C\eta\right)\tanh(x/2) \sech^2(x/2) < 0.
\end{equation*}
Thus, by symmetry, the repulsivity condition $x P_{2,\eta}'(x)\leq 0$ holds for all $x \in \bbR$, with strict inequality for $x\neq 0$.  We now prove the bound \eqref{phi4:main_estimate}.

\begin{proof}[Proof of \eqref{phi4:main_estimate}]
By symmetry, it suffices to fix $x>0$. From the explicit expression \eqref{phi4:first-potential}, we obtain the following tail bound
\begin{equation}\label{phi4:tail_bound_Q1eta}
    \left| \partial_{\eta}^j \partial_x^k Q_{1,\eta}(x) \right| \leq C \sech^2(x/2) \leq Ce^{-x}, \qquad 0 \leq j,k\leq 2.
\end{equation}
Since $Q_{1,\eta}<0$ and $Y_{1,\eta}>0$ everywhere, we obtain from \eqref{phi4:reduced-eigenvalue} that 
\begin{align*}
    Y_{1,\eta}'(x) + \kappa_\eta Y_{1,\eta}(x) &= -\int_x^\infty e^{-\kappa_\eta(y-x)} Q_{1,\eta}(y)Y_{1,\eta}(y) \,\ud y >0, \\
    Y_{1,\eta}'(x) - \kappa_\eta Y_{1,\eta}(x) &= \int_{-\infty}^x e^{\kappa_\eta(y-x)} Q_{1,\eta}(y)Y_{1,\eta}(y) \,\ud y <0. 
\end{align*}
It follows from $V_{1,\eta} = \frac{Y_{1,\eta}'}{Y_{1,\eta}}$ that 
$$- \kappa_\eta < V_{1,\eta}(x) < \kappa_\eta,$$ 
and for $y \geq x$ that
\begin{equation*}
    \frac{Y_{1,\eta}(y)}{Y_{1,\eta}(x)} = \exp\left( \int_x^y V_{1,\eta}(s) \,\ud s\right) \leq e^{\kappa_\eta(y-x)}.
\end{equation*}
Hence, using logarithmic differentiation and \eqref{phi4:tail_bound_Q1eta}, we obtain the bound 
\begin{align*}
    0 < V_{1,\eta}(x) + \kappa_\eta = - \int_x^\infty e^{-\kappa_\eta(y-x)}Q_{1,\eta}(y) \frac{Y_{1,\eta}(y)}{Y_{1,\eta}(x)} \,\ud y \leq - \int_x^\infty Q_{1,\eta}(y)\,\ud y \leq C e^{-x}.
\end{align*}
From \eqref{phi4:second-potential}, we have the integral formulation 
\begin{equation*}
    V_{1,\eta}(x) + \kappa_\eta = -\int_x^\infty e^{-2\kappa_\eta(y-x)} Q_{1,\eta}(y) \,\ud y + \int_x^\infty e^{-2\kappa_\eta(y-x)}\big( V_{1,\eta}(y) + \kappa_\eta\big)^2 \,\ud y.
\end{equation*}
Since $\kappa_\eta$ is analytic in $\eta$ for small $\eta$, we have $|\partial_\eta (\kappa_\eta)| + |\partial_\eta^2(\kappa_\eta)| \leq C$. Thus, we obtain from \eqref{phi4:tail_bound_Q1eta} that 
\begin{equation}\label{phi4:tail_bound_V1eta}
    \left| \partial_\eta^j \partial_x^k \big( V_{1,\eta}(x) +\kappa_\eta\big) \right| \leq C e^{-x}, \qquad 0 \leq j,k\leq 2.
\end{equation}
Rewriting \eqref{phi4:second-potential} as 
\begin{equation*}
    P_{2,\eta}-1 = -Q_{1,\eta} - 4\kappa_{\eta}(V_{1,\eta} + \kappa_\eta) + 2 (V_{1,\eta} + \kappa_\eta)^2,
\end{equation*}
we conclude from \eqref{phi4:tail_bound_Q1eta} and \eqref{phi4:tail_bound_V1eta} for all $x \in \bbR$ that 
\begin{equation*}
    \left| \partial_\eta^2 \partial_x^k \big( P_{2,\eta}(x) -1\big) \right| \leq C e^{-|x|}\leq C \sech^2(x/2), \qquad 0 \leq k\leq 2.
\end{equation*}
Note that $P_{2,0} \equiv 1$ and the first variation of $P_{2,\eta}$ is given by 
\begin{equation*}
\partial_{\eta} P_{2,\eta}|_{\eta=0} = \partial_{\eta} Q_{1,\eta}|_{\eta=0}-2\px(\partial_{\eta} V_{1,\eta}|_{\eta=0}).
\end{equation*}
From the explicit expression \eqref{phi4:first-potential}, we have 
\begin{equation*}
\partial_{\eta} Q_{1,\eta}|_{\eta=0} = \frac34 \sech^2(x/2) \big(2-\sech^2(x/2)\big).
\end{equation*}
Since $V_{1,\eta}$ satisfies the Riccati identity $V_{1,\eta}'+V_{1,\eta}^2 = \kappa_\eta^2 +Q_{1,\eta}$ and $V_{1,0}=-\tfrac12 \tanh(x/2)$, its first variation  $\sfv=(\partial_{\eta} V_{1,\eta})|_{\eta=0}$ solves the ODE
\begin{equation*}
    \sfv'(x) - \tanh(x/2) \sfv(x) = \frac32 \sech^2(x/2) - \frac34\sech^4(x/2) - \frac35, \qquad \sfv(0) =0.
\end{equation*}
The solution is given by 
\begin{equation}\label{phi4:V1eta_first_variation}
    \sfv(x) = - \frac{3}{10}\tanh(x/2)\sech^2(x/2) + \frac{3}{5}\tanh(x/2).
\end{equation}
Therefore, we have
\begin{equation*}
\partial_{\eta} P_{2,\eta}|_{\eta=0} = \partial_{\eta} Q_{1,\eta}|_{\eta=0} - 2\sfv'(x) =\frac{3}{20}\sech^2(x/2)\big(2+\sech^2(x/2)\big).
\end{equation*}
The preceding bounds for $P_{2,\eta}$ complete the proof of \eqref{phi4:main_estimate}.
\end{proof}

\subsection{Verifying the Fermi golden rule}\label{phi4:FGR}
From \eqref{phi4:internal mode-goal}, the continuum energy $4\lambda_\eta^2$ is expanded as
\begin{align}\label{phi4:FGR-energy}
4\lambda_\eta^2-1&=2+\frac{12}{5}\eta+O(\eta^2).
\end{align}
Hence
\begin{align*}
k_\eta=\sqrt{4\lambda_\eta^2-1}
 &=\sqrt2+\frac{3\sqrt2}{5}\eta+O(\eta^2).
\end{align*}
Let $g_\eta$ be the generalized eigenfunction satisfying
\begin{align}\label{phi4:continuum-normalization}
L_{0,\eta}g_\eta&=4\lambda_\eta^2g_\eta,\qquad
 g_\eta(0)=0,\qquad g_\eta'(0)=4\lambda_\eta^2,
\end{align}
and define
\begin{align}\label{phi4:FGR-amplitude}
\Gamma_\eta&:=-\frac12\int_{\mathbb R}
 W_\eta^{(3)}(H_\eta)Y_\eta^2g_\eta\,\ud x.
\end{align}
In contrast to the double sine--Gordon case, the unperturbed
coefficient is nonzero. It therefore suffices to compute
$\Gamma_0$ and prove $\Gamma_\eta=\Gamma_0+O(\eta)$.
See also \cite[Lemma G.2.1]{delort2022kink} for the Fermi golden rule constant for the $\phi^4$ model. In the following subsections, we compute \eqref{phi4:FGR-amplitude} by expanding each function in $\eta$.

\subsubsection{Expansion of the eigenfunction}\label{phi4:eigenfunction}
First, we normalize the eigenfunction $Y_\eta$ by
\begin{align}\label{phi4:internal-normalization}
\|Y_\eta\|_{L_x^2}&=1,\qquad Y_\eta(x)>0\quad(x>0).
\end{align}
We claim that 
\begin{align}\label{phi4:internal-eigenfunction}
\begin{aligned}
Y_\eta(x)&=\frac{\sqrt3}{2}\tanh(x/2)\sech(x/2)
                  [1+\eta B(x)]+R_\eta^Y(x),\\
B(x)&=-\frac32\tanh^2(x/2)+\frac65\log\cosh(x/2)
                  -\frac7{10}+\frac65\log2,
\end{aligned}
\end{align}
where the remainder satisfies 
\begin{align}\label{phi4:internal-remainder}
|\partial_x^mR_\eta^Y(x)|&\leq C_m\eta^2e^{-C|x|}.
\end{align}

Let us denote by $\dot Y_{\phi^4}:=\left.\partial_\eta Y_\eta\right|_{\eta=0}$ for the first variation of $Y_\eta$. Differentiating the eigenvalue equation $L_{0,\eta}Y_{\eta}=\lambda^2 Y_{\eta}$ with respect to $\eta$ gives 
\begin{align}\label{phi4:eigenfunction-variation}
\left(L_{\phi^4}-\frac34\right)\dot Y_{\phi^4}
       &=\left(\frac35-P\right)Y_{\phi^4}.
\end{align}
Recall that $Y_{\phi^4}(x)= \frac{\sqrt3}{2}\tanh (x/2)\sech (x/2)$. Since 
$\dot Y_{\phi^4}(x)=Y_{\phi^4}(x)B(x)$ and 
\begin{align*} 
\left(L_{\phi^4}-\frac34\right) \big( Y_{\phi^4}B\big)
       &=-Y_{\phi^4}B'' -2Y_{\phi^4}'B',
\end{align*}
the equation \eqref{phi4:eigenfunction-variation} reduces to
\begin{equation*}
    (Y_{\phi^4}^2 B')' = \left( P - \frac35\right) Y_{\phi^4}^2.
\end{equation*}
Integrating from zero gives 
\begin{equation*}
    B'(x) = \frac{1}{Y_{\phi^4}(x)}\int_0^x \left( P(y)- \frac35 \right) Y_{\phi^4}(y)^2 \,\ud y = \tanh(x/2) \left( \frac35 - \frac32 \sech^2(x/2)\right),
\end{equation*}
and integrating once more yields
\begin{equation*}
    B(x) = -\frac32 \tanh^2(x/2) + \frac65 \log\big(\cosh(x/2)\big)  + \mathrm{C},
\end{equation*}
for some constant $\mathrm{C}$. Moreover, differentiating $\|Y_{\eta}\|^2=1$ with respect to $\eta$ gives
\begin{equation}\label{phi4:orthogonal}
    \int_{\mathbb R}Y_{\phi^4}\dot Y_{\phi^4}\,\ud x =\int_{\mathbb R}Y_{\phi^4}^2B\,\ud x  =0.
\end{equation}
Using the fact that 
\begin{align*}
    \int_\bbR Y_{\phi^4}^2(x) \tanh^2(x/2)\,\ud x &= \frac35, \\
    \int_\bbR Y_{\phi^4}^2(x) \log\big(\cosh(x/2)\big)\,\ud x &= \frac43 - \log 2,
\end{align*}
we determine from \eqref{phi4:orthogonal} that $\mathrm{C}= - \frac{7}{10} + \frac{6}{5}\log 2$. It remains to bound the Taylor remainder $R_\eta^Y$. From $V_{1,\eta} = \log(Y_{1,\eta})'$, we obtain for $x>0$ that 
\begin{equation}\label{phi4:ratio}
    \frac{Y_{1,\eta}(x)}{Y_{1,\eta}(0)}= \exp \left[ - \kappa_\eta x + \int_0^x \big( V_{1,\eta}(y) + \kappa_\eta\big)\,\ud y\right].
\end{equation}
Using \eqref{phi4:tail_bound_Q1eta}, \eqref{phi4:tail_bound_V1eta}, and the identities
\begin{equation*}
    Y_{1,\eta}' = V_{1,\eta} Y_{1,\eta},\qquad Y_{1,\eta}'' = (\kappa_\eta^2 + Q_{1,\eta})Y_{1,\eta},
\end{equation*}
we obtain
\begin{align*}
    \left|\partial_\eta^j \px^m \left( \frac{Y_{1,\eta}(x)}{Y_{1,\eta}(0)}
\right)\right| \leq C\jx^j e^{-\kappa_\eta x}, \qquad 0 \leq j,m\leq 2.
\end{align*}
By symmetry, the exponential tail bounds extend to $x<0$. We reconstruct $Y_{\eta}$ using $U_{0,\eta}^* = -\px - V_{0,\eta}$ by 
\begin{equation*}
    Y_{\eta}(x) = \frac{U_{0,\eta}^*(Y_{1,\eta}(x)/Y_{1,\eta}(0))}{\|U_{0,\eta}^*(Y_{1,\eta}(\cdot)/Y_{1,\eta}(0))\|}.
\end{equation*}
The coefficients $V_{0,\eta}$ have uniformly bounded mixed derivatives, so we obtain that 
\begin{align*}
|\partial_\eta^j\partial_x^mY_\eta(x)|
 &\leq C_m\langle x\rangle^j e^{-C|x|},
 \qquad 0 \leq j\leq 2.
\end{align*}Taylor's formula therefore gives
\eqref{phi4:internal-remainder}. Applying $U_{0,\eta}$ to \eqref{phi4:internal-eigenfunction} gives
\begin{align*}
Y_{1,\eta}(x)=\frac{\sqrt3}{4}\sech(x/2)
       \left[1+\eta\left(-\frac3{10}\tanh^2(x/2)
          +\frac65\log\cosh(x/2)-\frac7{10}+\frac65\log2\right)\right]+O(\eta^2e^{-C|x|}),
\end{align*}
with similar mixed derivative bounds.

\subsubsection{Constructing the generalized eigenfunction}\label{phi4:continuum}

Let $h_\eta$ solve
\begin{align}\label{phi4:FGR-partner-IVP}
-h_\eta''+\eta P_{0,\eta}h_\eta&=k_\eta^2h_\eta,\qquad
 h_\eta(0)=0,\qquad h_\eta'(0)=1.
\end{align}
The function $h_\eta$ is odd. The two Darboux transformations give
\begin{align}\label{phi4:FGR-g-definition}
g_\eta&=-U_{0,\eta}^*U_{1,\eta}^*h_\eta.
\end{align}
To check the initial data, write
\begin{align*}
U_{0,\eta}^*U_{1,\eta}^*h_\eta
 &=h_\eta''+(V_{0,\eta}+V_{1,\eta})h_\eta'
       +(V_{1,\eta}'+V_{0,\eta} V_{1,\eta})h_\eta.
\end{align*}
Since $V_{0,\eta}(0)=V_{1,\eta}(0)=h_\eta''(0)=0$ and
$P_{2,\eta}(0)=-V_{0,\eta}'(0)-2V_{1,\eta}'(0)$, we obtain
$g_\eta(0)=0$ and
\begin{align*}
g_\eta'(0)
 &=k_\eta^2-P_{2,\eta}(0)+1-V_{0,\eta}'(0)-2V_{1,\eta}'(0)
  =4\lambda_\eta^2.
\end{align*}

For $x\geq0$, variation of constants gives
\begin{align}\label{phi4:FGR-Volterra}
h_\eta(x)&=\frac{\sin(k_\eta x)}{k_\eta}
 +\frac{\eta}{k_\eta}\int_0^x
     \sin(k_\eta(x-y))P_{0,\eta}(y)h_\eta(y)\,\ud y.
\end{align}
Gronwall's inequality and the
uniform $L^1$ bound for $P_{0,\eta}$ imply that 
\begin{align}\label{phi4:FGR-h-bounds}
\begin{aligned}
\|h_\eta\|_{L_x^\infty}+\|h_\eta'\|_{L_x^\infty}&\leq C,\\
\left|h_\eta(x)-\frac{\sin(k_\eta x)}{k_\eta}\right|
 +|h_\eta'(x)-\cos(k_\eta x)|&\leq C\eta.
\end{aligned}
\end{align}
Consequently, with $h_0(x)=\sin(\sqrt2x)/\sqrt2$, we obtain 
\begin{align}\label{phi4:FGR-h-comparison}
|h_\eta(x)-h_0(x)|+|h_\eta'(x)-h_0'(x)|
 &\leq C\eta\langle x\rangle.
\end{align}
These estimates for $h_\eta(x)$ extend to $x<0$ by symmetry. The equation for $h_\eta$
and \eqref{phi4:FGR-g-definition} also give uniform bounds for
$g_\eta$ and $g_\eta'$.

\subsubsection{Expansion of the localized forcing}

Set
\begin{align*}
F_\eta&:=W_\eta^{(3)}(H_\eta)Y_\eta^2,\qquad
G_\eta:=U_{1,\eta}U_{0,\eta}F_\eta.
\end{align*}
Taking adjoints yields
\begin{align}\label{phi4:FGR-transformed-amplitude}
\Gamma_\eta&=\frac12\int_{\mathbb R}G_\eta h_\eta\,\ud x.
\end{align}

The kink formula from \eqref{phi4:exact-kink} and the expansion for $Y_{\eta}$ from \eqref{phi4:internal-eigenfunction} give
\begin{align*}
F_0(x)=\frac94\tanh^3(x/2)\sech^2(x/2),\qquad
\sum_{j=0}^2|\partial_x^j(F_\eta-F_0)(x)|
 \leq C\eta e^{-C|x|}.
\end{align*}
Moreover, from \eqref{phi4:log-translation} and \eqref{phi4:V1eta_first_variation}, for small $\eta$, we have 
\begin{equation*}
    \| V_{0,\eta}-V_{0,0}\|_{W_x^{1,\infty}} + \| V_{1,\eta}-V_{1,0}\|_{L_x^\infty} \leq C |\eta|.
\end{equation*}
Thus, we have $G_\eta=G_0+R_\eta^G$ where the remainder satisfies $|R_\eta^G(x)|\leq C\eta e^{-C|x|}$, and where 
\begin{align}\label{phi4:FGR-G-expansion}
\begin{aligned}
G_0&=\left(\partial_x+\frac12\tanh(x/2)\right)
       \left(\partial_x+\tanh(x/2)\right)F_0\\
   &=\frac{27}{16}\tanh(x/2)\sech^4(x/2)
                    \bigl(5\sech^2(x/2)-3\bigr).
\end{aligned}
\end{align}

\subsubsection{Computation of the FGR constant}

Let $\Gamma_0=\frac12\int_{\mathbb R}G_0h_0\,\ud x$.
The preceding estimates imply
\begin{align*}
|\Gamma_\eta-\Gamma_0|
 &\leq\frac12\|G_\eta-G_0\|_{L_x^1}\|h_\eta\|_{L_x^\infty}
       +\frac12\int_{\mathbb R}|G_0|\,|h_\eta-h_0|\,\ud x\\
 &\leq C\eta.
\end{align*}
In Lemma~\ref{lemma:fourier-fgr-values} from
Appendix~\ref{app:fourier-integrals}, we compute that 
\begin{align*}
\Gamma_0
 &=\frac12\int_{\mathbb R}G_0(x)
                    \frac{\sin(\sqrt2x)}{\sqrt2}\,\ud x
  =\frac{27\sqrt2\pi}{8\sinh(\sqrt2\pi)}.
\end{align*}
This proves \eqref{phi4:FGR-main-result}. After decreasing $\eta_0$,
\begin{align*}
\Gamma_\eta&\geq\frac{27\sqrt2\pi}{16\sinh(\sqrt2\pi)}>0.
\end{align*}
Thus, we conclude the proof of Proposition~\ref{prop:phi4-main}.

\appendix
\section{Collection of technical lemmas}
Here we collect and recall some technical lemmas from \cite{kowalczyk2021sufficient}.

\begin{lemma}[Lemma~4.5 from \cite{kowalczyk2021sufficient}]\label{lemma45}
For any $T > 0$, let $w \in \mathcal{C}([0, T], H^1) \cap \mathcal{C}^1([0, T], L^2)$ and $K > 0$ be such that
\[
(t, x) \mapsto w(t, x) \cosh \left(\frac{\omega}{K} x\right) \in L^\infty([0, T] \times \mathbb{R}),
\]
and $\sigma \in \mathcal{C}^1([0, T], \mathbb{R})$ such that
\[
\sup_{t \in [0, T]} |\dot{\sigma}(t)| \le \frac{\omega}{8K}.
\]
Then, the following inequality holds:
\begin{align*}
&\int_0^T \int |w|^3 \cosh \left(\frac{\omega}{K} x\right) \ud x \, \ud t \\
&\quad \le \frac{24K}{\omega} \left\| w \cosh \left(\frac{\omega}{K} x\right) \right\|_{L^\infty([0, T] \times \mathbb{R})}\times \left[ \frac{K}{\omega} \int_0^T \int (\partial_x w + \sigma \dot{w})^2 \ud x \, \ud t + \sup_{t \in [0, T]} |\sigma(t)| \|w(t)\|_{L^2}^2 \right].    
\end{align*}
\end{lemma}

\begin{lemma}[Lemma~4.6 from \cite{kowalczyk2021sufficient}]\label{KMMvdB_lem4.6}
Let $J$ be any non-empty open interval. There exists $C = C(J) > 0$ such that, for any $K \ge 1$, $T > 0$, $w \in \mathcal{C}([0, T], H^1) \cap \mathcal{C}^1([0, T], L^2)$, and $\sigma \in \mathcal{C}^1([0, T], \mathbb{R})$ such that
\[
\sup_{t \in [0, T]} |\dot{\sigma}(t)| \le \frac{\omega}{4K},
\]
the following inequality holds
\begin{align*}
&\int_0^T \int w^2(t, x) \mathrm{sech} \left(\frac{\omega}{K} x\right) \ud x \, \ud t \\
&\quad \le C K^2 \int_0^T \int (\partial_x w + \sigma \dot{w})^2(t, x) \ud x \, dt  \\
&\qquad + C K \int_0^T \int_J w^2(t, x) \ud x \, \ud t + C K \sup_{t \in [0, T]} |\sigma(t)| \|w(t)\|_{L^2}^2.    
\end{align*}
\end{lemma}

\section{Evaluation of the Fourier integrals}\label{app:fourier-integrals}

We collect some standard Fourier identities used in
Sections~\ref{dsg:FGR} and~\ref{phi4:FGR}.
Related recurrence formulas have been used in the
Fermi golden rule computations in Section 6 of  Martel
\cite{martel2024asymptotic} and Section 3 of Cuccagna--Maeda
\cite{cuccagna2025fgr}.  We include the derivations for completeness and apply them to compute the Fermi golden rule constants for our 
two perturbation models. For integers $p\geq 1$ and $k\in\mathbb R$, we set
\begin{align*}
I_p(k):=\int_{\mathbb R}\sech^p x\cos(kx)\,\ud x,\qquad
M_p(k):=\int_{\mathbb R}\sech^p x\log(2\cosh x)\cos(kx)\,\ud x.
\end{align*}
In particular, we recall from \cite[Example~3.3]{SteinShakarchi_Complex} and \cite[Corollary~5.7]{luhrmann2023sinegordon} that
\begin{align}\label{app:fourier-base}
I_1(k)&=\pi\sech\left(\frac{\pi k}{2}\right),\qquad
I_2(k)=\frac{\pi k}{\sinh(\pi k/2)}.
\end{align}
In the following lemma, we prove recurrence formulas for $I_p$ and $M_p$.

\begin{lemma}\label{lemma:fourier-identities}
For all integers $p\geq1$ and real $k\in\mathbb R$, we have 
\begin{align}\label{app:fourier-recurrence}
I_{p+2}(k)&=\frac{p^2+k^2}{p(p+1)}I_p(k),\\
k I_p(k) &= p\int_{\bbR}\tanh (x)\sech^p (x)\sin(kx)\,\ud x
 ,\label{app:fourier-sine}\\
M_{p+2}(k)&=\frac{p^2+k^2}{p(p+1)}M_p(k)+\frac{(2p+1)k^2-p^2}{p(p+1)}I_p(k).
 \label{app:fourier-log-recurrence}
\end{align}
\end{lemma}

\begin{proof} The identities \eqref{app:fourier-recurrence} and \eqref{app:fourier-sine} have been derived in \cite[Lemma~A.1]{LL2} using integration by parts. Next, we note that 
\begin{equation*}
\px\big(\log(2\cosh(x))\big) = \tanh(x), \qquad \px^2\big(\log(2\cosh(x))\big) = \sech^2(x).
\end{equation*}
By the product rule and the hyperbolic identity $1-\tanh^2(x)=\sech^2(x)$, we obtain that 
\begin{align*}
    \px^2 \big( \sech^p(x) \log(2\cosh(x))\big) &= p^2 \sech^p(x) \log(2 \cosh(x)) - p(p+1)\sech^{p+2}(x)\log(2\cosh(x)) \\
    &\quad - 2p \sech^p(x) + (2p+1)\sech^{p+2}(x).
\end{align*}
The terms in the above identity have exponential decay and vanish at spatial infinities. Multiplying the above identity by $\cos(kx)$ and integrating by parts therefore yields 
\begin{equation*}
    -k^2M_p = p^2 M_p - p(p+1)M_{p+2} - 2p I_p + (2p+1) I_{p+2},
\end{equation*}
which is equivalent to \eqref{app:fourier-log-recurrence}.
\end{proof}
In the next lemma, we evaluate the leading integral terms in the Fermi golden rule constants. 
\begin{lemma}\label{lemma:fourier-fgr-values}
Let $\calG$ and $G_0$ be as in \eqref{dsg:FGR-H-explicit}
and~\eqref{phi4:FGR-G-expansion}, respectively. Then
\begin{align}\label{app:fourier-dsg-value}
\int_{\mathbb R}\calG(x)\cos(\sqrt3x)\,\ud x
 &=-16\pi\sech\left(\frac{\sqrt3\pi}{2}\right),\\
\frac12\int_{\mathbb R}G_0(x)\frac{\sin(\sqrt2x)}{\sqrt2}\,\ud x
 &=\frac{27\sqrt2\pi}{8\sinh(\sqrt2\pi)}.
 \label{app:fourier-phi4-value}
\end{align}
\end{lemma}

\begin{proof}
Recall that 
\begin{align*}
\calG(x)=&-\frac{32}{3}\sech (x)+\frac{526}{3}\sech^3(x)-\frac{1130}{3}\sech^5(x)+210\sech^7(x)\\
  &-32\big(\sech^3(x)-\sech^5(x)\big)\log(2\cosh x).
\end{align*}
Iterating \eqref{app:fourier-recurrence} gives
\begin{align*}
I_3\big(\sqrt3\big)&=I_5\big(\sqrt3\big)=2I_1\big(\sqrt3\big),\qquad
I_7\big(\sqrt3\big)=\frac{28}{15}I_1\big(\sqrt3\big).
\end{align*}
By \eqref{app:fourier-log-recurrence} and \eqref{app:fourier-recurrence},
\begin{align*}
M_3\big(\sqrt3\big)-M_5\big(\sqrt3\big)
 &=-\frac1{12}I_3\big(\sqrt3\big)=-\frac16I_1\big(\sqrt3\big).
\end{align*}
Substituting these identities and then using \eqref{app:fourier-base} yields
\begin{align*}
&\int_{\mathbb R}\calG(x)\cos(\sqrt3x)\,\ud x\\
 &=-\frac{32}{3}I_1\big(\sqrt3\big)+\frac{526}{3}I_3\big(\sqrt3\big)
       -\frac{1130}{3}I_5\big(\sqrt3\big)+210I_7\big(\sqrt3\big)-32\Big(M_3\big(\sqrt3\big)-M_5\big(\sqrt3\big) \Big)\\
 &=-\frac{48}{3}I_1\big(\sqrt3\big)
  =-16\pi\sech\left(\frac{\sqrt3\pi}{2}\right).
\end{align*}
Next, we recall that $G_0(x) = \frac{27}{16}\tanh(x/2)\sech^4(x/2)                   \bigl(5\sech^2(x/2)-3\bigr)$. 
From \eqref{app:fourier-recurrence}, we have 
\begin{align*}
I_4\big(2\sqrt2\big)&=2I_2\big(2\sqrt2\big),\qquad I_6\big(2\sqrt2\big)=\frac{12}{5}I_2\big(2\sqrt2\big).
\end{align*}
Using \eqref{app:fourier-sine}, the above identities, and then \eqref{app:fourier-base}, we obtain
\begin{align*}
\frac12\int_{\mathbb R}G_0(x)\frac{\sin(\sqrt2x)}{\sqrt2}\,\ud x
 =\frac{27}{8}\left(\frac56I_6\big(2\sqrt2\big)-\frac34I_4\big(2\sqrt2\big)\right)=\frac{27}{16}I_2\big(2\sqrt2\big)
  =\frac{27\sqrt2\pi}{8\sinh(\sqrt2\pi)},
\end{align*}as desired.
\end{proof}

\bibliographystyle{amsplain}
 \bibliography{ref}

@article{alejo2023sinegordon,
  author = {Alejo, M. A. and Mu{\~n}oz, C. and Palacios, J. M.},
  title = {On asymptotic stability of the {sine--Gordon} kink in the energy space},
  journal = {Comm. Math. Phys.},
  volume = {402},
  pages = {581--636},
  year = {2023}
}

@article{campbell1986kink,
  author = {Campbell, D. K. and Peyrard, M. and Sodano, P.},
  title = {Kink--antikink interactions in the double {sine--Gordon} equation},
  journal = {Phys. D},
  volume = {19},
  number = {2},
  pages = {165--205},
  year = {1986}
}

@misc{chen2024sinegordon,
  author = {Chen, G. and L{\"u}hrmann, J.},
  title = {Asymptotic stability of the {sine--Gordon} kink},
  year = {2024},
  eprint = {2411.07004},
  archiveprefix = {arXiv},
  url = {https://arxiv.org/abs/2411.07004},
  note = {Preprint}
}

@article{coles2016edge,
  author = {Coles, M. and Gustafson, S.},
  title = {A degenerate edge bifurcation in the 1{D} linearized nonlinear {Schr{\"o}dinger} equation},
  journal = {Discrete Contin. Dyn. Syst.},
  volume = {36},
  pages = {2991--3009},
  year = {2016}
}

@article{cuccagna2023kink,
  author = {Cuccagna, S. and Maeda, M.},
  title = {Asymptotic stability of kink with internal modes under odd perturbation},
  journal = {NoDEA Nonlinear Differential Equations Appl.},
  volume = {30},
  year = {2023},
  note = {Paper No. 1, 47 pp.}
}

@article{cuccagna2023stabilization,
  author = {Cuccagna, S. and Maeda, M. and Scrobogna, S.},
  title = {Small energy stabilization for 1{D} nonlinear {Klein--Gordon} equations},
  journal = {J. Differential Equations},
  volume = {350},
  pages = {52--88},
  year = {2023},
}

@article{cuccagna2025cubic,
  author = {Cuccagna, S. and Maeda, M.},
  title = {The asymptotic stability on the line of ground states of the pure power {NLS} with {$0<|p-3|\ll1$}},
  journal = {J. Funct. Anal.},
  volume = {288},
  year = {2025},
  note = {Paper No. 110861, 57 pp.}
}

@article{cuccagna2025quadratic,
  author = {Cuccagna, S. and Maeda, M.},
  title = {On the asymptotic stability on the line of ground states of the pure power {NLS} with {$0\leq2-p\ll1$}},
  journal = {J. Differential Equations},
  volume = {440},
  year = {2025},
  note = {Paper No. 113451, 46 pp.}
}

@book{delort2022kink,
  author = {Delort, J.-M. and Masmoudi, N.},
  title = {Long-time dispersive estimates for perturbations of a kink solution of one-dimensional cubic wave equations},
  series = {Mem. Eur. Math. Soc.},
  volume = {1},
  year = {2022},
  publisher = {EMS Press, Berlin}
}

@article{denzler1993nonpersistence,
  author = {Denzler, J.},
  title = {Nonpersistence of breather families for the perturbed {sine--Gordon} equation},
  journal = {Comm. Math. Phys.},
  volume = {158},
  pages = {397--430},
  year = {1993}
}

@article{germain2022quadratic,
  author = {Germain, P. and Pusateri, F.},
  title = {Quadratic {Klein--Gordon} equations with a potential in one dimension},
  journal = {Forum Math. Pi},
  volume = {10},
  year = {2022},
  note = {Paper No. e17, 172 pp.}
}

@article{germain2024review,
  author = {Germain, P.},
  title = {A review on asymptotic stability of solitary waves in nonlinear dispersive problems in dimension one},
  journal = {Journ{\'e}es {\'E}quations aux D{\'e}riv{\'e}es Partielles},
  pages = {1--22},
  year = {2024},
  note = {Exp. No. VI}
}

@article{henry1982stability,
  author = {Henry, D. B. and Perez, J. F. and Wreszinski, W. F.},
  title = {Stability theory for solitary-wave solutions of scalar field equations},
  journal = {Comm. Math. Phys.},
  volume = {85},
  pages = {351--361},
  year = {1982}
}

@article{kopylova2011kink,
  author = {Kopylova, E. and Komech, A. I.},
  title = {On asymptotic stability of kink for relativistic {Ginzburg--Landau} equations},
  journal = {Arch. Ration. Mech. Anal.},
  volume = {202},
  pages = {213--245},
  year = {2011}
}

@article{kopylova2011moving,
  author = {Kopylova, E. and Komech, A. I.},
  title = {On asymptotic stability of moving kink for relativistic {Ginzburg--Landau} equation},
  journal = {Comm. Math. Phys.},
  volume = {302},
  pages = {225--252},
  year = {2011}
}

@article{kowalczyk2017breathers,
  author = {Kowalczyk, M. and Martel, Y. and Mu{\~n}oz, C.},
  title = {Nonexistence of small, odd breathers for a class of nonlinear wave equations},
  journal = {Lett. Math. Phys.},
  volume = {107},
  pages = {921--931},
  year = {2017}
}

@article{kowalczyk2017kink,
  author = {Kowalczyk, M. and Martel, Y. and Mu{\~n}oz, C.},
  title = {Kink dynamics in the {$\phi^4$} model: asymptotic stability for odd perturbations in the energy space},
  journal = {J. Amer. Math. Soc.},
  volume = {30},
  number = {3},
  pages = {769--798},
  year = {2017},
}

@article{kowalczyk2021sufficient,
  author = {Kowalczyk, M. and Martel, Y. and Mu{\~n}oz, C. and {Van Den Bosch}, H.},
  title = {A sufficient condition for asymptotic stability of kinks in general $(1+1)$-scalar field models},
  journal = {Ann. PDE},
  volume = {7},
  number = {1},
  year = {2021},
  note = {Paper No. 10, 98 pp.}
}

@article{kowalczyk2022jems,
  author = {Kowalczyk, M. and Martel, Y. and Mu{\~n}oz, C.},
  title = {Soliton dynamics for the 1{D} {NLKG} equation with symmetry and in the absence of internal modes},
  journal = {J. Eur. Math. Soc. (JEMS)},
  volume = {24},
  number = {6},
  pages = {2133--2167},
  year = {2022}
}

@article{kowalczyk2022kink,
  author = {Kowalczyk, M. and Martel, Y.},
  title = {Kink dynamics under odd perturbations for $(1+1)$-scalar field models with one internal mode},
  journal = {Math. Res. Lett.},
  volume = {31},
  number = {3},
  pages = {795--832},
  year = {2024},
}

@article{kowalczyk2026breathers,
  author = {Kowalczyk, M. and Martel, Y.},
  title = {Asymptotic analysis of small energy breathers for the nonlinear {Klein--Gordon} equation},
  journal = {Ann. PDE},
  volume = {12},
  year = {2026},
  eprint = {2511.20553},
  archiveprefix = {arXiv},
  url = {https://arxiv.org/abs/2511.20553},
  note = {Paper No. 26}
}

@article{li2023quadratic,
  author = {Li, Y. and L{\"u}hrmann, J.},
  title = {Soliton dynamics for the 1{D} quadratic {Klein--Gordon} equation with symmetry},
  journal = {J. Differential Equations},
  volume = {344},
  pages = {172--202},
  year = {2023}
}

@article{lindblad2021asymptotics,
  author = {Lindblad, H. and L{\"u}hrmann, J. and Soffer, A.},
  title = {Asymptotics for 1{D} {Klein--Gordon} equations with variable coefficient quadratic nonlinearities},
  journal = {Arch. Ration. Mech. Anal.},
  volume = {241},
  pages = {1459--1527},
  year = {2021}
}

@article{lindblad2023modified,
  author = {Lindblad, H. and L{\"u}hrmann, J. and Schlag, W. and Soffer, A.},
  title = {On modified scattering for 1{D} quadratic {Klein--Gordon} equations with non-generic potentials},
  journal = {Int. Math. Res. Not. IMRN},
  volume = {2023},
  pages = {5118--5208},
  year = {2023}
}

@article{luhrmann2023sinegordon,
  author = {L{\"u}hrmann, J. and Schlag, W.},
  title = {Asymptotic stability of the {sine--Gordon} kink under odd perturbations},
  journal = {Duke Math. J.},
  volume = {172},
  pages = {2715--2820},
  year = {2023}
}

@article{martel2023asymptotic,
  author = {Martel, Y.},
  title = {Asymptotic stability of solitary waves for the 1{D} cubic--quintic {Schr{\"o}dinger} equation with no internal mode},
  journal = {Probab. Math. Phys.},
  volume = {3},
  number = {4},
  pages = {839--867},
  year = {2022},
}

@article{martel2024survey,
  author = {Martel, Y.},
  title = {Asymptotic stability of small solitons for one-dimensional nonlinear {Schr{\"o}dinger} equations},
  journal = {S{\'e}m. Laurent Schwartz EDP Appl.},
  pages = {1--22},
  year = {2024},
  note = {2023--2024, Exp. No. X}
}

@article{martel2024asymptotic,
  author = {Martel, Y.},
  title = {Asymptotic stability of small standing solitary waves of the one-dimensional cubic--quintic {Schr{\"o}dinger} equation},
  journal = {Invent. Math.},
  volume = {237},
  number = {3},
  pages = {1253--1328},
  year = {2024}
}

@article{rialland2024absence,
  author = {Rialland, G.},
  title = {Asymptotic stability of solitary waves for the 1{D} near-cubic nonlinear {Schr{\"o}dinger} equation in the absence of internal modes},
  journal = {Nonlinear Anal.},
  volume = {241},
  year = {2024},
  note = {Paper No. 113474, 30 pp.}
}

@article{rialland2025internal,
  author = {Rialland, G.},
  title = {Asymptotic stability of solitons for near-cubic {NLS} equation with an internal mode},
  journal = {Discrete Contin. Dyn. Syst.},
  volume = {45},
  pages = {1591--1652},
  year = {2025}
}

@article{segur1983wobbling,
  author = {Segur, H.},
  title = {Wobbling kinks in {$\phi^4$} and {sine--Gordon} theory},
  journal = {J. Math. Phys.},
  volume = {24},
  number = {6},
  pages = {1439--1443},
  year = {1983},
}

@article{shatah1985normal,
  author = {Shatah, J.},
  title = {Normal forms and quadratic {Klein--Gordon} equations},
  journal = {Comm. Pure Appl. Math.},
  volume = {38},
  pages = {685--696},
  year = {1985}
}

@article{sigal1993instability,
  author = {Sigal, I. M.},
  title = {Non-linear wave and {Schr{\"o}dinger} equations. {I}. {Instability} of periodic and quasiperiodic solutions},
  journal = {Comm. Math. Phys.},
  volume = {153},
  pages = {297--320},
  year = {1993}
}

@article{simon1976bound,
  author = {Simon, B.},
  title = {The bound state of weakly coupled {Schr{\"o}dinger} operators in one and two dimensions},
  journal = {Ann. Physics},
  volume = {97},
  pages = {279--288},
  year = {1976}
}

@article{soffer1999resonances,
  author = {Soffer, A. and Weinstein, M. I.},
  title = {Resonances, radiation damping and instability in {Hamiltonian} nonlinear wave equations},
  journal = {Invent. Math.},
  volume = {136},
  number = {1},
  pages = {9--74},
  year = {1999},
}

@article{Cuccagna_2024_CenterHypersurface,
  author   = {Cuccagna, S. and Maeda, M. and Murgante, F. and Scrobogna, S.},
  title    = {On asymptotic stability on a center hypersurface at the soliton for even solutions of the nonlinear {Klein--Gordon} equation when $2\ge p>\frac{5}{3}$},
  journal  = {SIAM J. Math. Anal.},
  volume   = {56},
  number   = {4},
  pages    = {5445--5473},
  year     = {2024},
}

@article{Snelson_2018_VariableSpeedKink,
  author   = {Snelson, S.},
  title    = {Asymptotic stability for odd perturbations of the stationary kink in the variable-speed {$\phi^4$} model},
  journal  = {Trans. Amer. Math. Soc.},
  volume   = {370},
  number   = {10},
  pages    = {7437--7460},
  year     = {2018}
}

@article{leger2021internal,
  title={{Internal modes and radiation damping for quadratic Klein-Gordon in 3D}},
  author={L{\'e}ger, T. and Pusateri, F.},
  journal={arXiv preprint arXiv:2112.13163},
  year={2021}
}

@article{palacios2026local,
  title={Local Energy control in the presence of a zero-energy resonance},
  author={Palacios, J. M. and Pusateri, F.},
  journal={Archive for Rational Mechanics and Analysis},
  volume={250},
  number={4},
  pages={66},
  year={2026},
  publisher={Springer}
}

@article{luhrmann2026asymptotic,
  title={{Asymptotic stability of the degree-one vortex in the abelian Yang-Mills-Higgs model: Linear Theory}},
  author={L{\"u}hrmann, J. and Palacios, J. M. and Pusateri, F. and Schlag, W. and Shahshahani, S.},
  journal={arXiv preprint arXiv:2608.10609},
  year={2026}
}

@article{luhrmann2026asymptoticstabilitydegreeonevortex,
title={{Asymptotic stability of the degree-one vortex in the abelian Yang-Mills-Higgs model: Spectral Theory and Numerics}}, 
author={L{\"u}hrmann, J. and Palacios, J. M. and Pusateri, F. and Schlag, W. and Shahshahani, S.},
      year={2026},
      eprint={2608.10608},
      archivePrefix={arXiv},
      primaryClass={math.AP},
      url={https://arxiv.org/abs/2608.10608}, 
}

@article{chen2026global,
  title={{Global dynamics for the 1D quartic Klein-Gordon equation with an internal mode}},
  author={Chen, G. and Diebou, G. Y. and Kairzhan, A. and Pusateri, F.},
  journal={arXiv preprint arXiv:2610.06456},
  year={2026}
}

@article{luhrmann2026asymptoticstabilitydegreeonevortexstability,
      title={{Asymptotic stability of the degree-one vortex in the abelian Yang-Mills-Higgs model}}, 
author={L{\"u}hrmann, J. and Palacios, J. M. and Pusateri, F. and Schlag, W. and Shahshahani, S.},
      year={2026},
      eprint={2608.10610},
      archivePrefix={arXiv},
      primaryClass={math.AP},
      url={https://arxiv.org/abs/2608.10610}, 
}

@book{Agmon1982Lectures,
  author    = {Agmon, S.},
  title     = {{Lectures on Exponential Decay of Solutions of Second-Order Elliptic Equations: Bounds on Eigenfunctions of N-Body Schr\"odinger Operators}},
  series    = {Mathematical Notes},
  number    = {29},
  publisher = {Princeton University Press},
  address   = {Princeton, N.J.},
  year      = {1982},
  isbn      = {9780691069166},
  doi       = {10.1515/9781400853076}
}

@article{cuccagna2025fgr,
  author = {Cuccagna, S. and Maeda, M.},
  title = {A note on the {Fermi Golden Rule} constant for the pure power {NLS}},
  journal = {Rend. Istit. Mat. Univ. Trieste},
  volume = {57},
  pages = {237--266},
  year = {2025}
}

@article {PoschlTeller,
       author = {{P{\"o}schl}, G. and {Teller}, E.},
        title = {{Bemerkungen zur Quantenmechanik des anharmonischen Oszillators}},
      journal = {Zeitschrift f{\"u}r Physik},
         year = {1933},
       volume = {83},
       number = {3-4},
        pages = {143-151},
}

@article{LL2, 
title={{Asymptotic stability of solitary waves for the 1D focusing cubic Schrödinger equation under even perturbations}}, volume={14}, DOI={10.1017/fms.2026.10268}, journal={Forum of Mathematics, Sigma}, author={Li, Y. and L{\"u}hrmann, J.}, year={2026}, pages={e125}}

@book {SteinShakarchi_Complex,
    AUTHOR = {Stein, E. M. and Shakarchi, R.},
     TITLE = {Complex analysis},
    SERIES = {Princeton Lectures in Analysis},
    VOLUME = {2},
 PUBLISHER = {Princeton University Press, Princeton, NJ},
      YEAR = {2003},
     PAGES = {xviii+379},
}

@article{kowalczyk2017asymptotic,
  title={On asymptotic stability of nonlinear waves},
  author={Kowalczyk, M. and Martel, Y. and Mu{\~n}oz, C.},
  journal={S{\'e}minaire Laurent Schwartz—EDP et applications},
  pages={1--27},
  year={2017}
}

@article{koch2026asymptotic,
  title={{Asymptotic stability of the sine-Gordon kinks under perturbations in weighted Sobolev norms}},
  author={Koch, H. and Yu, D.},
  journal={Communications in Partial Differential Equations},
  volume={51},
  number={2-3},
  pages={131--192},
  year={2026},
  publisher={Taylor \& Francis}
}

@article{chen2020long,
  title={{Long-time asymptotics and stability for the sine-Gordon equation}},
  author={Chen, G. and Liu, J. and Lu, B.},
  journal={arXiv preprint arXiv:2009.04260},
  year={2020}
}

@article{kowalczyk-martel26,
    author = {Kowalczyk, M. and Martel, Y.},
    title = {Asymptotic stability of moving kinks in the bag model},
    journal = {In preparation},
    year = {2026}
}

\end{document}